\documentclass[11pt]{article}

\usepackage{a4,amsfonts,exscale}

\usepackage{srcltx,hyperref}

\long\def\forget#1{}

\forget{
\pdfoutput=1

\usepackage[pdftex]{hyperref}
\hypersetup{%
 pdftitle     = {Period Spaces in Equal Characteristic},
 pdfauthor    = {Urs Hartl},
 pdfcreator = {pdfLaTeX},
 pdfproducer = {pdfLaTeX with hyperref}
}
}
\usepackage{amsmath}
\usepackage{amssymb}
\usepackage{amscd}
\usepackage{textcomp}

\usepackage{amsthm}
\theoremstyle{plain}
\newtheorem{theorem}{Theorem}[subsection]

\newtheorem{lemma}[theorem]{Lemma}

\newtheorem{corollary}[theorem]{Corollary}
\newtheorem{proposition}[theorem]{Proposition}

\theoremstyle{definition}
\newtheorem{definition}[theorem]{Definition}

\newtheorem{example}[theorem]{Example}
\newtheorem{remark}[theorem]{Remark}

\theoremstyle{remark}

\usepackage{xy}
\xyoption{all}

\newdir^{ (}{{}*!/-3pt/\dir^{(}}    
\newdir^{  }{{}*!/-3pt/\dir^{}}    
\newdir_{ (}{{}*!/-3pt/\dir_{(}}    
\newdir_{  }{{}*!/-3pt/\dir_{}}

\newcounter{zahl}

\def\theenumi{(\alph{enumi})}

\def\p@enumii{\theenumi}

\newcommand{\DS}{\displaystyle}
\newcommand{\TS}{\textstyle}
\newcommand{\SC}{\scriptstyle}
\newcommand{\SSC}{\scriptscriptstyle}

\DeclareMathOperator{\Aut}{Aut}

\DeclareMathOperator{\Cov}{Cov}

\DeclareMathOperator{\End}{End}

\DeclareMathOperator{\Ext}{Ext}

\DeclareMathOperator{\Gal}{Gal}
\DeclareMathOperator{\GL}{GL}
\DeclareMathOperator{\Koh}{H}

\DeclareMathOperator{\Hom}{Hom}
\newcommand{\CHom}{{\cal H}om}
\DeclareMathOperator{\Id}{Id}

\DeclareMathOperator{\PGL}{PGL}

\DeclareMathOperator{\Quot}{Frac}

\DeclareMathOperator{\Spm}{Sp}
\DeclareMathOperator{\Spec}{Spec}
\DeclareMathOperator{\Spf}{Spf}
\DeclareMathOperator{\Stab}{Stab}

\DeclareMathOperator{\Tor}{Tor}

\DeclareMathOperator{\Var}{V}

\newcommand{\ad}{{\rm ad}}
\newcommand{\alg}{{\rm alg}}
\newcommand{\an}{{\rm an}}

\DeclareMathOperator{\coker}{coker}

\DeclareMathOperator{\diag}{diag}
\newcommand{\et}{{\rm\acute{e}t}}

\DeclareMathOperator{\id}{\,id}
\DeclareMathOperator{\im}{im}

\renewcommand{\mod}{{\rm\,mod\,}}
\DeclareMathOperator{\ord}{ord}
\newcommand{\perf}{{\rm perf}}

\newcommand{\rig}{{\rm rig}}
\DeclareMathOperator{\rk}{rk}
\newcommand{\sep}{{\rm sep}}

\renewcommand{\phi}{\varphi}
\renewcommand{\epsilon}{\varepsilon}
\usepackage{amsfonts}
\newcommand{\BOne} {{\mathchoice{\hbox{\rm1\kern-2.7pt l\kern.9pt}}
                              {\hbox{\rm1\kern-2.7pt l\kern.9pt}}
                              {\hbox{\scriptsize\rm1\kern-2.3pt l\kern.4pt}}
                              {\hbox{\scriptsize\rm1\kern-2.4pt l\kern.5pt}}}}

\newcommand{\BA}{{\mathbb{A}}}

\newcommand{\BC}{{\mathbb{C}}}
\newcommand{\BD}{{\mathbb{D}}}

\newcommand{\BF}{{\mathbb{F}}}
\newcommand{\BG}{{\mathbb{G}}}
\newcommand{\BH}{{\mathbb{H}}}

\newcommand{\BN}{{\mathbb{N}}}

\newcommand{\BP}{{\mathbb{P}}}
\newcommand{\BQ}{{\mathbb{Q}}}
\newcommand{\BR}{{\mathbb{R}}}
\newcommand{\BS}{{\mathbb{S}}}

\newcommand{\BZ}{{\mathbb{Z}}}

\newcommand{\CC}{{\cal{C}}}

\newcommand{\CE}{{\cal{E}}}
\newcommand{\CF}{{\cal{F}}}

\newcommand{\CH}{{\cal{H}}}
\newcommand{\CI}{{\cal{I}}}
\newcommand{\CJ}{{\cal{J}}}
\newcommand{\CK}{{\cal{K}}}

\newcommand{\CM}{{\cal{M}}}
\newcommand{\CN}{{\cal{N}}}
\newcommand{\CO}{{\cal{O}}}
\newcommand{\CP}{{\cal{P}}}
\newcommand{\CQ}{{\cal{Q}}}
\newcommand{\CR}{{\cal{R}}}

\newcommand{\CV}{{\cal{V}}}

\newcommand{\FS}{{\mathfrak{S}}}

\newcommand{\Fa}{{\mathfrak{a}}}

\newcommand{\Fm}{{\mathfrak{m}}}
\newcommand{\Fn}{{\mathfrak{n}}}

\newcommand{\Fp}{{\mathfrak{p}}}
\newcommand{\Fq}{{\mathfrak{q}}}

\newcommand{\rinj}{ \mbox{\mathsurround=0pt \;\raisebox{0.63ex}{\small $\subset$} \hspace{-1.07em} $\longrightarrow$\;}}

\newcommand{\lbij}{\mbox{\mathsurround=0pt \;$\stackrel{\quad \sim}{\longleftarrow}\hspace{-0.75em}-$\;}}

\let\setminus\smallsetminus

\newcommand{\es}{\enspace}
\newcommand{\open}{^\circ}

\newcommand{\dual}{^{\SSC\lor}}

\newcommand{\kompl}{^{\SSC\land}}
\newcommand{\mal}{^{\SSC\times}}
\newcommand{\fdot}{{\,{\SSC\bullet}\,}}
\newcommand{\ul}[1]{{\underline{#1}}}
\newcommand{\ol}[1]{{\overline{#1}}}
\newcommand{\wh}[1]{{\widehat{#1}}}
\newcommand{\wt}[1]{{\widetilde{#1}}}

\usepackage{ifthen}

\newcommand{\invlim}[1][]{\ifthenelse{\equal{#1}{}}% falls Argument leer
{\DS \lim_{\longleftarrow}}%                         verwende niedrige Version
{\DS \lim_{\underset{#1}{\longleftarrow}}}%  sonst:  verwende Argument
}

\newcommand{\dirlim}[1][]{\ifthenelse{\equal{#1}{}}% falls Argument leer
{\DS \lim_{\longrightarrow}}%                        verwende niedrige Version
{\DS \lim_{\underset{#1}{\longrightarrow}}}% sonst:  verwende Argument
}

\newcommand{\dbl}{{\mathchoice{\mbox{\rm [\hspace{-0.15em}[}}
                              {\mbox{\rm [\hspace{-0.15em}[}}
                              {\mbox{\scriptsize\rm [\hspace{-0.15em}[}}
                              {\mbox{\tiny\rm [\hspace{-0.15em}[}}}}
\newcommand{\dbr}{{\mathchoice{\mbox{\rm ]\hspace{-0.15em}]}}
                              {\mbox{\rm ]\hspace{-0.15em}]}}
                              {\mbox{\scriptsize\rm ]\hspace{-0.15em}]}}
                              {\mbox{\tiny\rm ]\hspace{-0.15em}]}}}}
\newcommand{\dpl}{{\mathchoice{\mbox{\rm (\hspace{-0.15em}(}}
                              {\mbox{\rm (\hspace{-0.15em}(}}
                              {\mbox{\scriptsize\rm (\hspace{-0.15em}(}}
                              {\mbox{\tiny\rm (\hspace{-0.15em}(}}}}
\newcommand{\dpr}{{\mathchoice{\mbox{\rm )\hspace{-0.15em})}}
                              {\mbox{\rm )\hspace{-0.15em})}}
                              {\mbox{\scriptsize\rm )\hspace{-0.15em})}}
                              {\mbox{\tiny\rm )\hspace{-0.15em})}}}}

\newcommand{\BFZ}{{\BF_q\dpl\zeta\dpr}}
\newcommand{\ancon}[1][]{{\mathchoice
           {\TS\langle\frac{z}{\zeta^{#1}},z^{-1}\}}
           {\TS\langle\frac{z}{\zeta^{#1}},z^{-1}\}}
           {\SC\langle\frac{z}{\zeta^{#1}},z^{-1}\}}
           {\SSC\langle\frac{z}{\zeta^{#1}},z^{-1}\}}}}
\newcommand{\overcon}{\{z,z^{-1}\}}
\newcommand{\con}[1][]{{\mathchoice
           {\TS\langle\frac{z}{\zeta^{#1}}\rangle[z^{-1}]}
           {\TS\langle\frac{z}{\zeta^{#1}}\rangle[z^{-1}]}
           {\SC\langle\frac{z}{\zeta^{#1}}\rangle[z^{-1}]}
           {\SSC\langle\frac{z}{\zeta^{#1}}\rangle[z^{-1}]}}}

\newcommand{\BLoc}{\BF_q\dbl z\dbr\mbox{-}\ul{\rm Loc}}
\newcommand{\PLoc}{\BF_q\dpl z\dpr\mbox{-}\ul{\rm Loc}}
\DeclareMathOperator{\pair}{\BS}

\newcommand{\dotBD}{\vbox{\hbox{\kern2pt\bf.}\vskip-4.5pt\hbox{$\BD$}}}

\newcommand{\tminus}{t_{\SSC -}}
\newcommand{\tplus}{t_{\SSC +}}
\def\ulM{{\underline{M\!}\,}{}}
\def\ulN{{\underline{N\!}\,}{}}

\def\?{\ 

???\ \immediate\write16{}

\immediate\write16{Warning: There was still a question mark . . . }
\immediate\write16{}}

\def\longto{\longrightarrow}
\def\into{\rinj}

\DeclareMathOperator{\isoto}{\mbox{\hspace{1mm}\raisebox{+1.4mm}{$\SC\sim$}\hspace{-3.5mm}$\longrightarrow$}}

\def\leftmapsto{\gets\joinrel\shortmid}

\newbox\mybox
\def\arrover#1{\mathrel{
       \setbox\mybox=\hbox spread 1.4em{\hfil$\scriptstyle#1$\hfil}
       \vbox{\offinterlineskip\copy\mybox
             \hbox to\wd\mybox{\rightarrowfill}}}}

\begin{document}

%%%%%%%%%%%%%%%%%%%%%%%%%%%%%%%%%%%%%%%%%%%%%%%%%%%%%%%%%%%%%%%%%%%%%%

\author{Urs Hartl
\footnote{The author acknowledges support of the Deutsche Forschungsgemeinschaft in form of DFG-grant HA3006/2-1, SFB 478 and SFB 878}
}
\date{March 28, 2011}

\title{Period Spaces for Hodge Structures in Equal Characteristic}

\maketitle

\begin{abstract}
We develop the analog in equal positive characteristic of Fontaine's theory for crystalline %$p$-adic 
Galois representations of a $p$-adic field. In particular we describe the analog of Fontaine's functor which assigns to a crystalline Galois representation an isocrystal with a Hodge filtration. In equal characteristic the role of isocrystals and Hodge filtrations is played by $z$-isocrystals and Hodge-Pink structures. The latter were invented by Pink. Our first main result in this article is the analog of the Colmez-Fontaine Theorem that ``weakly admissible implies admissible''. Next we construct period spaces for Hodge-Pink structures on a fixed $z$-isocrystal. These period spaces are analogs of the Rapoport-Zink period spaces for Fontaine's filtered isocrystals in mixed characteristic and likewise are rigid analytic spaces. For our period spaces we prove the analog of a conjecture of Rapoport and Zink stating the existence of a ``universal local system'' on a Berkovich open subspace of the period space. As a consequence of ``weakly admissible implies admissible'' this Berkovich open subspace contains every classical rigid analytic point of the period space. 
As the principal tool to demonstrate these results we use the analog of Kedlaya's Slope Filtration Theorem which we also formulate and prove here.

\end{abstract}

\tableofcontents

\newcommand{\refAppRigFormal}{A.1}
\newcommand{\refAppBerkovich}{A.2}
\newcommand{\refAppEtaleSheaves}{A.3}
\newcommand{\refAppFundamentalGroups}{A.4}
\newcommand{\refAppApproxLemma}{A.5}
\newcommand{\refAppExtensionLemma}{A.6}

\bigskip

%%%%%%%%%%%%%%%%%%%%%%%%%%%%%%%%%%%%%%%%%%%%%%%%%%%%%%%%%%%%%%%%%%%%%%
%
%    Introduction
%
%%%%%%%%%%%%%%%%%%%%%%%%%%%%%%%%%%%%%%%%%%%%%%%%%%%%%%%%%%%%%%%%%%%%%%

\subsection*{Introduction}
\addcontentsline{toc}{section}{Introduction}

In \cite{Grothendieck} Grothendieck posed the problem to define a functor (the ``mysterious functor'') relating the $p$-adic \'etale cohomology $\Koh^i_\et(X_{\ol K},\BQ_p)$ of a proper smooth scheme $X$ over a $p$-adic field $K$ with good reduction to the crystalline cohomology $\Koh^i_{cris}(X)$ of $X$. This functor had before been constructed by Grothendieck~\cite{Grothendieck} and Tate~\cite{Tate2} in case $X$ is an abelian variety. The first step to solve the problem was Fontaine's~\cite{Fontaine2} definition of a functor from suitable $p$-adic representations of $\Gal(K^\sep/K)$, such as $\Koh^i_\et(X_{\ol K},\BQ_p)$, which Fontaine called \emph{crystalline} to filtered isocrystals, such as $\Koh^i_{cris}(X)$ with its Hodge-Tate filtration constructed by Deligne-Illusie~\cite{DI}. The problem was finally solved by Faltings~\cite{Faltings88} who showed that indeed $\Koh^i_\et(X_{\ol K},\BQ_p)$ is crystalline when $X$ is proper and smooth over $K$ with good reduction. See for example \cite{Illusie} for a survey of this theory which is called \emph{$p$-adic Hodge Theory}.

There remained various open questions. One of them was to describe which filtered isocrystals arise from crystalline Galois representations. These were called \emph{admissible}. There is a numerical criterion which is easily seen to be necessary and the filtered isocrystals satisfying this criterion were called \emph{weakly admissible}. Fontaine conjectured \cite[Conjecture 5.3]{Fontaine2} that weakly admissible implies admissible and this was proved by Colmez-Fontaine~\cite{CF}. Meanwhile there are many more ways to prove this conjecture \cite{Colmez1,Berger1,Kisin,GL10}. An important approach uses the relation with $p$-adic differential equations (also called $(\phi,\Gamma)$-modules over the Robba ring) introduced by Berger~\cite{Berger02}, and the Slope Filtration Theorem of Kedlaya~\cite{Kedlaya} which says that every $\phi$-module over the Robba ring admits a filtration with isoclinic factors of increasing slopes.

In this article we want to develop the analog for part of $p$-adic Hodge theory in the \emph{Arithmetic of Function Fields over Finite Fields}. We call this analog \emph{Hodge-Pink Theory} for reasons we give below. We should make clear from the beginning that we concentrate here on the analog of Fontaine's theory which relates crystalline Galois representations of $p$-adic fields to filtered isocrystals. We cannot say anything about the relation of our theory with the cohomology of varieties in positive characteristic. Let us explain the translation from $p$-adic Hodge theory to Hodge-Pink theory. $p$-adic Hodge theory is a sector of the arithmetic of the local field $\BQ_p$ and its ring of integers $\BZ_p$. In the Arithmetic of Function Fields the latter are replaced by the Laurent series field $\BF_q\dpl z\dpr$ over the finite field $\BF_q$ with $q$ elements and its ring of integers, the power series ring $\BF_q\dbl z\dbr$. The ring of \emph{Witt vectors} $W(k)$ over a perfect field $k$ containing $\BF_p$ is replaced by the power series ring $\ell\dbl z\dbr$ over a perfect field $\ell$ containing $\BF_q$. There even is no need to require $\ell$ to be perfect. The \emph{Frobenius lift} on $\ell\dbl z\dbr$ is the endomorphism
\[
\sigma:\es\ell\dbl z\dbr\longto \ell\dbl z\dbr\,,\es z\mapsto z\,,\es b\mapsto b^q\text{ for all }b\in\ell\,.
\]
Thus the role of the natural number $p$ in mixed characteristic is taken over by the indeterminate $z$ in equal characteristic. But $p$ enters in $p$-adic Hodge theory in a twofold way. Firstly it is the uniformizing parameter of $\BZ_p$ and  $W(k)$ and secondly it is an element of the base field $K$ over which the objects like varieties, $\Gal(K^\sep/K)$-representations, and filtered isocrystals are defined. The necessity to separate these two roles of $p$ in equal characteristic and to work with two indeterminates $z$ and $\zeta$ was first pointed out by Anderson~\cite{Anderson}. The moral reason is that the natural number $p$ never can act on a module or vector space as anything else than the scalar $p$ whereas there is no such restriction on $z$. Hence we use $z$ for the uniformizing parameter of the analogs of $\BZ_p$ and $W(k)$ and $\zeta$ for the element of the base field $L$ and we let $\sigma$ act as the $q$-Frobenius on $L$ and $\zeta$. (Strictly speaking this distinction between the two roles of $p$ is also searched for in Fontaine's theory where an object, called $[\tilde p]$ by Colmez~\cite{Colmez}, is constructed that behaves like our $\zeta$, whereas $p$ behaves like our $z$.)

Now the concept of \emph{$F$-isocrystal} over $k$ is translated into a finite dimensional $\ell\dpl z\dpr$-vector space $D$ together with an isomorphism $F_D:\sigma^\ast D\to D$ of $\ell\dpl z\dpr$-vector spaces, where we abbreviate $\sigma^\ast D:=D\otimes_{\ell\dpl z\dpr,\sigma}\ell\dpl z\dpr$. We call the pair $(D,F_D)$ a \emph{$z$-isocrystal} over $\ell$. The parallel between $F$-isocrystals and $z$-isocrystals reaches quite far, see \cite{Anderson2,Hartl,HartlDict}. The analog of \emph{filtered isocrystals} is defined as follows. Let $L$ be a field extension of $\BFZ$ which is complete with respect to an absolute value $|\,.\,|:L\to\BR_{\geq0}$ extending the absolute value on $\BFZ$. Let $\CO_L$ be its valuation ring and let $\ell$ be its residue field. Assume that there is a section $\ell\hookrightarrow\CO_L$ of the residue map $\CO_L\to\ell$ and fix one. This induces in particular a homomorphism $\ell\dpl z\dpr\to L\dbl z-\zeta\dbr$ into the power series ring over $L$ sending $z$ to $z=\zeta+(z-\zeta)$. We make the following

\bigskip

\noindent
{\bf Definition.} A \emph{Hodge-Pink structure} over $L$ on the $z$-isocrystal $(D,F_D)$ is an $L\dbl z-\zeta\dbr$-lattice $\Fq_D$ inside $\sigma^\ast D\otimes_{\ell\dpl z\dpr}L\dpl z-\zeta\dpr$.

\bigskip

We also pose $\Fp_D:=\sigma^\ast D\otimes_{\ell\dpl z\dpr}L\dbl z-\zeta\dbr$. Every Hodge-Pink structure determines in particular a \emph{Hodge-Pink filtration} $Fil^\bullet$ on the $L$-vector space $D_L:=\sigma^\ast D\otimes_{\ell\dpl z\dpr}L\dbl z-\zeta\dbr/(z-\zeta)=\sigma^\ast D\otimes_{\ell\dpl z\dpr,\,z\mapsto\zeta}L$ by letting $Fil^i D_L$ be the image of 
\[
\Fp_D\;\cap\;(z-\zeta)^i\Fq_D
\]
in $D_L$. Note that $z$ acts on $D_L$ as the scalar $\zeta$. However the fact that $z$ and $\zeta$ both play part of the role of $p$ makes it necessary to consider Hodge-Pink structures instead of only Hodge-Pink filtrations to get a reasonable category. This was first observed by Pink~\cite{Pink} who developed Hodge theory over $\BF_q\dpl z\dpr$ under this viewpoint and for this reason we propose the name ``Hodge-Pink Theory''. In Remark~\ref{Remark7.3} we discuss more reasons for the necessity to consider Hodge-Pink structures instead of Hodge-Pink filtrations.

With a $z$-isocrystal with Hodge-Pink structure $\ul D=(D,F_D,\Fq_D)$ of rank $n$ over $L$ one associates two numerical invariants, its \emph{Newton slope} $t_N(\ul D)$ and its \emph{Hodge slope} $t_H(\ul D)$. The former is the slope of $\wedge^n(D,F_D)$, that is $\ord_z(\det F_D)$ with respect to some (any) $\ell\dpl z\dpr$-basis of $D$. The latter is the integer $e$ such that $\wedge^n\Fq_D=(z-\zeta)^{-e}\wedge^n\Fp_D$. As in Fontaine's theory $\ul D$ is called \emph{weakly admissible} if $t_H(\ul D)=t_N(\ul D)$ and $t_H(\ul D')\leq t_N(\ul D')$ for any subobject $\ul D'\subset\ul D$, that is for any $F_D$-stable $\ell\dpl z\dpr$-subspace $D'\subset D$ and lattice $\Fq_{D'}\subset\Fq_D\cap\sigma^\ast D\otimes_{\ell\dpl z\dpr}L\dpl z-\zeta\dpr$.

To define the analog of \emph{Fontaine's functor} and the notion of \emph{admissibility} the guiding observation is the following. If $K$ is a finite extension of $\BQ_p$ with ring of integers $\CO_K$, Breuil~\cite[Th\'eor\`eme 1.4]{Breuil2} (see also Kisin~\cite[Theorem 0.3]{Kisin}) has shown that any crystalline representation of $\Gal(K^\sep/K)$ with all Hodge-Tate weights equal to $0$ or $1$ arises from a $p$-divisible group over $\CO_K$. Kisin's proof builds on Breuil's work and uses an intermediate category ${\rm BT}^\phi_{/\FS}\otimes\BQ_p$ whose analog in equal characteristic is as follows.

\bigskip

\noindent
{\bf Definition.} A \emph{local shtuka} over $\CO_L$ is a finite free $\CO_L\dbl z\dbr$-module $M$ together with an isomorphism
\[
\TS F_M:\es\sigma^\ast M[\frac{1}{z-\zeta}]\es\isoto\es M[\frac{1}{z-\zeta}]
\]
where the appendage $[\frac{1}{z-\zeta}]$ stands for inverting $z-\zeta$. A \emph{morphism} $f:M\to N$ of local shtukas is a morphism of $\CO_L\dbl z\dbr$-modules which satisfies $f\circ F_M=F_N\circ\sigma^\ast f$. An \emph{isogeny} of local shtukas is a morphism $f:M\to N$ such that for some morphism $g:N\to M$ the compositions $f\circ g$ and $g\circ f$ equal multiplication with a power of $z$.

\bigskip

With this terminology the category ${\rm BT}^\phi_{/\FS}\otimes\BQ_p$ is the analog of the isogeny category of local shtukas for which $F_M$ is actually a morphism $\sigma^\ast M\to M$ whose cokernel is killed by $z-\zeta$. Kisin~\cite[Theorem 0.1]{Kisin} goes on to show that the category of crystalline representations of $\Gal(K^\sep/K)$ is equivalent to a full subcategory of the category analogous to the isogeny category of all local shtukas; see also \cite[\S5.4]{HartlDict}. This observation has led A.\ Genestier and V.\ Lafforgue~\cite{GL} in case $L$ is discretely valued to construct a functor $\BH$ from local shtukas over $\CO_L$ to $z$-isocrystals with Hodge-Pink structure over $L$ and to call the objects in the essential image of $\BH$ \emph{admissible}. In Section~\ref{SectMysterious} we extend the construction of $\BH$ to arbitrary $L$ and investigate its properties. Also in Remark~\ref{RemLokShInsteadOfGalRep} we present some more reasons for this use of the term ``admissible'' which expresses the opinion that local shtukas are the appropriate analog of crystalline Galois representations and $\BH$ is the analog of the mysterious functor. On the other hand, with a local shtuka $M$ of rank $n$ over $\CO_L$ one can associate its \emph{Tate module}
\[
T_zM\es:=\es \{\,x\in M\otimes_{\CO_L\dbl z\dbr}L^\sep\dbl z\dbr:\es F_M(\sigma^\ast x)=x\,\}
\]
 which is a free $\BF_q\dbl z\dbr$-module of rank $n$ and defines a Galois representation
\[
\rho_M:\es\Gal(L^\sep/L)\es\longto\es\GL_n\bigl(\BF_q\dpl z\dpr\bigr)\,.
\]
This provides the link between local shtukas and Galois representations.

Now every admissible $z$-isocrystal with Hodge-Pink structure is weakly admissible. The main result of Chapter~\ref{ChaptHPTheory} is the following

\bigskip

\noindent
{\bf Theorem~\ref{Thm3.5}.} {\it
Let $\wt L$ be the closure of the compositum $\ell^\alg L$ inside a complete algebraically closed extension of $L$ and assume that $\wt L$ does not contain an element $x$ with $0<|x|<1$ such that all $q$-power roots of $x$ also lie in $\wt L$.
Then every weakly admissible $z$-isocrystal with Hodge-Pink structure over $L$ is already admissible.}

\bigskip

Note that the value groups of $L$ and $\wt L$ coincide, and therefore the hypothesis is satisfied if $L$ is discretely valued, or if the value group $\Gamma$ of $L$ is finitely generated, or even if $\Gamma$ does not contain a non-zero element which is arbitrarily often divisible by $q$. Theorem~\ref{Thm3.5} was previously obtained by Genestier-Lafforgue~\cite{GL} for discretely valued $L$ with perfect residue field. Our proof is entirely different and akin to the proofs of Berger~\cite{Berger1} and Kisin~\cite{Kisin} in mixed characteristic. For $L$ violating the assumption of Theorem~\ref{Thm3.5} we present weakly admissible Hodge-Pink structures which are not admissible in Example~\ref{Ex8.2}.

\medskip

In Chapter~\ref{ChaptPeriodSpaces} we construct period spaces for Hodge-Pink structures analogous to the period spaces for filtered isocrystals of Rapoport and Zink~\cite{RZ}. To any fixed $z$-isocrystal $(D,F_D)$ over $\BF_q^{\,\alg}$ and fixed numerical invariants the period space $\CH^{wa}$ parametrizes all weakly admissible Hodge-Pink structures on $(D,F_D)$. It is an admissible open rigid analytic subspace over $\BF_q^{\,\alg}\dpl\zeta\dpr$ of an algebraic variety. But contrary to mixed characteristic this algebraic variety is not a partial flag variety but a partial jet bundle 
over a partial flag variety. This corresponds to the fact that the Hodge-Pink structure contains more information than just the Hodge-Pink filtration; see Remark~\ref{Remark7.3}. To reveal their full information these period spaces are best considered as Berkovich spaces (see Appendix~\refAppBerkovich) instead of rigid analytic spaces. The reason is that one should not only consider their classical rigid analytic points (whose residue field is finite over $\BF_q^{\,\alg}\dpl\zeta\dpr$) but also their $L$-valued points for complete extensions $L$ of $\BF_q^{\,\alg}\dpl\zeta\dpr$ as above. Namely the $z$-isocrystals with Hodge-Pink structure at the former points are all admissible by Theorem~\ref{Thm3.5} whereas this is not true for the latter points. The Berkovich space associated with $\CH^{wa}$ naturally contains also the latter points. Our second main theorem is

\bigskip

\noindent
{\bf Theorem~\ref{Thm2.11}.} {\it
The set $\CH^a$ of admissible Hodge-Pink structures is Berkovich open inside the period space $\CH^{wa}$.}

\bigskip

By Theorem~\ref{Thm3.5} the set $\CH^a$ contains all classical rigid analytic points. We further show that both spaces are connected. In Example~\ref{Ex8.2} we show that $\CH^a$ may be strictly smaller than $\CH^{wa}$. This means that when viewed as rigid analytic spaces the morphism $\CH^a\to\CH^{wa}$ is \'etale and bijective on classical rigid analytic points but need not be an isomorphism of rigid analytic spaces. Then the two rigid analytic spaces have the same underlying set of points but differ by their Grothendieck topology. Our third main result says that not only every point of $\CH^a$ is admissible but that

\bigskip

\noindent
{\bf Theorem~\ref{Thm2.16}.} {\it
There exists an admissible formal scheme $X$ over $\Spf \BF_q^{\,\alg}\dbl\zeta\dbr$ in the sense of Raynaud~\cite{Raynaud} such that its associated rigid analytic space $X^\rig$ is an \'etale covering space of $\CH^a$, and there exists a local shtuka over $X$ giving rise to the universal $z$-isocrystal with Hodge-Pink structure on $\CH^a$.}

\bigskip

As a consequence the Tate module of this local shtuka over $X$ defines a local system $\CV$ of $\BF_q\dpl z\dpr$-vector spaces on $X^\rig$ which descends to $\CH^a$. For any geometric base point $\bar x$ of $\CH^a$ with underlying $L$-rational point $x$ this local system even gives rise to a representation of the \'etale fundamental group $\pi_1^\et(\CH^a,\bar x)$ 
\[
\rho_\CV:\es\pi_1^\et(\CH^a,\bar x)\es\longto\es\GL_n
\]
such that the composition with $\Gal(L^\sep/L)=\pi_1^\et(x,\bar x)\to\pi_1^\et(\CH^a,\bar x)$ yields the Galois representation $\rho_{M_x}$ associated with the local shtuka $M_x$ at $x$ from above.
This situation is analogous to a conjecture of Rapoport and Zink~\cite[p.\ 29]{RZ}, \cite[p.\ 429]{Rapoport94}. Namely, to a reductive algebraic group $G$ over $\BQ_p$ and fixed numerical invariants Rapoport and Zink consider the  period space $\CF^{wa}$ of weakly admissible filtered isocrystals. They conjecture the existence of a morphism $\CF'\to\CF^{wa}$ of rigid analytic spaces which is \'etale and bijective on classical rigid analytic points, and of a ``universal local system'' of $\BQ_p$-vector spaces on $\CF'$ which for any geometric base point $\bar x$ of $\CF'$ with underlying $K$-rational point $x$ induces a representation of the \'etale fundamental group $\rho:\pi_1^\et(\CF',\bar x)\to G$ such that the composition of $\rho$ with $\Gal(K^\sep/K)=\pi_1^\et(x,\bar x)\to\pi_1^\et(\CF',\bar x)$ yields the $\Gal(K^\sep/K)$-representation associated by Fontaine's inverse functor with the filtered isocrystal corresponding to $x$. 
Inspired by the results in the present article a proof for Rapoport's and Zink's conjecture in the case of Hodge-Tate weights $0$ and $1$ was given in Faltings~\cite{Faltings07} and Hartl~\cite{HartlCR-RZ,HartlRZ}.

\medskip

Our main tool to prove the results above is an analog of Kedlaya's Slope Filtration Theorem \cite{Kedlaya} for $\phi$-modules over the Robba ring. The latter comes into play in $p$-adic Hodge theory through Berger's construction assigning to a filtered isocrystal $(D,F_D,Fil^\bullet)$ over $K$ a $(\phi,\Gamma)$-module over the Robba ring which is a modification of $(D,F_D)$ tensored up to the Robba ring 
at the zeroes of $\log[\epsilon]$. These zeroes are precisely the $\sigma$-translates of the point corresponding to the maximal ideal of $B_{dR}^+$. (Using the notation of Colmez~\cite[\S5.6]{Colmez} we like to view $B_{dR}^+$ as the complete local ring at the point of $\Spec\wt A^+$ given by the ideal $(p-[\tilde p])$. This point of view is inspired by the situation in equal characteristic which we expose in this article.)
Correspondingly we consider the point $z=\zeta$ of $L\con$ and its $\sigma$-translates $z=\zeta^{q^i}$, $i\in\BN_0$. For a $z$-isocrystal with Hodge-Pink structure $(D,F_D,\Fq_D)$ over $L$ we modify the $\sigma$-module $\sigma^\ast(D,F_D)\otimes_{\ell\dpl z\dpr}L\con$ at $z=\zeta^{q^i}$ for $i\in\BN_0$ according to the data given by $\Fq_D$. Here $L\langle\frac{z}{\zeta}\rangle$ is the Tate algebra of convergent power series on the disc $\{|z|\le|\zeta|\}$.
We obtain a $\sigma$-module over the ring of rigid analytic functions on the punctured disc $\{0<|z|\leq|\zeta|\}$. We denote this ring by $L\ancon$ and view it as the analog of the Robba ring. 
Note that the elements of $L\con$ have at most poles at $z=0$ whereas the elements of $L\ancon$ may have essential singularities.
Our analog of the Slope Filtration Theorem is then the following

\bigskip

\noindent
{\bf Theorem~\ref{Thm6.13}.} {\it
Let $M$ be a $\sigma$-module over $L\ancon$. Then there exists a unique filtration $0=M_0\subset M_1\subset\ldots\subset M_\ell=M$ of $M$ by saturated $\sigma$-submodules with the following properties:
\begin{enumerate}
\item
For $i=1,\ldots,\ell$ the quotient $M_i/M_{i-1}$ is isoclinic of some slope $s_i$,
\item
$s_1<\ldots<s_\ell$.
\end{enumerate}
}

\medskip

As in Kedlaya's theory the isoclinic $\sigma$-modules descend as follows.

\medskip

\noindent
{\bf Theorem.} {\it (Special case of Theorem~\ref{Thm6.12}.)\\
Let $M$ be a $\sigma$-module over $L\ancon$ which is isoclinic.
Then there exists an isoclinic $\sigma$-module $M'$ over $L\con$ with $M'\otimes L\ancon\cong M$. It is unique up to canonical isomorphism.
}

\bigskip

The demonstration of these two theorems takes up the entire Chapter~\ref{ChaptSFT}. We may even prove the latter result over arbitrary affinoid $L$-algebras $B$ in place of $L$. This generalization to families is in fact the key to prove the analog of the Rapoport-Zink conjecture. Note that in $p$-adic Hodge theory it is as yet impossible to set up Theorem~\ref{Thm6.12} in a relative version. To be precise, there are currently two concepts for ``relative $p$-adic Hodge theory''. The one with relative coefficients is due to Berger, Colmez, Kisin, Liu, and others \cite{BergerColmez,Kisin08,Liu07}. It treats Galois representations still of a local $p$-adic field but on a module over an affinoid $\BQ_p$-algebra. It has the advantage that it makes an explicit separation between the field extensions and the coefficients of the Galois representation, which in the function field case is achieved via the distinction between $z$ and $\zeta$.

The other concept with relative base spaces is developed among others by Andreatta, Brinon, Faltings, Iovita, Tsuji~\cite{Andreatta,AI,AB,Brinon,Faltings88,Faltings02,Tsuji}. It treats representations of the \'etale fundamental group of \emph{certain} rigid analytic spaces on finite dimensional $\BQ_p$-vector spaces. This is the analog of our theory in Chapter~\ref{ChaptSFT}, that seems to be needed for proving the original conjecture of Rapoport and Zink in mixed characteristic. However so far the relative theory of the field of norms does exist only over certain formal $\BZ_p$-algebras (see \cite{Andreatta}) whose associated Berkovich spaces do not form a basis for the topology on a general smooth Berkovich space. Therefore the mixed characteristic relative version of Theorem~\ref{Thm6.12} is yet out of reach. In equal characteristic the relative version over $B$ of the Robba ring is just the ring $B\ancon$ of rigid analytic functions on the relative punctured disc $\{0<|z|\leq|\zeta|\}$ over $B$ and the theory of the field of norms is trivial. The situation is more transparent here since we can separate the two roles of $p$ into the two indeterminates $z$ and $\zeta$ and moreover, since the phenomenon of ramification of field extensions $K$ of $\BQ_p$ is not so dominant in equal characteristic. If these difficulties in mixed characteristic can be overcome, our arguments presented here could be translated to mixed characteristic to yield a proof of the original conjecture of Rapoport and Zink for arbitrary Hodge-Tate weights, extending \cite{Faltings07,HartlCR-RZ,HartlRZ}.

\medskip

Although rigid analytic geometry and even its variant  introduced by Berkovich~\cite{Berkovich1,Berkovich2} enter at various places we have tried to make this article accessible to readers without much background in rigid analytic geometry by working most of the time with modules over affinoid algebras rather than sheaves on rigid analytic spaces.
Also we review the facts we need from rigid analytic geometry in an appendix.

\medskip

{\it Acknowledgments.} I would like to thank V.\ Berkovich for his numerous explanations, A.\ Genestier, V.\ Lafforgue, and M.\ Rapoport for their interest in this work and many helpful discussions, and the unknown referee for his careful reading of the article.

%%%%%%%%%%%%%%%%%%%%%%%%%%%%%%%%%%%%%%%%%%%%%%%%%%%%%%%%%%%%%%%%%%%%%%
%
%    Chapter 1: The Slope Filtration Theorem
%
%%%%%%%%%%%%%%%%%%%%%%%%%%%%%%%%%%%%%%%%%%%%%%%%%%%%%%%%%%%%%%%%%%%%%%

\section{The Slope Filtration Theorem} \label{ChaptSFT}
\setcounter{equation}{0}

In this chapter we formulate and prove the equal characteristic analog of Kedlaya's~\cite{Kedlaya} Slope Filtration Theorem from $p$-adic Hodge theory. The latter is a powerful tool that has given rise to proofs of various famous conjectures, like the conjectures that ``weakly admissible implies admissible'' \cite{Berger1,Kisin}, that ``de Rham implies potentially semistable'' \cite{Berger02,Kedlaya}, or Crew's conjecture which is also called the $p$-adic local monodromy theorem \cite{Kedlaya2}. The results that we obtain in this chapter will enable us to prove the equal characteristic analogs of the fact that ``weakly admissible implies admissible'' in Chapter~\ref{ChaptHPTheory} and the conjecture of Rapoport and Zink mentioned in the introduction in Chapter~\ref{ChaptPeriodSpaces}.
We follow in large parts Kedlaya's \cite{Kedlaya} proof in mixed characteristic.
Nevertheless we present full proofs here for two reasons. First of all in the function field case various arguments can be somewhat simplified. 
And secondly we formulate and prove some analogs of Kedlaya's results in a relative situation. This is necessary for applications to the analog of the Rapoport-Zink conjecture.

%%%%%%%%%%%%%%%%%%%%%%%%%%%%%%%%%%%%%%%%%%%%%%%%%%%%%%%%%%%%%%%%%%%%%%
%
%    Notation
%
%%%%%%%%%%%%%%%%%%%%%%%%%%%%%%%%%%%%%%%%%%%%%%%%%%%%%%%%%%%%%%%%%%%%%%

\subsection{Notation} \label{SectNotation}

In this article we denote by
\begin{tabbing}
$B\con[n]$\es\= \kill
$\BN_0$ \> the set of non-negative integers, \\[2mm]
$\BF_q$ \> the finite field which has $q$ elements and characteristic $p$, \\[2mm]
$\BF_q\dbl\zeta\dbr$ \> the ring of formal power series over $\BF_q$ in the indeterminate $\zeta$, \\[2mm]
$\BFZ$\> its field of fractions, \\[2mm]
$R\supset\BF_q\dbl\zeta\dbr$ \> \parbox[t]{.79\textwidth}{a rank-$1$ valuation ring which is complete and separated with respect to the $\zeta$-adic topology (it is not assumed to be noetherian),} \\[2mm]
$L,\;\Fm_R$ \> the fraction field and the maximal ideal of $R$, \\[2mm]
$\ol L$ \> the completion of an algebraic closure of $L$. It is again algebraically closed. 
\end{tabbing}
We say a field $K$ is a \emph{complete extension} of $L$, if $K$ is a field extension of $L$ equipped with an absolute value $|\,.\,|:K\to\BR_{\geq0}$ which restricts on $L$ to the absolute value of $L$ and if $K$ is complete with respect to this absolute value. We further denote by
\begin{tabbing}
$B\con[n]$\es\= \kill
$B$ \> an affinoid $L$-algebra with an $L$-Banach norm $|\,.\,|$ (see Appendix~\refAppRigFormal), \\[2mm]
$B\dbl z\dbr$ \> the ring of formal power series over $B$ in the indeterminate $z$, \\[2mm]
$B\dbl z\dbr[z^{-1}]$ \> the ring of Laurent series over $B$ in $z$ with finite principal part.
\end{tabbing}

\noindent
For rational numbers $r$ and ${r'}$ with ${r'}\geq r>0$ we define the following $B$-algebras
\begin{tabbing}
$B\con[n]$\es\= \kill
$B\langle\frac{z}{\zeta^r}\rangle$ \>$\DS =\es\bigl\{\;\sum_{i=0}^\infty b_i z^i\in B\dbl z\dbr:\es |b_i|\,|\zeta|^{ri} \to 0 \es(i\to\infty)\;\bigr\}$, \\[2mm]
$B\con[r]$ \>$\DS =\es\bigl\{\;\sum_{i\gg-\infty}^\infty b_i z^i\in B\dbl z\dbr[z^{-1}]:\es |b_i|\,|\zeta|^{ri} \to 0 \es(i\to\infty)\;\bigr\}$, \\[2mm]
$B\langle\frac{\zeta^{r'}}{z}\rangle$\> $\DS =\es\bigl\{\;\sum_{i=-\infty}^0 b_i z^i:\es |b_i|\,|\zeta|^{{r'}i}\to0\es(i\to-\infty)\;\bigr\}$, \\[2mm]
$B\langle\frac{z}{\zeta^r},\frac{\zeta^{r'}}{z}\rangle$\> $\DS =\es\bigl\{\;\sum_{i=-\infty}^\infty b_i z^i:\es |b_i|\,|\zeta|^{ri} \to 0 \es(i\to\infty),\es |b_i|\,|\zeta|^{{r'}i}\to0\es(i\to-\infty)\;\bigr\}$, \\[2mm]
$B\ancon[r]$ \> $\DS =\es\bigcap_{r'\to\infty}\,B\langle{\TS \frac{z}{\zeta^r},\frac{\zeta^{r'}}{z}}\rangle\es =\es\bigl\{\;\sum_{i=-\infty}^\infty b_i z^i:\es |b_i|\,|\zeta|^{r'i}\to0\es(i\to\pm\infty)\text{ for all }r'\ge r\;\bigr\}$, \\[2mm]
$B\{z,\frac{\zeta^{r'}}{z}\rangle$ \> $\DS =\es\bigcap_{r\to0}\,B\langle{\TS \frac{z}{\zeta^r},\frac{\zeta^{r'}}{z}}\rangle\es =\es\bigl\{\;\sum_{i=-\infty}^\infty b_i z^i:\es |b_i|\,|\zeta|^{ri}\to0\es(i\to\pm\infty)\text{ for all }0<r\le r'\;\bigr\}$, \\[2mm]
$B\overcon$ \> $\DS =\es\bigcap_{r\to 0}\,B\ancon[r]\es =\es\bigl\{\;\sum_{i=-\infty}^\infty b_i z^i:\es |b_i|\,|\zeta|^{ri}\to0\es(i\to\pm\infty)\text{ for all }r>0\;\bigr\}$, \\[2mm]
$B\{z\}$\>$\DS =\es\bigcap_{r\to0}B\langle{\TS\frac{z}{\zeta^r}}\rangle\es =\es\bigl\{\;\sum_{i=0}^\infty b_i z^i:\es b_i\in B,\;|b_i\zeta^{ri}|\to0\es(i\to\infty) \es\text{for all }r>0\;\bigr\}$, \\[2mm]
$\|\,.\,\|_r$\> is the norm on $B\langle\frac{z}{\zeta^r},\frac{\zeta^r}{z}\rangle$ given by $\DS\Bigl\|\,\sum_{i\in\BZ} b_iz^i\,\Bigr\|_r\;:=\;\max_{i\in\BZ}|b_i|\,|\zeta|^{ri}$.
\end{tabbing}
Although the above definitions make sense for arbitrary $r$ we consider in this article \emph{only the case where $r$ and $r'$ are positive}. Observe that these $B$-algebras do not depend on the particular choice of the Banach norm on $B$ since all such norms induce the same topology on $B$.
We extend the Banach norm $|\,.\,|$ on $B$ to matrices $A=(A_{\mu\nu})\in M_{m\times n}(B)$ by setting $|A|:=\max |A_{\mu\nu}|$ and similarly for $\|\,.\,\|_r$~. We let

\begin{tabbing}
$\sigma:x\mapsto \sigma(x)=:x^\sigma$ \=\kill
$\sigma:x\mapsto \sigma(x)=:x^\sigma$ \> \parbox[t]{.79\textwidth}{be the ``endomorphisms'' of the above $B$-algebras which act as the identity on $z$ and as $b\mapsto b^q$ on the elements $b\in B$.}
\end{tabbing}

\smallskip

\noindent
Note that $\sigma$ is not quite an endomorphism of $B\langle\frac{z}{\zeta^r}\rangle$ but only a homomorphism $\sigma:B\langle\frac{z}{\zeta^r}\rangle\to B\langle\frac{z}{\zeta^{qr}}\rangle$ and similarly for $B\con[r]$ and $B\ancon[r]$. If $\CR$ is one of these $B$-algebras we let $\CR^\sigma$ \label{RSigma} be the target of $\sigma:\CR\to\CR^\sigma$. Explicitly $\bigl(B\langle\frac{z}{\zeta^r}\rangle\bigr)^\sigma=B\langle\frac{z}{\zeta^{qr}}\rangle$ as well as $\bigl(B\con[r]\bigr)^\sigma=B\con[qr]$ and $\bigl(B\ancon[r]\bigr)^\sigma=B\ancon[qr]$. Note that in addition there is also a natural inclusion of $B$-algebras
\[
\iota:\CR\to\CR^\sigma\,,\quad \TS \sum b_i z^i\mapsto \sum b_i z^i\,.
\]
Furthermore, in the case $\CR=B\langle\frac{z}{\zeta^r},\frac{\zeta^{r'}}{z}\rangle$ with $r'\geq qr$ we set $\CR^\sigma:=B\langle\frac{z}{\zeta^{qr}},\frac{\zeta^{r'}}{z}\rangle$ and in the cases $\CR=B\dbl z\dbr$, $B\dbl z\dbr[z^{-1}]$, $B\langle\frac{\zeta^{r'}}{z}\rangle$, $B\{z,\frac{\zeta^{r'}}{z}\rangle$, and $B\overcon$ we set $\CR^\sigma:=\CR$. Then the above pattern also holds in these cases.
Note the following rules for $\|\,.\,\|_r$ which are straightforward to prove:

\begin{lemma} \label{Lemma6.0}
\begin{enumerate}
\item
$a\in B\langle\frac{z}{\zeta^r},\frac{\zeta^{r'}}{z}\rangle$ if and only if $a^\sigma\in B\langle\frac{z}{\zeta^{qr}},\frac{\zeta^{qr'}}{z}\rangle$.
\item
For all $a\in B\langle\frac{z}{\zeta^r},\frac{\zeta^r}{z}\rangle$ we have $\|a^\sigma\|_{qr}\;=\;\|a\|_r^q$.
\item
In particular for $a=\sum_{i\leq0}a_i z^i\in B\langle\frac{\zeta^r}{z}\rangle$ we have 
\[
\|a^\sigma\|_r \es =\es \max_{i\leq0} |a_i^q|\,|\zeta|^{ir}\es\leq\es \max_{i\leq0}\bigl(|a_i|\,|\zeta|^{ir}\bigr)^q \es=\es \|a\|_r^q\,.
\]
\item
The norm $\|\,.\,\|_r$ is an $L$-Banach norm on $B\langle\frac{z}{\zeta^r},\frac{\zeta^r}{z}\rangle$, on $B\langle\frac{\zeta^r}{z}\rangle$ and on $B\langle\frac{z}{\zeta^r}\rangle$.
\item
For matrices $A,\wt A\in B\langle\frac{z}{\zeta^r},\frac{\zeta^r}{z}\rangle$ we have $\|A\wt A\|_r\leq\|A\|_r\cdot\|\wt A\|_r$.
\qed
\end{enumerate}
\end{lemma}

\medskip

The $B$-algebras defined above have an interpretation in terms of rigid analytic geometry (see Appendix~\refAppRigFormal). For positive rational numbers $r$ and $r'$ let $\BD(r)=\Spm L\langle\frac{z}{\zeta^r}\rangle$ be the disc with radius $|z|=|\zeta^r|$, and if $r'\geq r$ let $A(r,r')=\Spm L\langle\frac{z}{\zeta^r},\frac{\zeta^{r'}}{z}\rangle$ be the annulus with inner radius $|z|=|\zeta^{r'}|$ and outer radius $|z|=|\zeta^r|$. Moreover, let $\dotBD(r)$ be the admissible open subspace of $\BD(r)$ which is the complement of the point $z=0$, and let $\dotBD(0)\open$ be the union of $\dotBD(r)$ for all $r>0$. It is the open punctured unit disc on which $0<|z|<1$ and is also an admissible open subspace of $\BD(0)$. Then $B\langle\frac{z}{\zeta^r}\rangle$ is the ring of global rigid analytic functions on the relative disc $\Spm B\times_L\BD(r)$,\es $B\langle\frac{z}{\zeta^r},\frac{\zeta^{r'}}{z}\rangle$ is the ring of global rigid analytic functions on the relative annulus $\Spm B\times_L A(r,r')$, \es $B\ancon[r]$ is the ring of global rigid analytic functions on the relative punctured disc $\Spm B\times_L\dotBD(r)$, and $B\overcon$ is the ring of global rigid analytic functions on $\Spm B\times_L\dotBD(0)\open$. In particular if $B=L$ is a field then $L\dbl z\dbr$, $L\langle\frac{z}{\zeta^r}\rangle$, $L\con[r]$, and $L\langle\frac{z}{\zeta^r},\frac{\zeta^{r'}}{z}\rangle$ are principal ideal domains (by Lemma~\ref{LemmaPID} below), and $L\dpl z\dpr=L\dbl z\dbr[z^{-1}]$ is a field.

The reader should note that in case $B=L$ the rings $L\con$, $L\dpl z\dpr$, and $L\ancon$ are the equal characteristic analogs of the rings that Kedlaya~\cite{Kedlaya} calls $\Gamma_{\rm con}[\pi^{-1}]$, $\Gamma[\pi^{-1}]$, and $\Gamma_{\rm an,\,con}$, respectively, and Colmez~\cite{Colmez} calls $\mathbf{B}_L^\dagger$, $\mathbf{B}_L$, and $\mathbf{B}^\dagger_{\rig,L}$, respectively; compare \cite[\S\S2.3, 2.4, 2.8]{HartlDict}. Also note that it is clear how to define relative versions over affinoid algebras of these rings in equal characteristic, whereas the same task is not at all obvious in mixed characteristic to say the least.

The following lemma is proved in \cite[Proposition~4]{Lazard} and the subsequent corollary.

\begin{lemma} \label{LemmaPID}
The rings $L\langle\frac{z}{\zeta^r}\rangle$ and $L\langle\frac{z}{\zeta^r},\frac{\zeta^{r'}}{z}\rangle$ are principal ideal domains. Every element $f\in L\langle\frac{z}{\zeta^r}\rangle$ (respectively $f\in L\langle\frac{z}{\zeta^r},\frac{\zeta^{r'}}{z}\rangle$) can be written in the form $f=ug$ for a unit $u\in L\langle\frac{z}{\zeta^r}\rangle\mal$ (respectively $u\in L\langle\frac{z}{\zeta^r},\frac{\zeta^{r'}}{z}\rangle\mal$) and a polynomial $g\in L[z]$ whose zeroes $x\in\ol L$ all satisfy $|x|\leq|\zeta^r|$ (respectively $|\zeta^{r'}|\leq|x|\leq|\zeta^r|$). Moreover, $u$ and $g$ are uniquely determined up to a constant factor in $L\mal$.
\end{lemma}

We shall also need the following lemma about Galois descent.

\begin{lemma} \label{Lemma0.3}
Let $L'$ be a finite Galois extension of $L$ and extend the action of $\Gal(L'/L)$ to $L'\ancon[r]$ by the trivial action on $z$. Then
\begin{enumerate}
\item 
the fixed subring of $L'\ancon[r]$ under the action of $\Gal(L'/L)$ is equal to $L\ancon[r]$.
\item 
the fixed subring of the fraction field $\Quot\bigl( L'\ancon[r]\bigr)$ under the action of $\Gal(L'/L)$ is equal to $\Quot\bigl( L\ancon[r]\bigr)$.
\end{enumerate}
\end{lemma}

\medskip

\begin{proof}
Assertion (a) is straightforward to prove. To prove (b) let $x,y\in L'\ancon[r]$ such that $\frac{x}{y}$ is fixed under $G=\Gal(L'/L)$. Set $y'=\prod_{g\in G} {}^{g\!}y$ and $x'=y'x/y$. Then both $y'$ and $x'$ belong to $L'\ancon[r]$ and are fixed under $G$. By (a), $x',y'\in L\ancon[r]$ and $\frac{x}{y}=\frac{x'}{y'}\in\Quot L\ancon[r]$ as claimed.
\end{proof}

We end this section by recording the following two lemmas which are basic tools in estimating radii of convergence. 

\begin{lemma} \label{Lemma3b}
Let $A=\sum_{i=0}^\infty A_i z^i \in M_n\bigl(B\langle\frac{z}{\zeta^{qr}}\rangle\bigr)$ be a matrix and let
$x=\sum_{j=0}^\infty x_j z^j\in B\dbl z\dbr^n$ be a vector satisfying $x^\sigma=A\,x$.
Then $x$ belongs to $B\langle\frac{z}{\zeta^r}\rangle^n$.
\end{lemma}

\begin{proof}
Since $A\in M_n\bigl(B\langle\frac{z}{\zeta^{qr}}\rangle\bigr)$ there is a constant $c\ge 1$ with
$|A_i \zeta^{qri}|\leq c$ for all $i$. We expand the equation $x^\sigma=A\,x$ as
\[
\bigl(x_j \zeta^{rj}\bigr)^\sigma\es = \es \sum_{i=0}^j \bigl(A_{j-i} \zeta^{qr(j-i)}\bigr) \bigl(x_i \zeta^{ri}\bigr) \,\zeta^{ir(q-1)}\,.
\]
In view of $|\zeta^r|<1$ this implies the estimate
\[
|x_j\zeta^{rj}|^q\es \leq \es c\cdot \max\{\,|x_i\zeta^{ri}|:\,0\leq i\leq j\,\}
\]
from which induction yields $|x_j\zeta^{rj}|\leq c^{1/(q-1)}$ for all $j\ge0$. 
In particular $x\in B\langle\frac{z}{\zeta^{qr}}\rangle^n$. But now the equation $x^\sigma=A\,x$ shows that
$x^\sigma\in B\langle\frac{z}{\zeta^{qr}}\rangle^n$, hence $x\in B\langle\frac{z}{\zeta^r}\rangle^n$ as desired.
\end{proof}

\begin{lemma} \label{Lemma12}
Let $A=\sum_{i=0}^\infty A_i z^i \in M_n\bigl(B\langle\frac{z}{\zeta^q}\rangle\bigr)$ and let 
$x=\sum_{j=-\infty}^\infty x_j z^j\in B\langle\frac{z}{\zeta},\frac{\zeta^q}{z}\rangle^{\oplus n}$ with
$x-A\,x^\sigma \in B\con[q]^{\oplus n}$. Then $x\in B\con^{\oplus n}$.
\end{lemma}

\begin{proof}
There is a constant $c\geq1$ with $|A_i\zeta^{qi}|\leq c$ for all $i$. 
Assume that $x$ does not have finite principal part. Since $x\in B\langle\frac{z}{\zeta},\frac{\zeta}{z}\rangle^{\oplus n}$
we can find a negative integer $m$ with $x_m\neq0$, $x_m-\sum_{i=0}^\infty A_i x_{m-i}^\sigma=0$, $|x_m\zeta^m|=:d\leq c^{-1}$, and $|x_{m-i}\zeta^{m-i}|\leq d$ for all $i\geq0$.
We obtain
\[
|x_m\zeta^m|\es\leq\es|\zeta^{(1-q)m}|\,\max_{i\geq0}\bigl\{\,|A_i\zeta^{qi}|\,|x_{m-i}^\sigma\zeta^{qm-qi}|\,\bigr\}\es<\es c\,d^q\es\leq\es d\,,
\]
a contradiction.
\end{proof}

%%%%%%%%%%%%%%%%%%%%%%%%%%%%%%%%%%%%%%%%%%%%%%%%%%%%%%%%%%%%%%%%%%%%%%
%
%    $\sigma$-Modules
%
%%%%%%%%%%%%%%%%%%%%%%%%%%%%%%%%%%%%%%%%%%%%%%%%%%%%%%%%%%%%%%%%%%%%%%

\subsection{\texorpdfstring{$\sigma$}{Sigma}-Modules} \label{SectSigmaModules}

\begin{definition}\label{DefSigmaModule}
We let $\CR$ be one of the $B$-algebras $B\dbl z\dbr$, $B\dbl z\dbr[z^{-1}]$, $B\langle\frac{z}{\zeta^r}\rangle$, $B\con[r]$, $B\langle\frac{\zeta^{r'}}{z}\rangle$, or $B\langle\frac{z}{\zeta^r},\frac{\zeta^{r'}}{z}\rangle$. The associated $B$-algebra $\CR^\sigma$ was defined on page~\pageref{RSigma}. For an $\CR$-module $M$ we abbreviate $\sigma^\ast M:=M\otimes_{\CR,\sigma}\CR^\sigma$ and $\iota^\ast M:=M\otimes_{\CR,\iota}\CR^\sigma$. A \emph{$\sigma$-module} over $\CR$ is a pair $(M,F_M)$ consisting of a locally free $\CR$-module $M$ of finite rank together with an isomorphism of $\CR^\sigma$-modules $F_M:\sigma^\ast M\isoto \iota^\ast M$. We will usually abuse notation and write $F_M:\sigma^\ast M\isoto M$ instead. 
The \emph{rank} of a $\sigma$-module $M$ is by definition the rank of its underlying locally free module. It is denoted $\rk M$. The abelian group $\Hom_\sigma(M,N)$ of \emph{morphisms} between two $\sigma$-modules $M$ and $N$ consists of all homomorphisms $f:M\to N$ of $\CR$-modules which satisfy $F_N\circ\sigma^\ast f=\iota^\ast f\circ F_M$. 
\end{definition}

\noindent 
{\it Remark.} The definition of $\sigma$-modules over $B\dbl z\dbr$ and over $B\dbl z\dbr[z^{-1}]$ also makes sense if $B$ is an arbitrary $\BF_q$-algebra. The latter setting is of importance in studying the reduction modulo $\zeta$ of $\sigma$-modules and we will work in this framework in Chapters~\ref{ChaptHPTheory} and \ref{ChaptPeriodSpaces}. 

\medskip

To define $\sigma$-modules over $B\ancon[r]$, $B\{z,\frac{\zeta^{r'}}{z}\rangle$, $B\overcon$, and $B\{z\}$ we have to work with sheaves on the corresponding rigid analytic spaces mentioned above. This could be done also for the affinoid $B$-algebras $B\langle\frac{z}{\zeta^r}\rangle$ $B\langle\frac{\zeta^{r'}}{z}\rangle$, and $B\langle\frac{z}{\zeta^r},\frac{\zeta^{r'}}{z}\rangle$. In view of \cite[Proposition 7.3.2/3 and Theorem 9.4.3/3]{BGR} we get for these an equivalent definition as follows.

\begin{definition}\label{DefSigmaModule2}
We let $\CR$ be one of the $B$-algebras just mentioned. A \emph{$\sigma$-module} over $\CR$ is a pair $(M,F_M)$ consisting of a locally free sheaf $M$ of finite rank on the rigid analytic space corresponding to $\CR$ together with an isomorphism of locally free sheaves $F_M:\sigma^\ast M\isoto \iota^\ast M$ on the rigid analytic space corresponding to $\CR^\sigma$. 
The \emph{rank} $\rk M$ of a $\sigma$-module $M$ is by definition the rank of its underlying locally free sheaf. The abelian group $\Hom_\sigma(M,N)$ of \emph{morphisms} between two $\sigma$-modules $M$ and $N$ consists of all homomorphisms of the sheaves $f:M\to N$ which satisfy $F_N\circ\sigma^\ast f=\iota^\ast f\circ F_M$. 
\end{definition}

The reader who prefers the more down-to-earth formulation with rings and modules should note that by Proposition~\ref{Prop0.1} below the study of $\sigma$-modules over $B\ancon[r]$ is equivalent to the study of $\sigma$-modules over $B\langle\frac{z}{\zeta^r},\frac{\zeta^{qr}}{z}\rangle$. Moreover, $\sigma$-modules over $B\overcon$ are of minor, mainly technical interest in this article. So it is safe to work with Definition~\ref{DefSigmaModule} throughout.

\medskip

We want to introduce a few operations on $\sigma$-modules. 
Let $\CR$ be one of the above $B$-algebras and let $M$ and $N$ be $\sigma$-modules over $\CR$. The \emph{tensor product} $M\otimes N$ is defined in the obvious way as the tensor product of the underlying locally free modules $M\otimes_\CR N$ with the isomorphism $F_{M\otimes N}=F_M\otimes F_N$. Likewise \emph{symmetric} and \emph{alternating powers} are defined. The $\sigma$-module $(\CR,F=\sigma)$ is  a \emph{unit object} for the tensor product.

The \emph{inner hom} $\CHom(M,N)$ is defined as the inner hom of the underlying locally free modules $H=\CHom_\CR(M,N)$ together with $F_H:\sigma^\ast H\isoto H,\;f\mapsto F_N\circ f\circ F_M^{-1}$. In particular the \emph{dual} of a $\sigma$-module is defined as $M\dual=\CHom\bigl(M,(\CR,\sigma)\bigr)$. These definitions satisfy the usual compatibilities and make the category of $\sigma$-modules over $\CR$ into an additive rigid tensor category; see \cite{DM}. It is $\BF_q\dbl z\dbr$-linear and even $\BF_q\dpl z\dpr$-linear if $z$ is invertible in $\CR$. If in Definition~\ref{DefSigmaModule} we allow arbitrary finitely generated $\CR$-modules $M$ instead of only locally free ones and arbitrary morphisms $F_M:\sigma^\ast M\to \iota^\ast M$ instead of only isomorphisms, we obtain an abelian category. The category of $\sigma$-modules is a full subcategory of this abelian category. 

Let $n$ be a positive integer, then $(M,F_M)\mapsto(M,F_M^n:=F_M\circ\sigma^\ast F_M\circ\ldots\circ (\sigma^{n-1})^\ast F_M)$ gives a functor from $\sigma$-modules over $(\CR,\CR^\sigma)$ to $\sigma^n$-modules over $(\CR,\CR^{\sigma^n})$. We will use this functor as a technical tool. To avoid confusion we will indicate in the notation, like $\Hom_\sigma(M,N)$ versus $\Hom_{\sigma^n}(M,N)$, whether we talk about $\sigma$-modules or $\sigma^n$-modules.

\begin{definition}\label{DefShortExactSeq}
A \emph{short exact sequence} of $\sigma$-modules over $\CR$ 
\[
0\to M'\to M\to M''\to0
\]
is a sequence of $\sigma$-modules such that the underlying sequence of $\CR$-modules is exact. Note that this implies that the submodule $M'\subset M$ is saturated.
The category of $\sigma$-modules with these short exact sequences is an \emph{exact category} in the sense of Quillen~\cite[\S2]{Quillen}.
\end{definition}

\begin{definition}
Let $\CR\subset\CR'$ be an inclusion of two $B$-algebras of the above types (over the same $B$) and let $M$ be a $\sigma$-module over $\CR'$. An $\CR$-submodule $N\subset M$ with $N\otimes_\CR\CR'=M$ is called an \emph{$\CR$-lattice} in $M$. We say that $N$ is \emph{$F$-stable} if $F_M$ maps $\sigma^\ast N$ into $\iota^\ast N$. Clearly $\fdot\otimes_\CR\CR'$ is a faithful $\BF_q\dbl z\dbr$-linear exact tensor functor from $\sigma$-modules over $\CR$ to $\sigma$-modules over $\CR'$.
% NOT FULL !!!!!!!!!!!!!!!!!!!!!!!!!!!!!!!!!!!!!!!!!!!
\end{definition}

For later reference we record a special case of Lemma 3.6.2 from \cite{Kedlaya}, which applies to more general situations than the one considered above. Specifically, let $\CR\to\CR'$ be an inclusion of principal ideal domains containing $\BF_q$ equipped with compatible morphisms $\sigma:\CR\to\CR^\sigma$ and $\sigma:\CR'\to(\CR')^\sigma$. Let $(M,F_M)$ be a $\sigma$-module over $\CR$ (definition as  in~\ref{DefSigmaModule}) and let $N'$ be a saturated $\sigma$-submodule of $M\otimes_\CR \CR'$. We say that \emph{$N'$ descends to $\CR$} if there exists a saturated $\sigma$-submodule $N$ of $M$ with $N'=N\otimes_\CR \CR'$.

\begin{lemma} \label{LemmaKe3.6.2}
\cite[Lemma 3.6.2]{Kedlaya}
Let $M$ be a $\sigma$-module over $\CR$ and let $N'$ be a saturated $\sigma$-submodule of $M\otimes_\CR \CR'$. Put $d:=\rk_{\CR'}N'$. Then $N'$ descends to $\CR$, if and only if the $\sigma$-submodule $\wedge^dN'\subset(\wedge^dM)\otimes_\CR \CR'$ descends to $\CR$.
\end{lemma}

\medskip

Let $C$ be a $B$-algebra which is itself an affinoid $L$-algebra, let $\CR_B$ be one of the above $B$-algebras, and let $\CR_C$ be the corresponding $C$-algebra with the same convergence behavior. Then there is a natural morphism $\CR_B\to\CR_C$. It gives rise to a \emph{pullback functor} $M\mapsto M\otimes_{\CR_B}\CR_C$ from $\sigma$-modules over $\CR_B$ to $\sigma$-modules over $\CR_C$.

Regarding the local freeness of $\sigma$-modules the following useful result is proved in \cite{Lubo}. The precise statement which we shall need is best formulated using Berkovich spaces. Recall from Appendix~\refAppBerkovich{} the definition of the Berkovich space $\CM(B)$ associated with the affinoid $L$-algebra $B$.

\begin{lemma} \label{Lemma6.4} \label{Lemma6.5}
Let $M$ be a locally free $B\langle \frac{z}{\zeta^r}\rangle$-module (respectively $B\langle\frac{z}{\zeta^{r}},\frac{\zeta^{r'}}{z}\rangle$-module) of finite rank. Then there exists a finite collection of affinoid $L$-algebras $B_i\supset B$ defining affinoid subdomains, such that the pullback of $M$ to $B_i\langle \frac{z}{\zeta^r}\rangle$ (respectively to $B_i\langle\frac{z}{\zeta^{r}},\frac{\zeta^{r'}}{z}\rangle$) is free, and such that $\{\CM(B_i)\}_i$ is an affinoid covering of $\CM(B)$ in the sense of Definition~\ref{DefBerkovichSpaces}.
\end{lemma}

\begin{proof}
Without the last assertion that $\{\CM(B_i)\}_i$ is an affinoid covering of $\CM(B)$ the statement for $B\langle \frac{z}{\zeta^r}\rangle$ was proved by L\"utkebohmert~\cite[\S1, Satz 1]{Lubo}. We next treat $B\langle\frac{z}{\zeta^{r}},\frac{\zeta^{r'}}{z}\rangle$.
By \cite[\S1, Satz 2]{Lubo} the restriction of $M$ to $\Spm B\langle\frac{z}{\zeta^{r'}},\frac{\zeta^{r'}}{z}\rangle$ is free locally on the base $B$. We can thus locally glue $M$ with a free sheaf on $\Spm B\langle\frac{z}{\zeta^{r'}}\rangle$ to a locally free sheaf of finite rank on $\Spm B\langle\frac{z}{\zeta^r}\rangle$. Since the latter is free locally on the base by \cite[\S1, Satz 1]{Lubo}, the same holds for its restriction to $\Spm B\langle\frac{z}{\zeta^r},\frac{\zeta^{r'}}{z}\rangle$ which is $M$. We thus obtain a  collection of affinoid $L$-algebras $B_i\supset B$ defining affinoid subdomains and trivializing $M$.

It remains to show that $\{\CM(B_i)\}_i$ is an affinoid covering of $\CM(B)$. For this we may have to enlarge the $\CM(B_i)$. This is achieved by the following well known lemma applied to $C=L\langle\frac{z}{\zeta^r}\rangle$ and $A=B_i$.
\end{proof}

\begin{lemma} \label{LemmaGerritzen}
(Compare \cite[Lemma, p.\ 128]{Lubo}, or \cite[\S2, Lemma 4]{Gerritzen}.)
Let $B$ and $C$ be affinoid $L$-algebras and let $M$ be a locally free $B\wh\otimes_L C$-module of finite rank. Let $A\supset B$ be an affinoid $L$-algebra defining an affinoid subdomain, such that the pullback of $M$ to $A\wh\otimes_L C$ is free. Then there exists an affinoid $L$-algebra $A'$ with $B\subset A'\subset A$ such that the pullback of $M$ to $A'\wh\otimes_L C$ is likewise free and $\CM(A)$ lies in the open interior of $\CM(A')$ in $\CM(B)$.
\end{lemma}

\begin{proof}
This lemma is well known but not well documented in the literature. Therefore we include a proof. Since $B\subset A$ is dense, the same is true for $M\subset M'':=M\otimes_{B\wh\otimes C}A\wh\otimes_L C$. Let $f_1'',\ldots,f_n''$ be a basis of $M''$ over $A\wh\otimes_L C$ and approximate $f_i''$ by $f_i\in M$. Write $f_j=\sum_id_{ij}f_i''$ with $d_{ij}\in A\wh\otimes_L C$. If each $f_i$ is close enough to $f_i''$ the matrix $(d_{ij})$ is close to $\Id_n$, hence invertible. Thus we may assume that also the $f_i$ form a basis of $M''$ over $A\wh\otimes_L C$. Let $N\subset M$ be the $B\wh\otimes_L C$-submodule generated by the $f_i$. By \cite[Proposition~9.5.2/4]{BGR} the quotient $M/N$ is supported on a Zariski-closed subset $Z\subset\CM(B\wh\otimes_L C)$ which is disjoint from $\CM(A\wh\otimes_L C)$. Outside $Z$ the map $(B\wh\otimes_L C)^{\oplus n}\to M$ sending the basis vectors to the $f_i$ is an isomorphism (by Nakayama's Lemma applied to its kernel). Since $\CM(B\wh\otimes_L C)$ is compact, the image $Y$ of $Z$ in $\CM(B)$ is closed. Therefore $\CM(B)\setminus Y\hookrightarrow\CM(B)$ is an open immersion and a closed morphism in terms of \cite[Example 1.5.3(ii)]{Berkovich2} whose image contains from $\CM(A)$. By Temkin~\cite[Theorem 5.1]{Temkin} we find an affinoid $L$-algebra $B\subset A'\subset A$ with $\CM(A')\subset\CM(B)\setminus Y$ such that $\CM(A)$ lies in the open interior of $\CM(A')$ inside $\CM(B)$. Then $\CM(A'\wh\otimes_L C)\subset\CM(B\wh\otimes_L C)\setminus Z$ and this proves the lemma.
\end{proof}

\medskip

\begin{example}\label{Ex1.2.8}
Let us give the following fundamental examples of $\sigma$-modules. We denote the $\sigma$-module $(\CR,F=\sigma)$ by $\CO(0)$. Next assume that $z$ is invertible in $\CR$. For every integer $d$ we define the $\sigma$-module $\CO(d)$ over $\CR$ as $(\CR,F=z^{-d}\cdot\sigma)$. For more examples let $d$ and $n$ be relatively prime integers with $n>0$. Consider the matrix
\begin{equation}\label{EqMatrixA_dn}
A_{d,n} := \left( \raisebox{6.2ex}{$
\xymatrix @C=0.3pc @R=0.3pc {
0 \ar@{.}[rrr]\ar@{.}[drdrdrdr] & & & 0 & z^{-d}\\
1\ar@{.}[drdrdr]  & & &  & 0 \ar@{.}[ddd]\\
0 \ar@{.}[dd]\ar@{.}[drdr] & & & &  \\
& & & & \\
0 \ar@{.}[rr] & & 0 & 1 & 0\\
}$}
\right)\es\in\es\GL_n(\CR).
\end{equation}
Let $\CF_{d,n}=(\CR^n,F=A_{d,n}\cdot\sigma)$. It is a $\sigma$-module of rank $n$ over $\CR$. As a special case if $n=1$ we obtain $\CF_{d,1}=\CO(d)$. We will say that a basis of $\CF_{d,n}$ over $\CR$ is a \emph{standard basis} if $F$ is of the form $A_{d,n}\cdot\sigma$ with respect to this basis.

Since we didn't want to overload the notation, the ring $\CR$ over which these $\sigma$-modules are living is not visible in the symbols. We will compensate for this by always specifying over which $\CR$ we consider the $\sigma$-modules $\CO(d)$ and $\CF_{d,n}$. 
\end{example}

\medskip

Let $\ol k$ be an algebraically closed field containing $\BF_q$. As mentioned after definition~\ref{DefSigmaModule} we may consider $\sigma$-modules over $\ol k\dpl z\dpr$. The following theorem proved by Laumon \cite[Theorem 2.4.5]{Laumon} is the analog of the Dieudonn\'e-Manin classification \cite{Manin} of isocrystals over an algebraically closed field. 

\begin{theorem}\label{Thm2}
Every $\sigma$-module over $\ol k\dpl z\dpr$ is isomorphic to a direct sum $\bigoplus_{i=1}^\ell \CF_{d_i,n_i}$ where the $d_i$ and $n_i$ are integers with $n_i>0$ and $(d_i,n_i)=1$. The pairs $(d_i,n_i)$ are uniquely determined up to permutation. 
\end{theorem}

\begin{corollary}\label{Cor0.6b}
Every $\sigma$-module of rank $1$ over $\ol k\dpl z\dpr$ is isomorphic to $\CO(d)$ for a uniquely determined integer $d$.
\qed
\end{corollary}

\begin{proposition}\label{Prop0.7b}
The $\sigma$-modules $\CF_{d,n}$ over $\ol k\dpl z\dpr$ satisfy
$\DS\Hom_\sigma(\CF_{d,n},\CF_{d',n'})=0$ if $\frac{d}{n}\neq\frac{d'}{n'}$.
\end{proposition}

\begin{proof}
Let $f\in\Hom_\sigma(\CF_{d,n},\CF_{d',n'})$ with $f\neq0$. Then $f$ corresponds to a non-zero
matrix $A\in M_{n'\times n}\bigl(\ol k\dpl z\dpr\bigr)$ satisfying $A_{d',n'}\cdot A^\sigma=A\cdot A_{d,n}$. In particular
\[
z^{-d'n}A^{\sigma^{nn'}}\es=\es A_{d',n'}\cdot \ldots\cdot A_{d',n'}^{\sigma^{nn'-1}}\cdot A^{\sigma^{nn'}}\es=\es A\cdot A_{d,n}\cdot \ldots\cdot A_{d,n}^{\sigma^{nn'-1}}\es=\es z^{-dn'}A\,.
\]
Since $A$ has finite principal part this is only possible if $d'n=dn'$, that is if $\frac{d}{n}=\frac{d'}{n'}$.
\end{proof}

%%%%%%%%%%%%%%%%%%%%%%%%%%%%%%%%%%%%%%%%%%%%%%%%%%%%%%%%%%%%%%%%%%%%%%
%
%   $F$-Invariants
%
%%%%%%%%%%%%%%%%%%%%%%%%%%%%%%%%%%%%%%%%%%%%%%%%%%%%%%%%%%%%%%%%%%%%%%

\subsection{\texorpdfstring{$F$}{F}-Invariants} \label{SectFInvariants}

Let again $\CR$ be a $B$-algebra of the above type.

\begin{definition} \label{DefH^0}
Let $M$ be a $\sigma$-module over $\CR$. We define the \emph{set of $F$-invariants} of $M$ as
\[
M^F(B)\es:=\es\{\;x\in M:\es F_M(\sigma^\ast x)=x\;\}\,,
\]
where we denote by $\sigma^\ast x:=x\otimes1$ the image of $x$ under $M\to \sigma^\ast M=M\otimes_{\CR,\sigma}\CR^\sigma$.
\end{definition}

In fact $M^F(B)$ is a module over the ring $\CO(0)^F(B)\;=\;\{y\in\CR: y^\sigma=y\}$.

\begin{proposition}\label{PropH^0OfR}
Assume that $B$ has no non-trivial idempotents. Then
\begin{enumerate}
\item 
$\CO(0)^F(B)\es=\es\BF_q\dbl z\dbr$\quad if $z\in\CR, z^{-1}\notin\CR$.
\item 
$\CO(0)^F(B)\es=\es\BF_q\dpl z\dpr$\quad if $z,z^{-1}\in\CR$.
\item 
$\CO(0)^F(B)\es=\es\BF_q[z^{-1}]$\quad if $z^{-1}\in\CR, z\notin\CR$
\end{enumerate}
\end{proposition}

\begin{proof}
By definition $\CO(0)^F(B)$ consists of all Laurent series $y=\sum_{i=-\infty}^\infty b_i z^i\in\CR$ with $b_i^q=b_i$, that is, $b_i\in\BF_q$ since we have assumed that $B$ has no non-trivial idempotents. In particular $|b_i|=1$ if $b_i\neq0$. Observing the convergence condition for $\CR$ which are imposed on $y$ the claim follows.
\end{proof}

Again it would be more accurate to formulate this definition in terms of sheaves. One proceeds as follows. Denote by $\pi_\ast$ (respectively by $\wt\pi_\ast$) the functor from $\CR$-modules (respectively $\CR^\sigma$-modules) to $B$-modules coming from the inclusion $B\subset\CR$ (respectively $B\subset\CR^\sigma$). Denote by $\sigma_\ast$ the functor from $\CR^\sigma$-modules to $\CR$-modules coming from $\sigma:\CR\to\CR^\sigma$. It is right adjoint to $\sigma^\ast$. From the commutativity of the diagram
\[
\xymatrix @M+0.3pc {B \ar@{^{ (}->}[r] \ar[d]_{\sigma_B} & \CR\ar[d]^\sigma \\
B\ar@{^{ (}->}[r] & \CR^\sigma}
\]
with $\sigma_B:b\mapsto b^q$ we obtain $\pi_\ast\sigma_\ast=(\sigma_B)_\ast\wt\pi_\ast$. By adjunction the morphism $F_M:\sigma^\ast M\isoto\iota^\ast M$ induces a morphism $\sigma_\ast F_M:M\to\sigma_\ast\iota^\ast M$ of $\CR^\sigma$-modules and further a morphism $\pi_\ast\sigma_\ast F_M:\pi_\ast M\to\pi_\ast\sigma_\ast\iota^\ast M$ of $B$-modules. On the other hand the inclusion $\iota:\CR\to\CR^\sigma$ yields a morphism $i:\pi_\ast M\to\wt\pi_\ast\iota^\ast M$ of $B$-modules.

Now the $B$-module $\pi_\ast M$ defines a sheaf of abelian groups on the \'etale site $X_\et$ of the rigid analytic space $X=\Spm B$ (see Appendix~\refAppEtaleSheaves). Namely for any \'etale morphism $f:Y\to X$ we set $W(\pi_\ast M)(Y)\;=\;\Gamma(Y,f^\ast\pi_\ast M)$. In the same manner we obtain $W(\wt\pi_\ast\iota^\ast M)$ and $W({\sigma_B}_\ast\wt\pi_\ast\iota^\ast M)$. Note that the last two sheaves are equal as sheaves of abelian groups (since $\sigma_B$ only changes on them the scalar multiplication by $B$). Therefore we can consider the morphism of sheaves of abelian groups on $X_\et$
\[
W(\pi_\ast\sigma_\ast F_M)-W(i):\es W(\pi_\ast M)\to W(\wt\pi_\ast\iota^\ast M)\,.
\]

\begin{definition}\label{DefH^0b}
We define the \'etale sheaves $M^F$ and $M_F$ on $X_\et$ as the kernel and the cokernel of $W(\pi_\ast\sigma_\ast F_M)-W(i)$ in the category of \'etale sheaves of abelian groups on $X_\et$. They are sheaves of $\BF_q\dbl z\dbr$-modules if $z\in\CR$ (and even $\BF_q\dpl z\dpr$-modules if $z$ is invertible in $\CR$).
\end{definition}

Clearly $M^F(\Spm B)$ coincides with $M^F(B)$ from Definition~\ref{DefH^0} since as usual the ``presheaf kernel'' already is a sheaf. In the special case where $B=\ol L$ we find for the cokernel
\[
M_F(\ol L)\es=\es\iota^\ast M/\{\,F_M(\sigma^\ast x)-x:\;x\in M\,\}\,.
\]
Thus for any $\sigma$-module $M$ over $\ol L\overcon$ the $\BF_q\dpl z\dpr$-vector spaces $M^F(\ol L)$ and $M_F(\ol L)$ are precisely the ones which were denoted $\Koh^0(M)$ and $\Koh^1(M)$ in \cite{HP}.

By the snake lemma every short exact sequence $0\to M'\to M\to M''\to0$ of $\sigma$-modules over $\CR$ gives rise to a long exact sequence
\[
0\es\longto\es (M')^F \es\longto\es M^F \es\longto\es (M'')^F \es\longto\es M'_F \es\longto\es M_F \es\longto\es M''_F \es\longto\es 0\,.
\]
Indeed for $\CR=B\dbl z\dbr$, $B\dbl z\dbr[z^{-1}]$, $B\langle\frac{z}{\zeta^r}\rangle$, $B\con[r]$, $B\langle\frac{\zeta^{r'}}{z}\rangle$, $B\langle\frac{z}{\zeta^r},\frac{\zeta^{r'}}{z}\rangle$ the functors $\pi_\ast$ and $\wt\pi_\ast i^\ast$ are exact by Definition~\ref{DefSigmaModule}. Also for $\CR=B\ancon[r]$, $B\{z,\frac{\zeta^{r'}}{z}\rangle$, $B\overcon$, $B\{z\}$ the functors $\pi_\ast$ and $\wt\pi_\ast i^\ast$ are exact by \cite[Theorem A]{KiehlAB} because the rigid analytic spaces associated with these $B$-algebras are relative Stein spaces over $\Spm B$.

Let $M$ and $N$ be $\sigma$-modules over the $B$-algebra $\CR$. Then $\Hom_\sigma(M,N)\;=\;\bigl(\CHom(N,M)\bigr)^F(B)$ by definition. We further define 
\[
\Ext_\sigma(N,M)\es:=\es\bigl(\CHom(N,M)\bigr)_F(B)\,.
\]
The following proposition is proved by standard homological arguments as in \cite[Proposition 2.4]{HP}.

\begin{proposition}\label{PropExt}
Let $\CR$ be one of the above $B$-algebras for $B=\ol L$ and let $M$ and $N$ be $\sigma$-modules over $\CR$. Then the group $\Ext_\sigma(N,M)$ classifies extensions (short exact sequences) of $\sigma$-modules 
\[
0\to M\to E\to N\to 0
\]
up to isomorphisms of short exact sequences that are the identity on $M$ and $N$.
\qed
\end{proposition}

Let us remark that the proposition (and its proof) also hold for general $B$ but we will not need this here.

\medskip

\begin{proposition} \label{Prop13}
Let $M$ be a $\sigma$-module over $B\con$ which contains an $F$-stable $B\langle\frac{z}{\zeta}\rangle$-lattice $N$. Then the natural morphism
\begin{enumerate}
\item
$M^F\isoto\bigl(M\otimes B\ancon\bigr)^F$ is an isomorphism and
\item
$M_F\into\bigl(M\otimes B\ancon\bigr)_F$ is injective.
\end{enumerate}
If moreover $F:\sigma^\ast N\to N$ is an isomorphism then the natural morphism
\begin{enumerate}
\setcounter{enumi}{2}
\item 
$M^F\isoto\bigl(M\otimes B\dbl z\dbr[z^{-1}]\bigr)^F$ is an isomorphism.
\item 
The sheaves $\bigl(N\otimes B\dbl z\dbr/(z^m)\bigr)^F$ form a local system (\ref{DefLocalSystem}) of $\BF_q\dbl z\dbr$-lattices on the rigid analytic space $\Spm B$. The sheaf $\bigl(N\otimes B\dbl z\dbr\bigr)^F$ is the projective limit of $\Bigl(\bigl(N\otimes B\dbl z\dbr/(z^m)\bigr)^F\Bigr)_{m\in\BN_0}$.
\end{enumerate}
\end{proposition}

\begin{proof}
Working locally on $B$ we may by Lemma~\ref{Lemma6.5} choose a basis of $N$ and write $F$ with respect to this basis as a matrix $A\in M_n\bigl(B\langle\frac{z}{\zeta^q}\rangle\bigr)$.

To prove (a) note that the injectivity is obvious and the surjectivity follows from Lemma~\ref{Lemma12}.

To prove (b) it suffices to show that for any \'etale affinoid $B$-algebra $C$ and for any vector $v\in M\otimes C\con$ which satisfies $v=w-F(\sigma^\ast w)$ for some $w\in M\otimes C\ancon$, already $w\in M\otimes C\con$.
Again this follows from Lemma~\ref{Lemma12}.

In (c) the injectivity is obvious. To prove the surjectivity let $x\in B\dbl z\dbr^n$ satisfy $x=A\,x^\sigma$. Since $A\in\GL_n\bigl(B\langle\frac{z}{\zeta^q}\rangle\bigr)$ by assumption,
$x\in B\langle\frac{z}{\zeta}\rangle^n$ by Lemma~\ref{Lemma3b}.

To prove (d) write $A=\sum_{j\geq0}A_j z^j$. Then for any \'etale affinoid $B$-algebra $C$ we obtain
\begin{eqnarray*}
\bigl(N\otimes B\dbl z\dbr\bigr)^F(C) & = & \bigl\{\,\sum_{i=0}^\infty x_i z^i:\es x_i\in C^{\,n},\, x_i=\sum_{j=0}^i A_j x_{i-j}^\sigma\text{ for all }i\geq0\,\bigr\}\quad\text{and}\\
\bigl(N\otimes B\dbl z\dbr/(z^m)\bigr)^F(C) & = & \bigl\{\,\sum_{i=0}^{m-1} x_i z^i:\es x_i\in C^{\,n},\, x_i=\sum_{j=0}^i A_j x_{i-j}^\sigma\text{ for all }i=0,\ldots,m-1\,\bigr\}\,.
\end{eqnarray*}
Thus $\bigl(N\otimes B\dbl z\dbr\bigr)^F\;=\;\invlim[m]\bigl(N\otimes B\dbl z\dbr/(z^m)\bigr)^F$. The remaining assertions are a consequence of the following lemma.
\end{proof}

\begin{lemma} \label{Lemma6.10}
Assume that $B$ is connected. Let $m\in \BN_0$ and let $N$ be a locally free $B\dbl z\dbr/(z^m)$-module of rank $n$, equipped with an isomorphism $F:\sigma^\ast N\to N$.
Then there exists a finite \'etale connected $B$-algebra $B_m$ such that $N^F(B_m)$ is a free $\BF_q\dbl z\dbr/(z^m)$-module of rank $n$. Moreover there is a natural isomorphism
\[
N^F(B_m)\otimes_{\BF_q\dbl z\dbr/(z^m)}\;B_m\dbl z\dbr/(z^m) \isoto N\otimes_{B\dbl z\dbr/(z^m)}\;B_m\dbl z\dbr/(z^m)\,.
\]
\end{lemma}

\begin{proof}
This is proved in \cite[\S2]{BH}. It is the analog of the fact that a finite \'etale group scheme
is trivialized by a finite \'etale covering space. See \cite[Theorem 2.5]{BH} for details.
\end{proof}

Although the sheaves $M^F$ and $M_F$ will serve us well in this chapter, the above notion of $F$-invariants is not the most useful for $\sigma$-modules over $B\con$ or $B\dbl z\dbr[z^{-1}]$. The reason for this is that unless $B=\ol L$, the sheaf $M^F$ will in general be the zero sheaf as the following example shows.

Let $A=\sum_{j=0}^\infty A_jz^j\in\GL_n\bigl(B\dbl z\dbr[z^{-1}]\bigr)$ and $(M,F_M)=\bigl(B\dbl z\dbr[z^{-1}]^n,A\cdot\sigma\bigr)$. Then for any affinoid \'etale $B$-algebra $C$ 
\[
M^F(C)\es=\es \bigl\{\,\sum_{i=0}^\infty x_iz^i:\es x_i\in C^n\text{ with}\es x_i\;=\;\sum_{j=0}^i A_j x_{i-j}^\sigma\text{ for all }i\,\bigr\}\;\otimes_{\BF_q\dbl z\dbr}\;\BF_q\dpl z\dpr\,.
\]
Even if $A\in\GL_n\bigl(B\langle\frac{z}{\zeta^q}\rangle\bigr)$ and a matrix $\Phi\in\GL_n(C)$ exists with $\Phi=A_0\Phi^\sigma$ we have to solve infinitely many Artin-Schreier equations
\[
(\Phi^{-1}x_i)-(\Phi^{-1}x_i)^\sigma\es=\es\sum_{j=1}^i(\Phi^{-1}A_j\Phi^\sigma)(\Phi^{-1}x_{i-j})^\sigma\,.
\]
Unless $C=\ol L$ or $|\Phi^{-1}x_i|<1$ for $i\gg0$ we cannot hope to find all solutions $x_i\in C^n$. But the expected convergence condition for $\sum x_iz^i$ is $|x_i\zeta^i|\to0$, so we will find $|\Phi^{-1}x_i|<1$ for $i\gg0$ only in rare cases.

It seems more meaningful to work with the local system
 $T_zN:=\left(\bigr(N\otimes B\dbl z\dbr/(z^m)\bigr)^F\right)_{m\in\BN_0}$ of $\BF_q\dbl z\dbr$-lattices on $\Spm B$ from Proposition~\ref{Prop13} (d)
 instead. 
We call $T_zN$ the \emph{Tate module} of $N$. These local systems give rise to Galois representations exactly like the Tate module of a $p$-divisible group or abelian variety. See Section~\ref{SectLocalShtuka} and Appendix~\refAppFundamentalGroups{} for more on this subject. We want to end this discussion by noting the following analog of Katz's~\cite[Proposition~4.1.1]{Katz73}.

\begin{proposition}\label{PropTateModuleExact}
The functor $M\mapsto T_zM$ from the category of $\sigma$-modules over $B\dbl z\dbr$ to the category of local systems of $\BF_q\dbl z\dbr$-lattices on $\Spm B$ is an $\BF_q\dbl z\dbr$-linear exact tensor equivalence.
\end{proposition}

\begin{proof}
Linearity and compatibility with tensor products are clear. The full faithfulness follows from Lemma~\ref{Lemma6.10}. The inverse functor is constructed as follows. Let $\CF$ be a local system of $\BF_q\dbl z\dbr$-lattices on $\Spm B$. For every $m$, the sheaf $\CF/z^m\CF$ of $\BF_q\dbl z\dbr/(z^m)$-modules is trivialized by a finite \'etale Galois covering $\Spm B_m\to\Spm B$; compare Proposition~\ref{Prop2.13}. The coherent sheaf $\CF/z^m\CF\otimes_{\BF_q\dbl z\dbr/(z^m)} \CO_{\Spm B_m}$ is equipped with the Frobenius $F=\id\otimes\sigma$ and the Galois action $\gamma(f\otimes b)=\gamma^{-1}(f)\otimes\gamma^\ast(b)$ for $\gamma\in\Gal(B_m/B)$. It descends to a locally free $B\dbl z\dbr/(z^m)$-module $N_m$ with an isomorphism $F_m:\sigma^\ast N_m\isoto N_m$. Clearly the $(N_m,F_m)$ are compatible and their projective limit is the desired $\sigma$-module over $B\dbl z\dbr$.

To prove exactness let $0\to M'\to M\to M''\to0$ be an exact sequence of $\sigma$-modules over $B\dbl z\dbr$. We have to show that the associated sequence of sheaves on $(\Spm B)_\et$
\[
0\es\longto\es(M'\mod z^m)^F\es\longto\es(M\mod z^m)^F\es\longto\es(M''\mod z^m)^F\es\longto\es0
\]
is exact for all $m\in\BN_0$. Since this property is local on $B$ we may choose bases of $M'$, $M$, and $M''$ such that
\[
F_{M'}=A'\cdot\sigma\,,\es F_{M''}=A''\cdot\sigma\quad\text{and}\quad F_M=\left(\begin{array}{cc} A' & \wt A \\ 0 & A''\end{array}\right)\cdot\sigma
\]
for $A'\in\GL_{n'}\bigl(B\dbl z\dbr\bigr)$, $A''\in\GL_{n''}\bigl(B\dbl z\dbr\bigr)$ and $\wt A\in M_{n'\times n''}\bigl(B\dbl z\dbr\bigr)$. In this presentation only the exactness on the right is non-obvious. So let $C$ be an affinoid \'etale $B$-algebra and let $\sum_{i=0}^{m-1}y_iz^i\;\in(M''\mod z^m)^F(C)$. It suffices to find a finite \'etale $C$-algebra $\wt C$ and $x_i\in\wt C^{n'}$ with $\sum_{i=0}^{m-1}\binom{x_i}{y_i}z^i\in (M\mod z^m)^F(\wt C)$. This amounts to solving the finite \'etale equations 
\[
x_i-A'_0x_i^\sigma\es=\es\sum_{j=1}^i A'_jx_{i-j}^\sigma \;+\;\sum_{j=0}^i\wt A_j y_{i-j}^\sigma
\]
for $i=0,\ldots,m-1$ and clearly can be accomplished.
\end{proof}

%%%%%%%%%%%%%%%%%%%%%%%%%%%%%%%%%%%%%%%%%%%%%%%%%%%%%%%%%%%%%%%%%%%%%%
%
%    Results on $\sigma$-Bundles and Supplements
%
%%%%%%%%%%%%%%%%%%%%%%%%%%%%%%%%%%%%%%%%%%%%%%%%%%%%%%%%%%%%%%%%%%%%%%

\subsection{Results on \texorpdfstring{$\sigma$}{sigma}-Bundles and Supplements} \label{SectResultsAndSuppl}

In this section we review the main results from \cite{HP} where the term \emph{$\sigma$-bundle} was used for what here is called a $\sigma$-module over $\ol L\overcon$. Let $r$ be a positive rational number.

\begin{proposition}\label{Prop0.1}
\begin{enumerate}
\item 
The functor from the category of $\sigma$-modules over $B\ancon[r]$ to the category of $\sigma$-modules over $B\langle\frac{z}{\zeta^r},\frac{\zeta^{qr}}{z}\rangle$, $M\mapsto M\otimes B\langle\frac{z}{\zeta^r},\frac{\zeta^{qr}}{z}\rangle$ is an exact equivalence of rigid tensor categories.
\item 
The functor from the category of $\sigma$-modules over $\ol L\overcon$ to the category of $\sigma$-modules over $\ol L\langle\frac{z}{\zeta^r},\frac{\zeta^{qr}}{z}\rangle$, $M\mapsto M\otimes \ol L\langle\frac{z}{\zeta^r},\frac{\zeta^{qr}}{z}\rangle$ is an exact equivalence of rigid tensor categories.
\end{enumerate}
\end{proposition}

\begin{proof}
We first prove (a). To prove essential surjectivity, let $M$ be a $\sigma$-module over $B\langle\frac{z}{\zeta^r},\frac{\zeta^{qr}}{z}\rangle$. The isomorphism $F:\sigma^\ast M\otimes B\langle\frac{z}{\zeta^{qr}},\frac{\zeta^{qr}}{z}\rangle\isoto M\otimes B\langle\frac{z}{\zeta^{qr}},\frac{\zeta^{qr}}{z}\rangle$ allows us to glue the $B\langle\frac{z}{\zeta^r},\frac{\zeta^{qr}}{z}\rangle$-module $M$ with the $B\langle\frac{z}{\zeta^{qr}},\frac{\zeta^{q^2r}}{z}\rangle$-module $M\otimes_{B\langle\frac{z}{\zeta^r},\frac{\zeta^{qr}}{z}\rangle,\,\sigma}B\langle\frac{z}{\zeta^{qr}},\frac{\zeta^{q^2r}}{z}\rangle$ to a $\sigma$-module $M_1$ over $B\langle\frac{z}{\zeta^r},\frac{\zeta^{q^2r}}{z}\rangle$. We continue in this way and obtain a $\sigma$-module over $B\ancon[r]$.

To prove full faithfulness, let $M$ and $N$ be $\sigma$-modules over $B\ancon[r]$. Let $f:M\to N$ be a morphism of $\sigma$-modules, that is $f\circ F_M=F_N\circ \sigma^\ast f$. Therefore the restriction of $f$ to $B\langle\frac{z}{\zeta^{qr}},\frac{\zeta^{q^2r}}{z}\rangle$ is uniquely determined by the restriction of $\sigma^\ast f$ to $B\langle\frac{z}{\zeta^{qr}},\frac{\zeta^{q^2r}}{z}\rangle$, hence by the restriction of $f$ to $B\langle\frac{z}{\zeta^r},\frac{\zeta^{qr}}{z}\rangle$. An iteration of this argument yields full faithfulness. Clearly the functor is also compatible with tensor products and internal Hom's, that is, it is a tensor functor, and maps short exact sequences to short exact sequences.

To prove (b) we also must go in the direction $|z|\to1$. Note that $\sigma:\ol L\overcon\to\ol L\overcon$ is an isomorphism. So we can consider $(\sigma^{-1})^\ast M$ and 
\[
\TS(\sigma^{-1})^\ast F:\es M\otimes\ol L\langle\frac{z}{\zeta^r},\frac{\zeta^r}{z}\rangle\es\isoto\es (\sigma^{-1})^\ast M\otimes\ol L\langle\frac{z}{\zeta^r},\frac{\zeta^r}{z}\rangle\,.
\]
This allows us to glue the $\ol L\langle\frac{z}{\zeta^r},\frac{\zeta^{qr}}{z}\rangle$-module $M$ with the $\ol L\langle\frac{z}{\zeta^{r/q}},\frac{\zeta^r}{z}\rangle$-module $(\sigma^{-1})^\ast M$ to a $\sigma$-module $M'_1$ over $\ol L\langle\frac{z}{\zeta^{r/q}},\frac{\zeta^{qr}}{z}\rangle$. With this building stone the proof continues as the proof of (a) above.
\end{proof}

In view of (b) the results from \cite{HP} are as follows.

\begin{theorem} \label{Thm1a} 
\cite[Theorem 11.1 and Corollary 11.8]{HP} Every $\sigma$-module over $\ol L\ancon$ is isomorphic to a direct sum $\bigoplus_{i=1}^\ell \CF_{d_i,n_i}$ where the $d_i$ and $n_i$ are integers with $n_i>0$ and $(d_i,n_i)=1$ (and $\CF_{d_i,n_i}$ was defined in Example~\ref{Ex1.2.8}). The pairs $(d_i,n_i)$ are uniquely determined up to permutation.
\end{theorem}

Kedlaya pointed out the parallel between Theorems~\ref{Thm2} and \ref{Thm1a} in the case of mixed characteristic by giving a simultaneous proof for both; \cite[Theorem 4.5.7]{Kedlaya}. The same could be done in our situation. Note however that the proof of Theorem~\ref{Thm2} obtained in this way is much more complicated than the standard proof (\cite{Laumon}; or in mixed characteristic \cite{Manin}).

\begin{corollary}\label{Cor0.6}
Every $\sigma$-module of rank $1$ over $\ol L\ancon$ is isomorphic to $\CO(d)$ for a uniquely determined integer $d$.
\qed
\end{corollary}

\begin{proposition}\label{Prop0.5}
\cite[Propositions 3.1 and 5.1]{HP} For every $\alpha\in\ol L$ there exists an $f_\alpha\in\ol L\overcon$ satisfying $f_\alpha=z^{-1}\cdot(f_\alpha)^\sigma$ with simple zeroes at $\alpha^{q^\nu}$ for $\nu\in\BZ$ and no other zeroes. For the $\sigma$-module $\CO(d)$ over $\ol L\overcon$ we have
\[
\CO(d)^F(\ol L)\es=\es\left\{\begin{array}{ll}
(0)& \text{if }d<0\,, \\[1mm]
\BF_q\dpl z\dpr & \text{if } d=0\,, \\
\DS\Bigl\{\,\sum_{\nu\in\BZ}z^{-d\nu}\sum_{j=0}^{d-1}z^j u_j^{q^\nu}:\es u_j\in\ol L, \,|u_j|<1\,\Bigr\} \es= \\[6mm]
=\es\BF_q\dpl z\dpr\cdot\bigl\{\,\prod_{i=1}^d f_{\alpha_i}:\es\alpha_i\in\ol L\,,\,|\zeta^q|<|\alpha_i|\leq|\zeta|\,\bigr\} & \text{if }d>0\,.
\end{array}\right.
\]
\end{proposition}

\begin{proof}
The assertion on $\CO(d)^F(\ol L)$ is demonstrated in the proof of \cite[Theorem 5.4]{HP}.
\end{proof}

\begin{proposition}\label{Prop0.7}
\cite[Propositions 8.5, 8.6, and 8.8]{HP} The $\sigma$-modules $\CF_{d,n}$ over $\ol L\ancon$ satisfy
\begin{enumerate}
\item 
$\DS\dim_{\BF_q\dbl z\dbr}\Hom_\sigma(\CF_{d,n},\CF_{d',n'})\es=\es\left\{\begin{array}{ll}
\infty & \text{if }d/n < d'/n'\,, \\[1mm]
n^2 & \text{if }d/n = d'/n'\,, \\[1mm]
0 & \text{if }d/n > d'/n'\,.
\end{array}\right.$
\item 
$\DS\dim_{\BF_q\dbl z\dbr}\Ext_\sigma(\CF_{d,n},\CF_{d',n'})\es=\es\left\{\begin{array}{ll}
0 & \text{if }d/n \leq d'/n'\,, \\[1mm]
\infty & \text{if }d/n > d'/n'\,.
\end{array}\right.$
\item
$\End_\sigma(\CF_{d,n})$ is the central division algebra over $\BF_q\dpl z\dpr$ of dimension $n^2$ and Hasse invariant $-\frac{d}{n}\mod\BZ$. 
\end{enumerate}
\end{proposition}

%%%%%%%%%%%%%%%%%%%%%%%%%%%%%%%%%%%%%%%%%%%%%%%%%%%%%%%%%%%%%%%%%%%%%%
%
%    Slopes and Polygons
%
%%%%%%%%%%%%%%%%%%%%%%%%%%%%%%%%%%%%%%%%%%%%%%%%%%%%%%%%%%%%%%%%%%%%%%

\subsection{Slopes and Polygons} \label{SectSlopesAndPolygons}

As remarked above we may also consider $\sigma$-modules over $k\dpl z\dpr$ for arbitrary fields $k\supset\BF_q$. For such a field we let $\ol k\supset k$ be some algebraically closed extension. What follows does not depend on the choice of $\ol k$.

\begin{definition} \label{DefDegree}
Let $M$ be a $\sigma$-module of rank $1$ over $L\ancon$ (or over $k\dpl z\dpr$). We define the \emph{degree} $\deg M$ of $M$ as the unique integer $d$ from Corollary~\ref{Cor0.6} (respectively Corollary~\ref{Cor0.6b}) such that $M\otimes \ol L\ancon\cong\CO(d)$ over $\ol L\ancon$ (respectively $M\otimes\ol k\dpl z\dpr\cong\CO(d)$ over $\ol k\dpl z\dpr$). For a $\sigma$-module $M$ of rank $n$ over $L\ancon$ or over $k\dpl z\dpr$ we define $\deg M:=\deg\wedge^n M$. We define the \emph{slope} of $M$ as
\[
\lambda(M)\es:=\es-\,\frac{\deg M}{\rk M}\,.
\]
\end{definition}

For example the $\sigma$-modules $\CF_{d,n}$ over $L\ancon$ (and also over $k\dpl z\dpr$) from Example~\ref{Ex1.2.8} have degree $\deg\CF_{d,n}=d$ and slope $\lambda(\CF_{d,n})=-\frac{d}{n}$. The degree is additive in short exact sequences (see \cite[Proposition 6.1]{HP}). 

Note that compared to \cite{HP} the definitions of degree of a $\sigma$-module $M$ over $\ol L\ancon$ coincide whereas the slope is the negative of what was called the \emph{weight} of $M$ in \cite{HP}. Compared to \cite{Kedlaya} the definitions of slope are the same whereas the definitions of degree differ by their sign. This difference lies in the fact that the terminology of \cite{HP} is ``geometric'', expressing the analogy with the stability of vector bundles, and the terminology of \cite{Kedlaya} is ``arithmetic'', owing to the fact that the Dieudonn\'e module of a $p$-divisible group has slopes $0$ and $1$. We have chosen our maybe unnatural definitions here because on the one hand we wanted to be consistent with \cite{HP} through the term ``degree'', and on the other hand our applications through the notion of ``slope'' are of arithmetic nature.

\begin{definition} \label{DefStable}
A $\sigma$-module $M$ over $L\ancon$ or over $k\dpl z\dpr$ is called \emph{semistable} (respectively \emph{stable}) if $\lambda(N)\geq\lambda(M)$ (respectively $\lambda(N)>\lambda(M)$) for any non-zero $\sigma$-submodule $N$ of $M$. 
\end{definition}

Note that a direct sum of semistable $\sigma$-modules of the same slope is again semistable.
Using the additivity of degree and rank in short exact sequences one shows as usual that $M$ is \mbox{(semi-)}stable if and only if for any quotient $\sigma$-module $N$ of $M$ the inequality $\lambda(N)<\lambda(M)$ (respectively $\lambda(N)\leq\lambda(M)$) holds.

In order to define the Harder-Narasimhan polygon of a $\sigma$-module we first recall the following

\begin{definition} \label{DefSlopesAndPolygons}
To a multiset $S$ of $n$ real numbers one defines the \emph{Newton polygon} of $S$ to be the graph of the 
piecewise linear function on $[0,n]$ sending $0$ to $0$, whose slope on $[i-1,i]$ is the $i$-th smallest element of $S$.
One calls the point on the graph corresponding to the image of $n$ the \emph{endpoint} of the polygon.
Conversely for every such graph one defines its \emph{slope multiset} as the slopes of the piecewise linear function
on $[i-1,i]$ for $i=1,\ldots,n$.

One says that one Newton polygon $P$ \emph{lies above} another Newton polygon $P'$ if both have the same endpoint 
and if no vertex of $P$ lies below the polygon $P'$. In this case one writes $P\geq P'$.

Given Newton polygons $P_1\ldots,P_m$, one defines their \emph{sum} $P_1+\ldots+P_m$
as the Newton polygon whose slope multiset is the union of the slope multisets of $P_1,\ldots,P_m$.
\end{definition}

\begin{definition} \label{DefHNFiltration}
Let $M$ be a $\sigma$-module over $L\ancon$ or over $k\dpl z\dpr$. A \emph{semistable filtration} of $M$ is a filtration $0=M_0\subset M_1\subset\ldots\subset M_\ell=M$ of $M$ by saturated $\sigma$-submodules,
such that the successive quotients $M_i/M_{i-1}$ are semistable of some slope $s_i$. 
A \emph{Harder-Narasimhan filtration} (short \emph{HN-filtration}) of $M$ is a semistable filtration with $s_1<\ldots< s_\ell$.
\end{definition}

\begin{lemma}\label{Lemma1.5.4'}(Compare \cite[Proposition 4.2.2]{Andre08}.)
A HN-filtration is unique if it exists.
\end{lemma}

\begin{proof}
Let $0=M_0\subset M_1\subset\ldots\subset M_\ell=M$ and $0=M'_0\subset M'_1\subset\ldots\subset M'_{\ell'}=M$ be HN-filtrations. Assume that $\lambda(M_1)\le\lambda(M'_1)$ and let $1\le i\le \ell'$ be the minimal index for which $M_1\subset M'_i$. Then $M_1\not\subset M'_{i-1}$ and $M_1/(M_1\cap M'_{i-1})$ is both a non-zero quotient $\sigma$-module of $M_1$ and a $\sigma$-submodule of $M'_i/M'_{i-1}$. By semistability its slope $\lambda$ satisfies $\lambda(M_1)\ge\lambda\ge\lambda(M'_i/M'_{i-1})\ge \lambda(M'_1)$. Hence $\lambda(M_1)=\lambda(M'_1)$ and $i=1$, that is, $M_1\subset M'_1$. Now conversely $\lambda(M'_1)=\lambda(M_1)$ implies $M'_1\subset M_1$. Thus $M_1=M'_1$ and the lemma follows by induction.
\end{proof}

\begin{definition} \label{DefHNPolygon}
Let $M$ be a $\sigma$-module over $L\ancon$ or over $k\dpl z\dpr$. With any semistable filtration $0=M_0\subset M_1\subset\ldots\subset M_\ell=M$ of $M$
we associate the multiset consisting of the slope $\lambda(M_i/M_{i-1})$ with multiplicity $\rk M_i/M_{i-1}$ for all $i=1,\ldots,\ell$. 
We call this multiset the \emph{slope multiset} of the filtration, and we call the corresponding Newton polygon the \emph{slope polygon} of the filtration.
If $M$ has an HN-filtration, we call its slope multiset the \emph{Harder-Narasimhan slope multiset} (short \emph{HN-slope multiset}) of $M$, 
and its slope polygon the \emph{Harder-Narasimhan polygon} (short \emph{HN-polygon}) of $M$.

If $M$ is a $\sigma$-module over $L\con$ we call the HN-polygon of $M\otimes L\dpl z\dpr$ the \emph{generic HN-polygon} of $M$, and we call the HN-polygon of $M\otimes L\ancon$ the \emph{special HN-polygon} of $M$. The same naming is chosen for HN-filtrations and HN-slopes.
\end{definition}

The latter names were coined and motivated by Kedlaya~\cite[\S7.3]{Kedlaya} in mixed characteristic and adopted here to emphasize the analogy between equal and mixed characteristic.

The ultimate goal of this chapter, the Slope Filtration Theorem, characterizes the (special) HN-filtration of a $\sigma$-module over $L\ancon$ (Theorem~\ref{Thm6.13}) and says that the (special) HN-polygon does not change under passing from $L\ancon$ to $\ol L\ancon$ (see Corollary~\ref{Cor6.14}). As in \cite{Kedlaya} the relation between $\sigma$-modules over $L\ancon$ and $\ol L\ancon$ is quite complicated. Therefore we too will have to make a detour via $L\dpl z\dpr$ where the comparison of HN-polygons is much easier (see Corollary~\ref{Cor10}). Since we may not yet use the assertion of the Slope Filtration Theorem we temporarily make the following definition.

\begin{definition} \label{DefIsoclinic}
Let $M$ be a $\sigma$-module over $L\ancon$ or over $k\dpl z\dpr$. We call the HN-slopes and the HN-polygon of $M\otimes \ol{L}\ancon$ (respectively $M\otimes\ol k\dpl z\dpr$)
the  \emph{absolute HN-slopes} and the \emph{absolute HN-polygon} of $M$. We say that $M$ is \emph{isoclinic} of slope $s$
if the absolute HN-slopes of $M$ are all equal to $s$. This means that $M\otimes\ol L\ancon$ (respectively $M\otimes \ol k\dpl z\dpr$) is isomorphic to a direct sum of $\CF_{d,n}$ for $s=-\frac{d}{n}$ with $(d,n)=1$ by Theorem~\ref{Thm1a} (respectively Theorem~\ref{Thm2}). 
If $B$ is an affinoid $L$-algebra, we call a $\sigma$-module $M$ over $B\ancon$ \emph{isoclinic} if at every analytic point $x\in\CM(B)$ (see Definition~\ref{DefAnalyticPoint})
the $\sigma$-module $M\otimes\wh{\kappa(x)^\alg}\ancon$ is isoclinic of the same slope.

A $\sigma$-module $M$ over $B\con$ is called \emph{isoclinic} if at every analytic point $x\in\CM(B)$ the $\sigma$-module $M\otimes\kappa(x)\dpl z\dpr$
is isoclinic of the same slope.
\end{definition}

As a justification of the latter definition we will see in Proposition~\ref{Prop7} below that if a $\sigma$-module $M$ over $B\con$ is isoclinic then also $M\otimes B\ancon$ is isoclinic of the same slope.

\begin{lemma}\label{Lemma1.5.6'}
Any isoclinic $\sigma$-module $M$ over $L\ancon$ or over $k\dpl z\dpr$ is semi-stable.
\end{lemma}

\begin{proof}
Let $M\otimes\ol L\ancon\cong\CF_{d,n}^{\oplus m}$ (respectively $M\otimes \ol k\dpl z\dpr\cong\CF_{d,n}^{\oplus m}$) with $\lambda(M)=-\frac{d}{n}$ and let $N\subset M$ be a $\sigma$-submodule. Then $N\otimes\ol L\ancon\cong\bigoplus_i \CF_{d_i,n_i}$ by Theorem~\ref{Thm1a} (respectively $N\otimes \ol k\dpl z\dpr\cong\bigoplus_i \CF_{d_i,n_i}$ by Theorem~\ref{Thm2}). From Proposition~\ref{Prop0.7} (respectively Proposition~\ref{Prop0.7b}) we conclude $\frac{d_i}{n_i}\le\frac{d}{n}$ (respectively $\frac{d_i}{n_i}=\frac{d}{n}$) for all $i$. This implies $\deg N=\sum_i d_i\le\frac{d}{n}\cdot\sum_i n_i=-\lambda(M)\cdot\rk N$ as desired.
\end{proof}

\medskip

Let us gather a few facts about the generic and special HN-polygons. We start with the generic HN-polygon. Let $k\supset\BF_q$ be a field and let $\ol k$ be some algebraically closed extension of $k$. First of all there is no need to distinguish between generic and generic absolute HN-polygons due to the following proposition. As a preparation we need a lemma.

\begin{lemma} \label{Lemma10b}
Let $M$ be a $\sigma$-module over $k\dpl z\dpr$
such that all slopes of $M' = M \otimes \ol k\dpl z\dpr$ are non-negative.
Then $M$ contains an $F$-stable $k\dbl z\dbr$-lattice $N$.
%Moreover if all slopes of $M'$ are zero, then $N$ can be chosen such that $F:\sigma^\ast N \to N$ is an isomorphism.
\end{lemma}

\begin{proof}
By Theorem~\ref{Thm2}, $M'$ is a direct sum with summands of the form $\CF_{d,n}$ for $d\leq0$.
In this presentation the standard basis vectors of the $\CF_{d,n}$ generate an $F$-stable $\ol k\dbl z\dbr$-lattice
in $M'$. Its intersection with $M$ gives the desired lattice in $M$. (Compare \cite[5.1.2]{Kedlaya}.)
\end{proof}

\begin{proposition}\label{Prop10}
Let $M$ be a $\sigma$-module over $k\dpl z\dpr$. Then there exists a uniquely determined filtration $0=M_0\subset M_1\subset\ldots\subset M_\ell=M$ of $M$ by $\sigma$-submodules over $k\dpl z\dpr$ with the following properties:
\begin{enumerate}
\item 
For $i=1,\ldots,\ell$ the $\sigma$-module $M_i/M_{i-1}$ is isoclinic of some slope $s_i$,
\item 
$s_1<\ldots<s_\ell$.
\end{enumerate}
It is the HN-filtration of $M$.
Moreover, if $k$ is perfect this filtration splits canonically and yields the slope decomposition $M\cong\bigoplus_{i=1}^\ell M_i/M_{i-1}$.
\end{proposition}

\begin{proof}
Let $M\otimes_{k\dpl z\dpr}\ol k\dpl z\dpr\;\cong\;\bigoplus_{j=1}^\ell \CF_{d_j,n_j}^{\oplus m_j}$ be the decomposition from Theorem~\ref{Thm2}. We may assume that the slopes $s_i:=-\frac{d_i}{n_i}$ satisfy condition (b). Put $M_i:=\bigoplus_{j=1}^i \CF_{d_j,n_j}^{\oplus m_j}$. Then $M_i/M_{i-1}$ is isoclinic of slope $s_i$. By induction on $\ell$ it suffices to descend $M_1$ to $k\dpl z\dpr$.
Choose a $k\dpl z\dpr$-basis $e_1,\ldots,e_n$ of $M$ and let $A\in\GL_n\bigl(k\dpl z\dpr\bigr)$ be the matrix by which $F_M$ acts on this basis. 

\smallskip

1.\es We first claim that $M_1$ descends to $k^\sep\dpl z\dpr$ where $k^\sep$ is the separable closure of $k$ in $\ol k$. By Lemma~\ref{LemmaKe3.6.2} it suffices to treat the case where $\dim_{\ol k\dpl z\dpr}M_1=1$ and hence $s_1\in\BZ$. By twisting with $\CO(s_1)$ we may then assume that $s_1=0$. So by Lemma~\ref{Lemma10b}, $M$ contains an $F$-stable $k\dbl z\dbr$-lattice and we may assume that the basis $e_1,\ldots,e_n$ is in fact a $k\dbl z\dbr$-basis of that lattice. This implies $A=\sum_{\nu=0}^\infty A_\nu z^\nu\in M_n\bigl(k\dbl z\dbr\bigr)$. Let $v\in M\otimes_{k\dpl z\dpr}\ol k\dpl z\dpr$ be a generator of $M_1$ satisfying $F_M(\sigma^\ast v)=v$ and let $x\in\ol k\dpl z\dpr^n$ be its coordinate vector with respect to the above basis. After multiplication with a power of $z$  we have $x=\sum_{\mu=0}^\infty x_\mu z^\mu\in\ol k\dbl z\dbr^n$ with $x_\mu\in\ol k^n$. The equation $F_M(\sigma^\ast v)=v$ implies
\[
\sum_{\nu=1}^\mu A_\nu x_{\mu-\nu}^\sigma \es=\es x_\mu-A_0 x_\mu^\sigma \qquad\text{for all }\mu\ge0\,.
\]
This is an \'etale equation for the components of $x_\mu$ and by induction on $\mu$ we obtain $x_\mu\in (k^\sep)^n$. thus $M_1$ descends to $k^\sep\dpl z\dpr$ as claimed.

\smallskip

2.\es There are two things left to be proven. Firstly we have to show that the $\sigma$-submodule 
$M_1=\CF_{d_1,n_1}^{\oplus m_1}\subset M\otimes k^\sep\dpl z\dpr$ further descends to $k\dpl z\dpr$. Secondly, if $k$ is perfect and $\ol k:=k^\sep$ we must prove that the $\sigma$-submodules $\CF_{d_j,n_j}^{\oplus m_j}$ descend to $k\dpl z\dpr$. Thus let 
$N:=\CF_{d_j,n_j}^{\oplus m_j}\subset M\otimes_{k\dpl z\dpr}k^\sep\dpl z\dpr$. By Lemma~\ref{LemmaKe3.6.2} we may assume that $N$ has dimension $1$, that is, $m_j=n_j=1$, and $N$ is generated by a vector $v=\sum_i x_ie_i$ with $x_i\in k^\sep\dpl z\dpr$ satisfying $F(\sigma^\ast v)=z^{d_j}v$. Without loss of generality we have $x_1\ne0$ and we set $w:=x_1^{-1}v$. Now let $\gamma\in\Gal(k^\sep/k)$. Then $F(\sigma^\ast v^\gamma)=z^{d_j}v^\gamma$ and $v^\gamma$ belongs to $N$, that is $v^\gamma=\alpha\,v$ for some $\alpha\in k^\sep\dpl z\dpr$. In particular $w^\gamma=\alpha\,x_1/x_1^\gamma\cdot w$. Since the first coordinate of $w$ with respect to $e_1,\ldots,e_n$ is $1$ we infer $\alpha\,x_1/x_1^\gamma=1$ and $w^\gamma=w$. Thus $w\in M$ and $N$ descends to $k\dpl z\dpr$ as desired.
(Compare \cite[5.3.1]{Kedlaya}.)
\end{proof}

\begin{corollary} \label{Cor10}
Let $M$ be a $\sigma$-module over $k\dpl z\dpr$. Then the HN-filtration of $M$, tensored up to $\ol k\dpl z\dpr$ gives the HN-filtration of $M\otimes\ol k\dpl z\dpr$.
\end{corollary}

\begin{proof}
The characterization of the HN-filtration given in Proposition~\ref{Prop10} is stable under base change.
\end{proof}

Although we won't need it in the sequel the reader should note the following analog of the theorem of Grothendieck-Katz \cite{Katz} which is proved for instance in \cite[Theorem 3.6]{RapoportRichartz}.

\begin{theorem} \label{ThmGrothendieck-Katz}
Let $A$ be an $\BF_q$-algebra and let $M$ be a $\sigma$-module of rank $n$ over $A\dbl z\dbr[z^{-1}]$. Let $P$ be the graph of a continuous real valued function on $[0,n]$ which is linear between successive integers. Then the set of points in $\Spec A$ at which the HN-polygon of $M$ lies above $P$ is Zariski closed. 
\end{theorem}

The following tool will be useful for determining the generic HN-polygon.

\begin{lemma} \label{Lemma14}
Let $M$ be a $\sigma$-module over $\ol k\dpl z\dpr$. Suppose $M$ has a basis on which $F$ acts by a matrix 
$A\in\GL_n\bigl(\ol k\dpl z\dpr\bigr)$ which satisfies $AD^{-1}-\Id_n\,\in M_n\bigl(z\,\ol k\dbl z\dbr\bigr)$ for some diagonal matrix $D\in\GL_n\bigl(\ol k\dpl z\dpr)$.
Then the generic HN-slopes of $M$ are equal to $\ord_z D_{\nu\nu}$, where $D_{\nu\nu}$ are the diagonal entries of $D$.
\end{lemma}

\begin{proof}
There is no harm in multiplying both $A$ and $D$ by the same power of $z$. Thus we may assume $A,D\in M_n\bigl(\ol k\dbl z\dbr\bigr)$.
We construct a sequence of matrices $U_\ell\in\GL_n\bigl(\ol k\dbl z\dbr\bigr)$ with $U_0=\Id_n$, $U_{\ell+1}\equiv U_\ell\mod z^{\ell+1}$, and $U_\ell^{-1}AU_\ell^\sigma D^{-1}\equiv\Id_n\mod z^{\ell+1}$.
Then the limit $U$ of the $U_\ell$ satisfies $U^{-1}AU^\sigma=D$, proving the lemma.

By assumption the conditions are satisfied for $\ell=0$. Suppose $U_\ell$ has been constructed. Put $V=U_\ell^{-1}AU_\ell^\sigma D^{-1}-\Id_n\,\in M_n\bigl(z^{\ell+1}\ol k\dbl z\dbr\bigr)$. 
Define a matrix $W$ whose entry $W_{\mu\nu}$ for each $\mu$ and $\nu$ is a solution of the congruence
\[
W_{\mu\nu}-D_{\mu\mu}W_{\mu\nu}^\sigma D_{\nu\nu}^{-1} \es\equiv\es V_{\mu\nu} \mod z^{\ell+2}
\]
with $W_{\mu\nu}, D_{\mu\mu}W_{\mu\nu}^\sigma D_{\nu\nu}^{-1} \,\in z^{\ell+1}\ol k\dbl z\dbr$. Such a solution exists because $\ol k$ is algebraically closed and depending on $\ord_z(D_{\mu\mu}D_{\nu\nu}^{-1})$
one has to solve either a linear equation, an Artin-Schreier equation or extract a $q$-th root over $\ol k$. Then both $W$ and $DW^\sigma D^{-1}$ belong to $M_n\bigl(z^{\ell+1}\ol k\dbl z\dbr\bigr)$. 
Put $U_{\ell+1}= U_\ell(\Id_n+W)$. Then $U_{\ell+1}\in\GL_n\bigl(\ol k\dbl z\dbr\bigr)$ satisfies $U_{\ell+1}\equiv U_\ell\mod z^{\ell+1}$ and
\begin{eqnarray*}
U_{\ell+1}^{-1}AU_{\ell+1}^\sigma D^{-1} &=& (\Id_n+W)^{-1}U_\ell^{-1}AU_\ell^\sigma(\Id_n+W)^\sigma D^{-1} \\
&=& (\Id_n+W)^{-1}U_\ell^{-1}AU_\ell^\sigma D^{-1}(\Id_n+DW^\sigma D^{-1}) \\
&=& (\Id_n+W)^{-1}(\Id_n+V)(\Id_n+DW^\sigma D^{-1}) \\
&\equiv& \Id_n-W+V+DW^\sigma D^{-1} \\
&\equiv& \Id_n \mod z^{\ell+2}\,.
\end{eqnarray*}
Thus we obtain $U_{\ell+1}$ as desired and the lemma follows.
(Compare \cite[5.2.6]{Kedlaya}.)
\end{proof}

\bigskip

Next let us turn towards $\sigma$-modules over $L\ancon$ and the special HN-polygon.

\begin{proposition}\label{Prop0.8}
\cite[Propositions 6.3 and 8.2]{HP} The $\sigma$-modules $\CF_{d,n}$ over $L\ancon$ are stable. The $\sigma$-modules $\CF_{d,n}^{\oplus m}$ and in particular $\CO(d)^{\oplus m}$ over $L\ancon$ are semistable.
\end{proposition}

\begin{theorem}\label{Thm1}
\cite[Corollary 11.7]{HP} Let $M$ be a $\sigma$-module over $\ol L\ancon$ and let $M\cong\bigoplus_{i=1}^\ell\CF_{d_i,n_i}^{\oplus m_i}$ be its decomposition from Theorem~\ref{Thm1a}. Assume that $d_{i-1}/n_{i-1}>d_i/n_i$ for all $i$ and set $M_j:=\bigoplus_{i=1}^j\CF_{d_i,n_i}^{\oplus m_i}$. Then the HN-filtration of $M$ is $0=M_0\subset M_1\subset\ldots\subset M_\ell=M$ and the HN-slopes are the numbers $\frac{d_i}{n_i}$ with multiplicity $n_im_i$. In particular the HN-polygon of $M$ determines its isomorphy type.
\end{theorem}

\begin{proposition} \label{Prop8'}
Every $\sigma$-module $M$ over $L\ancon$ admits a Harder-Narasimhan filtration.
\end{proposition}

\begin{proof}
By Theorem~\ref{Thm1a} and Proposition~\ref{Prop0.8} the set of slopes of all nontrivial $\sigma$-submodules of $M$
has a smallest element $s_1$. 
Let $N_1$ and $N_2$ be two $\sigma$-submodules of $M$ of slope $s_1$. Then $N_1$ and $N_2$, and thus also $N_1\oplus N_2$ are semistable.
Since the sum $N_1+N_2$ of $N_1$ and $N_2$ inside $M$ is a quotient of $N_1\oplus N_2$ its slope is at most $s_1$. By minimality of $s_1$ it equals $s_1$.
Therefore there exists a unique maximal $\sigma$-submodule $M_1$ of $M$ of slope $s_1$.
Then $M_1$ is the first step in the HN-filtration of $M$. The remaining steps are obtained by replacing $M$ by $M/M_1$. (Compare \cite[4.2.5]{Kedlaya}.)
\end{proof}

\begin{proposition} \label{Prop8}
Let $M$ be a $\sigma$-module over $L\ancon$. Then the HN-polygon lies above the slope polygon of any semistable filtration of $M$.
\end{proposition}

\begin{proof}
Let $0=M_0\subset M_1\subset\ldots\subset M_\ell=M$ be the HN-filtration and let $0=M'_0\subset M'_1\subset\ldots\subset M'_m=M$
be a semistable filtration of $M$. It suffices to prove that for each $i=1,\ldots,\ell$ we can find $\rk M_i$ slopes from the 
slope multiset of the semistable filtration whose sum is less than or equal to $-\deg M_i$.

For $j=1,\ldots,m$ consider the $\sigma$-submodule $(M'_j\cap M_i)/(M'_{j-1}\cap M_i)$ of $M'_j/M'_{j-1}$. Note that the $M'_j\cap M_i$ are saturated $\sigma$-submodules of $M_i$. Since $M'_j/M'_{j-1}$ is semistable we have
\[
\lambda_j\es:=\es\lambda\bigl((M'_j\cap M_i)/(M'_{j-1}\cap M_i)\bigr) \es\geq\es\lambda(M'_j/M'_{j-1})\,.
\]
Set $r_j \;=\;\rk\,(M'_j\cap M_i)/(M'_{j-1}\cap M_i)\;\leq\;\rk M'_j/M'_{j-1}$. Since the $M'_j\cap M_i$ filter $M_i$ we have
\[
-\deg M_i\es =\es \sum_{j=1}^m r_j \lambda_j\es\ge\es\sum_{j=1}^m r_j\cdot\lambda(M'_j/M'_{j-1})\,.
\]
From this the desired inequality follows. (Compare \cite[3.5.4]{Kedlaya}.)
\end{proof}

\begin{proposition} \label{Prop16}
Let $0\to M_1\to M\to M_2\to 0$ be a short exact sequence of $\sigma$-modules over $L\ancon$
and let $P(M_1), P(M)$, and $P(M_2)$ be their absolute HN-polygons. Then
\begin{enumerate}
\item
$P(M)\geq P(M_1)+P(M_2)$,
\item
$P(M)= P(M_1)+P(M_2)$ if and only if the short exact sequence splits over $\ol{L}\ancon$.
\end{enumerate}
\end{proposition}

\begin{proof}
(The proof literally follows \cite[Proposition 4.7.2]{Kedlaya}.)
To prove (a) note that the HN-filtrations of $M_1\otimes \ol{L}\ancon$ and $M_2\otimes\ol{L}\ancon$
induce a semistable filtration of $M\otimes\ol{L}\ancon$ whose slope polygon is $P(M_1)+P(M_2)$. Then Proposition~\ref{Prop8} yields the claim.

To prove (b) note that the splitting of the sequence implies $P(M)=P(M_1)+P(M_2)$ by Theorem~\ref{Thm1}. Conversely suppose $P(M)=P(M_1)+P(M_2)$. 
We prove by induction on the rank of $M$ that the sequence splits in case $L=\ol{L}$.
The induction starts with the case where $M_1=\CF_{c,r}$ and $M_2=\CF_{d,s}$.
By Proposition~\ref{Prop0.7} the sequence splits if $\frac{d}{s}\leq\frac{c}{r}$. So now assume $\frac{d}{s}>\frac{c}{r}$.
By Theorem~\ref{Thm1} the assumption on $P(M)$ implies that $M\cong M_1\oplus M_2$, yielding a map $M_2\to M$.
Since $\End_\sigma(\CF_{d,s})$ is a division algebra by Proposition~\ref{Prop0.7} the composite
$M_2\to M\to M_2$ is either zero or an isomorphism. In the former case, the map factors through $M_2\to M_1$ by exactness.
But because $\frac{d}{s}>\frac{c}{r}$, this is impossible by Proposition~\ref{Prop0.7}.
Thus $M_2\to M\to M_2$ is an isomorphism splitting the sequence.

Next assume that $M_1$ decomposes as $M_1=N\oplus N'$. We compute
\[
\begin{array}{r@{\es}c@{\es}l@{\qquad}l}
P(M) & \geq & P(N)+P(M/N) & \text{[by (a)]} \\
& \geq & P(N)+P(M_1/N)+P(M_2) & \text{[by (a)]} \\
& = & P(M_1)+P(M_2) & \text{[because $N$ is a summand of $M_1$]} \\
& = & P(M) & \text{[by assumption]}
\end{array}
\]
and find that all the inequalities must be equalities. In particular $P(M/N)=P(M_1/N)+P(M_2)$.
Since $M/N$ is of smaller rank, the induction hypothesis implies that the sequence 
$0\to M_1/N\to M/N\to M_2\to 0$ splits. We obtain a map $M\to M_1/N\cong N'$ which splits the sequence
$0\to N'\to M\to M/N'\to 0$. It remains to consider the sequence $0\to N\to M/N'\to M_2\to 0$
which splits because $P(M/N')=P(N)+P(M_2)$ by a similar computation as above. Thus we have proved the claim whenever $M_1$ decomposes.

There remains the case, when only $M_2$ decomposes. By considering the dualized sequence, we reduce to the case in which $M_1$ decomposes. This proves the proposition.
\end{proof}

%%%%%%%%%%%%%%%%%%%%%%%%%%%%%%%%%%%%%%%%%%%%%%%%%%%%%%%%%%%%%%%%%%%%%%
%
%    Comparison of Harder-Narasimhan Polygons
%
%%%%%%%%%%%%%%%%%%%%%%%%%%%%%%%%%%%%%%%%%%%%%%%%%%%%%%%%%%%%%%%%%%%%%%

\subsection{Comparison of Harder-Narasimhan Polygons} \label{SectComparisonHNPolygons}

We will now study the relation between the generic and special HN-polygons. We need three preparatory lemmas. The first is proved in the same way as Lemma~\ref{Lemma10b}.

\begin{lemma} \label{Lemma4}
Let $M$ be a $\sigma$-module over $L\con$
such that all slopes of $M' = M \otimes \ol{L}\dpl z\dpr$ are non-negative.
Then $M$ contains an $F$-stable $L\langle\frac{z}{\zeta}\rangle$-lattice $N$.
Moreover if all slopes of $M'$ are zero, then $N$ can be chosen such that
$F:\sigma^\ast N \to N$ is an isomorphism.
\qed
\end{lemma}

\begin{lemma} \label{Lemma3}
Let $M$ be a $\sigma$-module over 
$\ol{L}\con$ such that all HN-slopes
of $M':=M\otimes \ol{L}\dpl z\dpr$ are non-positive. Let $v\in M'$ satisfy $F(\sigma^\ast v)=v$. Then $v\in M$.
\end{lemma}

\begin{proof}
By Lemma~\ref{Lemma4} we can find an $F$-stable $\ol{L}\langle\frac{z}{\zeta}\rangle$-lattice $N\dual$ in $M\dual$.
Therefore its dual lattice $N=\{w\in M:f(w)\in \ol L\langle\frac{z}{\zeta}\rangle \text{ for all }f\in N\dual\,\}$ in $M$ is stable under $F^{-1}:M\to\sigma^\ast M$. Let $e_1,\ldots,e_n$ be a basis of $N$ and let 
$A=(A_{\mu\nu})\in M_n\bigl(\ol{L}\langle\frac{z}{\zeta^q}\rangle\bigr)$ be the matrix with
$F^{-1}e_\nu=\sum_\mu A_{\mu\nu}\sigma^\ast(e_\mu)$. 

Let $x\in\ol{L}\dpl z\dpr^n$ be the coordinate vector of $v$ with respect to the basis $e_1,\ldots,e_n$. By multiplying $v$ with a power of $z$ we may assume that $x\in \ol{L}\dbl z\dbr^n$. Then $F(\sigma^\ast v)=v$ amounts to the equation $x^\sigma=A x$. We deduce from Lemma~\ref{Lemma3b} that $x\in \ol{L}\langle\frac{z}{\zeta}\rangle^n$
 and $v\in M$ as desired. (Compare \cite[5.4.1]{Kedlaya}.)
\end{proof}

\begin{lemma}\label{LemmaDescSemistFilt}
Let $M$ be a $\sigma$-module over $\ol L\con$ such that $M\otimes \ol L\dpl z\dpr\cong\CF_{d,n}^{\oplus m}$. Then $M$ itself is isomorphic to the $\sigma$-module $\CF_{d,n}^{\oplus m}$ over $\ol L\con$.
\end{lemma}

\begin{proof}
Let $e_1,\ldots,e_{mn}$ be the standard basis of the $\sigma$-module $\CF_{d,n}^{\oplus m}$ over $\ol L\dpl z\dpr$, that is, 
\[
\begin{array}{r@{\quad=\quad}l@{\es\qquad}l}
F(\sigma^\ast e_{\ell n+i}) & e_{\ell n+i+1} & \text{for all } i=1,\ldots,n-1 \quad\text{and }\ell=0,\ldots,m-1\,, \\[1mm]
F(\sigma^\ast e_{\ell n+n}) & z^{-d}e_{\ell n+1} & \text{for all } \ell=0,\ldots,m-1\,.
\end{array}
\]
Replace $F$ by $F':=z^d F^n:(\sigma^n)^\ast M\isoto M$. Then the vectors $e_i$ satisfy $F'(\sigma^\ast e_i)=e_i$. In particular $e_i\in M$ by Lemma~\ref{Lemma3}. Clearly $e_1,\ldots,e_{mn}$ form a basis of $M$ over $\ol L\con$. This proves the lemma.
\end{proof}

In order to compare the generic and special HN-polygons of a $\sigma$-module over $\ol L\con$ we need to descend the generic information from $\ol L\dpl z\dpr$ to $\ol L\con$. This is achieved by the reverse filtration.

\begin{definition} \label{Def5}
Let $M$ be a $\sigma$-module over $\ol{L}\dpl z\dpr$ with slope decomposition $P_1\oplus\ldots\oplus P_\ell$ labeled so that
$\lambda(P_1) >\ldots >\lambda(P_\ell)$. We call the filtration $0=M_0\subset M_1\subset\ldots\subset M_\ell=M$ 
with $M_i=P_1\oplus\ldots\oplus P_i$ for $i=0,\ldots,\ell$ the \emph{reverse filtration} of $M$.
By construction its slope polygon is the HN-polygon of $M$.
\end{definition}

\begin{proposition} \label{Prop6}
Let $M$ be a $\sigma$-module over $\ol{L}\con$, then the reverse filtration of $M\otimes\ol{L}\dpl z\dpr$ descends to $\ol{L}\con$.
\end{proposition}

\begin{proof}
Let $0=M'_0\subset M'_1\subset\ldots\subset M'_\ell =M'$ be the reverse filtration of $M'=M\otimes\ol{L}\dpl z\dpr$.
It suffices to show that $M'_1$ descends to $\ol{L}\con$. By Lemma~\ref{LemmaKe3.6.2} we may reduce to the case where $\rk M'_1=1$ by passing from $M$ to an exterior power.
By twisting we may then reduce to the case where $\lambda(M'_1)=0$. Thus $M'_1$ is isomorphic to $\CO(0)$ and therefore $M'_1$ is generated by a vector $v$ with $Fv=v$.
By Lemma~\ref{Lemma3}, $v$ belongs to $M$. Hence $M'_1$ descends to $\ol{L}\con$ proving the proposition.
(Compare \cite[5.4.3]{Kedlaya}.)
\end{proof}

\begin{proposition} \label{Prop7}
Let $M$ be a $\sigma$-module over $L\con$. Then the HN-polygon of $M\otimes \ol{L}\ancon$ lies above the HN-polygon of $M\otimes \ol{L}\dpl z\dpr$.
\end{proposition}

\begin{proof}
The HN-polygon of $M\otimes \ol{L}\dpl z\dpr$ equals the slope polygon of its (descended) reverse filtration. The factors of the latter are isomorphic to $\sigma$-modules $\CF_{d,n}^{\oplus m}$ over $\ol L\con$ for appropriate $d,n$ and $m$ by Lemma~\ref{LemmaDescSemistFilt}. This property is preserved by tensoring the reverse filtration up to $\ol L\ancon$. In this way we obtain a semistable filtration of $M\otimes \ol L\ancon$; use Proposition~\ref{Prop0.8}. The slope polygon of this filtration equals the HN-polygon of $M\otimes\ol L\dpl z\dpr$ and lies below the HN-polygon of $M\otimes\ol{L}\ancon$ by Proposition~\ref{Prop8}.
(Compare \cite[5.5.1]{Kedlaya}.)
\end{proof}

The case when both polygons coincide is particularly favorable.

\begin{theorem} \label{Thm9}
Let $M$ be a $\sigma$-module over $L\con$ whose generic and special absolute HN-polygons coincide. 
Then the generic and special absolute HN-filtrations of $M\otimes\ol{L}\dpl z\dpr$ and $M\otimes\ol{L}\ancon$, respectively, are both obtained from a filtration of $M$.
\end{theorem}

\begin{proof}
It is enough to show that the first steps of the generic and special absolute HN-filtrations descend and coincide. By Lemma~\ref{LemmaKe3.6.2} we may reduce to the case 
where the least slope of the common HN-polygon occurs with multiplicity one. Choose a basis $e_1,\ldots,e_n$ of $M$.
Let $v\in M\otimes\ol{L}\dpl z\dpr$ be a generator of the first step of the HN-filtration of $M\otimes\ol{L}\dpl z\dpr$ and write $v=a_1 e_1+\ldots +a_n e_n$ with $a_i\in\ol{L}\dpl z\dpr$. Without loss of generality, $a_1\neq0$ and we can assume $a_1=1$. 
By Corollary~\ref{Cor10} we have $a_i\in L\dpl z\dpr$ for all $i$.

Consider the descended reverse filtration $0=\wt M_0\subset\wt M_1\subset\ldots\subset\wt M_\ell$
of $\wt M_\ell:=M\otimes\ol{L}\con$ from Proposition~\ref{Prop6}. It satisfies $\wt M_\ell\otimes\ol L\dpl z\dpr=\wt M_{\ell-1}\otimes\ol L\dpl z\dpr\oplus\ol{L}\dpl z\dpr\cdot v$.
By Proposition~\ref{Prop7} the HN-polygons satisfy $P\bigr(\wt M_{\ell-1}\otimes\ol{L}\dpl z\dpr\bigr)\leq P\bigl(\wt M_{\ell-1}\otimes\ol{L}\ancon\bigr)$
and analogously for $\wt M_\ell/\wt M_{\ell-1}$. Since by assumption we have
\begin{eqnarray*}
P\bigl(\wt M_{\ell}\otimes\ol{L}\ancon\bigr) & =  & P\bigr(\wt M_{\ell}\otimes\ol{L}\dpl z\dpr\bigr) \\
& = & P\bigr(\wt M_{\ell-1}\otimes\ol{L}\dpl z\dpr\bigr) + P\bigr((\wt M_\ell/\wt M_{\ell-1})\otimes\ol{L}\dpl z\dpr\bigr)\\
& \le & P\bigr(\wt M_{\ell-1}\otimes\ol{L}\ancon\bigr) + P\bigr((\wt M_\ell/\wt M_{\ell-1})\otimes\ol{L}\ancon\bigr)\\
& \le & P\bigl(\wt M_\ell\otimes\ol L\ancon\bigr)\,,
\end{eqnarray*}
we deduce from Proposition~\ref{Prop16} that the sequence $0\to \wt M_{\ell-1}\to\wt M_\ell \to\wt M_\ell/\wt M_{\ell-1}\to0$ splits after tensoring up to $\ol{L}\ancon$.
Since all generic slopes of $\CHom(\wt M_\ell/\wt M_{\ell-1}\,,\,\wt M_{\ell-1})$ are positive, Lemma~\ref{Lemma4} allows us to apply Proposition~\ref{Prop13} to deduce 
that the sequence already splits over $\ol{L}\con$. Therefore there is a vector $\wt v=b_1e_1+\ldots+b_ne_n\in\wt M_\ell$, $b_i\in\ol L\con$, which generates a $\sigma$-submodule isomorphic to $\wt M_\ell/\wt M_{\ell-1}$.
Over $\ol L\dpl z\dpr$ the vectors $v$ and $\wt v\in\wt M_\ell\otimes\ol L\dpl z\dpr$ define splittings of the sequence
\[
0\es\longto\es \wt M_{\ell-1}\otimes\ol L\dpl z\dpr\es\longto\es\wt M_\ell\otimes\ol L\dpl z\dpr \es\longto\es(\wt M_\ell/\wt M_{\ell-1})\otimes\ol L\dpl z\dpr\es\longto\es0\,.
\] 
Since all generic slopes of $\wt M_\ell/\wt M_{\ell-1}$ are less than all generic slopes of $\wt M_{\ell-1}$, these splittings coincide by Proposition~\ref{Prop0.7b}. Hence $\wt v=c\,v$ for some $c\in\ol L\dpl z\dpr\mal$. By construction $c=b_1\in\ol L\con$. Since $c$ has only finitely many zeroes by Lemma~\ref{LemmaPID} there is a rational number $r\geq1$ with $c\in\ol L\con[r]\mal$. In particular
\[
a_i\es=\es c^{-1}b_i\es\in\es L\dpl z\dpr\cap\ol L\con[r]\es=\es L\con[r]\,.
\]
So the vector $v\in M\otimes L\con[r]$ defines a saturated $\sigma$-submodule $N'$ of $M\otimes L\con[r]$ with $N'\otimes\ol L\dpl z\dpr=\ol L\dpl z\dpr\cdot v$. By repeated application of Lemma~\ref{Lemma9b} below, $N'\otimes L\langle\frac{z}{\zeta^r},\frac{\zeta^{qr}}{z}\rangle$ extends to a saturated $\sigma$-submodule of $M\otimes L\langle\frac{z}{\zeta},\frac{\zeta^{qr}}{z}\rangle$ which we glue over $L\langle\frac{z}{\zeta^r},\frac{\zeta^{qr}}{z}\rangle$ with $N'$ to a saturated $\sigma$-submodule of $M$.
Thus the first step of the generic HN-filtration descends to $L\con$. Let $M_1$ be the corresponding $\sigma$-submodule of $M$.

Let $M'_1$ be the first step of the HN-filtration of $M\otimes\ol{L}\ancon$. So $M'_1$ is characterized as the maximal $\sigma$-submodule of minimal slope, which by assumption is $\lambda(M'_1)=\lambda(M_1)$.
In particular $M_1\otimes\ol{L}\ancon\subset M'_1$ and since their ranks and degrees are equal, $M_1\otimes\ol{L}\ancon$ and $M'_1$ coincide by \cite[Proposition 6.2]{HP}. 
Hence the first step of the special absolute HN-filtration also descends to $M_1$. This proves the theorem.
(Compare \cite[5.5.2]{Kedlaya}.)
\end{proof}

\begin{lemma}\label{Lemma9b}
Let $r$ and $r'$ be rational numbers with $r'\geq q^2r>0$. Let $M$ be a $\sigma$-module over $L\langle\frac{z}{\zeta^r},\frac{\zeta^{r'}}{z}\rangle$ and let $N'\subset M\otimes L\langle\frac{z}{\zeta^{qr}},\frac{\zeta^{r'}}{z}\rangle$ be a saturated $\sigma$-submodule. Then there exists a uniquely determined saturated $\sigma$-submodule $N$ of $M$ with $N'=N\otimes L\langle\frac{z}{\zeta^{qr}},\frac{\zeta^{r'}}{z}\rangle$.
\end{lemma}

\begin{proof}
Clearly $N$ is unique if it exists since it can be described as the intersection $M\cap N'$ inside $M\otimes L\langle\frac{z}{\zeta^{qr}},\frac{\zeta^{r'}}{z}\rangle$.

To prove existence let $d:=\rk N'$. By Lemma~\ref{LemmaKe3.6.2} it suffices to show that $\wedge^dN'\subset(\wedge^dM)\otimes L\langle\frac{z}{\zeta^{qr}},\frac{\zeta^{r'}}{z}\rangle$ descends to $L\langle\frac{z}{\zeta^r},\frac{\zeta^{r'}}{z}\rangle$. Thus we may assume that $\rk N'=1$. Choose a basis $e_1,\ldots,e_n$ of $M$ and let $\Phi\in \GL_n\bigl(L\langle\frac{z}{\zeta^{qr}},\frac{\zeta^{r'}}{z}\rangle\bigr)$ be the matrix by which $F_M$ acts on this basis. Let $v=a_1e_1+\ldots+a_ne_n$ with $a_i\in L\langle\frac{z}{\zeta^{qr}},\frac{\zeta^{r'}}{z}\rangle$ be a generator of $N'$. Without loss of generality $a_1\neq0$. By Lemma~\ref{LemmaPID} we may multiply $v$ by a unit and assume that $a_1\in L[z]$ is a monic polynomial whose zeroes $x$ all satisfy $|\zeta^{r'}|\leq|x|\leq|\zeta^{qr}|$. Since $N'$ is $F_M$-stable there exists a unit $\alpha\in L\langle\frac{z}{\zeta^{q^2r}},\frac{\zeta^{r'}}{z}\rangle\mal$ with $F_M(\sigma^\ast v)=\alpha v$. Let
\[
\left(\begin{array}{c} a'_1 \\ \vdots \\ a'_n \end{array}\right) \es:=\es
\Phi^{-1} \left(\begin{array}{c} a_1 \\ \vdots \\ a_n \end{array}\right) \es=\es
\alpha^{-1} \left(\begin{array}{c} a_1^\sigma \\ \vdots \\ a_n^\sigma\end{array}\right)
\quad\in\es \TS L\langle\frac{z}{\zeta^{qr}},\frac{\zeta^{r'}}{z}\rangle^{\oplus n}\,.
\]
Again by Lemma~\ref{LemmaPID} we write $a'_1=uf$ with a unit $u\in L\langle\frac{z}{\zeta^{qr}},\frac{\zeta^{r'}}{z}\rangle\mal$ and a monic polynomial $f\in L[z]$ whose zeroes $x$ satisfy $|\zeta^{r'}|\leq|x|\leq|\zeta^{qr}|$. By \cite[Proposition 2]{Lazard} we may split $f=g\tilde f$ with $g,\tilde f\in L[z]$ monic such that the zeroes $x$ of $g$ (respectively of $\tilde f$) all satisfy $|x|>|\zeta^{q^2r}|$ (respectively $|x|\leq|\zeta^{q^2r}|$). Similarly we split $a_1^\sigma=hb$ with $h,b\in L[z]$ monic such that the zeroes $x$ of $h$ (respectively of $b$) all satisfy $|x|<|\zeta^{r'}|$ (respectively $|\zeta^{r'}|\leq|x|\leq|\zeta^{q^2r}|$). Since $u,\alpha, g$ and $h$ are units in $L\langle\frac{z}{\zeta^{q^2r}},\frac{\zeta^{r'}}{z}\rangle$ the equation $ug\tilde f=a_1'=\alpha^{-1}hb$ implies $b=\tilde f$ by Lemma~\ref{LemmaPID}. Hence $ug=\alpha^{-1}h$. Now there exists a polynomial $\tilde g\in L[z]$ with $\tilde g^\sigma=g^q$, namely $\tilde g(z):=g(z^q)$. Consider the vector $\tilde g v\in N'$. We have
\[
(\tilde ga_i)^\sigma \es=\es g^q a_i^\sigma  \es=\es g^{q-1}hu^{-1} a_i' 
\TS\quad\in\es L\langle\frac{z}{\zeta^{qr}},\frac{\zeta^{r'}}{z}\rangle\,.
\]
Therefore $\tilde g a_i\in L\langle\frac{z}{\zeta^r},\frac{\zeta^{r'}}{z}\rangle$. After clearing a common factor of the $\tilde g a_i$ in $L\langle\frac{z}{\zeta^r},\frac{\zeta^{r'}}{z}\rangle$ if necessary, the vector $\tilde g v$ defines the desired saturated $\sigma$-submodule $N$ of $M$. It satisfies $N'=N\otimes L\langle\frac{z}{\zeta^{qr}},\frac{\zeta^{r'}}{z}\rangle$ since $\tilde g\in L\langle\frac{z}{\zeta^{qr}},\frac{\zeta^{r'}}{z}\rangle\mal$.
\end{proof}

%%%%%%%%%%%%%%%%%%%%%%%%%%%%%%%%%%%%%%%%%%%%%%%%%%%%%%%%%%%%%%%%%%%%%%
%
%    Descent of $\sigma$-Modules and HN-Filtrations
%
%%%%%%%%%%%%%%%%%%%%%%%%%%%%%%%%%%%%%%%%%%%%%%%%%%%%%%%%%%%%%%%%%%%%%%

\subsection{Descent of \texorpdfstring{$\sigma$}{sigma}-Modules and HN-Filtrations} \label{SectDescentHNFiltrations}

In this section we finally prove the Slope Filtration Theorem as well as the fact that isoclinic $\sigma$-modules over $B\ancon$ descend to $B\con$. For this purpose we will need to work locally on the Berkovich space $\CM(B)$ associated with $B$ (see Appendix~\refAppBerkovich). Also recall our definition of affinoid covering of $\CM(B)$ from Definition~\ref{DefBerkovichSpaces}.

\begin{lemma} \label{Lemma6.1}
Let $r$ be a positive rational number and let $D\in\GL_n\bigl(B\langle\frac{z}{\zeta^r},\frac{\zeta^r}{z}\rangle\bigr)$. Define the number $h:=\bigl(\|D\|_r\,\|D^{-1}\|_r\bigr)^{1/(q-1)}\,\geq\,1$
and let $A\in M_n\bigl(B\langle\frac{z}{\zeta^r},\frac{\zeta^r}{z}\rangle\bigr)$ satisfy $\|AD^{-1}-\Id_n\|_r<h^{-1}\le1$.
Then there exists a matrix $U\in\GL_n\bigl(B\langle\frac{\zeta^r}{z}\rangle\bigr)$ with $\|U-\Id_n\|_r<1$ and
\[
U^{-1}AU^\sigma D^{-1}-\Id_n \es\in\es M_n\bigl(zB{\TS\langle\frac{z}{\zeta^r}\rangle}\bigr)\qquad\text{and}\qquad\|U^{-1}AU^\sigma D^{-1}-\Id_n\|_r<1\,.
\]
\end{lemma}

\begin{proof}
We define sequences of matrices $U_0,U_1,\ldots$ and $A_0,A_1,\ldots$ as follows. We start with $U_0=\Id_n$. Given $U_\ell\in\GL_n\bigl(B\langle\frac{\zeta^r}{z}\rangle\bigr)$ we set $A_\ell=U_\ell^{-1}AU_\ell^\sigma$.
We define an additive function $f:B\langle\frac{z}{\zeta^r},\frac{\zeta^r}{z}\rangle\to B\langle\frac{\zeta^r}{z}\rangle$ by mapping $a=\sum_{i\in\BZ}a_i z^i$ to $\sum_{i\leq0}a_i z^i$ and we extend it to a function on matrices.
We have $\|f(a)\|_r\leq\|a\|_r$. Set $X_\ell=f(A_\ell D^{-1}-\Id_n)$ and $U_{\ell+1}=U_\ell(\Id_n+X_\ell)$. Put $c=\|AD^{-1}-\Id_n\|_r\cdot h<1$ and $c_\ell=\|X_\ell\|_r\cdot h$.

We claim that $c_\ell\leq c^{\ell+1}$ and $\|A_\ell D^{-1}-\Id_n\|_r\leq\|AD^{-1}-\Id_n\|_r<1$ for all $\ell\geq0$.

Clearly this holds for $\ell=0$. Assume the claim for $\ell$. Then $\|X_\ell\|_r<1$ and in particular $\Id_n+X_\ell$ and $U_{\ell+1}$ belong to $\GL_n\bigl(B\langle\frac{\zeta^r}{z}\rangle\bigr)$ and $A_\ell D^{-1}$ belongs to
$\GL_n\bigl(B\langle\frac{z}{\zeta^r},\frac{\zeta^r}{z}\rangle\bigr)$. We obtain
\[
A_{\ell+1}D^{-1}\es=\es (\Id_n+X_\ell)^{-1}A_\ell D^{-1}(\Id_n+DX_\ell^\sigma D^{-1})\,.
\]
Write
\begin{eqnarray*}
(\Id_n+X_\ell)^{-1}A_\ell D^{-1} &=& (\Id_n-X_\ell)A_\ell D^{-1}\,+\, (\Id_n+X_\ell)^{-1}X_\ell^2 A_\ell D^{-1} \\
&=& \Id_n+Y_\ell+V_\ell+W_\ell
\end{eqnarray*}
with
\begin{eqnarray*}
Y_\ell &=& A_\ell D^{-1}-\Id_n-X_\ell\,, \\
V_\ell &=& -X_\ell(A_\ell D^{-1}-\Id_n)\,, \\
W_\ell &=& X_\ell^2(\Id_n+X_\ell)^{-1}A_\ell D^{-1}\,.
\end{eqnarray*}
Note that $Y_\ell\in M_n\bigl(zB\langle\frac{z}{\zeta^r}\rangle\bigr)$, in particular $f(Y_\ell)=0$. We estimate using Lemma~\ref{Lemma6.0}
\begin{eqnarray*}
\bigl\| (\Id_n+X_\ell)^{-1}A_\ell D^{-1}\bigr\|_r &=& 1 \,, \\
\bigl\| V_\ell\bigr\|_r &\leq& c_\ell \,c\,h^{-2} \es\leq\es c^{\ell+2}\,h^{-2} \,, \\
\bigl\| W_\ell\bigr\|_r &\leq& (c_\ell\,h^{-1})^2 \es\leq\es c^{2\ell+2}\,h^{-2} \,, \\
\bigl\| Y_\ell\bigr\|_r &\leq& c\,h^{-1}  \,, \\
\bigl\| DX_\ell^\sigma D^{-1}\bigr\|_r &\leq& \|D\|_r\,\|D^{-1}\|_r\,\|X_\ell\|_r^q \es\leq\es h^{q-1}(c_\ell\,h^{-1})^q \es\leq\es c^{q(\ell+1)}\,h^{-1}\,.
\end{eqnarray*}
Putting everything together we obtain
\begin{eqnarray*}
\|A_{\ell+1}D^{-1}-\Id_n\|_r &=& \bigl\| Y_\ell+V_\ell+W_\ell+\,(\Id_n+X_\ell)^{-1}A_\ell D^{-1}(DX_\ell^\sigma D^{-1})\bigr\|_r \\
&\leq& c\,h^{-1}\es=\es \|AD^{-1}-\Id_n\|_r\,.
\end{eqnarray*}
Moreover,
\[
X_{\ell+1}\es=\es f(V_\ell)+f(W_\ell)+f\bigl((\Id_n+X_\ell)^{-1}A_\ell D^{-1}(DX_\ell^\sigma D^{-1})\bigr)\,.
\]
We compute $c_{\ell+1}\,h^{-1}\;=\;\|X_{\ell+1}\|_r\;\leq\; c^{\ell+2}\,h^{-1}$ and the claim follows.

Since $\|X_\ell\|_r\to0$ as $\ell\to\infty$ the sequence $U_\ell$ converges to a matrix $U\in\GL_n\bigl(B\langle\frac{\zeta^r}{z}\rangle\bigr)$ satisfying
$\|U^{-1}AU^\sigma D^{-1}-\Id_n\|_r\;\leq\;\|AD^{-1}-\Id_n\|_r\;<\;1$. Moreover, since $f$ is continuous with respect to the norm $\|\,.\,\|_r$ we obtain
\[
f(U^{-1}AU^\sigma D^{-1}-\Id_n)\es=\es \lim_{\ell\to\infty} f(A_\ell D^{-1}-\Id_n) \es=\es \lim_{\ell\to\infty} X_\ell \es=\es 0\,.
\]
Hence $U^{-1}AU^\sigma D^{-1}-\Id_n\;\in M_n\bigl(zB\langle\frac{z}{\zeta^r}\rangle\bigr)$. Since $U_{\ell+1}-\Id_n=(U_\ell-\Id_n)(\Id_n+X_\ell)+X_\ell$ and $\|X_\ell\|_r\le ch^{-1}$ we find $\|U_\ell-\Id_n\|_r\le ch^{-1}<1$ by induction. Hence $\|U-\Id_n\|_r<1$ as desired.
(Compare \cite[6.1.1]{Kedlaya}.)
\end{proof}

\begin{proposition} \label{Prop6.3}
Let $M$ be a $\sigma$-module over $B\ancon[r]$ and let $x\in\CM(B)$ be an analytic point. Then there exists a connected affinoid neighborhood $\CM(B')$ of $x$, a finite \'etale $B'$-algebra $C$,
a lift of $x$ to a point of $\CM(C)$, a positive integer $s$, and a $\sigma^s$-module $M'$ over $C\con[r]$ such that $M'\otimes C\ancon[r]\cong M\otimes C\ancon[r]$ as $\sigma^s$-modules, 
and such that for all analytic points $y\in\CM(C)$ the generic HN-polygon of $M'\otimes\kappa(y)\con[r]$ coincides with the special absolute HN-polygon of $M\otimes\kappa(x)\ancon[r]$. 
The integer $s$ can be chosen to be the smallest common denominator of all the slopes of $M\otimes\kappa(x)\ancon[r]$.
\end{proposition}

\begin{proof}
Note that we may assume that the slopes of $M\otimes \kappa(x)\ancon[r]$ are integers 
by replacing $F$ by $F^s$ for a suitable positive integer $s$. In the sequel we write again $q$ for $q^s$ and $\sigma$ for $\sigma^s$ to shorten notation.
By Lemma~\ref{Lemma6.4} we may assume that $M\otimes B\langle\frac{z}{\zeta^r},\frac{\zeta^{qr}}{z}\rangle$ is free. 
Let $e_1,\ldots,e_n$ be a basis and let $A\in\GL_n\bigl(B\langle\frac{z}{\zeta^{qr}},\frac{\zeta^{qr}}{z}\rangle\bigr)$ be the matrix by which $F$ acts on this basis.
Consider the base change of $M$ via $B\to\wh{\kappa(x)^\alg}=:\ol{K}$. Since all slopes of $M\otimes \ol{K}\ancon[r]$ are assumed to be integers, Theorem~\ref{Thm1a} provides us with a matrix $W\in\GL_n\bigl(\ol{K}\langle\frac{z}{\zeta^r},\frac{\zeta^{qr}}{z}\rangle\bigr)$ 
such that $D:=W^{-1}AW^{\sigma}$ is a diagonal matrix whose diagonal entries are powers of $z$.
Put $h:=\bigl(\|D\|_{qr}\,\|D^{-1}\|_{qr}\bigr)^{1/(q-1)}\geq 1$. By Lemma~\ref{Lemma6.7'} we find an affinoid neighborhood $\CM(B')$ of $x$, a finite \'etale $B'$-algebra $C$ such that $\CM(C)$ contains exactly one point $x'$ above $x$,
and a matrix $V\in M_n\bigl(C\langle\frac{z}{\zeta^r},\frac{\zeta^{qr}}{z}\rangle\bigr)$ with $\|V-W\|_{x',r'}<\|W^{-1}\|_{x',r'}^{-1}\,h^{-1}$ for all $r'\in[r,qr]$. Indeed it suffices to approximate finitely many coefficients 
in the Laurent series expansion of $W$. So we successively apply Lemma~\ref{Lemma6.7'} to each of their matrix entries and use Lemma~\ref{Lemma6.7''} whenever we need to further shrink $\CM(C)$.

We claim that after shrinking $\CM(B')$ the matrix $V$ belongs to $\GL_n\bigl(C\langle\frac{z}{\zeta^r},\frac{\zeta^{qr}}{z}\rangle\bigr)$.
Indeed, since $\|VW^{-1}-\Id_n\|_{x',r'}<1$ for all $r'\in[r,qr]$ we have $V\in\GL_n\bigl(\kappa(x')\langle\frac{z}{\zeta^r},\frac{\zeta^{qr}}{z}\rangle\bigr)$ at the point $x'$. By \cite[Lemma 9.7.1/1]{BGR} this is equivalent to the existence of an integer $m$ such that
$z^m\det V=\sum_{i\in\BZ}a_i z^i$ satisfies 
\begin{equation} \label{EqLemma6.3}
|a_i\zeta^{ri}|_{x'}\es<\es|a_0|_{x'} \quad\text{for all }i>0\qquad\text{and}\qquad |a_i\zeta^{qri}|_{x'}\es<\es|a_0|_{x'}\quad\text{for all }i<0\,.
\end{equation}
Clearly the equations (\ref{EqLemma6.3}) are satisfied on a whole affinoid neighborhood $\CM(C')$ of $x'$. Using Lemma~\ref{Lemma6.7''} the claim follows.

In $M_n\bigl(\ol{K}\langle\frac{z}{\zeta^{qr}},\frac{\zeta^{qr}}{z}\rangle\bigr)$ we write
\begin{eqnarray*}
V^{-1}AV^\sigma D^{-1} &=& (W^{-1}V)^{-1}D(W^{-1}V)^\sigma D^{-1} \\
&=& \bigl(W^{-1}(V-W)+\Id_n\bigr)^{-1} \,D\, \bigl(W^{-1}(V-W)+\Id_n\bigr)^\sigma D^{-1} \\
&=& \bigl(W^{-1}(V-W)+\Id_n\bigr)^{-1} \,D\, \bigl(W^{-1}(V-W)\bigr)^\sigma D^{-1} +\bigl(W^{-1}(V-W) +\Id_n\bigr)^{-1}\,.
\end{eqnarray*}
Using Lemma~\ref{Lemma6.0} we estimate
\begin{eqnarray*}
\bigl\|\bigl(W^{-1}(V-W)\bigr)^\sigma\bigr\|_{x',qr} &=&\|W^{-1}(V-W)\|_{x',r}^q\;<\; h^{-q}\,,\quad\text{and}\\[1mm]
\|\bigl(W^{-1}(V-W)+\Id_n\bigr)^{-1}-\Id_n\|_{x',qr}&\leq&\|W^{-1}(V-W)\|_{x',qr}\es<\es h^{-1}\,.
\end{eqnarray*}
It follows that $\|V^{-1}AV^\sigma D^{-1}-\Id_n\|_{x',qr}\;<\;h^{-1}$ at the point $x'$. Shrinking $\CM(B')$ further we may assume that
$\|V^{-1}AV^\sigma D^{-1}-\Id_n\|_{qr}\;<\;h^{-1}$ holds on all of $\CM(C)$. 
Then Lemma~\ref{Lemma6.1} 
yields a matrix $U\in\GL_n\bigl(C\langle\frac{\zeta^{qr}}{z}\rangle\bigr)$ with $\|U-\Id_n\|_{qr}<1$ for which 
$\wt A:=(VU)^{-1}A(VU)^\sigma$ satisfies $\wt AD^{-1}-\Id_n\in M_n\bigl(zC\langle\frac{z}{\zeta^{qr}}\rangle\bigr)$
and $\|\wt AD^{-1}-\Id_n\|_{qr}\;<\;1$. In particular $\wt AD^{-1}\in\GL_n\bigl(C\langle\frac{z}{\zeta^{qr}}\rangle\bigr)$ and $\wt A\in\GL_n\bigl(C\con[qr]\bigr)$. Consider the $\sigma$-module $M'$ over $C\con[r]$ with basis $v_1,\ldots,v_n$ on which Frobenius acts by the matrix $\wt A$. Then by what we have just shown
$M'\otimes C\langle\frac{z}{\zeta^r},\frac{\zeta^{qr}}{z}\rangle\cong M\otimes C\langle\frac{z}{\zeta^r},\frac{\zeta^{qr}}{z}\rangle$, hence $M'\otimes C\ancon[r]\cong M\otimes C\ancon[r]$ by Proposition~\ref{Prop0.1}.
By Lemma~\ref{Lemma14} for any point $y\in\CM(C)$ the generic HN-slopes of $M'\otimes\wh{\kappa(y)^\alg}\con[r]$ are the exponents of $z$ in the diagonal entries of $D$. By construction these are also the special absolute HN-slopes of $M'\otimes\wh{\kappa(x')^\alg}\con[r]$.
(Compare \cite[6.2.2]{Kedlaya}.)
\end{proof}

From the statement $\wt AD^{-1}\in\GL_n\bigl(C\langle\frac{z}{\zeta^{qr}}\rangle\bigr)$ in the proof we may also read off the following corollary.

\begin{corollary} \label{Cor6.3'}
Consider the situation of Proposition~\ref{Prop6.3}.
\begin{enumerate}
\item
If the absolute HN-slopes of $M\otimes\kappa(x)\ancon[r]$ are all non-negative then the $\sigma^s$-module $M'$ contains an $F^s$-stable $C\langle\frac{z}{\zeta^r}\rangle$-lattice $N$.
\item
If the absolute HN-slopes of $M\otimes\kappa(x)\ancon[r]$ are all equal to $\frac{d}{s}$ then the $\sigma^s$-module $M'$ contains a $C\langle\frac{z}{\zeta^r}\rangle$-lattice $N$
on which $z^{-d}F^s:(\sigma^s)^\ast N \to N$ is an isomorphism.
\qed
\end{enumerate}
\end{corollary}

In the situation where $B$ equals $L$ Proposition~\ref{Prop6.3} takes the following form.

\begin{corollary} \label{Cor6.3''}
Let $M$ be a $\sigma$-module over $L\ancon$. Then there exists a finite separable extension $L'$ of $L$, a positive integer $s$, and a $\sigma^s$-module $M'$ over $L'\con$ such that
$M'\otimes L'\ancon\cong M\otimes L'\ancon$ as $\sigma^s$-modules and the generic and special absolute HN-polygon of $M'$ coincide.
\qed
\end{corollary}

If $M$ is isoclinic we can even say more. For our definition of affinoid covering see Definition~\ref{DefBerkovichSpaces}.

\begin{theorem} \label{Thm6.12}
Let $M$ be a $\sigma$-module over $B\ancon$ which is isoclinic of slope $\frac{d}{s}$ with $d$ and $s$ relatively prime.
Then there exists a $\sigma$-module $\wt M$ over $B\con$, isoclinic of slope
$\frac{d}{s}$ with $\wt M\otimes B\ancon\cong M$. It is unique up to canonical isomorphism.

Moreover there exists an affinoid covering $\{\CM(B_i)\}_i$ of $\CM(B)$ by connected $\CM(B_i)$, and for each $i$ a finite \'etale $B_i$-algebra $C_i$ and a $C_i\langle\frac{z}{\zeta}\rangle$-lattice $\wt N_i$ in $\wt M\otimes C_i\con$, on which $z^{-d}F^s:(\sigma^s)^\ast \wt N_i\to \wt N_i$ is an isomorphism.
\end{theorem}

\begin{proof}
From Proposition~\ref{Prop6.3} we obtain an affinoid covering $\{\CM(B_i)\}_i$ of $\CM(B)$ with $\CM(B_i)$ connected,
and for each $i$ an isoclinic $\sigma^s$-module $M'_i$ over $C_i\con$ for a finite \'etale $B_i$-algebra $C_i$, and an isomorphism $\alpha_i:M'_i\otimes C_i\ancon\isoto M\otimes C_i\ancon$ of $\sigma^s$-modules. 
By Corollary~\ref{Cor6.3'} each $M'_i$ contains a $C_i\langle\frac{z}{\zeta}\rangle$-lattice $N'_i$ on which $z^{-d}F^s$ is an isomorphism.

We want to descend the $M'_i$ to a $\sigma$-module over $B\con$. Since $\CM(B)$ is compact there is a finite set $I$ such that the $\CM(B_i)$ for $i\in I$ cover $\CM(B)$. We let $C:=\bigoplus_{i\in I}C_i$ and we view $M':=\bigoplus_{i\in I}M'_i$ as a $\sigma^s$-module over $C\con$ containing the lattice $N':=\bigoplus_{i\in I}N'_i$. We set $C'':=C\otimes_B C$ and $C''':=C\otimes_B C\otimes_B C$. We let $p^\ast_\nu:C\to C''$ for $\nu=1,2$ be the homomorphism into the $\nu$-th factor and $q_{\mu\nu}^\ast:C''\to C'''$ for $\mu\nu\in\{12,13,23\}$ be the homomorphism into the $\mu$-th and $\nu$-th factor. Consider the $\sigma^s$-module $M'':=\CHom(p^\ast_1 M',p^\ast_2 M')$ over $C''\con$. The isomorphism $\alpha:=\bigoplus_{i\in I}\alpha_i:M'\otimes C\ancon\isoto M\otimes C\ancon$ yields an element $f:=p^\ast_2\alpha^{-1}\circ p^\ast_1\alpha$ of $\bigl(M''\otimes C''\ancon\bigr)^{F^s}(C'')$. Since $M''$ contains the $F^s$-stable $C''\langle\frac{z}{\zeta}\rangle$-lattice $\CHom(p^\ast_1 N',p^\ast_2 N')$, this element $f$ comes from an element $f\in(M'')^{F^s}(C'')$ by Proposition~\ref{Prop13}. The latter induces an isomorphism of $\sigma^s$-modules $p^\ast_1 M'\isoto p^\ast_2 M'$ over $C''\con$. Over $C'''\con$ this isomorphism satisfies the cocycle condition $q^\ast_{13}f=q^\ast_{23}f\circ q^\ast_{12}f$ by Proposition~\ref{Prop13} because this condition holds over $C'''\ancon$.

Now $B\con\to C\con$ is faithfully flat and quasi-compact. So Grothendieck's descent theory \cite[Theorem 6.1/4]{BLR} says that $M'$ comes from a $B\con$-module $\wt M$ together with an isomorphism $\beta:\wt M\otimes C\con\isoto M'$ and $f=p^\ast_2\beta\circ p^\ast_1\beta^{-1}$. Moreover $\wt M$ is a locally free $B\con$-module by \cite[II, Proposition 2.5.2]{EGA}. The isomorphism $\alpha\beta:\wt M\otimes C\ancon\isoto M\otimes C\ancon$ satisfies $p^\ast_1(\alpha\beta)=p^\ast_2\alpha\circ f\circ p^\ast_1\beta=p^\ast_2(\alpha\beta)$ over $C''\ancon$, hence comes from an isomorphism $\alpha\beta:\wt M\otimes B\ancon\isoto M$ by \cite[Proposition 6.1/1]{BLR}. To make $\wt M$ into a $\sigma$-module, consider the $\sigma^s$-module $H':=\CHom(\sigma^\ast M',M')$ over $C\con$. Then $H'\otimes C\ancon\isoto \CHom(\sigma^\ast M,M)\otimes C\ancon\,,\,h\mapsto \alpha\circ h\circ\sigma^\ast\alpha^{-1}$ is an isomorphism of $\sigma^s$-modules over $C\ancon$. The Frobenius $F$ on $M$ induces an element $F':=\alpha^{-1}\circ F\circ\sigma^\ast\alpha$ of $\bigl(H'\otimes C\ancon\bigr)^{F^s}(C)$. Since $H'$ contains the $F^s$-stable $C\langle\frac{z}{\zeta}\rangle$-lattice $\CHom(\sigma^\ast N',N')$, the element $F'$ lies in $(H')^{F^s}(C)$ by Proposition~\ref{Prop13} and gives an isomorphism $F':\sigma^\ast M'\isoto M'$. It satisfies 
\[
p^\ast_1 F'\es=\es p^\ast_1\alpha^{-1}\circ p^\ast_1 F\circ \sigma^\ast p^\ast_1\alpha\es=\es f^{-1}\circ p^\ast_2\alpha^{-1}\circ p^\ast_2 F\circ \sigma^\ast p^\ast_2\alpha\circ \sigma^\ast f\es=\es f^{-1}\circ p^\ast_2 F'\circ \sigma^\ast f
\]
in $\bigr(H'\otimes C''\ancon\bigr)^{F^s}(C'')$ and hence also in $\bigr(H'\otimes C''\con\bigr)^{F^s}(C'')$ by Proposition~\ref{Prop13}. Thus $F'$ descends to an isomorphism $\wt F:=\beta^{-1}\circ F'\circ\sigma^\ast\beta:\sigma^\ast\wt M\isoto M$ over $B\con$. This makes $(\wt M,\wt F)$ into a $\sigma$-module over $B\con$. By construction $\alpha\beta:\wt M\otimes B\ancon\isoto M$ is an isomorphism of $\sigma$-modules.
This proves the theorem. (Compare \cite[6.3.3]{Kedlaya}.)
\end{proof}

If $B=L$ is a field the theorem takes the following form:

\begin{corollary} \label{Cor6.12b}
Let $M$ be a $\sigma$-module over $L\ancon$ which is isoclinic of slope $\frac{d}{s}$ with $d$ and $s$ relatively prime. Then there exists a unique $\sigma$-module $M'$ over $L\con$, isoclinic of slope $\frac{d}{s}$ with $M'\otimes L\ancon\cong M$. Moreover $M'$ contains an $L\langle\frac{z}{\zeta}\rangle$-lattice on which $z^{-d}F^s$ is an isomorphism.
\qed
\end{corollary}

\begin{theorem} (Slope Filtration Theorem) \label{Thm6.13} \\
Let $M$ be a $\sigma$-module over $L\ancon$. Then there exists a unique filtration $0=M_0\subset M_1\subset\ldots\subset M_\ell=M$ of $M$ by saturated $\sigma$-submodules with the following properties:
\begin{enumerate}
\item
For $i=1,\ldots,\ell$ the quotient $M_i/M_{i-1}$ is isoclinic of some slope $s_i$ (Definition~\ref{DefIsoclinic}),
\item
$s_1<\ldots<s_\ell$.
\end{enumerate}
This filtration is a (the) Harder-Narasimhan filtration of $M$.
\end{theorem}

\begin{proof}
Since isoclinic $\sigma$-modules are semistable by Lemma~\ref{Lemma1.5.6'}, any filtration as in (a) and (b) is a Harder-Narasimhan filtration. In particular, by Lemma~\ref{Lemma1.5.4'} the filtration is unique if it exists.

To prove the existence, it suffices to show that the HN-filtration of $M\otimes \ol{L}\ancon$ from Theorem~\ref{Thm1} descends to $L\ancon$. By Corollary~\ref{Cor6.3''} there exists a finite separable extension $L'$ of $L$,
a positive integer $s$, and a $\sigma^s$-module $N$ over $L'\con$ such that $N\otimes L'\ancon\cong M\otimes L'\ancon$ as $\sigma^s$-modules, and such that the generic and special absolute HN-polygons of $N$ coincide. Clearly we may also assume that $L'$ is Galois over $L$.
By Theorem~\ref{Thm9} the HN-filtration of $M\otimes \ol{L}\ancon$ descends to a filtration of $N$ by $\sigma^s$-submodules, which in turn yields a filtration $0=M'_0\subset M'_1\subset\ldots\subset M'_\ell=M'$ of $M':=M\otimes L'\ancon$ by $\sigma^s$-submodules. 
Clearly the $M'_i$ are even $\sigma$-submodules since they are stable under $F$ after tensoring up to $\ol L\ancon$.

To prove that the filtration of $M'$ by the $M'_i$ further descends to $L\ancon$, 
it suffices to show that $M'_1$ descends to $L\ancon$. By Lemma~\ref{LemmaKe3.6.2} we only need to treat the case where $\rk M'_1=1$. In this case $M'_1$ is generated by a vector $v'\in M\otimes L'\ancon$. Let $e_1,\ldots,e_n$ be an $L\ancon$-basis of $M$.
Then $v'=a'_1 e_1+\ldots+a'_n e_n$ for $a'_i\in L'\ancon$ and we may assume $a'_1\ne0$. For each $g\in\Gal(L'/L)$, ${}^{g\!}v'$ is a scalar multiple of $v'$ since ${}^{g\!}v'\in M'_1$. Thus $\frac{a'_i}{a'_1}\in\Quot L'\ancon$ is $\Gal(L'/L)$-invariant for all $i$.
From Lemma~\ref{Lemma0.3} it follows that $\frac{a'_i}{a'_1}\in\Quot L\ancon$, hence $\frac{a'_i}{a'_1}=\frac{a_i}{a_1}$ for $a_i\in L\ancon$. Set $v=a_1e_1+\ldots a_ne_n\in M$ and let $M_1$ be the saturation in $M$ of the submodule generated by $v$. Them $M_1$ is a $\sigma$-submodule of $M$ over $L\ancon$ and $M_1\otimes L'\ancon$ is a saturated $\sigma$-submodule of $M'_1$ over $L'\ancon$, whence $M_1\otimes L'\ancon=M'_1$. Thus $M'_1$ descends to $L\ancon$ as claimed.
(Compare \cite[6.4.1]{Kedlaya}.)
\end{proof}

\begin{corollary} \label{Cor6.14}
For any complete extension $L'$ of $L$ and any $\sigma$-module $M$ over $L\ancon$, the HN-filtration of $M\otimes L'\ancon$ coincides with the result of tensoring the HN-filtration of $M$ up to $L'\ancon$.
\end{corollary}

\begin{proof}
The characterization of the HN-filtration given in Theorem~\ref{Thm6.13} is stable under base change.
\end{proof}

\begin{corollary}\label{Cor1.7.9}
Let $X$ be a Berkovich space (Definition~\ref{DefBerkovichSpaces}) and let $M$ be a $\sigma$-module over $\CO_X\ancon$. Let $P$ be a polygon in the sense of Definition~\ref{DefSlopesAndPolygons}. Then the set of points $x\in X$ such that the HN-Polygon of $M\otimes\kappa(x)\ancon$ lies above $P$ is open in $X$.
\end{corollary}

\begin{proof}
Let $x\in X$ be a point such that the HN-Polygon of $M\otimes\kappa(x)\ancon$ lies above $P$. By Proposition~\ref{Prop6.3} there exists a connected affinoid neighborhood $\CM(B')$ of $x$, a finite \'etale $B'$-algebra $C$, a lift of $x$ to a point of $\CM(C)$, a positive integer $s$, and a $\sigma^s$-module $M'$ over $C\con[r]$ such that $M'\otimes C\ancon[r]\cong M\otimes C\ancon[r]$ as $\sigma^s$-modules, and such that for all analytic points $y\in\CM(C)$ the generic HN-polygon of $M'\otimes\kappa(y)\con[r]$ coincides with the special absolute HN-polygon of the $\sigma^s$-module $M\otimes\kappa(x)\ancon[r]$. By Proposition~\ref{Prop8} this implies that the special HN-Polygon of the $\sigma^s$-module $M\otimes\kappa(y)\ancon$ lies above the one of $M\otimes\kappa(x)\ancon$. Hence also the special HN-Polygon of the $\sigma$-module $M\otimes\kappa(y)\ancon$ lies above $P$. By Corollary~\ref{Cor6.14} the same is then true for every analytic point $y$ in the neighborhood $\CM(B')$ of $x$.
\end{proof}

\vspace{2cm}%\pagebreak

%%%%%%%%%%%%%%%%%%%%%%%%%%%%%%%%%%%%%%%%%%%%%%%%%%%%%%%%%%%%%%%%%%%%%%
%
%    Hodge-Pink Theory
%
%%%%%%%%%%%%%%%%%%%%%%%%%%%%%%%%%%%%%%%%%%%%%%%%%%%%%%%%%%%%%%%%%%%%%%

\section{Hodge-Pink Theory} \label{ChaptHPTheory}
\setcounter{equation}{0}

In this chapter we develop the equal characteristic analog of $p$-adic Hodge theory. It is based on Pink's invention of Hodge theory over function fields \cite{Pink} and the concept of local shtukas which replaces the $p$-adic Galois representations from $p$-adic Hodge theory. In particular, we study in Section~\ref{SectMysterious} the analog of the mysterious functor relating crystalline Galois representations and filtered isocrystals. We define in Section~\ref{SectHPStructures} the analogs of (weakly) admissible filtered isocrystals and show in Section~\ref{SectWA=>A} that ``weakly admissible implies admissible'' if the (additively written) value group of the base field does not contain a non-zero element which is arbitrarily often divisible by $q$. The actual hypothesis on the base field is even weaker. As a universal tool we apply a criterion for (weak) admissibility in the spirit of Berger's~\cite{Berger1} criterion in $p$-adic Hodge theory using $p$-adic differential equations. We prove our criterion in Section~\ref{SectCriteriaWA}. The Slope Filtration Theorem and the other results from Chapter~\ref{ChaptSFT} enter at various places. We begin with an explanation of local shtukas.

%%%%%%%%%%%%%%%%%%%%%%%%%%%%%%%%%%%%%%%%%%%%%%%%%%%%%%%%%%%%%%%%%%%%%%
%
%    Local Shtukas
%
%%%%%%%%%%%%%%%%%%%%%%%%%%%%%%%%%%%%%%%%%%%%%%%%%%%%%%%%%%%%%%%%%%%%%%

\subsection{Local Shtukas} \label{SectLocalShtuka}

We retain the notation introduced in Section~\ref{SectNotation}.
For the notion of $\zeta$-adic formal scheme over $\Spf R$ we refer the reader to \cite[I$_{\rm new}$ \S10]{EGA}. In this article we are particularly interested in $\zeta$-adic formal schemes of the following two types:
\begin{itemize}
\item
schemes over $\Spec R$ on which $\zeta$ is locally nilpotent,
\item
admissible formal schemes over $\Spf R$ in the sense of Raynaud; see Appendix~\refAppRigFormal{}.
\end{itemize}

Let $X$ be a $\zeta$-adic formal scheme over $\Spf R$. We let $\CO_X\dbl z\dbr$ be the sheaf on $X$ of formal power series in $z$.
And we denote by $\sigma$ the endomorphism of $\CO_X\dbl z\dbr$ that acts as the identity on $z$ and as $b\mapsto b^q$ on local sections $b\in \CO_X$. For a sheaf $M$ of $\CO_X\dbl z\dbr$-modules on $X$ we abbreviate the sheaf $M\otimes_{\CO_X\dbl z\dbr,\sigma}\CO_X\dbl z\dbr$ by $\sigma^\ast M$.

\begin{definition} \label{Def1.1}
A \emph{local shtuka} of rank $n$ over $X$ is a pair $(M,F_M)$ consisting of 
\begin{itemize}
\item 
a sheaf $M$ of $\CO_X\dbl z\dbr$-modules on $X$, which Zariski-locally on $X$ is a free $\CO_X\dbl z\dbr$-module of rank $n$, and 
\item 
an isomorphism $F_M:\sigma^\ast M[\frac{1}{z-\zeta}] \isoto M[\frac{1}{z-\zeta}]$
(where we abbreviate $M\otimes_{\CO_X\dbl z\dbr}\CO_X\dbl z\dbr[\frac{1}{z-\zeta}]$ to $M[\frac{1}{z-\zeta}]$).
\end{itemize}
Let $\mu=(\mu_1\ge\ldots\ge\mu_n)$ be a decreasing sequence of integers. We say that $(M,F_M)$ is \emph{bounded by $\mu$} if
\begin{equation}\label{EqBounded}
 F_M\bigl(\wedge^i\sigma^\ast  M\bigr)\;\subset\;(z-\zeta)^{\mu_{n-i+1}+\ldots+\mu_n}\cdot\wedge^i M\qquad\text{for }1\le i\le n\text{ with equality for } i=n\,. 
\end{equation}
A \emph{morphism} of local shtukas $f:(M,F_M)\to(M',F_{M'})$ over $X$ is a morphism of the underlying sheaves $f:M\to M'$ which satisfies $F_{M'}\circ\sigma^\ast f = f\circ F_M$. An \emph{isogeny} is a morphism $f:(M,F_M)\to(M',F_{M'})$ such that locally on $X$ there is a morphism $g:(M',F_{M'})\to(M,F_M)$ such that the compositions $f\circ g$ and $g\circ f$ equal multiplication with a power of $z$.
\end{definition}

Note that the analog of Hilbert's Theorem 90 for $\CO_X\dbl z\dbr$-modules \cite[Proposition 2.3]{HV1} justifies working with the Zariski topology in the above definition. The notion of boundedness was thoroughly studied in \cite[Lemma 4.3 and \S3]{HV1}. We mention the following two results.

\begin{lemma}\label{LemmaBoundedLocFree}
If $(M,F_M)$ is bounded by $\mu=(\mu_1\ge\ldots\ge\mu_n)$ then for any $e\ge-\mu_n$ the quotient $(z-\zeta)^{-e}M/F_M(\sigma^\ast M)$ is a locally free $\CO_X$-module of finite rank.
\end{lemma}

\begin{proof}
First note that $F_M:\sigma^\ast M\to (z-\zeta)^{-e}M$ is injective. Let $x:\Spec\kappa(x)\to X$ be a point of $X$. Since $X$ is a $\zeta$-adic formal scheme we have $\zeta=0$ in $\kappa(x)$. Set $\CC=(z-\zeta)^{-e}M/F_M(\sigma^\ast M)$ and consider the exact sequence
\[
0\;\longto\;\Tor_1^{\CO_X}(\kappa(x),\CC)\;\longto\;x^\ast\sigma^\ast  M\;\xrightarrow{\es x^\ast F_M\;}\;x^\ast (z-\zeta)^{-e}M\;\longto\;x^\ast\CC\;\longto\;0\,.
\]
As $M$ is locally free of rank $n$ we have $x^\ast\sigma^\ast M\cong \kappa(x)\dbl z\dbr^{\oplus n}\cong x^\ast (z-\zeta)^{-e}M$. The boundedness condition (\ref{EqBounded}) says for $i=n$ that $\det(x^\ast F_M)$ equals $z^{\mu_1+\ldots+\mu_n}$ times a unit. Therefore the map $x^\ast F_M$ is injective. So $\CC=(z-\zeta)^{-e}M/F_M(\sigma^\ast M)$ is finitely presented by construction and flat over $X$ by Nakayama's Lemma; e.g.\ \cite[Exercise~6.2]{Eisenbud}, hence locally free of finite rank.
\end{proof}

In view of the following proposition we will only encounter bounded local shtukas in this article.

\begin{proposition}
Let $X$ be a quasi-compact $\zeta$-adic formal scheme such that for every open subset $\Spf B\subset X$ the ring $B$ is an integral domain. Then every local shtuka over $X$ is bounded by some $\mu$.
\end{proposition}

\begin{proof}
We need to find suitable integers $\mu_1\ge\ldots\ge\mu_n$. Since $X$ is quasi-compact we are allowed to find them locally on $X$. Hence we may assume that $X=\Spf B$ is affine, $M\cong B\dbl z\dbr^{\oplus n}$, and $F_M$ is given by a matrix $\Phi\in\GL_n\bigl(B\dbl z\dbr[\frac{1}{z-\zeta}]\bigr)$. By choosing $\mu_2,\ldots,\mu_n$ small enough we can achieve that (\ref{EqBounded}) is satisfied for all $1\le i<n$. We claim that $\det\Phi\in(z-\zeta)^r B\dbl z\dbr\mal$ for a (positive or negative) integer $r$. Then setting $\mu_1=r-\mu_2-\ldots-\mu_n$, condition (\ref{EqBounded}) is also satisfied for $i=n$.

To prove the claim note that there are positive integers $s,t$ with $(z-\zeta)^s\Phi, (z-\zeta)^t\Phi^{-1}\in M_n\bigl(B\dbl z\dbr\bigr)$. Set $f=(z-\zeta)^{sn}\det\Phi$ and $g=(z-\zeta)^{tn}\det\Phi^{-1}$. Then $fg=(z-\zeta)^{(s+t)n}\in (z-\zeta)B\dbl z\dbr$. Since $B$ is an integral domain, $(z-\zeta)B\dbl z\dbr$ is a prime ideal. So $f\in (z-\zeta)^u B\dbl z\dbr\mal$ for a positive integer $u$ and $\det\Phi\in (z-\zeta)^{u-sn}B\dbl z\dbr\mal$ proving our claim.
\end{proof}

Before we proceed further let us mention the main sources from which local shtukas arise. 

\begin{example} \label{ExGlobalObjects}
Let $C$ be a smooth projective geometrically irreducible curve over $\Spec\BF_q$ and let $\infty\in C$ be a closed point. Put $A:=\Gamma(C\setminus\{\infty\},\CO_C)$. Let $v\in C\setminus\{\infty\}$ be another closed point considered as a maximal ideal $ v\subset A$. Let $\wh\CO_{C,v}$ be the completed local ring and let $R$ be a complete extension of $\wh{\CO}_{C,v}$. Denote by $c:A\to R$ the inclusion morphism and by $\CI$ the kernel of $c\otimes\id:A\otimes R\to R$. We denote a uniformizing parameter of $\wh\CO_{C,v}$ by $z$ and its image in $R$ by $\zeta$. Then $\wh\CO_{C,v}\cong\BF_{q^e}\dbl z\dbr$ by Cohen's structure theorem, where $\kappa(v)=\BF_{q^e}$. 
Let $L$ be the fraction field of $R$.
For a $\zeta$-adic formal $R$-scheme $X$ consider the completion of $C_X:=C\times_{\BF_q}X$ along the closed subscheme $\{v\}\times X$. If $v\in C(\BF_q)$, and hence $e=1$, its structure sheaf is isomorphic to $\CO_X\dbl z\dbr$. If $e\ne1$ the situation is slightly more complicated and similar to \cite[before Proposition 8.5]{BH1}. For simplicity we now assume that $e=1$.
\begin{enumerate}
\item 
Let $X=\Spf R$ and let $(E,\phi)$ be a \emph{Drinfeld-$A$-module} \cite{Drinfeld} or more generally an \emph{abelian $t$-module} \cite{Anderson} over $L$ of dimension $d$. Let $(M_L,F)$ be the \emph{$t$-motive} associated with $(E,\phi)$ (more precisely the Drinfeld-Anderson-$A$-motive, see \cite[\S5]{Potemine}). $M_L$ is a locally free $A\otimes_{\BF_q}L$-module of finite rank equipped with an injective morphism $F:\sigma^\ast M_L\to M_L$, where 
\[
\sigma^\ast M_L\es:=\es M_L\otimes_{A\otimes L,\id\otimes Frob_q}A\otimes_{\BF_q}L\,, 
\]
such that $\coker F$ is an $L$-vector space of dimension $d$ and annihilated by $\CI^d$.

Assume that $(M_L,F)$ has \emph{good reduction}, that is, there exists a locally free $A\otimes_{\BF_q}R$-module $M$ of finite rank with $M_L=M\otimes_R L$ such that $F$ restricted to $M$ gives a morphism \mbox{$F:\sigma^\ast M\to M$} whose cokernel is a free $R$-module of rank $d$ and annihilated by $\CI^d$. Then the completion $(M\otimes R\dbl z\dbr,F\otimes\id)$ of $(M,F)$ along $\{v\}\times\Spf R$ is a local shtuka over $\Spf R$ which is bounded by $(d,0,\ldots,0)$. It can be viewed as a kind of Dieudonn\'e module of the $ v$-divisible group
\[
\dirlim[n] E[ v^n] \qquad\text{where}\quad E[ v^n] \es:=\es\bigcap_{a\in v^n}\ker\phi_a\,.
\]
of $(E,\phi)$; see \cite[\S\S6 and 7]{HartlDict,Hartl}. Recall that abelian $t$-modules are function field analogs of abelian varieties, so the $\dirlim E[ v^n]$ are analogs of $p$-divisible groups.
In that sense we are developing here the local theory of good reduction of $t$-motives at places which ``divide the residue characteristic''.
\item 
Let $(\CF_i,\Pi_i,\tau_i)$ be an \emph{abelian $\tau$-sheaf} of characteristic $c$ and dimension $d$ over $X$, see \cite[\S1]{Hartl}. Then the completions $\CF_i\otimes_{\CO_{C_X}}\CO_X\dbl z\dbr$ of $\CF_i$ along $\{v\}\times X$ are all isomorphic via $\Pi_i$. Let us denote their common value by $M$. Let us further denote the common value of the $\tau_i\otimes\id_{\CO_X\dbl z\dbr}$ by $F$. Then $(M,F)$ is a local shtuka over $X$ bounded by $(d,0,\ldots,0)$. We may even allow $v=\infty$ here. Then the formation of $(M,F)$ is slightly more complicated and was described in \cite[Construction 7.13]{Hartl}.
\item 
Let 
\[
\xymatrix @C=1pc @R=1pc {
\CE \ar[r]^j & \CE' & & & & \CE \\
\sigma^\ast\CE \ar[ur]_t& & & \Bigl(\text{respectively} & \CE' \ar[ur]^t \ar[r]_j & *!L(0.5) 
\objectbox{\;\sigma^\ast\CE\quad\Bigr)}
}
\]
be a \emph{right} (respectively \emph{left}) \emph{shtuka} \cite{Drinfeld3} over $X$. Assume that $\coker t$ is supported on $\{v\}\times X$ and that the support of $\coker j$ is disjoint from $\{v\}\times X$. Let $M:=\CE\otimes_{\CO_{C_X}}\CO_X\dbl z\dbr$ be the completion of $\CE$ at $\{v\}\times X$ and let $F:=j^{-1}\circ t:\sigma^\ast M[\frac{1}{z-\zeta}]\isoto M[\frac{1}{z-\zeta}]$. Then $(M,F)$ is a local shtuka over $X$ bounded by $(1,0,\ldots,0)$. The sense in which $(M,F)$ is \emph{local} is with respect to the \emph{coefficients}: $M$ lives over $\wh\CO_{C,v}$ as opposed to over all of $C$. This example gave rise to the name local shtuka.
\end{enumerate}
\end{example}

As mentioned in Example (a) local shtukas take up the role in the function field case that is played by $p$-divisible groups in number theory; see \cite[\S3]{HartlDict}. The concept of local shtukas is actually even more general since it is not restricted to weights $0$ and $1$; see Remark~\ref{Remark7.3} below. In fact if $X=\Spec k$ is the spectrum of a field in which $\zeta=0$, local shtukas over $X$ are precisely the analogs of $F$-crystals. Studying local shtukas up to isogeny leads to the notion of local isoshtukas.

\begin{definition} \label{Def1.2}
Let $X$ be a scheme over $\Spf R$ on which $\zeta$ is locally nilpotent. A \emph{local isoshtuka} of rank $n$ over $X$ is a pair $(D,F_D)$ consisting of a sheaf $D$ of $\CO_X\dbl z\dbr[z^{-1}]$-modules on $X$ which Zariski locally on $X$ is a free $\CO_X\dbl z\dbr[z^{-1}]$-module of rank $n$,
and an isomorphism $F_D:\sigma^\ast D\isoto D$.

Local isoshtukas over a field in which $\zeta$ is zero are also called \emph{$z$-isocrystals} since they behave very much like $F$-isocrystals in mixed characteristics; see \cite[\S3]{HartlDict}.

A \emph{morphism} of local isoshtukas over $X$ is a morphism of the underlying sheaves which is compatible with the $F$'s as in Definition~\ref{Def1.1}.
\end{definition}

For a morphism $\pi:X'\to X$ of formal schemes over $\Spf R$ we define the \emph{pullback functor} $(M,F_M)\mapsto \bigl(M\otimes_{\CO_X\dbl z\dbr}\CO_{X'}\dbl z\dbr,F_M\otimes\id\bigr)$ from local shtukas over $X$ to local shtukas over $X'$ (and similarly for local isoshtukas).

Note that a local isoshtuka over $\Spf A$, with $\zeta$ nilpotent in $A$, is nothing but a $\sigma$-module over $A\dbl z\dbr[z^{-1}]$ as defined in Section~\ref{SectSigmaModules}. We view a local shtuka over $X$ as an $\CO_X\dbl z\dbr$-lattice in the associated local isoshtuka.
We define the \emph{tensor product} of two local shtukas $(M,F_M)$ and $(N,F_N)$ over $X$ as the local shtuka
\[
\bigl(M\otimes_{\CO_X\dbl z\dbr}N\,,\,F_M\otimes F_N\bigr)\,.
\]
A \emph{unit object} for the tensor product is the local shtuka $(\CO_X\dbl z\dbr,F=\sigma)$, where the notation $F=\sigma$ means $F:\sigma^\ast\CO_X\dbl z\dbr=\CO_X\dbl z\dbr\otimes_{\CO_X\dbl z\dbr,\sigma}\CO_X\dbl z\dbr\to\CO_X\dbl z\dbr\,,\,a\otimes b\mapsto \sigma(a)\cdot b$. Also there is a natural definition of \emph{internal Hom's}.  
In particular the \emph{dual} $(M\dual, F_{M\dual})$ of a local shtuka $(M,F_M)$ over $X$ is defined as the sheaf $M\dual=\CHom_{\CO_X\dbl z\dbr}\bigl(M,\CO_X\dbl z\dbr\bigr)$ together with
\[
\TS F_{M\dual}\;=\;(\,.\,\circ F_M^{-1}):\;\sigma^\ast M\dual[\frac{1}{z-\zeta}]\isoto M\dual[\frac{1}{z-\zeta}]\,.
\]
Similar definitions apply to the category of local isoshtukas over $X$ making the category of local shtukas (respectively isoshtukas) over $X$ into an $\BF_q\dbl z\dbr$-linear (respectively $\BF_q\dpl z\dpr$-linear) additive rigid tensor category. If $X$ is the spectrum of a field then the category of local isoshtukas over $X$ is abelian. In general the category of local (iso-)shtukas is an exact category in the sense of Quillen~\cite[\S2]{Quillen} if one calls a short sequence of local shtukas \emph{exact} when the underlying sequence of sheaves of $\CO_X\dbl z\dbr$-modules is exact (compare Definition~\ref{DefShortExactSeq}).
Concerning boundedness we will need the following

\begin{lemma}\label{LemmaTensorBounded}
If $\ulM=(M,F_M)$ and $\ulN=(N,F_N)$ are bounded local shtukas over $X$ then also $\ulM\otimes\ulN$ and $\ulM\dual$ are bounded.
\end{lemma}

\begin{proof}
Let $\ulM$ be bounded by $(\mu_1\ge\ldots\ge\mu_m)$ and $\ulN$ be bounded by $(\nu_1\ge\ldots\ge\nu_n)$. We may even assume that $\mu_m\le0$, since otherwise we replace $\mu_m$ by $0$ and $\mu_1$ by $\mu_1+\mu_m$. Likewise we may assume $\nu_n\le0$. Now set $d=\mu_1+\ldots+\mu_m$ and $e=\nu_1+\ldots+\nu_n$. Then $F_M(\sigma^\ast M)\subset(z-\zeta)^{\mu_m}M$ and $F_N(\sigma^\ast N)\subset(z-\zeta)^{\nu_n}N$. Moreover, $F_M(\wedge^m\sigma^\ast M)=(z-\zeta)^d\wedge^m M$ and $F_N(\wedge^n\sigma^\ast N)=(z-\zeta)^e\wedge^n N$. Hence 
\[
(F_M\otimes F_N)\bigl(\wedge^i\sigma^\ast(M\otimes N)\bigr)\;\subset\; (z-\zeta)^{i(\mu_m+\nu_n)}\wedge^i (M\otimes N)
\]
first for $i=1$ and therefore also for $1\le i\le mn$. For $i=mn$ we even have 
\[
(F_M\otimes F_N)\bigl(\wedge^{mn}\sigma^\ast(M\otimes N)\bigr)\;=\; (z-\zeta)^{d+e}\wedge^{mn}(M\otimes N)
\]
because $\wedge^{mn}(M\otimes N)=\wedge^mM\otimes\wedge^nN$. This shows that the local shtuka $\ulM\otimes\ulN$ is bounded by $\bigl(d+e-(mn-1)(\mu_m+\nu_n)\,\ge\,\mu_m+\nu_n\ge\ldots\ge\mu_m+\nu_n\bigr)$. Note that this tuple is in decreasing order by our assumption $\mu_m,\nu_n\le0$.

To treat the dual of $\ulM$ we claim that $\ulM\dual$ is bounded by $(-\mu_n\ge\ldots\ge-\mu_1)$. In order to see this we choose locally on $X$ a basis of $M$ and write $F_M$ as an $n\times n$-matrix. Then Cramer's rule (e.g.~\cite[III.8.6, Formulas (21) and (22)]{BourbakiAlgebra}) gives a formula for the entries of the matrix $\wedge^jF_M^{-1}$ in terms of $(\det F_M)^{-1}$ times the entries of $\wedge^{n-j}F_M$, and this proves our claim.
\end{proof}

\medskip

\begin{remark} \label{RemarkGalReps}
Local shtukas give rise to \emph{Galois representations} as follows. Let $X$ be a quasi-paracompact admissible formal scheme over $\Spf R$ and let $X_L$ be the rigid analytic space over $L$ associated with $X$; see Appendix~\refAppBerkovich{}. Let $\bar x$ be a geometric base point of $X_L$ and let $\pi_1^\et(X_L,\bar x)$ be the \'etale fundamental group (\ref{DefFundamentalGroup}).
 For example take $X=\Spf R$\,,\, $X_L=\Spm L$\,,\, $\CO_X=R$\,,\, \mbox{$\CO_{X_L}=L$}\,,\linebreak\mbox{$(\bar x\!\to\! X_L)=(\ol L\!\hookleftarrow\! L)$}\,,\, $\pi_1^\et(X_L,\bar x)=\Gal(L^\sep/L)$. To any local shtuka $M$ over $X$ consider the $\sigma$-module $\bigl(M\otimes_{\CO_X\dbl z\dbr}\CO_{X_L}\langle\frac{z}{\zeta}\rangle\,,\,F_M\otimes\id\bigr)$ over $\CO_{X_L}\langle\frac{z}{\zeta}\rangle$. By Proposition~\ref{Prop13}
\[
T_zM\es:=\es\Bigl(\bigl(M\otimes_{\CO_X\dbl z\dbr}\CO_{X_L}\dbl z\dbr/(z^m)\bigr)^F\Bigr)_{m\in\BN_0}
\]
is a local system of $\BF_q\dbl z\dbr$-lattices on $X_L$. We call $T_zM$ the \emph{Tate module} of $M$. Its stalk $(T_zM)_{\bar x}$ at $\bar x$ defines a representation of the \'etale fundamental group $\rho_M:\pi_1^\et(X_L,\bar x)\to\GL_n\bigl(\BF_q\dbl z\dbr\bigr)$; see Proposition~\ref{Prop2.13}.
For a topological group $G$ denote by $\ul{\rm Rep}_{\BF_q\dbl z\dbr}G$ the category of continuous representations in finite free $\BF_q\dbl z\dbr$-modules.
\end{remark}

\begin{proposition}\label{Prop7.1}
Let $X_L$ be connected. Then the functor $M\mapsto\rho_M$ from the category of local shtukas over $X$ to $\ul{\rm Rep}_{\BF_q\dbl z\dbr}\pi_1^\et(X_L,\bar x)$ is a faithful $\BF_q\dbl z\dbr$-linear exact tensor functor.
\end{proposition}

\begin{proof}
Note that our functor is the composite of the $\BF_q\dbl z\dbr$-linear tensor functors
\[
M\es\longmapsto\es M\otimes_{\CO_X\dbl z\dbr}\CO_{X_L}\dbl z\dbr\es\longmapsto\es T_zM\es\longmapsto\es\rho_M\,.
\]
The first is clearly faithful and exact and the second is an exact equivalence by Proposition~\ref{PropTateModuleExact}. The third is an equivalence by Proposition~\ref{Prop2.13} and exact by \cite[Lemma 3.3.1]{JP} and \cite[0$_{\rm III}$, Proposition 13.2.2]{EGA}.
\end{proof}

Regarding full faithfulness there is the following

\begin{proposition}\label{PropKim}
Let $R$ (and $L$) be discretely valued. Then the functor $M\mapsto\rho_M$ from the category of local shtukas over $\Spf R$ to ${\rm Rep}_{\BF_q\dbl z\dbr}\Gal(L^\sep/L)$ is fully faithful.
\end{proposition}

\begin{proof}
This was first proved by Anderson~\cite[\S4.5, Theorem 1]{Anderson2} in case $F$ is topologically nilpotent (that is $\im F^n\subset zM+\Fm_R M$ for $n\gg0$). The general case is due to Kim~\cite[Theorem 5.2.3]{Kim}.
\end{proof}

We have not investigated the question when $M\mapsto\rho_M$ is fully faithful in general. Obviously it is not full if $X=\Spf R$ for $L$ algebraically closed, since then the non-isomorphic local shtukas
\[
\BOne(n)\es:=\es\bigl(M=R\dbl z\dbr\,,\es F_M:\sigma^\ast M\to M\,,\, 1\,\mapsto \,(z-\zeta)^n\bigr)
\]
for $n\in\BZ$ have isomorphic $\rho_{\BOne(n)}$. See also Remark~\ref{RemLokShInsteadOfGalRep}.

%%%%%%%%%%%%%%%%%%%%%%%%%%%%%%%%%%%%%%%%%%%%%%%%%%%%%%%%%%%%%%%%%%%%%%
%
%    Hodge-Pink Structures
%
%%%%%%%%%%%%%%%%%%%%%%%%%%%%%%%%%%%%%%%%%%%%%%%%%%%%%%%%%%%%%%%%%%%%%%

\subsection{Hodge-Pink Structures} \label{SectHPStructures}

We define the analogs of Fontaine's filtered isocrystals \cite{Fontaine}. The role of isocrystals is played in equal characteristic by the $z$-isocrystals defined in the previous section (Definition~\ref{Def1.2}). The Hodge filtration is replaced by a Hodge-Pink structure.
The theory of mixed Hodge-structures over $\BF_q\dpl z\dpr$ was developed by Pink~\cite{Pink} in analogy with the classical theory of Hodge-structures over local fields of characteristic zero. In the case of ``good reduction'' which we focus on here, we make the following definition which is a variant of Pink's definition.
Let $\CO_K\supset R$ be a rank-$1$ valuation ring which is complete and separated with respect to the $\zeta$-adic topology. Let $K$ be its fraction field, $\Fm_K$ its maximal ideal, and $k$ its residue field. 
Denote by $K\dbl z-\zeta\dbr$ the ring of formal power series over $K$ in the ``variable'' $z-\zeta$ and by $K\dpl z-\zeta\dpr$ its fraction field. We will always assume that there is a \emph{fixed section} $k\hookrightarrow\CO_K$ of the residue map $\CO_K\to k$. It yields a homomorphism $k\dpl z\dpr \to K\dbl z-\zeta\dbr$ by mapping $z$ to $z=\zeta+(z-\zeta)$.

\begin{definition} \label{Def1.3}
A \emph{$z$-isocrystal with Hodge-Pink structure} over $K$ is a triple $\ul D=(D,F_D,\Fq_D)$ consisting of a $z$-isocrystal $(D,F_D)$ over $k$ and a $K\dbl z-\zeta\dbr$-lattice $\Fq_D$ of full rank inside $\sigma^\ast D\otimes_{k\dpl z\dpr}K\dpl z-\zeta\dpr$.

A \emph{morphism} between 
$z$-isocrystals with Hodge-Pink structure $(D,F_D,\Fq_D)$ and $(D',F_{D'},\Fq_{D'})$ is a morphism $f:D\to D'$ of $z$-isocrystals which satisfies $(\sigma^\ast f\otimes\id)(\Fq_D)\subset \Fq_{D'}$.
Such a morphism is called \emph{strict} if 
\[
(\sigma^\ast f\otimes\id)(\Fq_D)\es=\es\Fq_{D'}\cap f(\sigma^\ast D)\otimes_{k\dpl z\dpr}K\dpl z-\zeta\dpr\,.
\]
\end{definition}

In addition to $\Fq_D$ there always is the tautological lattice $\Fp_D=\sigma^\ast D\otimes_{k\dpl z\dpr}K\dbl z-\zeta\dbr$. We define the \emph{Hodge-Pink weights} of $(D,F_D,\Fq_D)$ as the elementary divisors of $\Fq_D$ relative to $\Fp_D$. More precisely since $K\dbl z-\zeta\dbr$ is a principal ideal domain, there is a suitable $K\dbl z-\zeta\dbr$-basis $v_1,\ldots,v_n$ of $\Fp_D$ such that the lattice $\Fq_D$ has $K\dbl z-\zeta\dbr$-basis $\{(z-\zeta)^{w_i}v_i\}_i$. Then the Hodge-Pink weights are the integers $w_1,\ldots,w_n$, which we usually assume ordered $w_1\geq\ldots\geq w_m$. Alternatively if $e$ is large enough such that $\Fq_D\subset(z-\zeta)^{-e}\Fp_D$ or $(z-\zeta)^e\Fp_D\subset \Fq_D$ then the Hodge-Pink weights are characterized by
\begin{eqnarray*}
(z-\zeta)^{-e}\Fp_D/\Fq_D&\cong&\bigoplus_{i=1}^n K\dbl z-\zeta\dbr/(z-\zeta)^{e+w_i}\,,\\
\text{or}\quad\Fq_D/(z-\zeta)^e\Fp_D&\cong&\bigoplus_{i=1}^n K\dbl z-\zeta\dbr/(z-\zeta)^{e-w_i}
\end{eqnarray*}

\medskip

The category of $z$-isocrystals with Hodge-Pink structure possesses a \emph{tensor product} when we set
\[
\ul D\otimes\ul D'\es=\es \bigl(D\otimes_{k\dpl z\dpr} D', F_D\otimes F_{D'}, \Fq_D\otimes_{K\dbl z-\zeta\dbr}\Fq_{D'}\bigr)\,,
\]
\emph{duals} when we set
\[
\ul D\dual\es=\es \bigl(D\dual,F_{D\dual}, \Hom_{K\dbl z-\zeta\dbr}(\Fq_D,K\dbl z-\zeta\dbr)\bigr)\,,
\]
\emph{internal Hom's}, and the \emph{unit object} $\bigl(k\dpl z\dpr,F=\sigma,\Fq=\Fp\bigr)$. The endomorphism ring of the unit object is $\BF_q\dpl z\dpr$. A \emph{subobject} of $(D,F_D,\Fq_D)$ is given by $(D',F_{D'},\Fq_{D'})$ where $D'$ is an $F_D$-stable $k\dpl z\dpr$-subspace of $D$, $F_{D'}=F_D|_{D'}$, and $\Fq_{D'}\subset\Fq_D\cap\sigma^\ast D'\otimes_{k\dpl z\dpr}K\dpl z-\zeta\dpr$. If this last inclusion is an equality we say that $(D',F_{D'},\Fq_{D'})$ is a \emph{strict subobject} of $(D,F_D,\Fq_D)$. Clearly every subobject is contained in a uniquely determined strict subobject with the same underlying $z$-isocrystal. 
Dually a \emph{factor object} of $(D,F_D,\Fq_D)$ is given by $(D'',F_{D''},\Fq_{D''})$ where $(D'',F_{D''})$ is an $F$-equivariant quotient of $(D,F_D)$ and $\Fq_{D''}$ contains the corresponding quotient of $\Fq_D$. It is called a \emph{strict factor object} if moreover $\Fq_{D''}$ equals the quotient of $\Fq_D$.

The category of $z$-isocrystals with Hodge-Pink structure is an $\BF_q\dpl z\dpr$-linear additive rigid tensor category. It possesses kernels and cokernels, but is not abelian since the kernel of a morphism always is a strict subobject. Nevertheless, it is an exact category in the sense of Quillen~\cite[\S2]{Quillen} when we make the following

\begin{definition}\label{DefExSeqHP}
A short sequence of $z$-isocrystals with Hodge-Pink structure 
\[
0\longto \ul D'\xrightarrow{\es f\;}\ul D\xrightarrow{\es g\;}\ul D''\longto0
\]
is called \emph{(strict) exact} if the underlying sequence of $z$-isocrystals is exact and $f$ and $g$ are strict morphisms, that is, $(\sigma^\ast g\otimes\id)(\Fq_D)=\Fq_{D''}$ and $f$ identifies $\Fq_{D'}$ with $\Fq_D\cap\sigma^\ast D'\otimes_{k\dpl z\dpr}K\dpl z-\zeta\dpr$.
\end{definition}

\begin{remark} \label{Remark7.3}
In equal characteristic $z$-isocrystals with Hodge-Pink structure are the analogs of filtered isocrystals. Indeed the lattice $\Fq_D$ determines a decreasing filtration $Fil^\bullet$ on the $K$-vector space 
\[
D_K\es:=\es\sigma^\ast D\otimes_{k\dpl z\dpr}K\dbl z-\zeta\dbr/(z-\zeta)\es=\es\sigma^\ast D\otimes_{k\dpl z\dpr,\,z\mapsto\zeta}K
\]
by putting
\[
Fil^i D_K\es:=\es\bigl(\Fp_D\cap(z-\zeta)^i\Fq_D\bigr)/\bigl((z-\zeta)\Fp_D\cap(z-\zeta)^i\Fq_D\bigr)\,.
\]
This is called the \emph{Hodge-Pink filtration} of $\ul D$. It jumps precisely at the negatives of the Hodge-Pink weights.

However there is a fundamental difference to $p$-adic Hodge theory in that the lattice $\Fq_D$ contains more information than just the Hodge-Pink filtration. Let us explain the reasons for this. First of all, Pink~\cite{Pink} observed that the good definition of weak admissibility (Definition~\ref{Def1.4} below; Pink calls it \emph{semistability}) cannot be based on the Hodge-Pink filtration $Fil^\bullet$ alone. Instead one needs the finer information of the lattice $\Fq$ in order that the tensor product of two weakly admissible $z$-isocrystals with Hodge-Pink structure is again weakly admissible (Theorem~\ref{Thm1.4} below).

The second reason comes from a comparison with Fontaine's theory of $p$-adic Galois representations. The $\BF_q\dpl z\dpr$-algebra $K\dpl z-\zeta\dpr$ plays the role of Fontaine's filtered $\BQ_p$-algebra $B_{dR}$ which is a discretely valued field with valuation ring $B_{dR}^+$, uniformizing parameter $t$, and residue field $\BC_p$. Let us recall the situation in $p$-adic Hodge theory. If $K$ is a finite unramified extension of $\BQ_p$ and $(D,F_D,Fil^\bullet D_K)$ is a filtered isocrystal then we can construct a $\Gal(K^\sep/K)$-stable $B_{dR}^+$-lattice $\Fq_D$ in $\Fp_D\otimes_{B_{dR}^+}B_{dR}$ where $\Fp_D:=D_K\otimes_K B_{dR}^+$, by setting
\[
\Fq_D\es:=\es Fil^0(D_K\otimes_K B_{dR})\,.
\]
This lattice conversely determines the Hodge filtration by
\begin{equation}\label{EqFil}
Fil^i D_K\es=\es\Bigl(\bigl(\Fp_D\cap t^i\Fq_D\bigr)/\bigl(t\,\Fp_D\cap t^i\Fq_D\bigr)\Bigr)^{\Gal(K^\sep/K)}
\end{equation}
because $(B_{dR})^{\Gal(K^\sep/K)}=K$. It is easy to see that the compositum $Fil^\bullet D_K\mapsto \Fq_D\mapsto Fil^\bullet D_K$ of these two functors is the identity. The reason that also $\Fq_D\mapsto Fil^\bullet D_K\mapsto \Fq_D$ is the identity on $\Gal(K^\sep/K)$-stable lattices, is the fact that $B_{dR}^+$ is a successive extension of the Galois modules $\BC_p(0)$ by $\BC_p(1),\BC_p(2),\ldots$, which are all pairwise non isomorphic.
In this 1-1-correspondence between filtrations and lattices therefore only $\Gal(K^\sep/K)$-stable, and not all possible $B_{dR}^+$-lattices can arise. For example if 
\[
D_K=K^2\,,\quad Fil^0D_K\;=\; D_K\;\supset\; Fil^1D_K\;=\; Fil^2D_K\;=\; K\cdot v\;\supset\; Fil^3D_K\;=\;0
\]
for a vector $v\in D_K\setminus\{0\}$ then
\[
\Fq_D\es=\es \Fp_D\;+\;t^{-2}B_{dR}^+\cdot v  
\]
and the lattice $\Fq_D=\Fp_D+t^{-2}B_{dR}^{\SC +}\cdot (v+tv')$ can not occur if $v'$ is not a scalar multiple of $v$. 

On the other hand, in the theory of Hodge-Pink structures there is no such restriction. Namely, in equal characteristic, if $\BC$ denotes the completion of an algebraic closure of $K$, the $\Gal(K^\sep/K)$-modules $\BC(i)$ are all isomorphic as was observed by Anderson~\cite{Anderson2} and due to this pathology, Galois stability yields no restriction on $\Fq_D$ in equal characteristic. Moreover in equal characteristic the fixed field of $\Gal(K^\sep/K)$ inside $\BC$ is much larger than $K$. Namely by the Ax-Sen-Tate Theorem~\cite{Ax} it equals the completion $\wh{K^\perf}$ of the perfection of $K$. Hence, formula (\ref{EqFil}) does not produce a filtration of $D_K$ but only of $D_{\wh{K^\perf}}$. And for this reason we require that $\Fq_D$ is a lattice in $\sigma^\ast D\otimes_{k\dpl z\dpr}K\dpl z-\zeta\dpr$ instead of a $\Gal(K^\sep/K)$-stable lattice in $\sigma^\ast D\otimes_{k\dpl z\dpr}\BC\dpl z-\zeta\dpr$.
\end{remark}
\medskip

Let $\ul D=(D,F_D,\Fq_D)$ be a $z$-isocrystal of rank $n$ with Hodge-Pink structure over $K$. After choosing a basis of $D$, we let $\Phi_D\in\GL_n\bigl(k\dpl z\dpr\bigr)$ be the matrix by which $F_D$ acts on this basis. We define the \emph{Newton slope} $t_N(\ul D):=\ord_z(\det\Phi_D)$ as the valuation with respect to $z$ of $\det\Phi_D$. It does not depend on the chosen basis. 
On the other hand we define the \emph{Hodge slope} $t_H(\ul D)=-w_1-\ldots-w_n$ as the negative of the sum of the Hodge-Pink weights of $\ul D$. It also satisfies \mbox{$t_H(\ul D)=\sum_{i\in\BZ}i\cdot\dim_K Fil^i D_K/Fil^{i+1}D_K$} and $\wedge^n\Fq_D=(z-\zeta)^{-t_H(\ul D)}\wedge^n\Fp_D$.

\begin{definition} \label{Def1.4}
We say that a $z$-isocrystal with Hodge-Pink structure is \emph{weakly admissible} if $t_H(\ul D)=t_N(\ul D)$ and one of the following equivalent conditions holds:
\begin{enumerate}
\item 
$t_H(\ul D')\leq t_N(\ul D')$ for any (strict) subobject $\ul D'$ of $\ul D$,
\item 
$t_H(\ul D'')\geq t_N(\ul D'')$ for any (strict) factor object $\ul D''$ of $\ul D$.
\end{enumerate}
\end{definition}

\begin{proof}
The equivalence of (a) and (b) is standard, see for example \cite[Proposition 4.4]{Pink}.
\end{proof}

The arguments of \cite[\S\S4,5]{Pink} can easily  be adapted to our situation to yield the following result. Note that what we called weakly admissible here is called semistable in \cite{Pink}.

\begin{theorem} \label{Thm1.4}
The weakly admissible $z$-isocrystals with Hodge-Pink structure form a full subcategory of the category of $z$-isocrystals with Hodge-Pink structure, which has the following properties.
\begin{enumerate}
\item 
It is closed under the formation of tensor, symmetric, and exterior products, duals, extensions, kernels, and cokernels. Moreover it is abelian.
\item 
A direct sum of two $z$-isocrystals with Hodge-Pink structure is weakly admissible if and only if each summand is weakly admissible.
\item 
Any morphism between weakly admissible $z$-isocrystals with Hodge-Pink structure is automatically strict.
\qed
\end{enumerate}
\end{theorem}

\medskip

%%%%%%%%%%%%%%%%%%%%%%%%%%%%%%%%%%%%%%%%%%%%%%%%%%%%%%%%%%%%%%%%%%%%%%
%
%    The Mysterious Functor
%
%%%%%%%%%%%%%%%%%%%%%%%%%%%%%%%%%%%%%%%%%%%%%%%%%%%%%%%%%%%%%%%%%%%%%%

\subsection{The Mysterious Functor} \label{SectMysterious}

We come to the analog of Grothendieck's~\cite{Grothendieck} mysterious functor which was defined by Fontaine\ \cite{Fontaine2}. This analog was first constructed by A.\ Genestier and V.\ Lafforgue~\cite{GL} for discretely valued fields $K$. Let us recall its definition. For reasons why this functor should indeed be considered the analog of Fontaine's functor see Remark~\ref{RemLokShInsteadOfGalRep}. 

We keep the notation of the previous section.
Let $(M,F_M)$ be a local shtuka over $\Spf \CO_K$. It gives rise to a $z$-isocrystal $(D,F_D)=(M\otimes_{\CO_K\dbl z\dbr}k\dpl z\dpr, F_M\otimes\id)$ over $k$. 
Consider the  $\CO_K$-algebras
\begin{eqnarray*}\label{Page34}
\CO_K\dbl z,{\TS\frac{\zeta^r}{z}}\rangle[z^{-1}]&:=&\bigl\{\,\sum_{i=-\infty}^\infty b_iz^i:\es b_i\in\CO_K\,,\,|b_i|\,|\zeta|^{ri}\to0\;(i\to-\infty)\,\bigr\}\qquad\text{and}\\
\CO_K\dbl z,z^{-1}\}&:=&\bigl\{\,\sum_{i=-\infty}^\infty b_iz^i:\es b_i\in\CO_K\,,\,|b_i|\,|\zeta|^{ri}\to0\;(i\to-\infty) \text{ for all }r>0\,\bigr\}\,.
\end{eqnarray*}
The latter is contained in $K\overcon$.
We set $\DS \tminus:=\prod_{i\in\BN_0}(1-{\TS\frac{\zeta^{q^i}}{z}})\in\BF_q\dbl\zeta\dbr\dbl z,z^{-1}\}$. 
The following lemma is taken from \cite[Lemma 2.8 or Lemma 6.4]{GL}. For convenience of the reader we reproduce the proof from \cite{GL}. 

\begin{lemma} \label{LemmaGL}
If $K$ is discretely valued there is a unique functorial isomorphism
\[
\delta_M:\es M\otimes_{\CO_K\dbl z\dbr}\CO_K\dbl z,z^{-1}\}[\tminus^{-1}]\es\isoto\es D\otimes_{k\dpl z\dpr}\CO_K\dbl z,z^{-1}\}[\tminus^{-1}]
\]
which satisfies $\delta_M\circ F_M=F_D\circ\sigma^\ast\delta_M$ and which reduces to the identity modulo $\Fm_K$. 
\end{lemma}

\begin{proof}
Fix an $\CO_K\dbl z\dbr$-basis of $M$ and consider the induced $k\dpl z\dpr$-basis of $D$. Denote the matrices by which $F_M$ and $F_D$ act on these bases by $A\in\GL_n\bigl(\CO_K\dbl z\dbr[\frac{1}{z-\zeta}]\bigr)$ and $B=A\mod\Fm_K\in\GL_n\bigl(k\dpl z\dpr\bigr)$. 
To prove the existence of $\delta_M$ consider first the case where $A^{-1}\in M_n\bigl(\CO_K\dbl z\dbr\bigr)$.
Via the section $k\hookrightarrow\CO_K$ we view $B$ as an element of $\GL_n\bigl(\CO_K\dbl z\dbr[z^{-1}]\bigr)$ and we put
\[
C_m\es:=\es B\cdot B^\sigma\cdot\ldots\cdot B^{\sigma^m}\cdot (A^{\sigma^m})^{-1}\cdot\ldots\cdot (A^\sigma)^{-1}\cdot A^{-1} \quad\in\es M_n\bigl(\CO_K\dbl z\dbr[z^{-1}]\bigr)\,.
\]
Let $\pi$ be a uniformizer of $\CO_K$ and let $d\in\BN$ be a constant with $A\in M_n\bigl((z-\zeta)^{-d}\CO_K\dbl z\dbr\bigr)$. Then $B\in M_n\bigl(z^{-d}\CO_K\dbl z\dbr\bigr)$ and $A^{-1}-B^{-1}\in M_n\bigl(\pi\CO_K\dbl z\dbr\bigr)$. It follows that $C_m A=B C_{m-1}^\sigma$ and
\[
C_m-C_{m-1}\es=\es B\cdot\ldots\cdot B^{\sigma^m}\cdot(A^{-1}-B^{-1})^{\sigma^m}\cdot (A^{\sigma^{m-1}})^{-1}\cdot\ldots\cdot A^{-1} \quad\in\es M_n\bigl({\TS\frac{\pi^{q^m}}{z^{(m+1)d}}}\,\CO_K\dbl z\dbr\bigr)\,.
\]
Thus the sequence $C_m$ converges in $M_n\bigl(\CO_K\dbl z,z^{-1}\}\bigr)$ to a matrix $C$ which satisfies $CA=BC^\sigma$ and $C\equiv\Id_n\mod\Fm_K$. 

Since $(z-\zeta)^dA\in M_n\bigl(\CO_K\dbl z\dbr\bigr)$ and hence $(z^dB)^{-1}\in M_n\bigl(z^{-d}\CO_K\dbl z\dbr\bigr)$ the same reasoning shows that 
\[
\wt C_m\es:=\es (z-\zeta)^dA\cdot\ldots\cdot\bigl((z-\zeta)^dA\bigr)^{\sigma^m}\cdot(z^{-d}B^{-1})^{\sigma^m}\cdot\ldots\cdot(z^{-d}B^{-1})
\]
converges in $M_n\bigl(\CO_K\dbl z,z^{-1}\}\bigr)$ to a matrix $\wt C$ with $\wt C\,z^dB=(z-\zeta)^dA\wt C^\sigma$ and $\wt C\equiv\Id_n\mod\Fm_K$. In particular 
\[
C_m\wt C_m\es=\es\prod_{i=0}^m(1-{\TS\frac{\zeta^{q^i}}{z}})^d\cdot\Id_n\es=\es\wt C_m C_m
\]
converges to $C\wt C=\wt CC=\tminus^d\Id_n$. This shows that $C\in\GL_n\bigl(\CO_K\dbl z,z^{-1}\}[\tminus^{-1}]\bigr)$ and $C$ gives the desired isomorphism $\delta_M$.

It remains to show that $\delta_M$ is uniquely determined and functorial.
We establish both assertions simultaneously by proving the following claim: If $f:(M_1,F_{M_1})\to (M_2,F_{M_2})$ is a morphism of local shtukas over $\CO_K$ and 
\[
\delta_i:\es M_i\otimes_{\CO_K\dbl z\dbr}\CO_K\dbl z,z^{-1}\}[\tminus^{-1}]\es\isoto\es D_i\otimes_{k\dpl z\dpr}\CO_K\dbl z,z^{-1}\}[\tminus^{-1}]
\]
are isomorphisms with $\delta_i\circ F_{M_i}=F_{D_i}\circ\sigma^\ast\delta_i$ and $\delta_i\equiv\id\mod\Fm_K$ then $\delta_2\circ f=(f\mod\Fm_K)\circ\delta_1$.

Fix bases of $M_i$ and reduce them to bases of $D_i=M_i\otimes_{\CO_K\dbl z\dbr}k\dpl z\dpr$. Let $A_i,B_i$ and $C_i$ be the matrices by which $F_{M_i},F_{D_i}$ and $\delta_i$ act on these bases. Let $H$ be the matrix corresponding to $f$ and set $\ol H:=H\mod\Fm_K$. We must show that $C_2H=\ol H C_1$. By construction $\wt H:=C_2H-\ol H C_1$ satisfies $\wt H\in M_{n_2\times n_1}\bigl(\pi\CO_K\dbl z,\frac{\zeta^{1/q}}{z}\rangle[z^{-1}]\bigr)$ and $\wt H A_1=B_2\wt H^\sigma$. Since $\CO_K\dbl z\dbr[\frac{1}{z-\zeta}]\subset\CO_K\dbl z,\frac{\zeta^{1/q}}{z}\rangle[z^{-1}]$ we have $A_1\in \GL_{n_1}\bigl(\CO_K\dbl z,\frac{\zeta^{1/q}}{z}\rangle[z^{-1}]\bigr)$. So in fact $\wt H\in M_{n_2\times n_1}\bigl(\pi^q\CO_K\dbl z,\frac{\zeta^{1/q}}{z}\rangle[z^{-1}]\bigr)$ and iterating this argument shows that $\wt H=0$ as desired.
\end{proof}

\noindent 
{\it Remark.} Note that the analogous argument in mixed characteristic goes back to Dwork and is commonly called ``Dwork's trick''. A good account of it is de~Jong's proof~\cite[Lemma 6.3]{dJ98} of the triviality of a connection on a Frobenius-module on an open analytic disc over a $p$-adic field. Namely, in de Jong's situation the columns of the matrix we called $C$ form horizontal sections for the connection.

\bigskip

As a consequence of this lemma $\sigma^\ast\delta_M$ induces an isomorphism
\[
\sigma^\ast\delta_M:\;\sigma^\ast M\otimes_{\CO_K\dbl z\dbr}K\dbl z-\zeta\dbr\es\isoto\es\sigma^\ast D\otimes_{k\dpl z\dpr}K\dbl z-\zeta\dbr\es=\es\Fp_D\,.
\]
The isomorphism $F_M$ extends to an isomorphism 
\[
F_M:\es\sigma^\ast M\otimes_{\CO_K\dbl z\dbr}K\dpl z-\zeta\dpr\es\isoto\es M\otimes_{\CO_K\dbl z\dbr}K\dpl z-\zeta\dpr\,.
\]
We put $\Fq_D:=\sigma^\ast\delta_M\circ F_M^{-1}\bigl(M\otimes_{\CO_K\dbl z\dbr}K\dbl z-\zeta\dbr\bigr)\subset\Fp_D[\frac{1}{z-\zeta}]$.
This defines a functor $\BH$ which associates with the local shtuka $(M,F_M)$ over $\CO_K$ a $z$-isocrystal with Hodge-Pink structure $(D,F_D,\Fq_D)$ over $K$. By construction it transforms isogenies into isomorphisms and hence factors through the category of local shtukas over $\CO_K$ up to isogeny.

\bigskip

The author had initially hoped that the functor $\BH$ might also exist if $K$ is not discretely valued but this is not the case. The situation is as follows. Let $K$ be arbitrary and consider the category of \emph{rigidified local shtukas} over $\CO_K$, that is, triples $(M,F_M,\delta_M)$ where $(M,F_M)$ is a local shtuka over $\CO_K$ and 
\[
\delta_M:\es M\otimes_{\CO_K\dbl z\dbr}\CO_K\dbl z,z^{-1}\}[\tminus^{-1}]\es\isoto\es D\otimes_{k\dpl z\dpr}\CO_K\dbl z,z^{-1}\}[\tminus^{-1}]
\]
is an isomorphism with $\delta_M\circ F_M=F_D\circ\sigma^\ast\delta_M$ and $\delta_M\equiv\id\mod\Fm_K$, and where as above $(D,F_D)=(M,F_M)\mod\Fm_K$. A morphism in this category is a morphism $f:(M,F_M)\to (M',F_{M'})$ of local shtukas which satisfies $\delta_{M'}\circ f=(f\mod\Fm_K)\circ\delta_M$. On the category of rigidified local shtukas one defines the functor $\BH: (M,F_M,\delta_M)\mapsto (D,F_D,\Fq_D)$ exactly as above. This yields a diagram of functors
\[
\xymatrix { & **{!R(0.3) =<8pc,2pc>} \objectbox{\bigl\{\,\text{local shtukas over }\CO_K\text{ up to  isogeny}\,\bigr\}} &\\
\bigl\{\,\text{rigidified local shtukas over }\CO_K\text{ up to  isogeny}\,\bigr\}\ar[ur]_{\TS\alpha} \ar[dr]^{\TS\BH} & &\\
& **{!R(0.3) =<8pc,2pc>} \objectbox{\bigl\{\,\text{$z$-isocrystals with Hodge-Pink structure over }K\,\bigr\}} &
}
\]
where $\alpha$ is the forgetful functor. Lemma~\ref{LemmaGL} shows that $\alpha$ is an equivalence of categories if $K$ is discretely valued. We will further show the following.
\begin{itemize}
\item 
$\alpha$ is always faithful (this is obvious).
\item 
Fix $\lambda\in\BQ$. Then $\alpha$ restricted to the categories where the $\sigma$-module $M\otimes_{\CO_K\dbl z\dbr}K\{z,\frac{\zeta^{1/q}}{z}\rangle$ is isoclinic of slope $\lambda$ (Definition~\ref{DefIsoclinic}) is an equivalence of categories (Proposition~\ref{Prop4.1}).
\item 
If $K$ is discretely valued, $\alpha$ is an equivalence of categories (Genestier-Lafforgue, Lemma~\ref{LemmaGL}).
\item 
If $K$ is algebraically closed, $\alpha$ is essentially surjective but not full (Proposition~\ref{Prop4.2}).
\item 
Consider the following condition on the field $K$:
\begin{equation} \label{CondDoubleStar}
\text{\parbox{0.75\textwidth}{
Let $\wt K$ be the closure of the compositum $k^\alg K$ inside $\ol K$ and assume that $\wt K$ does not contain an element $a$ with $0<|a|<1$ such that all the $q$-power roots of $a$ also lie in $\wt K$.
}}\qquad\mbox{ }
\end{equation}
If $K$ satisfies (\ref{CondDoubleStar}) then $\alpha$ is full but need not be essentially surjective (Lemma~\ref{Lemma4.5} and Example~\ref{Ex4.3}).
\item 
The functors from the categories of the global objects from Example~\ref{ExGlobalObjects} to the category of local shtukas all factor canonically through the category of rigidified local shtukas (Remark~\ref{Remark4.4}).
\end{itemize}

The above condition (\ref{CondDoubleStar}) will also appear in Theorem~\ref{Thm3.5} as the minimal requirement on $K$ under which every weakly admissible $z$-isocrystal with Hodge-Pink structure arises from a rigidified local shtuka. So let us remark that (\ref{CondDoubleStar}) is satisfied whenever the value group of $K$ does not contain a non-zero element which is arbitrarily often divisible by $q$. This is due to the fact that the value groups of $K$ and $\wt K$ coincide. In particular, (\ref{CondDoubleStar}) is satisfied if $K$ is discretely valued, or if its value group is finitely generated.

However, before we prove these properties of the functor $\alpha$ let us study the functor $\BH$.

\begin{proposition}\label{PropTensorFunctor1}
The functor $\BH$ from the category of rigidified local shtukas over $\Spf\CO_K$ up to isogeny to the category of $z$-isocrystals with Hodge-Pink structure over $K$ is an $\BF_q\dpl z\dpr$-linear exact tensor functor.
\end{proposition}
We will see below that it is fully faithful (Proposition~\ref{PropHFaithfull}) and factors through the full subcategory of weakly admissible $z$-isocrystals with Hodge-Pink structure (Corollary~\ref{Cor2.9}).

\begin{proof}
One easily checks that the functor is $\BF_q\dpl z\dpr$-linear and compatible with tensor products and internal Hom's, thus a tensor functor. It remains to show that it preserves short exact sequences. So let $0\to M'\to M\to M''\to0$ be an exact sequence of rigidified local shtukas over $\Spf\CO_K$ and let $0\to \ul D'\xrightarrow{f}\ul D\xrightarrow{g}\ul D''\to0$ be the resulting sequence of $z$-isocrystals with Hodge-Pink structure over $K$. Clearly the underlying sequence of $z$-isocrystals over $k$ is exact. It remains to show that $f$ and $g$ are strict. This follows from the facts that $M\otimes_{\CO\dbl z\dbr}K\dbl z-\zeta\dbr\to M''\otimes_{\CO\dbl z\dbr}K\dbl z-\zeta\dbr$ is surjective and that
\[
M'\otimes_{\CO\dbl z\dbr}K\dbl z-\zeta\dbr\es=\es M\otimes_{\CO\dbl z\dbr}K\dbl z-\zeta\dbr\;\cap\; M'\otimes_{\CO\dbl z\dbr}K\dpl z-\zeta\dpr
\]
inside $M\otimes_{\CO\dbl z\dbr}K\dpl z-\zeta\dpr$.
\end{proof}

\begin{definition} \label{Def1.5}
A $z$-isocrystal with Hodge-Pink structure over $K$ is called \emph{admissible} if it lies in the essential image of the tensor functor $\BH$.
\end{definition}

As an immediate consequence of Proposition~\ref{PropTensorFunctor1} we obtain:

\begin{theorem} \label{Thm1.5b}
The category of admissible $z$-isocrystals with Hodge-Pink structure over $K$ is closed under the formation of tensor products, duals and direct sums.
\qed
\end{theorem}

\begin{example}\label{ExTateObjects}
As an example we define the \emph{Tate objects} $\BOne(n)$ for $n\in\BZ$ in the category of rigidified local shtukas over $\Spf\BF_q\dbl\zeta\dbr$ as the rigidified local shtuka
\[
\BOne(n)\es:=\es\bigl(M=\BF_q\dbl\zeta\dbr\dbl z\dbr\,,\es F_M:1\,\mapsto \,(z-\zeta)^n\,,\es\delta_M:1\mapsto \tminus^{-n}\bigr)\,.
\]
They give rise to the \emph{Tate objects} in the category of $z$-isocrystals with Hodge-Pink structure over $\BFZ$
\[
\BOne(n)\es:=\es\Bigl(D=\BF_q\dpl z\dpr\,,\,F_D=z^n\cdot\sigma\,,\,\Fq_D=(z-\zeta)^{-n}\Fp_D\Bigr)\,. 
\]
Clearly $\BOne(n)=\BOne(1)^{\otimes n}=\BOne(-n)\dual$.
\end{example}

\begin{remark}\label{RemLokShInsteadOfGalRep}
We will explain our reasons for viewing rigidified local shtukas as the appropriate analogs in equal characteristic for the crystalline Galois representations from $p$-adic Hodge theory. Let $K\supset\BFZ$ be a finite extension and let $\ol K$ be the completion of an algebraic closure of $K$. Let $\CO_{\ol K}$ be the ring of integers of $\ol K$. The first approach for an analog of Fontaine's functor would be to search for a tensor functor $\CF$ from a full subcategory of $\ul{\rm Rep}_{\BF_q\dpl z\dpr}\Gal(K^\sep/K)$ to the category of $z$-isocrystals with Hodge-Pink structure over $K$. However there are two objections to this plan.

First of all the fixed field in $\ol K$ for the continuous $\Gal(K^\sep/K)$-action induced by the density of $K^\sep$ in $\ol K$ is larger than $K$. Indeed by the Ax-Sen-Tate Theorem~\cite{Ax} it equals the closure of the perfection of $K$. Thus if one imitates Fontaine's construction of $B_{cris}$ and considers as a candidate for $\CF$ the functor
\[
\Bigl(\,\rho:\Gal(K^\sep/K)\to\GL(V)\,\Bigr) \es\longmapsto\es\bigl(V\otimes_{\BF_q\dpl z\dpr}B_{cris}\bigr)^{\Gal(K^\sep/K)}
\]
 in equal characteristic, $B_{cris}$ will be an $\CO_{\ol K}$-algebra and the functor does not yield finite dimensional $k\dpl z\dpr$-vector spaces.

Secondly rigidified local shtukas should be considered as analogs of the crystals associated with $p$-divisible groups. Hence if $M$ is a rigidified local shtuka over $\CO_K$ it is desirable that the $z$-isocrystal with Hodge-Pink structure obtained by this hypothetical functor $\CF$ from the Galois representation $\rho_M$ of Proposition~\ref{Prop7.1} is isomorphic to $\BH(M)$. So consider the commutative diagram of categories and functors
\[
\xymatrix {
\bigl\{\,\text{rigidified local shtukas over $\CO_K$ up to isogeny}\,\bigr\} \ar[rr]^{\TS \qquad\qquad\qquad M\mapsto\rho_M}\ar[dr]_{\TS\BH} && \es\CC \ar@{-->}[dl]^{\TS \CF} \\
&**{!R(0.3) =<8pc,2pc>} \objectbox{\bigl\{\,\text{weakly admissible $z$-isocrystals with Hodge-Pink structure over $K$}\,\bigr\}}
}
\]
in which the suitable full subcategory $\CC$ of $\ul{\rm Rep}_{\BF_q\dpl z\dpr}\Gal(K^\sep/K)$ and the functor $\CF$ are yet to be defined.

If $K$ satisfies condition (\ref{CondDoubleStar}) from page~\pageref{CondDoubleStar}, we will see in Theorem~\ref{Thm3.5} that $\BH$ is an equivalence of categories. So the existence of the functor $\CF$ together with its expected property to be fully faithful is equivalent to $M\mapsto\rho_M$ being an equivalence of categories. Thus we may work as well with rigidified local shtukas as with Galois representations.

On the other hand if $K$ is algebraically closed, the notion of Galois representation is vacuous, whereas rigidified local shtukas over $\CO_K$ still are interesting objects. Actually the Galois representations associated with the rigidified local shtukas $\BOne(n)$ over $\CO_K$ from 
Example~\ref{ExTateObjects}
are all isomorphic, whereas there are only the zero morphisms between the $z$-isocrystals with Hodge-Pink structure associated to $\BOne(m)$ and $\BOne(n)$ for $m\neq n$. Because $\BH$ is fully faithful the functor $\CF$ cannot exist in this case.

Since the category of rigidified local shtukas is large enough such that $\BH$ is an equivalence of categories if $K$ satisfies condition (\ref{CondDoubleStar}), we propose to view rigidified local shtukas as the appropriate analogs in Hodge-Pink theory of the crystalline Galois representations from $p$-adic Hodge theory. This point of view is further supported by the following three facts. Firstly one may even replace $\CO_K$ by an arbitrary admissible formal $\Spf R$-scheme $X$ and define the functor $\BH$ for rigidified local shtukas over $X$ with constant $z$-isocrystal (see the remark after Proposition~\ref{PropHFaithfull}). Secondly in this setting an analog of the conjecture of Rapoport and Zink mentioned in the introduction can be formulated and proved in Section~\ref{SectConjRZ}. And thirdly by the results of Kisin~\cite{Kisin} one could also in mixed characteristic replace crystalline Galois representations by the analogs of rigidified local shtukas as we have explained in the introduction and take this point of view.
\end{remark}

\medskip

It remains to establish the announced properties of the forgetful functor $\alpha$. We start with the following

\begin{lemma} \label{Lemma4.0}
Let $K$ be arbitrary. Let $N\in\BN_0$ and let $A,\wt A,C\in M_n\bigl(K\{z,\frac{\zeta^r}{z}\rangle\bigr)$ for some $r$ satisfy $C^\sigma=z^{-N}\wt ACA$. 
\begin{enumerate}
\item 
If $A,\wt A\in M_n\bigl(\CO_K\dbl z\dbr\bigr)$ then $C\in M_n\bigl(\CO_K\dbl z,\frac{\zeta^r}{z}\rangle[z^{-1}]\bigr)$.
\item 
If $A,\wt A\in\GL_n\bigl(\CO_K\dbl z\dbr[\frac{1}{z(z-\zeta)}]\bigr)$ then $C\in M_n\bigl(\CO_K\dbl z,z^{-1}\}[\tminus^{-1}]\bigr)$; see page~\pageref{Page34} for the notation.
\end{enumerate}
\end{lemma}

\begin{proof}
(a) \es The conditions on $C$ imply that for all $m\ge1$
\begin{equation}\label{EqForC}
C^{\sigma^m}\es=\es z^{-mN}\wt A^{\sigma^{m-1}}\cdot\ldots\cdot\wt A\cdot C\cdot A\cdot \ldots\cdot A^{\sigma^{m-1}}  \,.
\end{equation}
We write $C=\sum_{i=-\infty}^\infty C_iz^i$ with $C_i\in M_n(K)$. Let $d:=\|C\|_r=\max\{\,|C_i|\,|\zeta|^{ri}:i\in\BZ\,\}$. Then $|C_i|\le d\,|\zeta|^{-ri}$. Expanding (\ref{EqForC}) in powers of $z$ we get an expression for $C_i^{\sigma^m}$ which involves only terms from $\CO_K$ and the $C_j$ for $j\le i+mN$. Thus
\[
|C_i^{\sigma^m}|\es\le\es\max\{\,|C_j|:\es j\le i+mN\,\}\es\le \es d\,|\zeta|^{-ri-rmN}\,.
\]
In particular $|C_i|^{q^m}|\zeta|^{rmN}\le d\,|\zeta|^{-ri}$. When $m$ tends to infinity this implies that $|C_i|\le1$ proving part (a).

\smallskip
\noindent
(b) \es There exist integers $d,e\ge0$ with 
\[
z^d(z-\zeta)^dA\,,\es z^d(z-\zeta)^d\wt A\,,\es z^e(z-\zeta)^eA^{-1}\,,\es z^e(z-\zeta)^e\wt A^{-1}\es\in\es M_n\bigl(\CO_K\dbl z\dbr\bigr)\,.
\]
Then 
\[
(\tminus^{-2d}C)^\sigma\es=\es z^{-N-4d}\cdot z^d(z-\zeta)^d\wt A\cdot(\tminus^{-2d}C)\cdot z^d(z-\zeta)^dA
\]
and by (a) we have $\tminus^{-2d}C\;\in\; M_n\bigl(\CO_K\dbl z,\frac{\zeta^r}{z}\rangle[z^{-1}]\bigr)$ and hence also $\tminus^{2e}C=\tminus^{2(e+d)}\cdot \tminus^{-2d}C\;\in\; M_n\bigl(\CO_K\dbl z,\frac{\zeta^r}{z}\rangle[z^{-1}]\bigr)$. Now the equation 
\[
\tminus^{2e}C\es=\es z^{N-4e}\cdot z^e(z-\zeta)^e\wt A^{-1}\cdot(\tminus^{2e}C)^\sigma\cdot z^e(z-\zeta)^e A^{-1}
\]
shows that $\tminus^{2e}C$ converges on all of $0<|z|<1$. Therefore $\tminus^{2e}C\in M_n\bigl(\CO_K\dbl z,z^{-1}\}\bigr)$ and $C\in M_n\bigl(\CO_K\dbl z,z^{-1}\}[\tminus^{-1}]\bigr)$ as claimed.
\end{proof}

\begin{lemma}\label{Lemma4.5}
\begin{enumerate}
\item 
If $K$ satisfies condition (\ref{CondDoubleStar}) from page~\pageref{CondDoubleStar} then $\alpha$ is fully faithful.
\item 
Let $\lambda\in\BQ$. Then the restriction of $\alpha$ to the categories where $M\otimes_{\CO_K\dbl z\dbr}K\{z,\frac{\zeta^{1/q}}{z}\rangle$ is isoclinic of slope $\lambda$ (Definition~\ref{DefIsoclinic}) is fully faithful.
\end{enumerate}
\end{lemma}

\begin{proof}
Let $(M,F_M,\delta_M)$ and $(M',F_{M'},\delta_{M'})$ be rigidified local shtukas over $\CO_K$. Consider a morphism $f:(M,F_M)\to(M',F_{M'})$ of local shtukas and let $\bar f:=f\mod\Fm_K:(D,F_D)\to(D',F_{D'})$ be the induced morphism of $z$-isocrystals. We must show that $f$ respects the rigidifications, that is $\delta_{M'}\circ f=\bar f\circ\delta_M$.

By Theorem~\ref{Thm1a} there is an isomorphism of $\sigma$-modules over $\ol K\{z,\frac{\zeta^{1/q}}{z}\rangle$
\[
\phi:\es M\otimes_{\CO_K\dbl z\dbr}\ol K\{z,{\TS\frac{\zeta^{1/q}}{z}}\rangle\es\isoto\es\bigoplus_i\CF_{d_i,n_i}\,.
\]
By Lemma~\ref{Lemma4.0}, $\phi$ is in fact an isomorphism of $\sigma$-modules
\[
\phi:\es M\otimes_{\CO_K\dbl z\dbr}\CO_{\ol K}\dbl z,z^{-1}\}[\tminus^{-1}]\es\isoto\es\bigoplus_i\CF_{d_i,n_i}
\]
over $\CO_{\ol K}\dbl z,z^{-1}\}[\tminus^{-1}]$. We consider the reduction 
\[
\bar\phi\es:=\es\phi\mod\Fm_{\ol K}:\es D\otimes_{k\dpl z\dpr}k^\alg\dpl z\dpr\es\isoto\es\bigoplus_i\CF_{d_i,n_i}
\]
of $\phi$ over $k^\alg\dpl z\dpr$. Let $\wt K$ be the closure of the compositum $k^\alg K$ inside $\ol K$. We view $\bar\phi$ via the inclusion $k^\alg\dpl z\dpr\subset \CO_{\wt K}\dbl z,z^{-1}\}[\tminus^{-1}]$ as an isomorphism 
\[
\bar\phi:\es D\otimes_{k\dpl z\dpr}\CO_{\wt K}\dbl z,z^{-1}\}[\tminus^{-1}]\es\isoto\es\bigoplus_i\CF_{d_i,n_i}
\]
over $\CO_{\wt K}\dbl z,z^{-1}\}[\tminus^{-1}]$. We replace $\phi$ by $\bar\phi\circ\delta_M$. Since $\delta_M\equiv\id\mod\Fm_K$ this does not alter $\bar\phi$. In this way $\phi$ becomes an isomorphism
\[
\phi:\es M\otimes_{\CO_K\dbl z\dbr}\CO_{\wt K}\dbl z,z^{-1}\}[\tminus^{-1}]\es\isoto\es\bigoplus_i\CF_{d_i,n_i}
\]
of $\sigma$-modules over $\CO_{\wt K}\dbl z,z^{-1}\}[\tminus^{-1}]$ with $\phi=\bar\phi\circ\delta_M$. 
Note that we have lowered the field over which $\phi$ is defined from $\ol K$ to $\wt K$. 
We make the same construction for $M'$ obtaining the isomorphism $\phi'$ with $\phi'=\bar\phi'\circ\delta_{M'}$, and we set $g:=\phi'\circ f\circ\phi^{-1}:\bigoplus_i\CF_{d_i,n_i}\to\bigoplus_j\CF_{d'_j,n'_j}$ and $\bar g:=\bar\phi'\circ\bar f\circ\bar\phi^{-1}$. Then $\delta_{M'}\circ f=\bar f\circ\delta_M$ is equivalent to $g=\bar g$. To show that the latter is satisfied we break up $g$ into components $g_{ij}:\CF_{d_i,n_i}\to\CF_{d'_j,n'_j}$ for all $i,j$.

Fix a pair $i,j$. If $\frac{d_i}{n_i}=\frac{d'_j}{n'_j}$ then $g_{ij} \in\End_\sigma(\CF_{d_i,n_i})$. From the explicit computation of this endomorphism ring in \cite[Proposition 8.6]{HP} one sees that $g_{ij}\in M_{n_i}\bigl(\BF_{q^{n_i}}\dpl z\dpr\bigr)$. This implies that $g_{ij}=\bar g_{ij}$, and in particular, (b) follows.

If $\frac{d_i}{n_i}>\frac{d'_j}{n'_j}$ then $g_{ij}=0$ by Proposition~\ref{Prop0.7} and therefore $g_{ij}=\bar g_{ij}$.

If $\frac{d_i}{n_i}<\frac{d'_j}{n'_j}$ we use the assumption of (a) on $K$. Let $n$ be the least common multiple of $n_i$ and $n'_j$ and set $d:=n(\frac{d'_j}{n'_j}-\frac{d_i}{n_i})>0$. Then 
\[
g_{ij}\;\in\;\Hom_\sigma(\CF_{d_i,n_i}\,,\,\CF_{d'_j,n'_j})\;\subset\; \Hom_{\sigma^n}\bigl({\TS\CO(\frac{nd_i}{n_i})^{\oplus n_i}\,,\,\CO(\frac{nd'_j}{n'_j}})^{\oplus n'_j}\bigr)\;=\;M_{n'_j\times n_i}\bigl(\CO(d)^{F^n}(\wt K)\bigr)\;=\;(0)\,.
\]
Here the inclusion comes from viewing the $\sigma$-modules $\CF_{d_i,n_i}$ as $\sigma^n$-modules, as which they are isomorphic to $\CO(\frac{nd_i}{n_i})^{\oplus n_i}$. The first equality is due to $\CO(\frac{nd_i}{n_i})\dual\otimes\CO(\frac{nd'_j}{n'_j})\cong\CO(d)$. Finally, the last equality follows from Proposition~\ref{Prop0.5} since there is no $u\in\wt K$, $u\ne0$ with $u^{q^{n\nu}}\in\wt K$ for all $\nu\in\BZ$. Thus we have $g_{ij}=0=\bar g_{ij}$ also in this case and (a) is proved.
\end{proof}

\begin{proposition}\label{Prop4.1}
Let $d,s\in\BZ$ be relatively prime with $s>0$. Then the functor $\alpha$ restricted to the categories where $M\otimes_{\CO_K\dbl z\dbr}K\{z,\frac{\zeta^{1/q}}{z}\rangle$ is isoclinic of slope $\frac{d}{s}$ is an equivalence of categories.
\end{proposition}

\begin{proof}
Full faithfulness was established in Lemma~\ref{Lemma4.5}. We prove essential surjectivity. Let $(M,F_M)$ be a local shtuka over $\CO_K$ with $M\otimes_{\CO_K\dbl z\dbr}K\{z,\frac{\zeta^{1/q}}{z}\rangle$ isoclinic of slope $\frac{d}{s}$. Choose a basis of $M$ and let $A\in\GL_n\bigl(\CO_K\dbl z\dbr[\frac{1}{z-\zeta}]\bigr)$ be the matrix by which $F_M$ acts on this basis. Set $B:=A\mod\Fm_K$\,, $A':=AA^\sigma\cdot\ldots\cdot A^{\sigma^{s-1}}$ and $B':=A'\mod\Fm_K$. 

A close examination of the proof of Proposition~\ref{Prop6.3} shows that we can find a finite extension $K'$ of $K$ and a matrix $VU\in\GL_n\bigl(\CO_{K'}\dbl z,\frac{\zeta^{1/q}}{z}\rangle[z^{-1}]\bigr)$ such that $\wt A':=(VU)^{-1}A'(VU)^{\sigma^s}$ satisfies $z^{-d}\wt A'\in\GL_n\bigl(\CO_{K'}\dbl z\dbr\bigr)$. Namely we use Proposition~\ref{Prop6.3} for $r=q^{-s-1}$ and the $\sigma$-module $M\otimes_{\CO_K\dbl z\dbr}K\langle\frac{z}{\zeta^{q^{-s-1}}},\frac{\zeta^{1/q}}{z}\rangle$ over $K\langle\frac{z}{\zeta^{q^{-s-1}}},\frac{\zeta^{1/q}}{z}\rangle$. In the notation of the proof of Proposition~\ref{Prop6.3} we have $D=z^d\Id_n$ and by Lemma~\ref{Lemma4.0} the matrix $W$ with $W^{-1}A'W^{\sigma^s}=z^d\id_n$ satisfies $W\in\GL_n(\CO_{\ol K}\dbl z,\frac{\zeta^{1/q}}{z}\rangle[z^{-1}]\bigr)$. Hence we may choose $V\in M_n\bigl(\CO_{K'}[z,z^{-1}]\bigr)$ close to $W$ for a finite extension $K'$ of $K$. Then $V\in\GL_n\bigl(\CO_{K'}\dbl z,\frac{\zeta^{1/q}}{z}\rangle[z^{-1}]\bigr)$ and the same is true for $V^{-1}A'V^{\sigma^s}$. Now Lemma~\ref{Lemma6.1} produces a matrix $U\in\GL_n\bigl(\CO_{K'}\langle\frac{\zeta^{1/q}}{z}\rangle\bigr)$ with $z^{-d}\wt A'-\Id_n\in M_n\bigl(z\CO_{K'}\dbl z\dbr\bigr)$ where $\wt A':=(VU)^{-1}A'(VU)^{\sigma^s}$. In particular $z^{-d}\wt A'\in\GL_n\bigl(\CO_{K'}\dbl z\dbr\bigr)$.

After enlarging $K'$ we may assume that the section $k\hookrightarrow\CO_K$ extends to $k'\hookrightarrow\CO_{K'}$, where $k'$ is the residue field of $\CO_{K'}$. Let the matrices $\wt B':=\wt A'\mod\Fm_{K'}$ and $\ol{VU}:=VU\mod\Fm_{K'}$ be viewed as elements of $\GL_n\bigl(\CO_{K'}\dbl z,z^{-1}\}\bigr)$ via this section. Let
\begin{eqnarray*}
\wt C_m & := & \wt B'\cdot\ldots\cdot\wt B'{}^{\sigma^{sm}}\cdot (\wt A'{}^{\sigma^{sm}})^{-1}\cdot\ldots\cdot(\wt A')^{-1} \\[2mm]
& = & (z^{-d}\wt B')\cdot\ldots\cdot(z^{-d}\wt B')^{\sigma^{sm}}\cdot (z^{-d}\wt A'{}^{\sigma^{sm}})^{-1}\cdot\ldots\cdot(z^{-d}\wt A')^{-1}\quad\in\es\GL_n\bigl(\CO_{K'}\dbl z\dbr\bigr)\,.
\end{eqnarray*}
Now the arguments given in the proof of Lemma~\ref{LemmaGL} show that the sequence $\wt C_m$ converges to a matrix $\wt C\in\GL_n\bigl(\CO_{K'}\dbl z\dbr\bigr)$ which satisfies $\wt C\equiv\Id_n\mod\Fm_{K'}$ and $\wt C\wt A'=\wt B'\wt C^{\sigma^s}$. Set $C:=\ol{VU}\,\wt C\,(VU)^{-1}\in\GL_n\bigl(\CO_{K'}\dbl z,\frac{\zeta^{1/q}}{z}\rangle[z^{-1}]\bigr)$. Then $C\equiv\Id_n\mod\Fm_{K'}$ and $CA'=B'C^{\sigma^s}$. Moreover, we obtain the equality
$(BC^\sigma A^{-1})A'=B'(BC^\sigma A^{-1})^{\sigma^s}$ and since by Lemma~\ref{Lemma4.5} the isomorphism $\delta_M$ is uniquely determined we conclude that $CA=BC^\sigma$. 

It remains to show that indeed $C\in\GL_n\bigl(\CO_K\dbl z,z^{-1}\}[\tminus^{-1}]\bigr)$. Repeated application of the equation $C=BC^\sigma A^{-1}$ shows that $C\in\GL_n\bigl(\CO_{K''}\dbl z,\frac{\zeta^{1/q}}{z}\rangle[z^{-1}]\bigr)$ where $K''$ is the separable closure of $K$ inside $K'$. We may replace the extension $K''/K$ by its Galois closure. If $\phi\in\Gal(K''/K)$ we find $\phi(C)A=B\phi(C)^\sigma$. So again by Lemma~\ref{Lemma4.5} we must have $\phi(C)=C$, whence $C\in\GL_n\bigl(\CO_K\dbl z,\frac{\zeta^{1/q}}{z}\rangle[z^{-1}]\bigr)$. Now the equation $C=BC^\sigma A^{-1}$ shows that $C\in\GL_n\bigl(\CO_K\dbl z,z^{-1}\}[\tminus^{-1}]\bigr)$ as desired. Thus $C$ defines a rigidification $\delta_M$ and this proves essential surjectivity.
\end{proof}

\begin{proposition} \label{Prop4.2}
If $K$ is algebraically closed $\alpha$ is essentially surjective but not full.
\end{proposition}

\begin{proof}
To prove essential surjectivity let $(M,F_M)$ be a local shtuka over $\CO_K$. By Theorem~\ref{Thm1a} there exists an isomorphism $\phi:M\otimes_{\CO_K\dbl z\dbr}K\{ z,\frac{\zeta^{1/q}}{z}\rangle\isoto\bigoplus_i\CF_{d_i,n_i}$. By Lemma~\ref{Lemma4.0}, $\phi$ is in fact an isomorphism
\[
\phi:\es M\otimes_{\CO_K\dbl z\dbr}\CO_K\dbl z,z^{-1}\}[\tminus^{-1}]\es\isoto\es\bigoplus_i\CF_{d_i,n_i}
\]
of $\sigma$-modules over $\CO_K\dbl z,z^{-1}\}[\tminus^{-1}]$. Let
\[
\bar\phi:=\phi\mod\Fm_K:\es D\otimes_{k\dpl z\dpr}\CO_K\dbl z,z^{-1}\}[\tminus^{-1}]\es\isoto\es\bigoplus_i\CF_{d_i,n_i}\,.
\]
Then $\delta_M:=\bar\phi^{-1}\phi$ is a rigidification of $(M,F_M)$.

\smallskip

The easiest example showing that $\alpha$ is not full is the following. Let $\tplus=\sum_{i=0}^\infty t_iz^i$ with $t_i\in K$ satisfy $\tplus^\sigma=(z-\zeta)\tplus$, that is, $t_i^q+\zeta t_i=t_{i-1}$. This implies $|t_i|=|\zeta|^{q^{-i}/(q-1)}<1$. Consider the rigidified local shtukas $M=\bigl(\CO_K\dbl z\dbr,(z-\zeta)\cdot\sigma,\tminus^{-1}\bigr)$ and $M'=\bigl(\CO_K\dbl z\dbr,\sigma,1\bigr)$ over $\CO_K$ and the morphism of local shtukas $f:M\to M', 1\mapsto \tplus$. It satisfies $\bar f:=f\mod\Fm_K=0$ and $f\ne\bar f\,\tminus^{-1}$. So $f$ does not lie in the image of $\alpha$.
\end{proof}

\begin{example}\label{Ex4.3}
We show that in general $\alpha$ need not be essentially surjective.
Let $K$ be the $\zeta$-adic completion of the field $\BF_q^{\,\alg}\dpl\zeta\dpr(\zeta^{(q+1)^{-i}}:i\in\BN)$. Then the value group of $K$ is isomorphic to $\bigl(\BZ[\frac{1}{q+1}],+\bigr)$ when we map $\zeta$ to $1$. Set $a=\sum_{i=0}^\infty \zeta^{(q+1)^{-i}}z^i$ and let \mbox{$(M,F_M)=\bigl(\CO_K\dbl z\dbr^{\oplus 2},A\cdot\sigma\bigr)$} where \mbox{$A=\left(\begin{array}{cc} 1 & a \\ 0 & z-\zeta \end{array}\right)$}. We obtain $(D,F_D)=\bigl(\BF_q^{\,\alg}\dpl z\dpr^{\oplus 2},B\cdot\sigma\bigr)$ with $B=\left(\begin{array}{cc} 1 & 0 \\ 0 & z \end{array}\right)$. Assume that there exists a rigidification $\delta_M$ on the local shtuka $(M,F_M)$, that is a matrix $C=\left(\begin{array}{cc} u & v \\ w & x \end{array}\right)\in\GL_2\bigl(\CO_K\dbl z,z^{-1}\}[\tminus^{-1}]\bigr)$ with $C\equiv\Id_2\mod\Fm_K$ and $CA=BC^\sigma$. The last equation amounts to
\[
\left(\begin{array}{cc} u & ua+(z-\zeta)v \\ w & wa+(z-\zeta)x \end{array}\right) \es=\es\left(\begin{array}{rr} u^\sigma & v^\sigma \\ z\,w^\sigma & z\,x^\sigma \end{array}\right)\,.
\]
In particular $w\in\CO(-1)^F(K)=(0)$ and $u,x\tminus\in\CO(0)^F(K)=\BF_q\dpl z\dpr$ by Proposition~\ref{Prop0.5}. Since $u,x\equiv 1\mod\zeta$ we must have $u=1,x=\tminus^{-1}$ and thus $a+(z-\zeta)v=v^\sigma$. Expand $v=\sum_{i=-\infty}^\infty v_iz^i$ with $v_i\in\CO_K$. Then 
\begin{equation}\label{Eq2.19}
\zeta^{(q+1)^{-i}}+v_{i-1}\es=\es v_i^q+\zeta v_i\qquad \text{for all }i\ge0\,.
\end{equation}
If for some $i\ge0$ the absolute value $|v_{i-1}|>|\zeta|^{(q+1)^{-i}}>|\zeta|^{q/(q-1)}$ then $|v_i|^q=|v_{i-1}|>|\zeta v_i|$. Indeed, $|v_i|^q\le|\zeta v_i|$ implies $|v_i|\le|\zeta|^{1/(q-1)}$ and $|v_i^q+\zeta v_i|\le|\zeta|^{q/(q-1)}$. If we assume $|v_{j-1}|>|\zeta|^{(q+1)^{-j}}$ for all $j\ge i$ this would imply $|v_{i-1}|^{q^{i-j}}=|v_{j-1}|>|\zeta|^{(q+1)^{-j}}$, whence $q^{i-j}\cdot\log_{|\zeta|}|v_{i-1}|<(q+1)^{-j}$. However, the right hand side decreases faster than the left hand side. So we eventually must find some $j>0$ with $|v_{j-1}|\le|\zeta|^{(q+1)^{-j}}$. This implies  $|v_j|\le|\zeta|^{(q+1)^{-j}/q}<|\zeta|^{(q+1)^{-j-1}}$ since the opposite inequality $|v_j|>|\zeta|^{(q+1)^{-j}/q}>|\zeta|^{1/(q-1)}$ yields $|v_j|^q >|\zeta v_j|$. Thus for $i=j+1$ the left hand side of Equation~(\ref{Eq2.19}) has absolute value equal to $|\zeta|^{(q+1)^{-j-1}}$. Reasoning as before we find $|v_{j+1}|=|\zeta|^{(q+1)^{-j-1}/q}$. In particular 
\[
\TS\log_{|\zeta|}|v_{j+1}|\es=\es\frac{1}{q(q+1)^{j+1}}\es\notin\es\BZ[{\TS\frac{1}{q+1}}]\,.
\]
This is a contradiction and therefore the isomorphism $\delta_M$ cannot exist in this example.\qed
\end{example}

\begin{remark}\label{Remark4.4}
In contrast to the previous example, local shtukas which arise from global objects as in Example~\ref{ExGlobalObjects} carry a canonical functorial rigidification regardless of the properties of $K$. We indicate this in a special case and leave the general case to the interested reader. Let $A=\BF_q[z]$ and fix $c:A\to\CO_K$ with $c(z)=\zeta\in\Fm_K$. Consider an Anderson $A$-motive over $K$ with good reduction, that is, a finite locally free $\CO_K[z]$-module $M$ together with an injective morphism $F_M:\sigma^\ast M\to M$ such that $\coker F_M$ is a free $\CO_K$-module and annihilated by a power of $z-\zeta$. Fix an $\CO_K[z]$-basis of $M$ and let $A\in\GL_n\bigl(\CO_K[z][\frac{1}{z-\zeta}]\bigr)$ be the matrix by which $F_M$ acts on this basis. Put $B:=A\mod\Fm_K$. Now the local shtuka associated with $(M,F_M)$ is $\wh M=\bigl(M\otimes_{\CO_K[z]}\CO_K\dbl z\dbr,F_M\otimes\id\bigr)$.
 Since $A-B\in M_n\bigl(\pi\CO_K[z][\frac{1}{z-\zeta}]\bigr)$ for some element $\pi\in\Fm_K$ the argument given in the proof of Lemma~\ref{LemmaGL} shows the existence of a canonical functorial rigidification on $\wh M$.
\end{remark}

%%%%%%%%%%%%%%%%%%%%%%%%%%%%%%%%%%%%%%%%%%%%%%%%%%%%%%%%%%%%%%%%%%%%%%
%
%    Criteria for Admissibility
%
%%%%%%%%%%%%%%%%%%%%%%%%%%%%%%%%%%%%%%%%%%%%%%%%%%%%%%%%%%%%%%%%%%%%%%

\subsection{Criteria for Admissibility} \label{SectCriteriaWA}

We want to relate Hodge-Pink structures to the $\sigma$-modules we studied in Chapter~\ref{ChaptSFT}. This will lead to a description of Hodge-Pink structures analogous to the description in mixed characteristic of filtered isocrystals by $p$-adic differential equations; see~\cite{Berger1}. Here in equal characteristic it is even possible to formulate this in a relative situation, what we now want to do. This is already also a preparation for the analog of the Rapoport-Zink conjecture. To be precise we make the following definition. Let $k$ be a field containing $\BF_q$ and let the complete local ring $R$ be an extension of $k\dbl\zeta\dbr$. Let $X$ be a quasi-paracompact admissible formal scheme over $\Spf R$ and let $X_L$ be the rigid analytic space over $L$ associated with $X$; see Appendix~\refAppBerkovich{}. The main example to keep in mind is $X=\Spf B\open$ for an admissible formal $R$-algebra $B\open$ (Appendix~\refAppRigFormal) and $X_L=\Spm B$ where $B:=B\open\otimes_R L$. Then $\CO_X=B\open$ and $\CO_{X_L}=B$.

\begin{definition}\label{DefFamilyOfHPStr}
Let $(D,F_D)$ be a $z$-isocrystal of rank $n$ over $k$ and consider the sheaf $\Fp:=\sigma^\ast D\otimes_{k\dpl z\dpr}\CO_{X_L}\dbl z-\zeta\dbr$ on $X_L$. Fix Hodge-Pink weights $\ul w=(w_1\geq\ldots\geq w_n)$. A \emph{family of Hodge-Pink structures over $X_L$ with Hodge-Pink weights $\ul w$} on the (constant) $z$-isocrystal $(D,F_D)$ is a sheaf of $\CO_{X_L}\dbl z-\zeta\dbr$-lattices $\Fq$ in $\Fp[\frac{1}{z-\zeta}]\;:=\;\Fp\otimes_{\CO_{X_L}\dbl z-\zeta\dbr}\CO_{X_L}\dbl z-\zeta\dbr[\frac{1}{z-\zeta}]$ which is a direct summand as $\CO_{X_L}$-module, such that for every point $x\in X_L$ and for every integer $e\geq -w_n$
\[
(z-\zeta)^{-e}\Fp_x/\Fq_x\es\cong\es\bigoplus_{i=1}^n \kappa(x)\dbl z-\zeta\dbr/(z-\zeta)^{e+w_i}\,.
\]
\end{definition}
With the same definitions of tensor products, etc., as in the absolute case we obtain an $\BF_q\dpl z\dpr$-linear exact rigid tensor category.

In the situation of Definition~\ref{DefFamilyOfHPStr} we let $\ul D=(D,F_D,\Fq)$. We define the $\sigma$-module
\[
\CP\es:=\es\sigma^\ast D\otimes_{k\dpl z\dpr}\CO_{X_L}\ancon\,,\quad F_\CP\es:=\es\sigma^\ast F_D\otimes\id
\]
over $\CO_{X_L}\ancon$.
Then $\CP\otimes \CO_{X_L}\dbl z-\zeta\dbr=\Fp$ and the Hodge-Pink structure $\Fq\subset\CP\otimes \CO_{X_L}\dbl z-\zeta\dbr[\frac{1}{z-\zeta}]$ defines a $\sigma$-module $\CQ$ over $\CO_{X_L}\ancon$ which is a modification of $\CP$ at $z=\zeta^{q^i}$ for $i\in\BN_0$ as follows.
Consider the isomorphism $\eta_i\;=\;\bigl(F_\CP\circ\ldots\circ(\sigma^{i-1})^\ast F_\CP\bigr)\otimes\id$
\[
\eta_i:\es(\sigma^i)^\ast\bigl(\Fp[{\TS\frac{1}{z-\zeta}}]\bigr)\es=\es \bigl((\sigma^i)^\ast\CP\bigr)\otimes \CO_{X_L}\dbl z-\zeta^{q^i}\dbr[{\TS\frac{1}{z-\zeta^{q^i}}}]\es\isoto\es\CP\otimes \CO_{X_L}\dbl z-\zeta^{q^i}\dbr[{\TS\frac{1}{z-\zeta^{q^i}}}]
\]
and set $\DS \tminus=\prod_{i\in\BN_0}\bigl(1-{\TS\frac{\zeta^{q^i}}{z}}\bigr)\in \CO_{X_L}\ancon$.
We define $\CQ$ as the $\CO_{X_L}\ancon$-submodule of $\CP[\tminus^{-1}]$ which coincides with $\CP$ outside $z=\zeta^{q^i}$ for $i\in\BN_0$ and at $z=\zeta^{q^i}$ satisfies $\CQ\otimes \CO_{X_L}\dbl z-\zeta^{q^i}\dbr\;=\;\eta_i(\sigma^{i\ast}\Fq)$. 
By construction $F_\CP$ induces on $\CQ$ the structure of a $\sigma$-module over $\CO_{X_L}\ancon$.

\begin{definition} \label{Def1.6}
The pair $(\CP,\CQ)$ just constructed is called the \emph{pair of $\sigma$-modules} over $\CO_{X_L}\ancon$ associated with $\ul D$. We denote the functor that maps $\ul D$ to the pair $(\CP,\CQ)$ by $\pair$.
\end{definition}

The category of pairs of $\sigma$-modules over $\CO_{X_L}\ancon$ has as objects all pairs $(\CP,\CQ)$ of $\sigma$-modules such that $\CP[\tminus^{-1}]=\CQ[\tminus^{-1}]$ and $F_\CP=F_\CQ$. The morphisms between $(\CP,\CQ)$ and $(\CP',\CQ')$ are the morphisms of $\sigma$-modules $\CP\to\CP'$ which at the same time are morphisms $\CQ\to\CQ'$.

\begin{theorem} \label{Thm1.7}
The functor $\pair$ is a faithful $\BF_q\dpl z\dpr$-linear exact tensor functor.
% NOT FULL !!!!!!!!!!!!!!!!!!!!!!!!!!!!!!!!!!!!!!!!!!!!!!!!!!!!!!!!!!!!!!
\end{theorem}

\begin{proof}
The $\BF_q\dpl z\dpr$-linearity and the compatibility with tensor products and duals follows directly from the construction.
To prove exactness let $0\to\ul D'\xrightarrow{f}\ul D\xrightarrow{g}\ul D''\to0$ be an exact sequence of (constant) $z$-isocrystals with Hodge-Pink structure over $X_L$. Then the exactness of the underlying sequence of $z$-isocrystals implies that the sequence of $\sigma$-modules $0\to\CP'\xrightarrow{f}\CP\xrightarrow{g}\CP''\to0$ over $\CO_{X_L}\ancon$ is exact. Moreover the sequence 
\begin{equation}\label{EqQ's}
0\longto\CQ'\xrightarrow{\es f\;}\CQ\xrightarrow{\es g\;}\CQ''\longto0
\end{equation}
of sheaves on $X_L\times_L\dotBD(1)$ is exact outside $z=\zeta^{q^i}$ for $i\in\BN_0$ since there it coincides with $0\to\CP'\xrightarrow{f}\CP\xrightarrow{g}\CP''\to0$. At $z=\zeta$ the sequence of the completed stalks equals $0\to\Fq'\xrightarrow{f}\Fq\xrightarrow{g}\Fq''\to0$ and is exact by the strictness of $f$ and $g$. Hence (\ref{EqQ's}) is exact at $z=\zeta$ and the commutative diagram
\[
\xymatrix {
0\ar[r] & \sigma^\ast\CQ'\ar[r]\ar[d]^{F_{\CQ'}} & \sigma^\ast\CQ\ar[r]\ar[d]^{F_{\CQ}} & \sigma^\ast\CQ''\ar[r]\ar[d]^{F_{\CQ''}} & 0 \\
0\ar[r] & \CQ'\ar[r] & \CQ\ar[r] & \CQ''\ar[r] & 0
}
\]
shows that (\ref{EqQ's}) is exact everywhere. Finally the faithfulness of $\pair$ results from the faithfulness of the functor $\sigma^\ast D\mapsto\CP=\sigma^\ast D\otimes_{k\dpl z\dpr}\CO_{X_L}\ancon$.
\end{proof}

We determine the composition $\pair\BH$ in the situation of the previous section. Recall that there are two morphisms
\[
\alpha:\es\CO_K\dbl z\dbr \xrightarrow{\mod\Fm_K} k\dbl z\dbr \into K\ancon \qquad\text{and}\qquad \beta:\es\CO_K\dbl z\dbr\subset K\ancon
\]
mapping $z$ to $z$.

\begin{proposition}\label{Prop7.5}
There exists a natural isomorphism of tensor functors from the category of rigidified local shtukas over $\CO_K$ to the category of pairs of $\sigma$-modules over $K\ancon$ between $\pair\BH$ and the functor $(M,F_M,\delta_M)\mapsto(\CP,\CQ)$ where
\begin{eqnarray*}
\CP&=&\sigma^\ast M\otimes_{\CO_K\dbl z\dbr,\alpha}K\ancon\es=\es\bigl(\sigma^\ast M\otimes_{\CO_K\dbl z\dbr}k\dpl z\dpr\bigr)\otimes_{k\dpl z\dpr}K\ancon\qquad\text{and}\\
\CQ&=& M\otimes_{\CO_K\dbl z\dbr,\beta}K\ancon
\end{eqnarray*}
and where the morphism $\sigma^\ast\delta_M\circ F_M^{-1}:\CQ[\tminus^{-1}]\isoto\CP[\tminus^{-1}]$ identifies $\CQ$ with a submodule of $\CP[\tminus^{-1}]$.
\end{proposition}

\begin{proof}
Let $\pair\BH(M)=(\CP',\CQ')$. Clearly $\CP=\CP'$. The $\sigma$-module $\CQ$ is isomorphic to its image $\CQ''$ inside $\CP[\tminus^{-1}]$ under $\sigma^\ast\delta_M\circ F_M^{-1}$
\[
\xymatrix @C+6pc {
\sigma^\ast M\otimes_{\CO_K\dbl z\dbr,\beta}K\ancon{}[\tminus^{-1}] \ar[r]^{\TS\qquad\qquad\sigma^\ast\delta_M} & \CP[\tminus^{-1}]\es . \\
M\otimes_{\CO_K\dbl z\dbr,\beta}K\ancon \ar[u]^{\TS F_M^{-1}}\ar[ur]_{\TS\sigma^\ast\delta_M\circ F_M^{-1}}
}
\]
It satisfies $\CQ''\otimes K\dbl z-\zeta\dbr\;=\;\Fq\;\subset\;\CP\otimes K\dpl z-\zeta\dpr$. Moreover, the commutativity of the diagram
\[
\begin{CD}
\sigma^{i\ast} M\otimes K\dbl z-\zeta^{q^i}\dbr & @>{\;\sigma^{i\ast}(\sigma^\ast\delta_M\circ F_M^{-1})\;}>> & 
\sigma^{i\ast}\CQ''\otimes K\dbl z-\zeta^{q^i}\dbr \es= \es \sigma^{i\ast}\Fq 
\es\subset\es & \sigma^{i\ast}\CP\otimes K\dpl z-\zeta^{q^i}\dpr \\[1mm]
@VV{F_M\circ\ldots\circ(\sigma^{i-1})^\ast F_M}V & & & @V{\eta_i}VV \\[1mm]
M\otimes K\dbl z-\zeta^{q^i}\dbr & @>{\quad\sigma^\ast\delta_M\circ F_M^{-1}\quad}>> & \CQ''\otimes K\dbl z-\zeta^{q^i}\dbr \es=\es \eta_i(\sigma^{i\ast}\Fq) 
\es\subset\es & \CP\otimes K\dpl z-\zeta^{q^i}\dpr \quad
\end{CD}
\]
in which all arrows are isomorphisms, shows that $\CQ''\otimes K\dbl z-\zeta^{q^i}\dbr\;=\;\eta_i(\sigma^{i\ast}\Fq)\;\subset\;\CP\otimes K\dpl z-\zeta^{q^i}\dpr$. So $M\otimes_{\CO_K\dbl z\dbr,\beta}K\ancon\cong\CQ''=\CQ'$ as desired.
\end{proof}

As an example let us compute the effect of $\pair$ on the Tate object $\BOne(n)$ from Example~\ref{ExTateObjects} in the category of rigidified local shtukas with Hodge-Pink structure over $\BFZ$. We find $\pair\bigl(\BOne(n)\bigr)=\bigl(\CO(-n),\tminus^{-n}\CO(-n)\bigr)$. Since the $\sigma$-module $\tminus^{-n}\CO(0)$ over $\BFZ\ancon$ equals $\bigl(\BFZ\ancon,F=(1-\frac{\zeta}{z})^n\cdot\sigma\bigr)$ we obtain
\[
\tminus^{-n}\CO(-n)\es=\es \tminus^{-n}\CO(0)\otimes\CO(-n)\es=\es\bigl(\BFZ\ancon, F=(z-\zeta)^n\cdot\sigma\bigr)
\]
and this also illustrates Proposition~\ref{Prop7.5}.

\medskip

Next we want to derive a criterion for (weak) admissibility in terms of the pair $(\CP,\CQ)$. As a consequence we shall see that ``admissible implies weakly admissible''. In the next section the criterion will allow us to prove that also the converse is true, that is, ``weakly admissible implies admissible'', provided $K$ satisfies condition (\ref{CondDoubleStar}) from page~\pageref{CondDoubleStar}.

For an integer $e$ we consider the $\sigma$-submodule $\tminus^{-e}\CP$ of $\CP[\tminus^{-1}]$ over $K\ancon$. It satisfies $F_{\tminus^{-e}\CP}=(1-\frac{\zeta}{z})^e\cdot F_\CP$ and is isomorphic to $\CP\otimes \tminus^{-e}\CO(0)$. We want to construct an isomorphism between $\tminus^{-e}\CO(0)$ and $\CO(e)$ over $\ol K\ancon$. For this we need to find a solution $\tplus=\sum_{i\geq0}t_i z^i\in\ol K\ancon\mal$ for the equation $z^{-1}\tplus^\sigma=(1-\frac{\zeta}{z})\,\tplus$ which we expand to $t_i^q+\zeta t_i=t_{i-1}$. Since $\ol K$ is algebraically closed we may indeed find solutions $t_i\in\ol K$ of the last equation. One easily sees that $|t_i|=|\zeta|^{q^{-i}/(q-1)}$, hence $\tplus\in\ol K\ancon\mal$. Multiplication with $\tplus^e$ defines an isomorphism $\tminus^{-e}\CO(0)\otimes\ol K\ancon\isoto\CO(e)$ over $\ol K\ancon$. Hence the $\sigma$-modules $\tminus^{-e}\CP$ and $\CP\otimes\CO(e)$ become isomorphic over $\ol K\ancon$. In particular $\deg \tminus^{-e}\CP=e\cdot\rk\CP+\deg\CP$. This discussion corresponds to the fact that $t:=\tminus\tplus$ satisfies $z^{-1}\sigma^\ast t=t$, that is $t\in\CO(1)^F(\ol K)$. So multiplication with $t^e$ yields an embedding $\CO(0)\to\CO(e)$ of $\sigma$-modules over $\ol K\ancon$. Note that $t$ is an analog of the complex period $2\pi i$; see \cite[\S2.7]{HartlDict}.

\begin{lemma} \label{Lemma2.7}
Let $\ul D$ be a $z$-isocrystal with Hodge-Pink structure over $K$. Consider the pair $(\CP,\CQ)=\pair(\ul D)$ of $\sigma$-modules over $K\ancon$. Then
\[
t_N(\ul D)\es=\es -\deg\CP\,,\qquad\text{and}\qquad t_H(\ul D)\es=\es\deg\CQ-\deg\CP\,.
\]
\end{lemma}

\begin{proof}
To prove the first assertion let $n$ be the rank of $D$ and let $d=t_N(\ul D)$. Then the $\sigma$-module $\wedge^n(D,F_D)\otimes_{k\dpl z\dpr}\ol K\dpl z\dpr$ over $\ol K\dpl z\dpr$ is isomorphic to $\CO(-d)$ by Theorem~\ref{Thm2} and Lemma~\ref{Lemma14}. Therefore already the $\sigma$-module $\wedge^n(D,F_D)\otimes_{k\dpl z\dpr} \ol K\con$ over $\ol K\con$ is isomorphic to $\CO(-d)$ by Lemma~\ref{LemmaDescSemistFilt}. So $\wedge^n\CP\otimes \ol K\ancon\;=\;\wedge^n(D,F_D)\otimes_{k\dpl z\dpr}\ol K\ancon\;\cong\;\CO(-d)$ and this implies that $\deg\CP=-d$ proving the claim. Alternatively one could also argue that $\wedge^n(D,F_D)\cong\CO(-d)$ as $z$-isocrystals over $k^\alg$ and then transport this isomorphism to $\ol K\ancon$ via some (any) embedding $k^\alg\dbl z\dbr\hookrightarrow\ol K\ancon$.

To prove the second assertion let $w_1\geq\ldots\geq w_n$ be the Hodge-Pink weights of $\ul D$ and let $e$ be an integer such that $\CQ\subset \tminus^{-e}\CP$. Then over $\ol K\ancon$ we have $\CO(\deg\CQ)\;\cong\;\wedge^n\CQ\;\subset\;\wedge^n\tminus^{-e}\CP\;\cong\;\CO(ne+\deg\CP)$. Since $\CQ$ differs from $\CP$ only at $z=\zeta^{q^i}$ for $i\geq0$,
the $\ol K\dbl z-\zeta\dbr$-module
\[
\wedge^n(z-\zeta)^{-e}\Fp_D\,/\,\wedge^n\Fq_D\es\cong\es \Bigl(\CO(ne+\deg\CP)\,/\,\CO(\deg\CQ)\Bigr)\otimes\ol K\dbl z-\zeta\dbr
\]
has length $ne+\deg\CP-\deg\CQ$. On the other hand it has length $\sum_{i=1}^n(e+w_i)$. From this the second assertion follows. 
\end{proof}

\begin{proposition} \label{Prop2.8}
Let $\ul D=(D,F_D,\Fq_D)$ be a $z$-isocrystal with Hodge-Pink structure over $K$ and let $(\CP,\CQ)=\pair(\ul D)$ be the associated pair of $\sigma$-modules over $K\ancon$. Then $\ul D$ is weakly admissible if and only if the following conditions hold:
\begin{enumerate}
\item
$\deg\CQ=0$ and
\item
for any $F_D$-stable $k\dpl z\dpr$-subspace $D'\subset D$, the saturation $\CQ'$ of $\bigl(D'\otimes_{k\dpl z\dpr} K\ancon, F_{D'}\bigr)\cap\CQ$ inside $\CQ$ satisfies $\deg\CQ'\leq0$.
\end{enumerate}
\end{proposition}

\begin{proof}
Note that $\bigl((D'\otimes K\ancon, F_{D'}),\CQ'\bigr)$ is the pair of $\sigma$-modules associated with the strict subobject $(D',F_{D'},\Fq_{D'})$ of $\ul D$ with $\Fq_{D'}=\Fq_D\cap \sigma^\ast D'\otimes_{k\dpl z\dpr}K\dpl z-\zeta\dpr$. Thus the claim follows from Lemma~\ref{Lemma2.7}.
\end{proof}

\begin{theorem} \label{Thm2.6}
A $z$-isocrystal with Hodge-Pink structure $\ul D$ over $K$ is admissible if and only if the associated pair of $\sigma$-modules $(\CP,\CQ)$ over $K\ancon$ satisfies the condition that $\CQ$ is isoclinic of slope zero (Definition~\ref{DefIsoclinic}).
\end{theorem}

%This is the analog of Berger's criterion for admissibility \cite[Proposition IV.2.2]{Berger1}. 
%THIS IS NOT QUITE TRUE !!!!!!!!!!!!!!!!!!!!!!!!!!!!!!!!
As an immediate consequence we obtain

\begin{corollary} \label{Cor2.9}
Every admissible $z$-isocrystal with Hodge-Pink structure is weakly admissible.
\end{corollary}

\begin{proof}
Let $\ul D$ be an admissible $z$-isocrystal with Hodge-Pink structure over $K$ and let $(\CP,\CQ)=\pair(\ul D)$. Then $\CQ$ is isoclinic of slope zero by Theorem~\ref{Thm2.6}, hence semistable (Definition~\ref{DefStable}) by Lemma~\ref{Lemma1.5.6'}. 
Therefore any non-zero $\sigma$-submodule $\CQ'$ of $\CQ$ has $\deg\CQ'\leq\deg\CQ=0$ and the assertion follows from Proposition~\ref{Prop2.8}.
\end{proof}

\begin{proof}[Proof of Theorem~\ref{Thm2.6}]
If $\ul D$ is admissible let $M$ be a rigidified local shtuka over $\CO_K$ which induces $\ul D$. We have $\CQ\cong M\otimes_{\CO_K\dbl z\dbr}K\ancon$ by Proposition~\ref{Prop7.5}. Fix a basis of $M$ and let $\Phi=\sum_{i\geq0}\Phi_i z^i\,\in\GL_n\bigl(\CO_K\dbl z\dbr[\frac{1}{z-\zeta}]\bigr)\subset\GL_n\bigl(K\dbl z\dbr\bigr)$
be the matrix by which $F_M$ acts on this basis. Let $\ol K$ be the completion of an algebraic closure of $K$. Since $\Phi_0\in\GL_n(K)$ there exists a matrix $S\in\GL_n(\ol K)$ such that $S^{-1}\Phi_0S=\Id_n$ by Lang's theorem \cite{Lang57}. Hence
\[
S^{-1}\Phi S^\sigma\es\in\es\Id_n+M_n\bigl(z\ol K\dbl z\dbr\bigl)\,.
\]
and $M\otimes\ol K\dpl z\dpr$ is isoclinic of slope zero by Lemma~\ref{Lemma14}. By Proposition~\ref{Prop7}, $\CQ\otimes\ol K\ancon\cong M\otimes\ol K\ancon$ is isoclinic of slope zero. By its very definition this means that $\CQ$ is isoclinic of slope zero.
Alternatively one could consider the equation $U=\Phi U^\sigma$ and explicitly show that it has a solution with $U\in\GL_n(\ol K\langle\frac{z}{\zeta}\rangle)$. Then $\CQ\otimes\ol K\ancon\cong\CO(0)^{\oplus n}$.

Using Corollary~\ref{Cor6.12b} the converse follows from the following Proposition~\ref{Prop2.6b} with $L=K, R=\CO_K, X=\Spf\CO_K, X_L=\Spm K, X_0=\Spec k$.
\end{proof}

Proposition~\ref{Prop2.6b} below gives a criterion for a family of Hodge-Pink structures on a constant $z$-isocrystal to arise from a bounded rigidified local shtuka. Although we need this result for the proof of Theorem~\ref{Thm2.6} only in case $L=K$ and $X_L=\Spm K$, we formulate it here and use it in Theorem~\ref{Thm2.16} in a relative setting over quasi-paracompact quasi-separated rigid analytic spaces (see Appendix~\refAppBerkovich). These are precisely the spaces which admit quasi-paracompact admissible formal models; see Theorem~\ref{ThmFormalRigBerkovich}. For any such formal $\Spf R$-scheme $X$ we denote the special fiber $X\otimes_R k$ by $X_0$. Recall that for an admissible formal $R$-algebra $B\open$ one sets 
\begin{tabbing}
$B\con[n]$\es\= \kill
$B\open\langle z\rangle$\> $\DS =\es\bigl\{\,\sum_{i=0}^\infty b_iz^i:\es b_i\in B\open,\;|b_i|\to0 \;(i\to\infty)\;\bigr\}$ \quad and\\[2mm]
$B\open\langle z,z^{-1}\rangle$\> $\DS =\es\bigl\{\,\sum_{i=-\infty}^\infty b_iz^i:\es b_i\in B\open,\;|b_i|\to0 \;(i\to\pm\infty)\;\bigr\}$\,.
\end{tabbing}
These are again admissible formal $R$-algebras. We further denote the passage to the associated affinoid $L$-algebra $B:=B\open\otimes_R L$ by dropping the superscript $\open$.

\begin{proposition} \label{Prop2.6b}
Let $(D,F_D)$ be a $z$-isocrystal over $k$. Let $X_L$ be a quasi-paracompact quasi-separated rigid analytic space over $L$. Let $\Fq$ be a family of Hodge-Pink structures on $(D,F_D)$ over $X_L$. Consider the associated pair $(\CP,\CQ)$ of $\sigma$-modules over $\CO_{X_L}\ancon$. 
Assume that $\CQ$ contains an $\CO_{X_L}\langle\frac{z}{\zeta}\rangle$-lattice $\CQ'$ on which $F_\CQ$ is an isomorphism. Then there exists a quasi-paracompact admissible formal model $X$ of $X_L$ over $\Spf R$, and a bounded rigidified local shtuka $(M,F_M,\delta_M)$ over $X$
with $M\otimes_{\CO_X\dbl z\dbr}\CO_{X_0}\dbl z\dbr[z^{-1}] \;\cong\; D\otimes_{k\dpl z\dpr}\CO_{X_0}\dbl z\dbr[z^{-1}]$ and such that the $\CO_{X_L}\dbl z-\zeta\dbr$-lattices $\Fq$ and $\sigma^\ast\delta_M\circ F_M^{-1}\bigl(M\otimes_{\CO_X\dbl z\dbr}\CO_{X_L}\dbl z-\zeta\dbr\bigr)$ coincide inside
\[
\TS \sigma^\ast\delta_M\left(\sigma^\ast M\otimes_{\CO_X\dbl z\dbr}\CO_{X_L}\dbl z-\zeta\dbr[\frac{1}{z-\zeta}]\right)\es=\es\sigma^\ast D\otimes_{k\dpl z\dpr}\CO_{X_L}\dbl z-\zeta\dbr[\frac{1}{z-\zeta}]\,.
\]
\end{proposition}

\begin{proof}
Let $\{X_L^i\}_i$ be an affinoid covering of $X_L$ which is locally finite. For each $i$ we will construct a quasi-compact admissible formal $R$-scheme $X^i$ and a rigidified local shtuka over $X^i$ as in the assertion which satisfies a maximality property by which it is uniquely defined. Then by \cite[Theorem 2.8/3]{Bosch} there exist a quasi-paracompact admissible formal $R$-model $X$ of $X_L$, open formal subschemes $U^i\subset X$, and admissible formal blowing-ups $U^i\to X^i$. By their maximality property the pullbacks to $U^i$ of the rigidified local shtukas glue yielding the desired rigidified local shtuka $M$ on $X$.

Now fix an $i$ and an admissible formal $R$-model $X^i=\Spf B\open$ of $X_L^i$.
During the proof we will replace $\Spf B\open$  by various admissible blowing-ups $X'$ of $\Spf B\open$. To ease notation we will then work on an affine open subset of $X'$. The reader should note however, that all our operations and the further admissible blowing-ups we need are globally defined on $X'$. 

Using Lemma~\ref{LemmaTensorBounded} we may tensor with a suitable Tate object and therefore assume that $F_D\in M_n\bigl(k\dbl z\dbr\bigr)\cdot\sigma$ and $\CP\subset\CQ$. We fix a finite free $k[z]$-submodule $E$ of $\sigma^\ast D$ with $\sigma^\ast D=E\otimes_{k[z]}k\dpl z\dpr$ and consider the sheaves $\CN:=E\otimes_{k[z]}B\open\langle z,z^{-1}\rangle$ on $\Spf B\open\langle z,z^{-1}\rangle$ and $\ol\CP:=E\otimes_{k[z]}B\langle z\rangle$ on $\Spm B\langle z\rangle$. Then $\CP=\ol\CP\otimes_{B\langle z\rangle}B\ancon$.
Consider the $B\langle\frac{z}{\zeta}\rangle$-lattice $\CQ'\subset\CQ$. By construction $\CQ'$ and $\ol\CP$ coincide on $\Spm B\langle\frac{z}{\zeta},\frac{\zeta}{z}\rangle\setminus\Var(z-\zeta)$ and are related at $z=\zeta$ by
\[
\CQ'\otimes B\dbl z-\zeta\dbr\es=\es\Fq\es\supset\es\Fp\es=\es\ol\CP\otimes B\dbl z-\zeta\dbr\,.
\]
We glue $\CQ'$ and $\ol\CP$ over $\Spm B\langle\frac{z}{\zeta},\frac{\zeta}{z}\rangle\setminus\Var(z-\zeta)$ to a locally free sheaf $\CF_L$ of finite rank on $\Spm B\langle z\rangle$ which coincides with $\CQ'$ on $\Spm B\langle\frac{z}{\zeta}\rangle$ and with $\ol\CP$ on $\Spm B\langle z,\frac{\zeta}{z}\rangle\setminus\Var(z-\zeta)$. The isomorphisms $\sigma^\ast F_D$ and $ F_\CQ$ define a morphism $ F_{\CF_L}:\sigma^\ast\CF_L\otimes_{B\langle z\rangle}B\{z\}\to\CF_L\otimes_{B\langle z\rangle}B\{z\}$ which is an isomorphism outside $z=\zeta$ and satisfies $\coker F_{\CF_L}=\Fq/\Fp$. We can think of $\CF_L\otimes_{B\langle z\rangle}B\{z\}$ as the rigid analytic part of the desired rigidified local shtuka $M$. Our task is now to actually construct the formal model $M$ of $\CF_L\otimes_{B\langle z\rangle}B\{z\}$.

Note that $\CF_L\otimes_{B\langle z\rangle}B\langle z,z^{-1}\rangle=\CN\otimes_{B\open\langle z,z^{-1}\rangle}B\langle z,z^{-1}\rangle$.
By Lemma~\ref{Lemma2.17} we may replace $\Spf B\open$ by an admissible blowing-up and assume that $\CF_L$ and $\CN$ both come from a locally free sheaf $\CF$ of finite rank on $\Spf B\open\langle z\rangle$ with the maximality property that for any morphism $\beta:\Spf C\open\to \Spf B\open$ of admissible formal $\Spf R$-schemes 
\[
\TS\Gamma\bigl(\Spf C\open\langle z\rangle\,,\,\beta^\ast\CF\bigr)\es=\es\Gamma\bigl(\Spm (C\open\otimes_R L)\langle z\rangle\,,\,\beta^\ast\CF_L\bigr)\;\cap\;\Gamma\bigl(\Spf C\open\langle z,z^{-1}\rangle\,,\,\beta^\ast\CN\bigr)
\]
inside $\Gamma\bigl(\Spm (C\open\otimes_R L)\langle z,z^{-1}\rangle\,,\,\beta^\ast\CF_L\bigr)$. This property uniquely determines $\CF$.
We put $\CM:=\CF\otimes_{B\open\langle z\rangle}B\open\dbl z\dbr$. 
The morphism $\sigma^\ast F_D$ induces an isomorphism 
\[
F_\CN\es:=\es\sigma^\ast F_D\otimes\id:\es\sigma^\ast\CN\otimes_{B\open\langle z,z^{-1}\rangle}B\open\dbl z,z^{-1}\rangle \es\longto\es\CN\otimes_{B\open\langle z,z^{-1}\rangle}B\open\dbl z,z^{-1}\rangle\,,
\]
where \quad $\DS\es B\open\dbl z,z^{-1}\rangle\; :=\;\bigl\{\;\sum_{i=-\infty}^\infty b_i z^i:\es b_i\in B\open,\; |b_i|\to0\es(i\to -\infty)\;\bigr\}$\,.

\noindent
Then $F_{\CF_L}$ and $F_\CN$ extend by Lemma~\ref{Lemma2.18} to a morphism $F_\CM:\sigma^\ast\CM\to\CM$. We determine its cokernel. Since $\CF$ is locally free, $\CM$ and $\sigma^\ast\CM$ are locally on $\Spf B\open$ isomorphic to $B\open\dbl z\dbr^n$. After we fix bases of $\CM$ and $\sigma^\ast\CM$ the morphism $F_\CM$ corresponds to a matrix $\Phi\in M_n\bigl(B\open\dbl z\dbr\bigr)$. Consider the commutative diagram with exact rows
\[
\xymatrix {
0\ar[r] & *+<0.5pc,1pc> \objectbox{\sigma^\ast\CM} \ar[r]^{F_\CM}\ar@{^{ (}->}[d] & *+<0.5pc,1pc> \objectbox{\CM} \ar[r]\ar@{^{ (}->}[d] & V \ar[r] \ar[d]& 0 \\
0\ar[r] & **{!U(0.15) +<0.2pc,0.5pc>} \objectbox{\sigma^\ast\CM\otimes_{B\open\dbl z\dbr}B\{z\} \ar[r]^{F_{\CF_L}}} & **{!U(0.15) +<0.2pc,0.5pc>} \objectbox{\CM\otimes_{B\open\dbl z\dbr}B\{z\}} \ar[r] & W \ar[r] & 0 \,.
}
\]
Since $W=\coker F_{\CF_L}$ is killed by $(z-\zeta)^e$ for some $e\in\BN_0$, we have $(z-\zeta)^e\CM\otimes_{B\open\dbl z\dbr}B\{z\}\subset\im F_{\CF_L}$. In particular there is a matrix $\Psi\in M_n\bigl(B\{z\}\bigr)$ representing the morphism $F_{\CF_L}^{-1}\circ(z-\zeta)^e$ with $\Phi\Psi=(z-\zeta)^e$. Now $\Phi\in\GL_n\bigl(B\open\dbl z,z^{-1}\rangle\bigr)$ since $\CM\otimes_{B\open\dbl z\dbr}B\open\dbl z,z^{-1}\rangle\;=\;\sigma^\ast D\otimes_{k\dpl z\dpr}B\open\dbl z,z^{-1}\rangle$ and $\sigma^\ast F_D$ is an isomorphism. Therefore $\Psi$ has its coefficients in $B\{z\}\cap B\open\dbl z,z^{-1}\rangle\;=\;B\open\dbl z\dbr$. We deduce that $(z-\zeta)^e\CM$ lies in the image of $F_\CM$ and $V$ is killed by $(z-\zeta)^e$. Consider the exact sequence of $B\open\langle\TS\frac{z}{\zeta}\rangle$-modules
\[
\xymatrix {
0\ar[r] & \sigma^\ast\CM\otimes_{B\open\dbl z\dbr}B\open\langle\TS\frac{z}{\zeta}\rangle \ar[r]^{F_\CM} & \CM\otimes_{B\open\dbl z\dbr}B\open\langle\TS\frac{z}{\zeta}\rangle \ar[r] & V\ar[r] & 0\,.
}
\]
By \cite[Theorem 4.1]{FRG2} there exists an admissible blowing-up $X''$ of $\Spf B\open$ such that on any affine open subset $\Spf C\open$ of $X''$ the $C\open$-module $V\otimes_{B\open}C\open/(\zeta\text{-torsion})$ is locally free of finite rank. Since $V\otimes_{B\open\langle\frac{z}{\zeta}\rangle}C\open\langle\frac{z}{\zeta}\rangle=V\otimes_{B\open\langle\frac{z}{\zeta}\rangle/(z-\zeta)^e}C\open\langle\frac{z}{\zeta}\rangle/(z-\zeta)^e=V\otimes_{B\open}C\open$ we have an exact sequence of $C\open\langle\TS\frac{z}{\zeta}\rangle$-modules
\[
\xymatrix {
0\ar[r] & \Tor_1^{B\open\langle\frac{z}{\zeta}\rangle}(V,C\open\langle\frac{z}{\zeta}\rangle)\ar[r] & \sigma^\ast\CM\otimes_{B\open\dbl z\dbr}C\open\langle\TS\frac{z}{\zeta}\rangle\ar[r] & \CM\otimes_{B\open\dbl z\dbr}C\open\langle\TS\frac{z}{\zeta}\rangle\ar[r] & V\otimes_{B\open}C\open \ar[r] & 0\,.
}
\]
The $C\open\langle\frac{z}{\zeta}\rangle$-module $\Tor_1^{B\open\langle\frac{z}{\zeta}\rangle}(V,C\open\langle\frac{z}{\zeta}\rangle)$ is killed by $(z-\zeta)^e$. So it must vanish since $\sigma^\ast\CM\otimes_{B\open\dbl z\dbr}C\open\langle\TS\frac{z}{\zeta}\rangle$ is free. Also by \cite[Lemma 4.5]{FRG2}, $V\otimes_{B\open}C\open$ has no $\zeta$-torsion and is therefore a locally free $C\open$-module of finite rank. Replacing $B\open$ by $C\open$ we conclude that $(\CM,F_\CM:\sigma^\ast\CM\to\CM)$ is a local shtuka over $X''$. It is bounded by $(d,0,\ldots,0)$ for $d=\rk_{B\open}V$. 
To construct a rigidification on $\CM$ we define
\begin{tabbing}
$B\con[n]$\es\= \kill
$B\langle z,z^{-1}\}$\> $\DS =\es\bigl\{\,\sum_{i=-\infty}^\infty b_iz^i:\es
b_i\in B,\;|b_i|\to0\;(i\to\infty)\;,|b_i|\,|\zeta|^{ri}\to0\;(i\to-\infty)\;\forall r>0\,\bigr\}$ \quad and\\[2mm]
$B\open\langle z,z^{-1}\}$\> $\DS =\es\bigl\{\,\sum_{i=-\infty}^\infty b_iz^i:\es b_i\in B\open,\;|b_i|\to0 \;(i\to\infty)\;,|b_i|\,|\zeta|^{ri}\to0\;(i\to-\infty)\;\forall r>0\,\bigr\}$\,.
\end{tabbing}
Then there is an exact sequence of $B\open\langle z,z^{-1}\}$-modules
\begin{equation} \label{Eq2.4.9}
\xymatrix @R=0pc {0\ar[r] & B\open\langle z,z^{-1}\} \ar[r] & B\open\langle z,z^{-1}\rangle \oplus B\langle z,z^{-1}\} \ar[r] & B\langle z,z^{-1}\rangle\\
& {\qquad f\qquad} \ar@{|->}[r] & {\quad(f,f)\es,\qquad (g,h)\quad}\ar@{|->}[r]& {\quad g-h\quad}
}
\end{equation}
Moreover, $\CF\otimes_{B\open\langle z\rangle}B\langle z,z^{-1}\}=\CF_L\otimes_{B\langle z\rangle}B\langle z,z^{-1}\}$, and $E\otimes_{k[z]}B\langle z,z^{-1}\}=\ol{\CP}\otimes_{B\langle z\rangle}B\langle z,z^{-1}\}$, and $\CF\otimes_{B\open\langle z\rangle}B\open\langle z,z^{-1}\rangle =\CN=E\otimes_{k[z]}B\open\langle z,z^{-1}\rangle$. The inclusions $\tminus^e\CQ\subset\CP\subset\CQ$ induce inclusions
\[
\tminus^e\cdot\bigl(\CF_L\otimes_{B\langle z\rangle}B\langle z,z^{-1}\}\bigr)\;\subset\;\ol\CP\otimes_{B\langle z\rangle}B\langle z,z^{-1}\}\;\subset\CF_L\otimes_{B\langle z\rangle}B\langle z,z^{-1}\}\,.
\]
Tensoring (\ref{Eq2.4.9}) with the locally free $B\open\langle z,z^{-1}\}$-modules $\CF\otimes_{B\open\langle z\rangle}B\open\langle z,z^{-1}\}$ and $E\otimes_{k[z]}B\open\langle z,z^{-1}\}$ yields
\[
\xymatrix @C=1.5pc {
0\ar[r]&\tminus^e\!\cdot\!\bigl(\CF\otimes_{B\open\langle z\rangle}B\open\langle z,z^{-1}\}\bigr) \ar@{^{ (}-->}[d] \ar[r] & \tminus^e\!\cdot\!\CN\oplus\tminus^e\!\cdot\!\bigl(\CF_L\otimes_{B\langle z\rangle}B\langle z,z^{-1}\}\bigr) \ar@{^{ (}->}[d] \ar[r] & \tminus^e\!\cdot\!\bigl(\CF_L\otimes_{B\langle z\rangle}B\langle z,z^{-1}\rangle\bigr)\ar@{=}[d] \\
0\ar[r] & E\otimes_{k[z]}B\open\langle z,z^{-1}\} \ar@{^{ (}-->}[d] \ar[r] & \CN \oplus\ol\CP\otimes_{B\langle z\rangle}B\langle z,z^{-1}\} \ar@{^{ (}->}[d] \ar[r] & \CN\otimes_{B\open\langle z,z^{-1}\rangle}B\langle z,z^{-1}\rangle \ar@{=}[d]\\
0\ar[r]&\CF\otimes_{B\open\langle z\rangle}B\open\langle z,z^{-1}\} \ar[r] & \CN \oplus \CF_L\otimes_{B\langle z\rangle}B\langle z,z^{-1}\} \ar[r] & \CF_L\otimes_{B\langle z\rangle}B\langle z,z^{-1}\rangle
}
\]
Tensoring the left column further with $B\open\dbl z,z^{-1}\}[\tminus^{-1}]$ yields the equality in the upper row of
\[
\xymatrix {
\CM\otimes_{\CO_{X''}\dbl z\dbr}\CO_{X''}\dbl z,z^{-1}\}[\tminus^{-1}] \ar@{=}[r] \ar[dr]_{\delta_\CM} &
\sigma^\ast D\otimes_{k\dpl z\dpr}\CO_{X''}\dbl z,z^{-1}\}[\tminus^{-1}] \ar[d]^{F_D\otimes\id}\\
 & D\otimes_{k\dpl z\dpr}\CO_{X''}\dbl z,z^{-1}\}[\tminus^{-1}]
}
\]
This is the desired rigidification on $\CM$. It satisfies
\[
\CM\otimes_{\CO_{X''}\dbl z\dbr}\CO_{X''_0}\dbl z\dbr[z^{-1}]\es=\es\sigma^\ast D\otimes_{k\dpl z\dpr}\CO_{X''_0}\dbl z\dbr[z^{-1}]\es\xrightarrow{{}_\sim\es\; F_D\otimes\id\;}\es D\otimes_{k\dpl z\dpr}\CO_{X''_0}\dbl z\dbr[z^{-1}]
\]
and $\sigma^\ast\delta_\CM\circ F_\CM^{-1}\bigl(\CM\otimes_{\CO_{X''}\dbl z\dbr}B\dbl z-\zeta\dbr\bigr)\;=\;\Fq$ inside $\sigma^\ast D\otimes_{k\dpl z\dpr}B\dbl z-\zeta\dbr[\frac{1}{z-\zeta}]$.
\end{proof}

By similar arguments one can also prove the following

\begin{proposition}\label{PropHFaithfull}
The functor $\BH:(M,F_M,\delta_M)\mapsto (D,F_D,\Fq)$ between the category of rigidified local shtukas over $\CO_K$ up to isogeny and the category of $z$-isocrystals with Hodge-Pink structure over $K$ is fully faithful.
\end{proposition}

\begin{proof}
The faithfulness follows from the faithfulness of $\pair\BH:M\mapsto(\CP,\CQ)$ which is a consequence of Proposition~\ref{Prop7.5}. To prove fullness let $(M_i,F_i,\delta_i)$ for $i=1,2$ be rigidified local shtukas over $\CO_K$. Let $\ul D_i:=\BH(M_i,F_i,\delta_i)$ be their associated $z$-isocrystals with Hodge-Pink structure over $K$ and let $(\CP_i,\CQ_i):=\pair(\ul D_i)\cong\pair\BH(M_i,F_i,\delta_i)$ be their associated pairs of $\sigma$-modules over $K\ancon$. Consider the $\sigma$-module $\CQ'_i:=M_i\otimes_{\CO_K\dbl z\dbr}K\langle\frac{z}{\zeta}\rangle$ over $K\langle\frac{z}{\zeta}\rangle$. It satisfies $\CQ_i\cong\CQ'_i\otimes K\ancon$. Let $f:\ul D_1\to\ul D_2$ be a morphism and let $f:\CQ_1\to\CQ_2$ be the induced morphism. Then by Proposition~\ref{Prop13}(a)
\[
f\in\CHom_\sigma(\CQ_1,\CQ_2)^F(K)\es=\es\CHom_\sigma\bigl(\CQ'_1\otimes K\con,\CQ'_2\otimes K\con\bigr)^F(K)\,.
\]
Therefore we may fix an integer $n$ such that $z^nf$ induces a morphism of $\sigma$-modules $z^nf:\CQ'_1\to\CQ'_2$.

Fix a rational number $r$ with $0<r<1$ and let $\ol\CP_i:=\sigma^\ast D_i\otimes_{k\dpl z\dpr}K\langle\frac{z}{\zeta^{r/q}},\frac{\zeta}{z}\rangle$. The $K\langle\frac{z}{\zeta^{r/q}}\rangle$-module $M_i\otimes_{\CO_K\dbl z\dbr}K\langle\frac{z}{\zeta^{r/q}}\rangle$ is obtained as the gluing of $\ol\CP_i$ with $\CQ'_i$ via $F_i\circ(\sigma^\ast\delta_i)^{-1}$ over $\Spm K\langle\frac{z}{\zeta},\frac{\zeta}{z}\rangle\setminus\Var(z-\zeta)$.
Moreover, since $M_i$ is a (locally) free $\CO_K\dbl z\dbr$-module we can recover $M_i$ as the intersection
\[
\TS M_i\es=\es M_i\otimes_{\CO_K\dbl z\dbr}K\langle\frac{z}{\zeta^{r/q}}\rangle\es\cap\es M_i\otimes_{\CO_K\dbl z\dbr}\CO_K\dbl z,\frac{\zeta^r}{z}\rangle[z^{-1}]
\]
inside $M_i\otimes_{\CO_K\dbl z\dbr}K\langle\frac{z}{\zeta^{r/q}},\frac{\zeta^r}{z}\rangle$, where we set
\[
\CO_K{\TS\dbl z,\frac{\zeta^r}{z}\rangle[z^{-1}]}\es:=\es\bigl\{\,\sum_{i=-\infty}^\infty b_iz^i:\es b_i\in\CO_K\,,\,|b_i\zeta^{ri}|\to0\;(i\to-\infty)\,\bigr\}\,.
\]
By what was said above the morphism $z^nf:\ul D_1\to \ul D_2$ induces the two compatible morphisms $z^nf:M_1\otimes_{\CO_K\dbl z\dbr}K\langle\frac{z}{\zeta^{r/q}}\rangle\to M_2\otimes_{\CO_K\dbl z\dbr}K\langle\frac{z}{\zeta^{r/q}}\rangle$ and
\[
\xymatrix {
M_1\otimes_{\CO_K\dbl z\dbr}\CO_K\dbl z,\frac{\zeta^r}{z}\rangle[z^{-1}] \ar[r]^{z^nf} \ar[d]^{\delta_1} &
M_2\otimes_{\CO_K\dbl z\dbr}\CO_K\dbl z,\frac{\zeta^r}{z}\rangle[z^{-1}] \ar[d]^{\delta_2} \\
D_1\otimes_{k\dpl z\dpr}\CO_K\dbl z,\frac{\zeta^r}{z}\rangle[z^{-1}] \ar[r]^{z^nf} &
D_2\otimes_{k\dpl z\dpr}\CO_K\dbl z,\frac{\zeta^r}{z}\rangle[z^{-1}]
}
\]
and therefore the desired morphism $z^nf:M_1\to M_2$ of rigidified local shtukas.
\end{proof}

\noindent
{\it Remark.}
Let $X$ be an arbitrary admissible formal scheme over $\Spf R$ and let $(D,F_D)$ be a $z$-isocrystal over $k$. Consider the category of bounded rigidified local shtukas $(M,F_M,\delta_M)$ over $X$ with constant $z$-isocrystal $(M,F_M)\otimes_{\CO_X\dbl z\dbr}\CO_{X_0}\dbl z\dbr[z^{-1}] \;\cong\; (D,F_D)\otimes_{k\dpl z\dpr}\CO_{X_0}\dbl z\dbr[z^{-1}]$. Here the rigidification consists of an isomorphism
\[
\delta_M:\es M\otimes_{\CO_{X}\dbl z\dbr}\CO_{X}\dbl z,z^{-1}\}[\tminus^{-1}]\es\isoto\es D\otimes_{k\dpl z\dpr}\CO_{X}\dbl z,z^{-1}\}[\tminus^{-1}]
\]
with $\delta_M\circ F_M=F_D\circ\sigma^\ast\delta_M$ and $\delta_M\otimes\id_{\CO_{X_0}\dbl z\dbr[z^{-1}]}=\id_D$. Then one can define the functor $\BH_X:(M,F_M,\delta_M)\mapsto (D,F_D,\Fq)$ from rigidified local shtukas $(M,F_M,\delta_M)$
 to families of Hodge-Pink structures on $(D,F_D)$ like in case $X=\Spf\CO_K$ and the above proof shows the full faithfulness also for $\BH_X$ over $X$. Genestier and Lafforgue \cite[\S2]{GL} have worked out this point of view in detail.

%%%%%%%%%%%%%%%%%%%%%%%%%%%%%%%%%%%%%%%%%%%%%%%%%%%%%%%%%%%%%%%%%%%%%%
%
%    Weakly Admissible Implies Admissible
%
%%%%%%%%%%%%%%%%%%%%%%%%%%%%%%%%%%%%%%%%%%%%%%%%%%%%%%%%%%%%%%%%%%%%%%

\subsection{Weakly Admissible Implies Admissible} \label{SectWA=>A}

As before assume that there is a fixed section $k\hookrightarrow\CO_K$ of the residue map $\CO_K\to k$.
Assume that $K$ satisfies condition (\ref{CondDoubleStar}) from page~\pageref{CondDoubleStar}.
Let $\ul D$ be a $z$-isocrystal with Hodge-Pink structure over $K$. 
We want to show in this section that if $\ul D$ is weakly admissible then it is already admissible. This fact was first proved by Genestier and Lafforgue \cite{GL} under the assumption that $K$ is discretely valued with perfect residue field $k$. Here we give an entirely different proof which instead parallels the proofs of Berger~\cite{Berger1} and Kisin~\cite{Kisin} of the analogous result in mixed characteristic.
As a preparation we need the following lemma which holds even without the above assumptions on $K$.

\begin{lemma} \label{Lemma3.2}
Let $\ol M$ be a $\sigma$-module over $K\overcon$ and let $N\subset \ol M\otimes K\ancon$ be a saturated $\sigma$-submodule over $K\ancon$. Then there exists a unique saturated $\sigma$-submodule $\ol N\subset \ol M$ over $K\overcon$ with $\ol N\otimes K\ancon=N$. 
\end{lemma}

\begin{proof}
Repeated application of Lemma~\ref{Lemma9b} allows us to extend $N\otimes K\langle\frac{z}{\zeta},\frac{\zeta^q}{z}\rangle$ to a saturated $\sigma$-submodule $N_i$ of $\ol M\otimes K\langle\frac{z}{\zeta^{q^{-i}}},\frac{\zeta^q}{z}\rangle$ for all $i\in\BN_0$. Gluing the $N_i$ with $N$ over $K\langle\frac{z}{\zeta},\frac{\zeta^q}{z}\rangle$ defines a saturated subsheaf $\ol N$ of $\ol M$ which is $F_{\ol M}$-invariant. By definition $\ol N$ is the desired saturated $\sigma$-submodule of $\ol M$.
\end{proof}

The key to showing that weakly admissible implies admissible is the following proposition.

\begin{proposition} \label{Prop3.3}
Assume that $K$ satisfies condition (\ref{CondDoubleStar}) from page~\pageref{CondDoubleStar}.
Let $\ul D$ be a Hodge-Pink structure  over $K$. Consider the associated pair $(\CP,\CQ)=\pair(\ul D)$ of $\sigma$-modules over $K\ancon$. Let $\CQ'\subset\CQ$ be a saturated $\sigma$-submodule over $K\ancon$ and let $\CP'$ be the saturation inside $\CP$ of $\CP\cap\CQ'$. Then there exists a subobject $\ul D'$ of $\ul D$ with $(\CP',\CQ')=\pair(\ul D')$.
\end{proposition}

\begin{proof}
If the subobject $\ul D'=(D',F_{D'},\Fq_{D'})$ of $\ul D=(D,F_D,\Fq_D)$ exists we necessarily must have $\Fq_{D'}=\Fq_D\cap\sigma^\ast D'\otimes_{k\dpl z\dpr}K\dpl z-\zeta\dpr$ since $\CQ'\subset\CQ$ is saturated. Conversely it suffices to establish the existence of an $F_D$-stable subspace $D'\subset D$ with $\sigma^\ast(D',F_D|_{D'})\otimes_{k\dpl z\dpr}K\ancon\; =\;\CP'$ as $\sigma$-submodules of $\CP$, since then $\Fq_D\cap\sigma^\ast D'\otimes_{k\dpl z\dpr}K\dpl z-\zeta\dpr$ is the desired Hodge-Pink structure.

\medskip

(a) \es First we claim that it suffices to prove the existence of $D'$ under the additional assumption that $k$ is algebraically closed. Namely, let $k^\alg$ be the algebraic closure of $k$ inside $\ol K$ and let $\wt K$ be the closure inside $\ol K$ of the compositum $k^\alg\cdot K$. Assume we are given an $F_D$-stable $k^\alg\dpl z\dpr$-subspace $D'$ such that 
\[
\sigma^\ast D'\otimes_{k^\alg\dpl z\dpr}\wt K\ancon\es=\es\CP'\otimes\wt K\ancon
\]
as submodules of $\sigma^\ast D\otimes_{k\dpl z\dpr}\wt K\ancon$. Set $\Fq_{D'}:=\Fq_D\otimes\wt K\dbl z-\zeta\dbr\cap \sigma^\ast D'\otimes_{k^\alg\dpl z\dpr}\wt K\dpl z-\zeta\dpr$ and assume that $\pair(D',F_{D'},\Fq_{D'})=(\CP',\CQ')\otimes\wt K\ancon$. We must show that $D'$ descends to $k\dpl z\dpr$.

By Lemma~\ref{LemmaKe3.6.2} we may assume that $\dim_{k^\alg\dpl z\dpr}D'\;=\;\rk_{K\ancon}\CP'\;=\;1$. Choose a basis $e_1,\ldots,e_n$ of $D$ over $k\dpl z\dpr$. Then $D'$ is generated over $k^\alg\dpl z\dpr$ by a vector $v=a_1 e_1+\ldots+a_n e_n$ with $a_i\in k^\alg\dpl z\dpr$. Without loss of generality, $a_1\neq0$ and we can assume $a_1=1$. On the other hand $\CP'$ is generated over $K\ancon$ by a vector $w=b_1\,\sigma^\ast e_1+\ldots+b_n\,\sigma^\ast e_n$ with $b_i\in K\ancon$. Since $\sigma^\ast D'\otimes\wt K\ancon=\CP'\otimes\wt K\ancon$ there exists a unit $a\in\wt K\ancon\mal$ with $w=a\,\sigma^\ast v$. We find that $b_1=a\in\wt K\ancon\mal\cap K\ancon=K\ancon\mal$ (for the last equality use \cite[Lemma 9.7.1/1]{BGR}). This yields $a_i^\sigma=a^{-1}b_i\in k^\alg\dpl z\dpr\cap K\ancon$. By expanding in powers of $z$ we conclude that $a_i^\sigma\in k\dpl z\dpr$ because $k^\alg\cap K=k$.
% REASON: $0\neq x\in k^\alg\cap K$ implies $|x|=1$ implies $x\in\CO_K$ implies $\exists \bar x\in k$ with $x-\bar x\in\Fm_K$ implies $x-\bar x\in k^\alg\cap\Fm_K=(0)$.
So $\sigma^\ast D'$ descends to $k\dpl z\dpr$. Since $F_D$ maps $\sigma^\ast D'$ isomorphically onto $D'$ also $D'$ descends to $k\dpl z\dpr$ as desired.
% REASON: $F_D(\sigma^\ast v)=\alpha v=\alpha e_1+\ldots$ implies $\alpha\in k\dpl z\dpr\mal$.

\medskip

(b) \es From now on we work over $\wt K$. By what was said above it suffices to assume that $\rk\CP'=1$ and to prove that $\CP'\otimes\wt K\ancon$ descends to a $\sigma$-submodule $\sigma^\ast D'\subset\sigma^\ast D\otimes k^\alg\dpl z\dpr$. Then $D'=(\sigma^{-1})^\ast\sigma^\ast D'\subset D\otimes k^\alg\dpl z\dpr$. Lemma~\ref{Lemma3.2} allows us to extend $\CP'\otimes \wt K\ancon$ to a saturated $\sigma$-submodule $\ol\CP{}'$ of $\ol\CP:=\sigma^\ast D\otimes_{k\dpl z\dpr}\wt K\overcon$. Now consider the  slope decomposition
\[
(D\otimes k^\alg\dpl z\dpr,F_D)\es\cong\es\bigoplus_{i=1}^\ell\CF_{d_i,n_i}^{\oplus m_i}
\]
from Theorem~\ref{Thm2} and set $D_\lambda:=\CF_{d_i,n_i}^{\oplus m_i}$ for $\lambda=\frac{d_i}{n_i}$. It induces the slope decomposition $\ol\CP\cong\bigoplus_\lambda\ol\CP_\lambda$ with $\ol\CP_\lambda:=\sigma^\ast D_\lambda\otimes_{k^\alg\dpl z\dpr}\wt K\overcon$. Let $d:=\deg \ol\CP{}'$. We claim that $\ol\CP{}'\subset\ol\CP_{d/1}$.

To prove the claim consider the induced morphisms $f_\lambda:\ol\CP{}'\hookrightarrow\ol\CP\to\ol\CP_\lambda$. By Proposition~\ref{Prop0.7}, $f_\lambda=0$ if $\lambda<d$. Thus it suffices to show that $f_\lambda=0$ also for $\lambda>d$. Twisting with $\CO(-d)$ allows us to assume that $d=0$. Set $\frac{e}{s}:=\lambda>0$ with positive integers $e$ and $s$ and assume $f_\lambda\neq0$. Replacing $\sigma$ by $\sigma^s$ (and $q$ by $q^s$) we may assume that $\ol\CP_\lambda=\CO(e)^{\oplus n}$, that is $F_{\ol\CP_\lambda}=z^{-e}\cdot\sigma$. Projecting onto a suitable summand we find a proper inclusion $\ol\CP{}'\subset\CO(e)$. There exists a rational number $r$ with $1\leq r<q$ such that the cokernel of $\ol\CP{}'\subset\CO(e)$ is supported on $\bigcup_{\nu\in\BZ}\Spm \wt K\langle\frac{z}{\zeta^{q^{-\nu}}},\frac{\zeta^{rq^{-\nu}}}{z}\rangle$. We tensor the inclusion $\ol\CP{}'\subset\CO(e)$ with the principal ideal domain $\wt K\langle\frac{z}{\zeta},\frac{\zeta^r}{z}\rangle$ and we let $\sum_{\ell=0}^e a_\ell z^\ell\in \wt K[z]$ be the corresponding elementary divisor. Note that the length of the cokernel $\bigl(\CO(e)/\ol\CP{}'\bigr)\otimes \wt K\langle\frac{z}{\zeta},\frac{\zeta^r}{z}\rangle$ is $e$ by Proposition~\ref{Prop0.5}.
We may assume $a_e=1$. Then the $a_\ell$ are uniquely determined by Lemma~\ref{LemmaPID}. In particular $0<|a_0|<1$ since otherwise the length of the cokernel would be less than $e$. The $\sigma$-invariance of the inclusion $\ol\CP{}'\subset\CO(e)$ implies that we find $\sum_{\ell=0}^e a_\ell^{1/q}z^\ell\in \wt K[z]$ as the elementary divisor after tensoring with $\wt K\langle\frac{z}{\zeta^{1/q}},\frac{\zeta^{r/q}}{z}\rangle$. Proceeding in this manner we obtain $a_0^{1/q^{\nu}}\in \wt K$ for arbitrarily large $\nu$. But this is impossible since $K$ satisfies condition (\ref{CondDoubleStar}) from page~\pageref{CondDoubleStar}. Thus $f_\lambda=0$ for all $\lambda\neq d$ and so $\ol\CP{}'\subset\ol\CP_{d/1}$.

\medskip

(c) \es After twisting $D$ with $\CO(-d)$ we obtain $\CP'\subset\CP_\lambda$ for $\lambda=0$. We let $e_1,\ldots,e_n$ be a basis of $D_\lambda$ on which $F_D$ acts as $\sigma$. Then $F_{\CP_\lambda}=\sigma$ and one sees that $\sigma^\ast D_\lambda=(\CP_\lambda)^F(\wt K)\otimes_{\BF_q\dpl z\dpr}k^\alg\dpl z\dpr$. We claim that 
\[
\sigma^\ast D'\es:=\es(\CP')^F(\wt K)\otimes_{\BF_q\dpl z\dpr}k^\alg\dpl z\dpr
\]
is the desired $\sigma$-submodule of $\sigma^\ast D\otimes k^\alg\dpl z\dpr$. Let $v=a_1\,\sigma^\ast e_1+\ldots+a_n\,\sigma^\ast e_n$ with $a_i\in\wt K\ancon$ be a generator of $\CP'$. Since $\deg\CP'=0$ there exists an isomorphism $\CP'\otimes\ol K\ancon\cong\CO(0)$ of $\sigma$-modules over $\ol K\ancon$. Thus we find a unit $b\in\ol K\ancon\mal$ with $F_{\CP_\lambda}(b^\sigma\sigma^\ast v)=bv$, whence $(ba_i)^\sigma=ba_i\in\BF_q\dpl z\dpr$. It follows that $b\in\wt K\ancon\mal$. Replacing $v$ by $bv$ yields $v\in\sigma^\ast D_\lambda$ and so $\CP'$ descends to a $\sigma$-submodule $\sigma^\ast D'$ with $\sigma^\ast D'\otimes_{k^\alg\dpl z\dpr}\wt K\ancon\cong\CP'\otimes\wt K\ancon$. This proves the proposition.
\end{proof}

\begin{theorem} \label{Thm3.5}
Assume that $K$ satisfies condition (\ref{CondDoubleStar}) from page~\pageref{CondDoubleStar}.
Then the functor $\BH$ is a tensor equivalence between the category of rigidified local shtukas over $\CO_K$ up to isogeny and the category of $z$-isocrystals with Hodge-Pink structure over $K$ which are weakly admissible. In other words, weakly admissible implies admissible.
\end{theorem}

\noindent
{\it Remark.} 
We will show in Example~\ref{Ex8.2} that weakly admissible need not imply admissible if $K$ does not satisfy condition (\ref{CondDoubleStar}).

\medskip

\noindent
{\it Remark.} 
To define a Hodge-Pink structure on a fixed $z$-isocrystal one needs only finitely many elements from $K$. Therefore every $z$-isocrystal with Hodge-Pink structure can be defined over a field $K$ which is topologically finitely generated over its residue field (which it contains by assumption). It might at first glance seem surprising that such a field could fail to satisfy condition (\ref{CondDoubleStar}). See however Example~\ref{ExAnalyticPoints}\ref{ExAnalyticPointsD} which shows that such a field can even be algebraically closed.

\begin{proof}
Let $(\CP,\CQ)=\pair(\ul D)$. By Theorem~\ref{Thm2.6} we must show that $\CQ$ is isoclinic of slope zero. By Proposition~\ref{Prop2.8} the weak admissibility implies that $\deg\CQ=0$. Assume that $\CQ$ is not isoclinic. Then the Slope Filtration Theorem~\ref{Thm6.13} yields a saturated $\sigma$-submodule $\CQ'\subset\CQ$ over $K\ancon$ which is isoclinic of negative slope, hence has positive degree. Let $\CP'$ be the saturation inside $\CP$ of $\CP\cap\CQ'$. From Proposition~\ref{Prop3.3} we obtain a subobject $\ul D'$ of $\ul D$ with $\pair(\ul D')=(\CP',\CQ')$. Then $t_H(\ul D')-t_N(\ul D')=\deg\CQ'>0$ by Lemma~\ref{Lemma2.7} contradicting the weak admissibility. This proves the theorem.
\end{proof}

\vspace{2cm}%\pagebreak

%%%%%%%%%%%%%%%%%%%%%%%%%%%%%%%%%%%%%%%%%%%%%%%%%%%%%%%%%%%%%%%%%%%%%%
%
%    Period Spaces
%
%%%%%%%%%%%%%%%%%%%%%%%%%%%%%%%%%%%%%%%%%%%%%%%%%%%%%%%%%%%%%%%%%%%%%%

\section{Period Spaces} \label{ChaptPeriodSpaces}
\setcounter{equation}{0}

This chapter is devoted to the investigation of period spaces for Hodge-Pink structures.
In Sections~\ref{SectPeriodSpaces} and \ref{SectAdmLocusOpen} we construct these period spaces, similarly to the period spaces Rapoport and Zink~\cite{RZ} constructed for Fontaine's filtered isocrystals, and we show that the admissible locus is Berkovich open. We illustrate this theory by describing some examples for period spaces in Section~\ref{SectExamplesPeriodSp}. Finally in Section~\ref{SectConjRZ} we formulate and prove a precise analog of the conjecture of Rapoport and Zink mentioned in the introduction. 

%%%%%%%%%%%%%%%%%%%%%%%%%%%%%%%%%%%%%%%%%%%%%%%%%%%%%%%%%%%%%%%%%%%%%%
%
%    Period Spaces for Hodge-Pink Structures
%
%%%%%%%%%%%%%%%%%%%%%%%%%%%%%%%%%%%%%%%%%%%%%%%%%%%%%%%%%%%%%%%%%%%%%%

\subsection{Period Spaces for Hodge-Pink Structures} \label{SectPeriodSpaces}

In this section our aim is to study period spaces for Hodge-Pink structures which arise similarly to the period spaces of Rapoport and Zink \cite{RZ} in mixed characteristic.
Let $k$ be an extension of $\BF_q$. Let $n$ be a positive integer and let $b\in\GL_n\bigl(k\dpl z\dpr\bigr)$. Let $\ul{\rm Rep}_{\BF_q\dpl z\dpr}\GL_n$ be the category of representations of $\GL_n$ in finite dimensional $\BF_q\dpl z\dpr$-vector spaces. With any object $V$ of $\ul{\rm Rep}_{\BF_q\dpl z\dpr}\GL_n$ we associate a $z$-isocrystal over $k$
\[
(D,F_D)\;=\;\bigl(V\otimes k\dpl z\dpr,\, b\cdot(\id\otimes\sigma)\bigr)\,.
\]
This defines a faithful $\BF_q\dpl z\dpr$-linear exact tensor functor from $\ul{\rm Rep}_{\BF_q\dpl z\dpr}\GL_n$ to the category of $z$-isocrystals over $k$.
% NOT FULL SINCE COMMUTATION WITH $b$ DOES NOT IMPLY COMMUTATION WITH ALL OF $\GL_n$ !!!!!!!!!!!!!!!!!!!!!!!!!!!!!!!!!!!!!!!!!!!!!!!!!!!!!!!!!!
Let $g\in\GL_n\bigl(k\dpl z\dpr\bigr)$ and $b'=gb(g^\sigma)^{-1}$. In this case we say that $b$ and $b'$ are \emph{$\sigma$-conjugate}. Then multiplication with $g$ defines an isomorphism between the tensor functor associated with $b$ and the tensor functor associated with $b'$.

\begin{definition} \label{Def2.2}
A $\sigma$-conjugacy class $\bar b$ in $\GL_n(k\dpl z\dpr)$ is called \emph{decent} if there exists an element $b\in\bar b$ and a positive integer $s$ such that
\[
(b\cdot\sigma)^s\es=\es b\cdot b^\sigma\cdot\ldots\cdot b^{\sigma^{s-1}}\cdot\sigma^s\es=\es \diag(z^{d_1},\ldots,z^{d_n})\cdot\sigma^s
\]
where $d_i/s$ are the slopes of the $z$-isocrystal $(D,F_D)$ over $k\dpl z\dpr$. We call $b$ a \emph{decent element} in $\bar b$ and the above equation a \emph{decency equation for $b$ with the integer $s$}.
\end{definition}

The following result is proved in the same manner as \cite[Corollaries 1.9 and 1.10]{RZ}.

\begin{proposition} \label{Prop2.3}
\begin{enumerate}
\item
If $\bar b$ is decent and $b$ and $s$ are as in Definition~\ref{Def2.2}. Then $b\in\GL_n\bigl(k'\dpl z\dpr\bigr)$ where $k'=\BF_{q^s}\cap k$.
\item
If $b_1,b_2\in\bar b$ satisfy a decency equation with the same integer $s$, then $b_1$ and $b_2$ are $\sigma$-conjugate in $\GL_n\bigl(k'\dpl z\dpr\bigr)$.
\qed
\end{enumerate}
\end{proposition}

\begin{proposition} \label{Prop2.4}
If $k$ is algebraically closed then any $\sigma$-conjugacy class in $\GL_n\bigl(k\dpl z\dpr\bigr)$ is decent.
\end{proposition}

\begin{proof}
This follows directly from Theorem~\ref{Thm2}.
\end{proof}

\begin{lemma} \label{LemmaDecentAndStable}
If $b\in\bar b$ satisfies a decency equation with the integer $s$ then every $b\cdot(\id\otimes\sigma)$-stable subspace of $V\otimes\BF_q^{\,\alg}\dpl z\dpr$ is obtained from a unique $b\cdot(\id\otimes\sigma)$-stable subspace of $V\otimes\BF_{q^s}\dpl z\dpr$ by base change to $\BF_q^{\,\alg}\dpl z\dpr$.
\end{lemma}

\begin{proof}
Let $D'\subset V\otimes\BF_q^{\,\alg}\dpl z\dpr$ be a $b\cdot(\id\otimes\sigma)$-stable subspace. By Lemma~\ref{LemmaKe3.6.2} it suffices to show that $D'$ descents to $\BF_{q^s}\dpl z\dpr$ provided that $D'$ has dimension one. Choose a basis $e_1,\ldots,e_n$ of $V$ and a generator $v=a_1e_1+\ldots+a_ne_n$ of $D'$ with $a_i\in\BF_q^{\,\alg}\dpl z\dpr$. Without loss of generality $a_1\neq0$ and we may assume $a_1=1$. There exists a unit $\alpha\in\BF_q^{\,\alg}\dpl z\dpr\mal$ with $b\,\sigma^\ast v=\alpha\,v$. By the decency condition this implies
\[
\diag(z^{d_1},\ldots,z^{d_n})\cdot
\left(\begin{array}{c} a_1^{\sigma^s} \\  \vdots\\ a_n^{\sigma^s} \end{array}\right)
\es=\es
\alpha\,\alpha^\sigma\cdot\ldots\cdot\alpha^{\sigma^{s-1}}
\left(\begin{array}{c} a_1 \\  \vdots\\ a_n \end{array}\right)\,.
\]
Hence $z^{d_1}=\alpha\,\alpha^\sigma\cdot\ldots\cdot\alpha^{\sigma^{s-1}}$ and as desired all $a_i$ belong to $\BF_{q^s}\dpl z\dpr$ (and are zero if $d_i\neq d_1$).
\end{proof}

Now let $K$ be a complete extension of $k\dpl\zeta\dpr$ and let $b\in\GL_n\bigl(k\dpl z\dpr\bigr)$ and $\gamma\in\GL_n\bigl(K\dpl z-\zeta\dpr\bigr)$. With any $\BF_q\dpl z\dpr$-rational representation $\rho:\GL_n\to\GL(V)$ we associate the $z$-isocrystal with Hodge-Pink structure
\[
H(V)\es:=\es H_{b,\gamma}(V)\es:=\es\Bigl(V\otimes k\dpl z\dpr\,,\;\rho(b)\cdot(\id\otimes\sigma)\,,\;\rho(\gamma)\cdot\bigl(V\otimes K\dbl z-\zeta\dbr\bigr)\Bigr)\,.
\]
Consider the composition $\det_V\circ\rho$. It factors through $\det:\GL_n\to\BG_m$ in the form $\det_V\circ\rho=(\det)^{m_{\SSC V}}$ for an integer $m_{\SSC V}$. We obtain
\[
t_N\bigl(H_{b,\gamma}(V)\bigr)\es=\es m_{\SSC V}\ord_z(\det b)\;,\qquad t_H\bigl(H_{b,\gamma}(V)\bigr)\es=\es -m_{\SSC V}\ord_{z-\zeta}(\det\gamma)\,.
\]

The following equivalence of assertions parallels the case of mixed characteristic. Note however that its proof is more subtle since the category $\ul{\rm Rep}_{\BF_q\dpl z\dpr}\GL_n$ is not semi-simple as opposed to $\ul{\rm Rep}_{\BQ_p}\!\!\GL_n$.

\begin{definition}\label{DefAdmPair}
We call the pair $(b,\gamma)$ \emph{admissible} if $\ord_z(\det b)=-\ord_{z-\zeta}(\det\gamma)$, and one of the following equivalent conditions holds:
\begin{enumerate}
\item 
for any $\BF_q\dpl z\dpr$-rational representation $V$ of $\GL_n$ the $z$-isocrystal with Hodge-Pink structure $H_{b,\gamma}(V)$ is admissible,
\item 
there exists a faithful $\BF_q\dpl z\dpr$-rational representation $V$ of $\GL_n$ for which the $z$-isocrystal with Hodge-Pink structure $H_{b,\gamma}(V)$ is admissible.
\end{enumerate}
We make the same definition for \emph{weakly admissible}.
\end{definition}

\begin{proof}
Clearly (a) implies (b). For the converse fix a faithful representation $V$. Then any $\BF_q\dpl z\dpr$-rational representation $W$ of $\GL_n$ is a subquotient of $U:=\bigoplus_{i=1}^r V^{\otimes m_i}\otimes (V\dual)^{\otimes n_i}$ for suitable $r$, $m_i$ and $n_i$. Since $U$ is admissible by Theorem~\ref{Thm1.5b} (respectively weakly admissible by Theorem~\ref{Thm1.4}) if $V$ is, it suffices to show that (weak) admissibility is preserved under passage to sub-representations and quotient representations.

The condition $\ord_z(\det b)=-\ord_{z-\zeta}(\det\gamma)$ implies
\[
t_N\bigl(H_{b,\gamma}(W)\bigr)\es=\es t_H\bigl(H_{b,\gamma}(W)\bigr)
\]
for all representations $W$. Now let $W$ be a representation such that $H_{b,\gamma}(W)$ is weakly admissible. Then the equivalent conditions from Definition~\ref{Def1.4} show that for every sub-representation or quotient representation $W'$ of $W$ also $H_{b,\gamma}(W')$ is weakly admissible.

If $W$ is actually admissible let $(\CP,\CQ)=\pair\bigl(H_{b,\gamma}(W)\bigr)$. Then $\CQ$ is isoclinic of slope zero by Theorem~\ref{Thm2.6}. If $W'\subset W$ is a sub-representation set $(\CP',\CQ')=\pair\bigl(H_{b,\gamma}(W')\bigr)$. Over an algebraically closed complete extension $\BC$ of $K$ the $\sigma$-submodule $\CQ'$ of $\CQ$ is isomorphic to $\bigoplus_i\CF_{d_i,n_i}$ by Theorem~\ref{Thm1a}. Since $\CQ$ is isoclinic of slope zero, all $d_i\leq0$ by Proposition~\ref{Prop0.7}. Since
\[
\sum_i d_i\es=\es\deg\CQ'\es=\es t_H\bigl(H_{b,\gamma}(W')\bigr)- t_N\bigl(H_{b,\gamma}(W')\bigr)\es=\es0
\]
by Lemma~\ref{Lemma2.7}, all $d_i$ must be zero and $\CQ'$ is isoclinic of slope zero. By Theorem~\ref{Thm2.6}, $H_{b,\gamma}(W')$ is admissible.
Dually if $W''$ is a quotient representation of $W$ and $(\CP'',\CQ'')=\pair\bigl(H_{b,\gamma}(W'')\bigr)$ with $\CQ''\otimes\BC\ancon\cong\bigoplus_i\CF_{d_i,n_i}$\,, the fact that $\CQ$ is isoclinic implies $d_i\geq0$ and $\deg\CQ''=0$ implies $d_i=0$. So again $\CQ''$ is isoclinic of slope zero and $H_{b,\gamma}(W'')$ is admissible.
\end{proof}

We consider the (faithful) standard representation $V=\BF_q\dpl z\dpr^{\oplus n}$ in $\ul{\rm Rep}_{\BF_q\dpl z\dpr}\GL_n$. We also fix Hodge-Pink weights $\ul w=(w_1\geq\ldots\geq w_n)$.

\begin{definition}\label{DefPeriodSpace}
The \emph{space $\CH=\CH_\ul{w}$ of Hodge-Pink structures of weights $\ul w$} is defined as follows. For any $\BFZ$-algebra $B$ we set $\Fp_B=V\otimes_{\BF_q\dpl z\dpr}B\dbl z-\zeta\dbr$ and we let $\CH(B)$ be the set
\begin{eqnarray*}
& \Bigl\{ \text{$B\dbl z-\zeta\dbr$-lattices $\Fq$ in $\Fp_B[\frac{1}{z-\zeta}]$ which are direct summands as $B$-modules}\\ 
& \text{such that for every point $x\in \Spec B$ and every integer $e\geq -w_n$}\\
& (z-\zeta)^{-e}\Fp_x/\Fq_x\es\cong\es\bigoplus_{i=1}^n \kappa(x)\dbl z-\zeta\dbr/(z-\zeta)^{e+w_i} \quad\Bigr\}\,.
\end{eqnarray*}
\end{definition}

\begin{remark}
As was noted in Remark~\ref{Remark7.3} the lattice $\Fq$ contains more information than just the Hodge-Pink filtration. We have to allow for this in the definition of $\CH$. It is not possible here to define $\CH$ solely through the Hodge cocharacter as is done in mixed characteristic \cite[1.31]{RZ}. Correspondingly $\CH$ will not be a partial flag variety but a partial jet bundle over a partial flag variety (see below).
 The jets precisely take account of the fact mentioned in \ref{Remark7.3} that more general lattices can occur in equal characteristic as opposed to $p$-adic Hodge theory.
\end{remark}

Now let 
\[
\Fq_0\es:=\es\bigoplus_{i=1}^n(z-\zeta)^{w_i}\BFZ\dbl z-\zeta\dbr\es\subset\es(\Fp_\BFZ)[\TS\frac{1}{z-\zeta}]\,.
\]
Then $\Fq_0$ has Hodge-Pink weights $\ul w$. Fix an integer $e$ with $e\geq w_1$ and $w_n\geq -e$. Then 
\[
(z-\zeta)^e\Fp\es\subset\es\Fq_0\es\subset\es(z-\zeta)^{-e}\Fp\,.
\]
If $K\supset\BFZ$ is a field and $\Fq\in\CH(K)$ then by the theorem on elementary divisors, $\Fq$ is conjugate to $\Fq_0$ under $\GL_n\bigl(K\dbl z-\zeta\dbr\bigr)$. Therefore $\CH_\ul{w}$ is represented by the homogeneous space
\[
\GL_n\bigl((z-\zeta)^{-e}\Fp\bigr)/\Stab(\Fq_0)\es=\es 
\GL_n\Bigl((z-\zeta)^{-e}\Fp/(z-\zeta)^e\Fp\Bigr)/\Stab\Bigl(\Fq_0/(z-\zeta)^e\Fp\Bigr)\,.
\]
Let $\wt G$ be the Weil restriction
\[
\wt G\es:=\es Restr_{\bigl(\BFZ\dbl z-\zeta\dbr/(z-\zeta)^{2e}\bigr)/\BFZ} \GL_n\,.
\]
It is a smooth connected algebraic group over $\BFZ$ and isomorphic to $\GL_n\bigl((z-\zeta)^{-e}\Fp/(z-\zeta)^e\Fp\bigr)$. 
The 
stabilizer $\Stab\bigl(\Fq_0/(z-\zeta)^e\Fp\bigr)$ is the 
smooth closed algebraic subgroup $S\subset\wt G$ defined over $\BFZ$, whose $B$-valued points are
\begin{eqnarray*}
S(B)&=&\Bigl\{\,(a_{ij})_{i,j}\;\in\;\wt G(B)\,=\,\GL_n\bigl(B\dbl z-\zeta\dbr/(z-\zeta)^{2e}\bigr):\\
& & \qquad a_{ij}\;\in\; (z-\zeta)^{w_i-w_j}B\dbl z-\zeta\dbr/(z-\zeta)^{2e}\quad\text{for }i<j\,\Bigr\}\,.
\end{eqnarray*}
Thus the homogeneous space $\CH\cong\wt G/S$ is a smooth algebraic variety over $\BFZ$ and the quotient morphism $\wt G\to\CH$ is smooth.
$\CH$ is quasi-projective, since giving $\Fq_x$ is equivalent to giving the subspace
\[
\Fq_x/(z-\zeta)^e\Fp_x\es\subset\es(z-\zeta)^{-e}\Fp_x/(z-\zeta)^e\Fp_x
\]
and this defines an embedding of $\CH$ into some Grassmanian. The group $S$ is contained in the parahoric subgroup $\wt S\subset\wt G$ with
\begin{eqnarray*}
\wt S(B)&=&\Bigl\{\,(a_{ij})_{i,j}\;\in\;\wt G(B)\,=\,\GL_n\bigl(B\dbl z-\zeta\dbr/(z-\zeta)^{2e}\bigr):\\
& & \qquad a_{ij}\;\in\; (z-\zeta)B\dbl z-\zeta\dbr/(z-\zeta)^{2e}\quad\text{whenever }w_j<w_i\,\Bigr\}
\end{eqnarray*}
and the quotient $\wt G/\wt S$ is a partial flag variety over $\BFZ$ onto which $\CH$ naturally projects. The fibers of this projection are isomorphic to $\wt S/S$, that is they are affine spaces.
Altogether we have proved the following

\begin{proposition}\label{PropPeriodSpSmooth}
The space $\CH_\ul{w}$ is a smooth quasi-projective variety over $\BFZ$. More precisely, it is a partial jet bundle over a partial flag variety. 
\qed
\end{proposition}

\begin{example}
Let $n=3$ and $w_0=1>w_1=0>w_2=-1$. Then $e=1$ and $\CH$ equals the quotient of $\wt G= Restr_{\left(\BFZ\dbl z-\zeta\dbr/(z-\zeta)^{2}\right)/\BFZ} \GL_3$ by
\[
S\es=\es\Bigl\{ \left(\begin{array}{ccc} 
\ast & a\cdot(z-\zeta) & 0 \\ 
\ast & \ast & b\cdot(z-\zeta) \\
\ast & \ast & \ast
\end{array} \right)\Bigr\}\,.
\]
Here 
\[
\wt S\es=\es\Bigl\{ \left(\begin{array}{ccc} 
\ast & a\cdot(z-\zeta) & c\cdot(z-\zeta) \\ 
\ast & \ast & b\cdot(z-\zeta) \\
\ast & \ast & \ast
\end{array} \right)\Bigr\}
\]
and so $\CH$ is a line bundle over the full flag variety.
\end{example}

Next we want to construct analytic subspaces of $\CH_\ul{w}$ defined over $\BF_q^{\,\alg}\dpl\zeta\dpr$. We fix an element $b\in\GL_n\bigl(\BF_q^{\,\alg}\dpl z\dpr\bigr)$. Let $x\in\CH_\ul{w}(K)$ be a point for a complete extension $K$ of $\BF_q^{\,\alg}\dpl\zeta\dpr$ and let $\gamma_x$ be an element of $\GL_n\bigl(K\dpl z-\zeta\dpr\bigr)$ with $\Fq_x=\gamma_x\cdot\Fp_x$. The point $x$ is called (\emph{weakly}) \emph{admissible} with respect to $b$ if the pair $(b,\gamma_x)$ is (weakly) admissible, that is if the $z$-isocrystal with Hodge-Pink structure 
\[
H_{b,\gamma_x}(V)\es=\es\Bigl(V\otimes\BF_q^{\,\alg}\dpl z\dpr,\,b\cdot(\id\otimes\sigma),\,\Fq_x\Bigr)
\]
is (weakly) admissible. This definition does not depend on the choice of $\gamma_x$.

We will show in Theorem~\ref{ThmParacompact} below that the subset of all weakly admissible points is a rigid analytic subspace of $\CH_{\ul w}\otimes\BF_q^{\,\alg}\dpl\zeta\dpr$. However, in order to define the subspace of all admissible points it is best to consider the variety $\CH_{\ul w}\otimes\BF_q^{\,\alg}\dpl\zeta\dpr$ as a Berkovich space (see Appendix~\refAppBerkovich) instead of rigid analytic space. The reason is that one should not only consider its classical rigid analytic points (whose residue field is finite over $\BF_q^{\,\alg}\dpl\zeta\dpr$) but also their $K$-valued points for arbitrary complete extensions $K$ of $\BF_q^{\,\alg}\dpl\zeta\dpr$. Namely by Theorem~\ref{Thm3.5} the $z$-isocrystals with Hodge-Pink structure at the former points are already admissible as soon as they are weakly admissible, whereas this is not true for the latter points. The difference between the subspaces of admissible respectively weakly admissible points can thus only be revealed if one considers also the latter points. The Berkovich space associated with $\CH_{\ul w}\otimes\BF_q^{\,\alg}\dpl\zeta\dpr$ naturally contains those. 

So we let $\CH_\ul{w}^\an$ be the Berkovich space over $\BF_q^{\,\alg}\dpl\zeta\dpr$ associated with the variety $\CH_{\ul w}\otimes \BF_q^{\,\alg}\dpl\zeta\dpr$ (see Appendix~\refAppBerkovich) and we denote by $\CH^{wa}_b=\CH^{wa}_{\ul{w},b}$ (respectively $\CH^{a}_b=\CH^{a}_{\ul{w},b}$) the subset of weakly admissible (respectively admissible) points of $\CH_\ul{w}^\an$. Of course these sets may be empty if the numerical invariants $\ul w$ are wrong. In the next section we shall see that  they are open analytic subspaces of $\CH_\ul{w}^\an$. 

\begin{definition}
The space $\CH_{\ul w,b}^{wa}$ is called the \emph{period space} of Hodge-Pink structures of weights $\ul w$ on the (constant) $z$-isocrystal $\bigl(V\otimes k\dpl z\dpr,\, b\cdot(\id\otimes\sigma)\bigr)$.
\end{definition}

\begin{proposition}\label{Prop2.11b}
Tensoring with the Tate object $\BOne(e)$, that is replacing $b$ by $z^eb$ and $\Fq$ by $(z-\zeta)^{-e}\cdot\Fq$, defines an isomorphism $\CH_\ul{w}\to\CH_{\ul w'}$, where $\ul w'=(w_1-e\geq\ldots\geq w_n-e)$. This isomorphism maps $\CH_{\ul w,b}^a$ onto $\CH_{\ul w',z^eb}^a$ and $\CH_{\ul w,b}^{wa}$ onto $\CH_{\ul w',z^eb}^{wa}$.
\end{proposition}

\begin{proof}
This follows from the admissibility of the Tate object $\BOne(e)$ and Theorems~\ref{Thm1.4} and \ref{Thm1.5b}.
\end{proof}

The space $\CH_\ul{w}$ is equipped with an action of the automorphism group $J_b\bigl(\BF_q\dpl z\dpr\bigr)$ of the $z$-isocrystal $(D,F_D)=\bigl(V\otimes k\dpl z\dpr,\, b\cdot(\id\otimes\sigma)\bigr)$. This action stabilizes the subsets $\CH_{\ul w,b}^a$ and $\CH_{\ul w,b}^{wa}$. The group $J_b$ can be described analogously to the situation in mixed characteristic as follows.

\begin{proposition}[{Compare \cite[Proposition 1.12]{RZ}}]\label{PropLeviGp}
Let $b\in\GL_n\bigl(\BF_q^{\,\alg}\dpl z\dpr\bigr)$ then the functor on the category of $\BF_q\dpl z\dpr$-algebras $A$
\[
J_b(A)\es:=\es\bigl\{\,g\in\GL_n\bigl(A\otimes_{\BF_q\dpl z\dpr}\BF_q^{\,\alg}\dpl z\dpr\bigr):\es g\,b\;=\;b\,g^\sigma\,\bigr\}
\]
is representable by a smooth affine group scheme over $\BF_q\dpl z\dpr$.
\qed
\end{proposition}

In particular the automorphism group $\Aut(D,F_D)$ equals $J_b\bigl(\BF_q\dpl z\dpr\bigr)$.

\begin{proposition}[{Compare \cite[Corollary 1.14]{RZ}}]\label{PropLeviGpDecent}
Assume that $b\in\GL_n\bigl(\BF_{q^s}\dpl z\dpr\bigr)$ satisfies a decency equation of the form
\[
(b\cdot\sigma)^s\es=\es \diag(z^{d_1},\ldots,z^{d_n})\cdot\sigma^s
\]
and let $s\nu:\BG_m\to\GL_n$ be the cocharacter over $\BF_q\dpl z\dpr$ with $s\nu(z)=\diag(z^{d_1},\ldots,z^{d_n})$. Then $J_b$ is an inner form of the centralizer of $s\nu$. It is a Levi subgroup of $\GL_n\otimes\BF_{q^s}\dpl z\dpr$. In particular $\Aut(D,F_D)\subset \GL_n\bigl(\BF_{q^s}\dpl z\dpr\bigr)$.
\qed
\end{proposition}

%%%%%%%%%%%%%%%%%%%%%%%%%%%%%%%%%%%%%%%%%%%%%%%%%%%%%%%%%%%%%%%%%%%%%%
%
%    Openness of the Admissible Locus
%
%%%%%%%%%%%%%%%%%%%%%%%%%%%%%%%%%%%%%%%%%%%%%%%%%%%%%%%%%%%%%%%%%%%%%%

\subsection{Openness of the Admissible Locus} \label{SectAdmLocusOpen}

The key to showing that the admissible locus is Berkovich open is the following lemma. Recall the definition of the analytic spectrum $\CM(B)$ of an affinoid $\BFZ$-algebra $B$ from Definition~\ref{DefAnalyticPoint}.

\begin{lemma} \label{Lemma2.10}
Let $U=\CM\bigl(\BFZ\langle\frac{u}{\zeta}\rangle\bigr)$ be the closed disc with radius $|\zeta|$ over $\BFZ$. Let $j,d,e$ be integers with $0\le j< d$, $0<e$. For a point $x\in U$ consider the element
\[
\sum_{\nu\in\BZ} z^{-d\nu+j}\,u^{q^\nu}\mod(z-\zeta)^e\es=\es \sum_{\rho=0}^{e-1}y_\rho (z-\zeta)^\rho
\]
in $\wh{\kappa(x)^\alg}\dbl z-\zeta\dbr/(z-\zeta)^e$. Let $\BA^e_\BFZ=\CM\bigl(\BFZ[Y_0,\ldots,Y_{e-1}]\bigr)$ be the Berkovich space associated with affine $e$-space over $\BFZ$ and let $\alpha(x)$ be the point of $\BA^e_\BFZ$ given by
\[
\BFZ[Y_0,\ldots,Y_{e-1}]\es\longto\es\wh{\kappa(x)^\alg}\,,\quad Y_\rho\es\mapsto\es y_\rho\,.
\]
Then the map $\alpha:U\to\BA^e_\BFZ,\,x\mapsto\alpha(x)$ is a continuous map of topological Hausdorff spaces.
\end{lemma}

\begin{proof}
Note the following explicit formulas for the $y_\rho$:
\[
y_\rho\es=\es \raisebox{1ex}{``\,}{\frac{1}{\rho !}\,\frac{d^\rho}{dz^\rho}}\raisebox{1ex}{\,''\;}\Bigl(\sum_{\nu\in\BZ} z^{-d\nu+j}\,u^{q^\nu}\Bigr)\Bigr|_{z=\zeta}\es=\es\sum_{\nu\in\BZ} \binom{-d\nu+j}{\rho} \zeta^{-d\nu+j-\rho} u^{q^\nu}
\,.
\]
Since $|u|_x\leq|\zeta|$ for all $x\in U$ we see that the $y_\rho(x)$ are bounded. Hence the image of $\alpha$ lands in a polydisc $\BD(\theta)^e=\CM\bigl(\BFZ\langle\frac{y_0}{\theta},\ldots,\frac{y_{e-1}}{\theta}\rangle\bigr)$ with radii $(|\theta|,\ldots,|\theta|)$ for some $\theta\in\BFZ$.
By definition of the topology on $\BA^e_\BFZ$ (Definition~\ref{DefAnalyticPoint}), $\alpha$ is continuous if and only if for any $f\in\BFZ[Y_0,\ldots,Y_{e-1}]$ and any open interval $I\subset\BR_{\geq0}$ the preimage under $\alpha$ of the open set
\[
\{\,x\in\BA^e_\BFZ:\es|f|_x\in I\,\}
\]
is open in $U$. Write $f=\sum_{\ul\mu}a_{\ul\mu}\ul Y^{\ul\mu}$ where $a_{\ul\mu}\in\BFZ$ for every multi-index $\ul\mu\in\BN_0^d$. If we set $Y_\rho=Y'_\rho+Y''_\rho$ we obtain from the Taylor expansion of $f$
in powers of $\ul Y''$ a bound $c$ such that the condition $|Y''_\rho|_x\leq c$ for all $\rho$ implies
\[
|f(\ul Y)|_x\in I\quad\text{and}\quad |\ul Y|_x\leq |\theta|\qquad\Longleftrightarrow\qquad |f(\ul Y')|_x\in I\quad\text{and}\quad|\ul Y'|_x\leq |\theta|\,.
\]
Now we fix a negative integer $m$ and set 
\[
y'_\rho\es:=\es\sum_{\nu\geq m}\binom{-d\nu+j}{\rho}\zeta^{-d\nu+j-\rho}u^{q^\nu}\qquad\text{and}\qquad
y''_\rho\es:=\es\sum_{\nu<  m}\binom{-d\nu+j}{\rho}\zeta^{-d\nu+j-\rho}u^{q^\nu}\,.
\]
Since $|u^{q^\nu}|_x<1$ for any $\nu\in\BZ$ and any $x\in U$ we can find a small enough $m$ such that $|y''_\rho|_x\leq c$ for all $\rho$ and for all $x\in U$.
Thus we are reduced to showing that
\[
W\es=\es\{\,x\in U:\es\bigl|f\bigl(y'_0(u),\ldots,y'_{e-1}(u)\bigr)\bigr|_x\in I\,\}
\]
is open in $U$. For this purpose consider the morphism $\beta:\CM\bigl(\BFZ\langle\frac{v}{\zeta^{q^m}}\rangle\bigr)\to U$ given by \mbox{$\beta^\ast u= v^{q^{-m}}$} (note that $m<0$). It is a homeomorphism since it is a continuous map between compact spaces,
which is bijective because any multiplicative semi-norm $|\,.\,|$ on $\BFZ\langle\frac{u}{\zeta}\rangle$ uniquely lifts to a multiplicative semi-norm on $\BFZ\langle\frac{v}{\zeta^{q^m}}\rangle$ with $|\frac{v}{\zeta^{qm}}|:=|\frac{u}{\zeta}|^{q^m}$. 
We have $y'_\rho\in \BFZ\langle\frac{v}{\zeta^{q^m}}\rangle$ for all $\rho$, and since the morphism
\[
\TS \CM\bigl(\BFZ\langle\frac{v}{\zeta^{q^m}}\rangle\bigr)\to\CM\bigl(\BA^e_\BFZ\bigr)\,,\quad v\mapsto \bigl(y'_0(v),\ldots,y'_{e-1}(v)\bigr)
\]
is continuous, the set $W$ is open as desired.
\end{proof}

\begin{theorem} \label{Thm2.11}
Let $\bar b$ be a $\sigma$-conjugacy class in $\GL_n\bigl(\BF_q^{\,\alg}\dpl z\dpr\bigr)$ and let $b\in\GL_n\bigl(\BF_{q^s}\dpl z\dpr\bigr)$ be a decent element of $\bar b$ satisfying a decency equation with the integer $s>0$ (Definition~\ref{Def2.2}). Let $\CH_s^\an$ be the Berkovich space over $\BF_{q^s}\dpl\zeta\dpr$ associated with the variety $\CH_\ul{w}\otimes_{\BFZ}\BF_{q^s}\dpl\zeta\dpr$.
Then the set
\[
\CH^a_b\es:=\es\CH^a_{\ul w,b}\es:=\es\{\,x\in\CH_s^\an:\es(D,F_D,\Fq_x) \es\text{is admissible}\,\}
\]
is Berkovich open in $\CH_s^\an$. If $b'\in\bar b$ satisfies a decency equation with the same integer $s$ then any $g\in\GL_n\bigl(\BF_{q^s}\dpl z\dpr\bigr)$ with $b'=gb(g^{-1})^\sigma$ 
induces an isomorphism $\Fq_x\mapsto g(\Fq_x)=:\Fq_{g(x)}$ between $\CH^a_b$ and $\CH^a_{b'}$.
\end{theorem}

\begin{proof}
(a) \es The last statement is clear. We prove the first. By Proposition~\ref{Prop2.11b} we may tensor with a suitable Tate object, and assume that all Hodge-Pink weights are non-negative and that $b^{-1}\in M_n\bigl(k\dbl z\dbr\bigr)$. Fix an integer $e$ bigger than all Hodge-Pink weights. So $(z-\zeta)^e\Fp\subset\Fq\subset\Fp$. Consider the pair $(\CP,\CQ)$ of $\sigma$-modules over $\CO_{\CH_s^\an}\ancon$ defined by the universal Hodge-Pink structure (Definition~\ref{Def1.6}).
By Theorem~\ref{Thm2.6} a point $x\in\CH_s^\an$ belongs to $\CH^a_b$ if and only if $\CQ_x$ is isoclinic of slope zero, hence if and only if its HN-polygon is constant zero. Thus the openness of $\CH^a_b$ follows from Corollary~\ref{Cor1.7.9}.

However, we want to give a second proof for the openness which gives a kind of parameterization of its complement and whose ideas we will again use in Theorem~\ref{Thm2.11'}.
By Theorem~\ref{Thm1a} and Proposition~\ref{Prop0.7} a point $x\in\CH_s^\an$ belongs to the complement of $\CH^a_b$ if and only if for some algebraically closed complete extension $\BC$ of $\kappa(x)$ there exists a non-trivial homomorphism of $\sigma$-modules over $\BC\ancon$ from $\CF_{1,n}$ to $\CQ_x\otimes\BC\ancon$.
Let $s'$ be the least common multiple of $s$ and $n$ and replace $s$ by $s'$. From now on we view all $\sigma$-modules as $\sigma^s$-modules by replacing $F$ by $F^s$. As a $\sigma^s$-module $\CF_{1,n}$ is isomorphic to $\CO(\frac{s}{n})^{\oplus n}$.
So if $x$ is not admissible then for some algebraically closed complete extension $\BC$ of $\kappa(x)$ the $\sigma^s$-module $\CQ_x\otimes \BC\ancon$ contains a $\sigma^s$-module $\CO(1)$. Conversely the latter implies that $\CQ_x\otimes\BC\ancon\not\cong\CO(0)$ as $\sigma$-modules and $x$ is not admissible. Now there is an isomorphism
\[
\Hom_{\sigma^s}\bigl(\CO(1)\,,\,\CQ_x\otimes\BC\ancon\bigr)\es\isoto\es\CQ_x(-1)^{F^s}(\BC)\,,\quad\phi\es\mapsto\es\phi(1)\,,
\]
where we identify the $K\ancon$-modules underlying $\CQ_x$ and $\CQ_x(-1)$, that is, we view $\CQ_x(-1)$ as the $\sigma^s$-module \mbox{$(\CQ_x,F^s=z\cdot F_{\CQ_x}^s)$}. By construction of $\CQ_x$ in Definition~\ref{Def1.6}
\begin{eqnarray*}
\Hom_{\sigma^s}\bigl(\CO(1)\,,\,\CQ_x\otimes\BC\ancon\bigr)&\cong&\Bigl\{\,f\in\bigl(\CP(-1)\otimes\BC\ancon\bigr)^{F^s}(\BC): \\
& &  \quad f\in\eta_r(\sigma^{r\ast}\Fq_x)\otimes\BC\dbl z-\zeta^{q^r}\dbr\es\text{for all }r=0,\ldots,s-1\,\Bigr\}\,,
\end{eqnarray*}
where $\eta_r=F_\CP\circ\ldots\circ(\sigma^{r-1})^\ast F_\CP$. Let $b_r:=b\,b^\sigma\cdot\ldots\cdot b^{\sigma^{r-1}}$ be the matrix by which $\eta_r$ acts on the standard basis of $\CP$.
Recall the algebraic groups $\wt G$ and $S\subset\wt G$ from the previous section and the isomorphism $\CH\cong\wt G/S$. Let $\wt G_s^\an$ be the Berkovich space over $\BF_{q^s}\dpl\zeta\dpr$ associated with the group scheme $\wt G\otimes\BF_{q^s}\dpl\zeta\dpr$. Thus $x\in\CH_s^\an$ is not admissible if and only if there is a point $g\in\wt G_s^\an$ mapping to $x$, an algebraically closed complete extension $\BC$ of $\kappa(x)$, and an $f\in\bigl(\CP(-1)\otimes\BC\ancon\bigr)^{F^s}(\BC), \;f\neq0$ with
\begin{equation} \label{EqCondOnF}
(g^{-1})^{\sigma^r}\,b_r^{-1}\,f\es\in\es(\sigma^r)^\ast\Fq_0\otimes\BC\dbl z-\zeta^{q^r}\dbr\qquad \text{for all }r=0,\ldots,s-1\,.
\end{equation}
(Note that $\Fq_x=g\cdot\Fq_0$.) A priori $(g^{-1})^{\sigma^r}b_r^{-1}f\in(\sigma^r)^\ast\Fp\otimes\BC\dbl z-\zeta^{q^r}\dbr$, so the latter condition is Zariski-closed.

\medskip

(b) \es For $\theta\in\BF_{q^s}\dpl\zeta\dpr$ consider the $sne$-dimensional polydisc $\BD(\theta)^{sne}$ over $\BF_{q^s}\dpl\zeta\dpr$ with radii $(|\theta|,\ldots,|\theta|)$. We represent the $K$-valued points of $\BD(\theta)^{sne}$ in the form $(h_r)_{r=0}^{s-1}$ with $h_r=\sum_{\rho=0}^{e-1}y_{\rho r}(z-\zeta^{q^r})^\rho$ where $y_{\rho r}\in\BD(\theta)^n(K)\;=\;\{\,y\in K^n:|y|\leq|\theta|\,\}$ for all $j,r$.
We shall exhibit in (c) below a constant $\theta\in\BF_{q^s}\dpl\zeta\dpr$ and a compact subset $C$ of the polydisc $\BD(\theta)^{sne}$ such that if $\BC$ is an algebraically closed complete extension of $\BFZ$ and $f\in\bigl(\CP(-1)\otimes\BC\ancon\bigr)^{F^s}(\BC),\; f\neq0$ then for some integer $N$,
\[
\Bigl(z^N f\mod(z-\zeta^{q^r})^{e}\Bigr)_{r=0}^{s-1}
\]
is a $\BC$-valued point of $C$, and $C$ consists precisely of those points. Note that $z\in\BC\dbl z-\zeta\dbr\mal$ is a unit.
Now we view $\Fp/(z-\zeta)^e\Fp$ as affine $ne$-space $\BA^{ne}$ over $\BF_{q^s}\dpl\zeta\dpr$ by interpreting, for any $\BF_{q^s}\dpl\zeta\dpr$-algebra $B$, an element $\sum_{i=0}^{e-1}f_i(z-\zeta)^i\;\in\;\Fp/(z-\zeta)^e\Fp\otimes B$ with $f_i\in B^n$ as the point $(f_i)_{i=0}^{e-1}$ of $\BA^{ne}$. We consider the morphism
\begin{eqnarray*}
\beta:\es\wt G_s^\an\times_{\BF_{q^s}\dpl\zeta\dpr}\BD(\theta)^{sne} & \longto & \prod_{r=0}^{s-1}(\sigma^r)^\ast\bigl(\Fp\,/\,(z-\zeta)^e\Fp\bigr)\qquad=\es\BA^{sne} \\
\Bigl(\,g\,,\,(h_r)_{r=0}^{s-1}\,\Bigr) & \mapsto & \Bigl(\,(g^{-1})^{\sigma^r}\!\cdot b_r^{-1}\cdot h_r \mod (z-\zeta^{q^r})^e\,\Bigr)_{r=0}^{s-1}
\end{eqnarray*}
Note that $g^{-1}\mapsto(g^{-1})^{\sigma^r}$ is the $q^r$-Frobenius morphism $\wt G_s^\an\mapsto(\sigma^r)^\ast\wt G_s^\an$. Let $W$ be the closed subset of $\wt G_s^\an\times_{\BF_{q^s}\dpl\zeta\dpr}\BD(\theta)^{sne}$ defined by (the pullback of) condition (\ref{EqCondOnF}) and the condition that $(h_r)_{r=0}^{s-1}$ belongs to $C$.
Furthermore, we consider the projection map 
\[
pr_1:\es\wt G_s^\an\times\BD(\theta)^{sne}\es\to\es\wt G_s^\an
\]
onto the first factor and the canonical map $\gamma:\wt G_s^\an\to\CH_s^\an$ coming from the isomorphism $\CH\cong\wt G/S$. By what we have said above, $x\in\CH_s^\an$ is not admissible if and only if $x\in\gamma\circ pr_1(W)$.

Since $\BD(\theta)$ is quasi-compact the projection $pr_1$ is a proper map of topological spaces by \cite[Proposition 3.3.2]{Berkovich1}, thus in particular it is closed. Hence $pr_1 (W)\subset\wt G_s^\an$ is closed.
Note that $\CH_s^\an$ carries the quotient topology under $\gamma$ since $\gamma$ is a smooth morphism of schemes, hence open by \cite[Proposition 3.5.8 and Corollary 3.7.4]{Berkovich2}. Since by construction $pr_1 (W)\;=\;\gamma^{-1}\bigl(\gamma\circ pr_1 (W)\bigr)$ 
we conclude that $\CH_s^\an\setminus\CH^a_b\;=\;\gamma\circ pr_1 (W)$ is closed in $\CH_s^\an$ as desired.

\medskip

(c) \es It remains to construct the compact set $C$. Due to the decency of $b$ we can write the $\sigma^s$-module $\CP(-1)$ as $\bigoplus_{i=1}^n\CO(d_i)$. Assume that in the beginning the Tate object was chosen such that $d_i\geq1$ for all $i$. Then for any algebraically closed complete extension $\BC$ of $\BFZ$
\begin{eqnarray*}
\bigl(\CP(-1)\otimes\BC\ancon\bigr)^{F^s}(\BC) &=& \bigoplus_{i=1}^n\CO(d_i)^{F^s}(\BC) \\
&=& \bigoplus_{i=1}^n\;\Bigl\{\,f\es=\es\sum_{\nu\in\BZ} z^{-d_i\nu}\,\sum_{j=0}^{d_i-1} z^j u_{i,j}^{q^{ s\nu}}:\es u_{i,j}\in\BC,\,|u_{i,j}|<1\,\Bigr\}
\end{eqnarray*}
by Proposition~\ref{Prop0.5}.
Consider the compact set
\[
\TS U\es:=\es\CM\Bigl(\BF_{q^s}\dpl\zeta\dpr\langle\frac{u_{i,j}}{\zeta}:\;i=1,\ldots,n,\,j=0,\ldots,d_i-1\rangle\Bigr)
\]
and for $1\leq k\leq n$ and $0\leq\ell\leq d_k-1$ the following compact subsets of $U$
\begin{eqnarray*}
\TS U_{k,\ell}&:=&\TS\CM\Bigl(\BF_{q^s}\dpl\zeta\dpr\langle\frac{u_{i,j}}{\zeta}:\;i=1,\ldots,n,\,j=0,\ldots,d_i-1;\es\frac{\zeta^{q^s}}{u_{k,\ell}}\rangle\Bigr)\\[2mm]
&= & \bigl\{\,|u_{i,j}|\leq|\zeta|\text{ for all }i,j\quad\text{and}\quad|\zeta|^{q^s}\leq|u_{k,\ell}|\,\bigr\}\,.
\end{eqnarray*}
For $x\in U$ we have $u_{i,j}\in\kappa(x)$ and we define $y_{\rho r}=(y_{\rho r,i})_{i=1}^n\in\wh{\kappa(x)^\alg}^{\oplus n}$ by
\[
\sum_{\rho=0}^{e-1}y_{\rho r}(z-\zeta^{q^r})^\rho\es=\es\bigoplus_{i=1}^n\;\sum_{\nu\in\BZ}z^{-d_i \nu}\,\sum_{j=0}^{d_i-1}z^j u_{i,j}^{q^{ s\nu}}\es\mod(z-\zeta^{q^r})^{e}\es=\es f\mod(z-\zeta^{q^r})^e
\]
in $\wh{\kappa(x)^\alg}\dbl z-\zeta^{q^r}\dbr\,/\,(z-\zeta^{q^r})^{e}$ for all $r=0,\ldots,s-1$. We let
\[
\BA^{sne}_{\BF_{q^s}\dpl\zeta\dpr}\es:=\es\CM\bigl(\BF_{q^s}\dpl\zeta\dpr[Y_{\rho r,i}:\;0\le\rho\le e-1\,,\,0\le r\le s-1\,,\,1\le i\le n]\bigr)
\]
be the Berkovich space associated with affine $sne$-space over $\BF_{q^s}\dpl\zeta\dpr$ and
we obtain a map of sets \mbox{$\alpha:U\to\BA^{sne}_{\BF_{q^s}\dpl\zeta\dpr}$} mapping the $\kappa(x)$-valued point $x$ given by $u_{i,j}\in\kappa(x)$ 
to the $\wh{\kappa(x)^\alg}$-valued point given by $Y_{\rho r,i}\mapsto y_{\rho r,i}\in\wh{\kappa(x)^\alg}$. By the same argument as in Lemma~\ref{Lemma2.10} one proves that $\alpha$ is continuous (compared to Lemma~\ref{Lemma2.10} here only the number of variables has blown up).
Thus $C:=\alpha\bigl(\bigcup_{k,\ell}U_{k,\ell}\bigr)\;\subset\;\BA^{sne}_{\BF_{q^s}}$ is compact and contained in a polydisc $\BD(\theta)^{sne}=\CM\bigl(\BFZ\langle\frac{Y_{\rho r,i}}{\theta}\rangle\bigr)$ for large enough $|\theta|$.
Note that multiplying 
\[
f\es=\es\sum_{\nu\in\BZ}z^{-d_i\nu}\,\sum_{j=0}^{d_i-1}z^j u_{i,j}^{q^{ s\nu}}
\]
with $z$ changes $u_{i,j}$ to $u_{i,j-1}$ for $1\leq j\leq d_i-1$ and $u_{i,0}$ to $u_{i,d_i-1}^{q^s}$. Thus we can adjust the $u_{ij}$ such that $|u_{ij}|\le|\zeta|$ for all $i,j$, and $|u_{k\ell}|\ge|\zeta|^{q^s}$ for some $k,\ell$. This shows that $C$ has the desired property formulated in (b). Putting everything together the theorem is proved.
\end{proof}

\begin{corollary} \label{Cor2.12}
There exists an affinoid covering of $\CH^a_b$ (in the sense of Definition~\ref{DefBerkovichSpaces}) by connected $\CM(B_i)$ and for each $i$ a finite \'etale $B_i$-algebra $C_i$ and a $C_i\langle\frac{z}{\zeta}\rangle$-lattice $N_i$ in $\CQ\otimes C_i\ancon$ on which $F:\sigma^\ast N_i\to N_i$ is an isomorphism.
Moreover, the $N_i\otimes C_i\con$ descend and glue to a canonical $\sigma$-module over $\CO_{\CH^a_b}\con$.
\end{corollary}

\begin{proof}
By Definition~\ref{DefIsoclinic} the $\sigma$-module $\CQ$ over $\CO_{\CH_b^a}\ancon$ is isoclinic of slope zero.
Choose an affinoid covering of $\CH^a_b$ by sets $\CM(B)$. Then the assertion follows from Theorem~\ref{Thm6.12}. Note that the $N_i\otimes C_i\con$ descend and are uniquely determined by Theorem~\ref{Thm6.12}. Thus they glue canonically.
\end{proof}

The strategy for the proof of Theorem~\ref{Thm2.11} also allows us to show that the weakly admissible locus is Berkovich open.

\begin{theorem}\label{Thm2.11'}
Let $\bar b$ be a $\sigma$-conjugacy class in $\GL_n\bigl(\BF_q^{\,\alg}\dpl z\dpr\bigr)$ and let $b\in\GL_n\bigl(\BF_{q^s}\dpl z\dpr\bigr)$ be a decent element of $\bar b$ satisfying a decency equation with the integer $s>0$ (Definition~\ref{Def2.2}). Let $\CH_s^\an$ be the Berkovich space over $\BF_{q^s}\dpl\zeta\dpr$ associated with the variety $\CH_\ul{w}\otimes_{\BFZ}\BF_{q^s}\dpl\zeta\dpr$.
Then the set
\[
\CH^{wa}_b\es:=\es\CH^{wa}_{\ul w,b}\es:=\es\{\,x\in\CH_s^\an:\es(D,F_D,\Fq_x) \es\text{is weakly admissible}\,\}
\]
is Berkovich open in $\CH_s^\an$. If $b'\in\bar b$ satisfies a decency equation with the same integer $s$ then any $g\in\GL_n\bigl(\BF_{q^s}\dpl z\dpr\bigr)$ with $b'=gb(g^{-1})^\sigma$ 
induces an isomorphism $\Fq_x\mapsto g(\Fq_x)=:\Fq_{g(x)}$ between $\CH^{wa}_b$ and $\CH^{wa}_{b'}$.
\end{theorem}

\begin{proof}
By Proposition~\ref{Prop2.11b} we may tensor with a suitable Tate object and assume that all Hodge-Pink weights are non-negative and thus $\Fq\subset\Fp$. 
If $-w_1-\ldots-w_n\neq t_N(D,F_D)$ then $\CH_b^{wa}$ is empty. So we now assume equality. Then at every point of $\CH_s^\an$ we have $\deg \CQ_x=0$. By Proposition~\ref{Prop2.8} and Lemma~\ref{LemmaDecentAndStable} a point $x\in\CH_s^\an$ is not weakly admissible if and only if there exists an $F_D$-stable $\BF_{q^s}\dpl z\dpr$-subspace $D'\subset D$ such that $\pair(D',F_{D'},\Fq'_x)=(\CP'_x,\CQ'_x)$ satisfies $\deg\CQ'_x>0$ where $\Fq'_x=\Fq_x\cap\sigma^\ast D'\otimes_{\BF_{q^s}\dpl z\dpr}\kappa(x)\dpl z-\zeta\dpr$. Let $n':=\dim_{\BF_{q^s}\dpl z\dpr}D'$ and consider the family of Hodge-Pink structures $\wedge^{n'}\Fq$ on the exterior power $\wedge^{n'}D$ over $\CH_s^\an$. Then the subspace $\wedge^{n'}D'$ is one-dimensional. Let $\lambda'\in\BZ$ be its slope, that is, $\wedge^{n'}D'\otimes\BF_q^{\,\alg}\dpl z\dpr\cong\CO(-\lambda')$. Consider the slope decomposition of the $z$-isocrystal $\wedge^{n'}D=\bigoplus_\lambda D_\lambda$ from Proposition~\ref{Prop10} and the induced decomposition of $\wedge^{n'}\CP=\bigoplus_\lambda \CP_\lambda$ with $\CP_\lambda:=\sigma^\ast D_\lambda\otimes_{\BF_{q^s}\dpl z\dpr}\CO_{\CH_s^\an}\ancon$. We have $\wedge^{n'}D'\subset D_{\lambda'}$ because after the faithful extension to $\BF_q^{\,\alg}\dpl z\dpr$ all morphisms $\wedge^{n'}D'\hookrightarrow\wedge^{n'}D\to D_\lambda$ are zero for $\lambda\neq\lambda'$ by Proposition~\ref{Prop0.7b}.
At the non-weakly admissible point $x$ there exists an algebraically closed complete extension $\BC$ of $\kappa(x)$ and an isomorphism $\wedge^{n'}\CQ'\otimes\BC\ancon\cong\CO(c)$ of $\sigma$-modules over $\BC\ancon$ by Corollary~\ref{Cor0.6} where $c:=\deg\CQ'>0$.
The saturation of $\wedge^{n'}\CQ'$ inside $\wedge^{n'}\CP$ is $\wedge^{n'}\CP'=\sigma^\ast\!\wedge^{n'}\!D'\otimes_{\BF_{q^s}\dpl z\dpr}\kappa(x)\ancon$. Since $\wedge^{n'}\CP$ and $\wedge^{n'}\CQ$ differ only at $z=\zeta^{q^\nu}$ for $\nu\in\BN_0$ the same is true for $\wedge^{n'}\CP'$ and $\wedge^{n'}\CQ'$. Considering their degrees we see that $\wedge^{n'}\Fq'_x=(z-\zeta)^{-\lambda'-c}\cdot\wedge^{n'}\Fp'_x$. By the decency assumption $(F_{D_{\lambda'}})^s=z^{s\lambda'}\cdot\sigma^s$, hence $\sigma^\ast\!\wedge^{n'}\!D'$ can be recovered as
\[
\sigma^\ast\!\wedge^{n'}\!D'\es=\es \bigl(\wedge^{n'}\CP'\otimes\CO(\lambda')\bigr)^{F^s}(\BC)\es\subset\es\bigl(\CP_{\lambda'}\otimes\CO(\lambda')\bigr)^{F^s}(\BC)\es=\es \sigma^\ast\!\wedge^{n'}\!D_{\lambda'}\,.
\]
Let $f\in\bigl(\CP_{\lambda'}\otimes\CO(\lambda')\bigr)^{F^s}(\BC)$ be a generator of $\sigma^\ast\!\wedge^{n'}\!D'$ over $\BF_{q^s}\dpl z\dpr$, respectively of $\wedge^{n'}\CP'$ over $\kappa(x)\ancon$. In particular $(z-\zeta)^{-\lambda'-1}f\in\wedge^{n'}\Fq'_x$. To summarize, if a point $x\in\CH_s^\an$ is not weakly admissible then the following holds:
\begin{equation}\label{EqCondNotWA}
\text{\parbox{0.9\textwidth}{there exists an integer $n'$ with $1\leq n'\leq n$, an isoclinic summand $\CP_{\lambda'}$ of $\wedge^{n'}\CP$ of slope $\lambda'\in\BZ$, an algebraically closed complete extension $\BC$ of $\kappa(x)$, and a non-zero element $f\in\bigl(\CP_{\lambda'}\otimes\CO(\lambda')\bigr)^{F^s}(\BC)$, such that $(z-\zeta)^{-\lambda'-1}f\in\wedge^{n'}\Fq_x$, and $f$ and $F_{D_{\lambda'}}(\sigma^\ast f)$ are linearly dependent over $\BC\dpl z-\zeta\dpr$. }}
\end{equation}
Note that the elements $f$ and $F_{D_{\lambda'}}(\sigma^\ast f)$ of the $\BF_{q^s}\dpl z\dpr$-vector space $\sigma^\ast\!\wedge^{n'}\!D_{\lambda'}$ are linearly dependent if and only if they are over $\BC\dpl z-\zeta\dpr$.
So conversely (\ref{EqCondNotWA}) implies that $f$ defines a subobject of rank one of $(\wedge^{n'}D,\wedge^{n'}F_D,\wedge^{n'}\Fq_x)$ which contradicts weak admissibility and hence $x\notin\CH_b^{wa}$ by Theorem~\ref{Thm1.4}. Thus it suffices to fix an integer $n'$ with $1\leq n'\leq n$ and a slope $\lambda'\in\BZ$, and to show that the set $V_{n',\lambda'}$ of the points $x\in\CH_s^\an$ with property (\ref{EqCondNotWA}) is closed.

Before Lemma~\ref{Lemma2.7} we have constructed a $\sigma$-module $\tminus^{-1}\CO(0)$ over $\BFZ\ancon$ where $\tminus:=\prod_{i\in\BN_0}\bigl(1-{\TS\frac{\zeta^{q^i}}{z}}\bigr)\in \BFZ\ancon$, and an isomorphism of $\sigma$-modules $\tminus^{-1}\CO(0)\isoto\CO(1)$, $1\mapsto \tplus$ over $\ol{\BFZ}\ancon$ for a suitable unit $\tplus\in\ol{\BFZ}\ancon\mal$. Consider the $\sigma$-module $\wt\CP:=\CP_{\lambda'}\otimes\CO(\lambda')\otimes \tminus^{-1}\CO(0)$ over $\CO_{\CH_s^\an}\ancon$ and the inclusion $\CP_{\lambda'}\otimes\CO(\lambda')\subset\wt\CP$ coming from the inclusion $\CO(0)\subset \tminus^{-1}\CO(0)$. It identifies $\wedge^{n'}\Fp$ with $(z-\zeta)\wt\CP\otimes\CO_{\CH_s^\an}\dbl z-\zeta\dbr$. By the decency assumption, $\CP_{\lambda'}\otimes\CO(\lambda')$ is as $\sigma^s$-module equal to $\CO(0)^{\oplus m}$ where $m=\rk\CP_{\lambda'}$. Hence multiplication with $\tplus$ is an isomorphism of $\sigma^s$-modules $\wt\CP\isoto\CO(s)^{\oplus m}$ over $\ol{\BFZ}\ancon$. Then $\bigl(\CP_{\lambda'}\otimes\CO(\lambda')\bigr)^{F^s}(\BC)$ equals the set
\begin{eqnarray*}
& \Bigl\{\,f\in\wt\CP^{F^s}(\BC):\es f\in(z-\zeta^{q^r})\wt\CP\otimes\BC\dbl z-\zeta^{q^r}\dbr\text{ for all }0\leq r\leq s-1\,\Bigr\} \\[2mm]
& \cong \es \Bigl\{\,\tplus f\in\bigl(\CO(s)^{\oplus m}\bigr)^{F^s}(\BC):\es f\in(\sigma^r)^\ast\wedge^{n'}\Fp\otimes\BC\dbl z-\zeta^{q^r}\dbr\text{ for all }0\leq r\leq s-1\,\Bigr\}\,.
\end{eqnarray*}
Recall the algebraic groups $\wt G$ and $S\subset\wt G$ defined in the previous section and the isomorphism $\CH\cong\wt G/S$. Let $\wt G_s^\an$ be the Berkovich space over $\BF_{q^s}\dpl\zeta\dpr$ associated with the group scheme $\wt G\otimes\BF_{q^s}\dpl z\dpr$. Then $x\in\CH_s^\an$ belongs to $V_{n',\lambda'}$ if and only if there exists an algebraically closed complete extension $\BC$ of $\kappa(x)$, a $\BC$-valued point $g\in\wt G_s^\an$ mapping to $x$, and a non-zero element $\tilde f=\tplus f\in\bigl(\CO(s)^{\oplus m}\bigr)^{F^s}(\BC)$ with
\begin{eqnarray} \label{EqCondition2.4}
f \text{ and } F_{D_{\lambda'}}(\sigma^\ast f) \hspace{-1em} & & \hspace{-1em}\text{linearly dependent over }\BC\dpl z-\zeta\dpr\,, \nonumber \\[2mm]
f & \in & (\sigma^r)^\ast\wedge^{n'}\Fp\otimes\BC\dbl z-\zeta^{q^r}\dbr \qquad\text{and}\\[2mm]
(z-\zeta^{q^r})^{-\lambda'-1}\cdot(g^{-1})^{\sigma^r}F_{D_{\lambda'}}^{-r}(f) & \in & (\sigma^r)^\ast\wedge^{n'}\Fq_0\otimes\BC\dbl z-\zeta^{q^r}\dbr\nonumber
\end{eqnarray}
for all $0\leq r\leq s-1$. Note that if $f$ and $F_{D_{\lambda'}}(\sigma^\ast f)$ are linearly dependent then it suffices to verify the conditions of lines 2 and 3 of (\ref{EqCondition2.4}) only for $r=0$. A priori $f\,,\,g^{-1}f\;\in\;(z-\zeta)^{-1}\wedge^{n'}\!\Fp\otimes\BC\dbl z-\zeta\dbr$, so condition (\ref{EqCondition2.4}) is closed. This property is similar to the property formulated at the end of part (a) of the proof of Theorem~\ref{Thm2.11} and from here on we proceed completely analogous to parts (b) and (c) of that proof. We leave the details to the reader.
\end{proof}

We have seen in the preceding theorems that the sets $\CH_b^a$ and $\CH_b^{wa}$ are open subspaces of the Berkovich space $\CH_s^\an$. Now there is a functor
\[
X\es\longmapsto\es X^\rig\es:=\es\{\,x\in X:\es\kappa(x)\text{ is a finite extension of }L\,\}
\]
from Berkovich spaces over $L$ to rigid analytic spaces over $L$. It restricts to an equivalence of categories on the category of paracompact Berkovich spaces (see Appendix~\refAppBerkovich{} for a discussion of the terminology and this equivalence). Therefore the following result is of interest. It will allow us to view $\CH^{wa}_b$ and $\CH^a_b$ as Berkovich spaces or rigid analytic spaces whichever is more convenient. 

\begin{theorem}\label{ThmParacompact}
In the situation of Theorems~\ref{Thm2.11} and \ref{Thm2.11'} the Berkovich spaces $\CH_b^{a}$ and $\CH_b^{wa}$ are arcwise connected and paracompact. Moreover the open immersion $\CH_b^{wa}\hookrightarrow\CH_s^\an$ identifies $(\CH_b^{wa})^\rig$ with an admissible open subset of $(\CH_s^\an)^\rig$.
\end{theorem}

\begin{proof}
(a)\es By construction $\CH_s^\an$ admits a countable covering  by (strictly) $\BF_{q^s}\dpl\zeta\dpr$-affinoid subsets (see Appendix~\refAppBerkovich). Then by Lemma~\ref{LemmaParacompact} the open subsets $\CH_b^{a}$ and $\CH_b^{wa}$ are paracompact.

We show that the subset $(\CH_b^{wa})^\rig\subset(\CH_s^\an)^\rig$ is admissible open. By Lemma~\ref{LemmaParacompact} we can choose for every analytic point $x\in\CH_b^{wa}$ an affinoid neighborhood $U_x\subset\CH_b^{wa}$. Then $(\CH_b^{wa})^\rig=\bigcup_{x\in\CH_b^{wa}}(U_x)^\rig$. By \cite[Proposition 9.1.4/2]{BGR} it suffices to show that for any morphism $\phi:Y\to\CH_s^\an$ from an affinoid $Y$ with $\phi(Y^\rig)\subset(\CH_b^{wa})^\rig$ the covering $\{\phi^{-1}U_x^\rig\}_{x\in\CH_b^{wa}}$ of $Y^\rig$ is admissible. For that it suffices by \cite[Lemma 1.6.2]{Berkovich2} to show that for every analytic point $y\in Y$ the point $\phi(y)$ belongs to $\CH_b^{wa}$. Namely then $\phi^{-1}U_{\phi(y)}$ is a neighborhood of $y$ and the hypotheses of loc.\ cit.\ are fulfilled.

\medskip
\noindent
(b)\es So now assume that there exists a morphism $\phi:Y\to\CH_s^\an$ with $\phi(Y^\rig)\subset(\CH_b^{wa})^\rig$ and an analytic point $y\in Y$ such that $x=\phi(y)$ is not weakly admissible. Then we have shown in condition (\ref{EqCondNotWA}) in the proof of Theorem~\ref{Thm2.11'} that there exists an integer $n'$ with $1\leq n'\leq n$, an isoclinic summand $D_{\lambda'}$ of $\wedge^{n'}D$ of some slope $\lambda'\in\BZ$, and a non-zero element $f\in\sigma^\ast\!D_{\lambda'}$ generating an $F_D$-stable subspace of dimension one, such that $(z-\zeta)^{-\lambda'-1}f\in\wedge^{n'}\Fq_x$. Now the set of points
\[
\{\,y'\in Y:\es (z-\zeta)^{-\lambda'-1}f\in\wedge^{n'}\Fq_{\phi(y')}\,\}
\]
is non-empty and Zariski-closed in $Y$ because $\Fp/\Fq$ is locally free of finite rank over $\CH_s^\an$.
In particular we can find a point $y'$ in this set whose residue field $\kappa(y')$ is a finite extension of $\BF_{q^s}\dpl\zeta\dpr$. But then $\phi(y')\in \phi(Y^\rig)\setminus(\CH_b^{wa})^\rig$ which is a contradiction. 

\medskip
\noindent
(c)\es To prove the connectedness of $\CH^{wa}_b$ and $\CH^a_b$ note that every Berkovich space is locally arcwise connected by \cite[Theorem 3.2.1]{Berkovich1}. By Proposition~\ref{PropPeriodSpSmooth} the variety $\CH_{\ul w}$ has a finite covering by affine spaces. Hence the space $\CH^\an_s$ has a countable affinoid covering by polydiscs $U\cong\CM\bigl(\BF_{q^s}\dpl\zeta\dpr\langle\frac{x_1}{\theta},\ldots,\frac{x_n}{\theta}\rangle\bigr)$ for varying $\theta\in\BF_{q^s}\dpl\zeta\dpr$. Since the points $x\in\CH^{wa}_b$ with $\kappa(x)$ finite over $\BF_{q^s}\dpl\zeta\dpr$ lie dense in $\CH^{wa}_b$ 
by \cite[Proposition 2.1.15]{Berkovich1}
it suffices to exhibit for every such point $x\in\CH^{wa}_b\cap U$ a continuous map $\alpha$ from the compact interval $[0,1]$ into $\CH^{wa}_b\cap U$ such that $\alpha(0)=x$ and $\alpha(1)$ is the point corresponding to the supremum norm on $U$; see Example~\ref{ExAnalyticPoints}\ref{ExAnalyticPointsB}

So let $x\in\CH^{wa}_b\cap U$ with $\kappa(x)$ finite over $\BF_{q^s}\dpl\zeta\dpr$ be the point with coordinates $x_i=b_i\in\kappa(x)$, $|b_i|\le|\theta|$. We let $\alpha$ map $u\in[0,1]$ to the point
\[
P\bigl(\,(b_1,\ldots,b_n)\,,\,(u|\theta|,\ldots,u|\theta|)\,\bigr)\,;
\]
see Example~\ref{ExAnalyticPoints}\ref{ExAnalyticPointsC}. Then $\alpha$ is easily seen to be continuous. It remains to show that $\alpha(u)$ lies in $\CH^{wa}_b$. For $u>0$ the prime ideal $\ker|\,.\,|_{\alpha(u)}$ corresponding to $\alpha(u)$ is the zero ideal. Since we have shown in (b) above that the set $U\setminus\CH^{wa}_b$ is a union of Zariski closed subsets of $U$ and since it does not contain the point $\alpha(0)$, we obtain $\alpha(u)\in\CH^{wa}_b$.

At last consider the set $\CH^a_b$. If $x\in\CH^a_b\cap U$ is a point with $\kappa(x)$ finite over $\BF_{q^s}\dpl\zeta\dpr$ then $\alpha\bigl([0,1]\bigr)\subset\CH^a_b\cap U$ by Theorem~\ref{Thm3.5}. Namely, the value group of $\kappa\bigl(\alpha(u)\bigr)$ is generated by $\bigl|\BF_{q^s}\dpl\zeta\dpr(b_1,\ldots,b_n)\mal\bigr|$ and $u|\theta|$ and hence is finitely generated. So $\kappa\bigl(\alpha(u)\bigr)$ satisfies condition (\ref{CondDoubleStar}) from page~\pageref{CondDoubleStar}. This establishes the connectedness of $\CH^{wa}_b$ and $\CH^a_b$.
\end{proof}

We should also note what happens over $\BF_q^\alg$.

\begin{theorem}
\begin{enumerate}
\item 
Let $b\in\GL_n\bigl(\BF_q^{\,\alg}\dpl z\dpr\bigr)$. Then the sets $\CH^{wa}_b$ and $\CH^a_b$ of weakly admissible, respectively admissible, points inside the Berkovich space $\bigl(\CH_{\ul w}\otimes\BF_q^{\,\alg}\dpl\zeta\dpr\bigr)^\an$ are open, arcwise connected, and paracompact and $(\CH^{wa}_b)^\rig$ is an admissible open rigid analytic subspace of $\bigl(\CH_{\ul w}\otimes\BF_q^{\,\alg}\dpl\zeta\dpr\bigr)^\rig$.
\item 
If $b'=g\,b\,(g^{-1})^\sigma$ for some $g\in\GL_n\bigl(\BF_q^{\,\alg}\dpl z\dpr\bigr)$ then the automorphism $\Fq_x\mapsto g(\Fq_x)=:\Fq_{g(x)}$ of $\CH_{\ul w}\otimes\BF_q^{\,\alg}\dpl\zeta\dpr$ maps $\CH^{wa}_b$ (respectively $\CH^a_b$) isomorphically onto $\CH^{wa}_{b'}$ (respectively $\CH^a_{b'}$).
\item
If $b$ satisfies a decency equation with the integer $s$ then the spaces from (a) are obtained by base change from the corresponding spaces over $\BF_{q^s}\dpl\zeta\dpr$ constructed in Theorem~\ref{ThmParacompact}.
\end{enumerate}
\end{theorem}

\begin{proof}
Assertion (b) is obvious and (c) follows from Lemma~\ref{LemmaDecentAndStable}. Both (b) and (c) imply (a) in view of Proposition~\ref{Prop2.4}.
\end{proof}

%%%%%%%%%%%%%%%%%%%%%%%%%%%%%%%%%%%%%%%%%%%%%%%%%%%%%%%%%%%%%%%%%%%%%%
%
%    Examples of Period Spaces
%
%%%%%%%%%%%%%%%%%%%%%%%%%%%%%%%%%%%%%%%%%%%%%%%%%%%%%%%%%%%%%%%%%%%%%%

\subsection{Examples of Period Spaces} \label{SectExamplesPeriodSp}

We want to illustrate the theory developed in the previous sections.

\begin{example}[The Hopkins-Gross Period Space] \mbox{ }\\\label{Ex8.1}
Let $n$ be a positive integer and let $D=\BF_q\dpl z\dpr^{\oplus n}$. Fix a basis $e_1,\ldots,e_n$ of $D$ and let
\[
b\es=\es \left( \raisebox{6.2ex}{$
\xymatrix @C=0.3pc @R=0.3pc {
0 \ar@{.}[rrr]\ar@{.}[drdrdrdr] & & & 0 & z\\
1\ar@{.}[drdrdr]  & & &  & 0 \ar@{.}[ddd]\\
0 \ar@{.}[dd]\ar@{.}[drdr] & & & &  \\
& & & & \\
0 \ar@{.}[rr] & & 0 & 1 & 0\\
}$}
\right) \es\in\es\GL_n\bigl(\BF_q\dpl z\dpr\bigr)
\]
Then $b$ satisfies a decency equation with the integer $s=n$. Fix the Hodge-Pink weights $\ul w=(w_1=\ldots=w_{n-1}=0>w_n=-1)$. Let $\Fq_0$ be the span
\[
\Fq_0\es=\es\bigl\langle\,e_1,\ldots,e_{n-1}\,,\,(z-\zeta)^{-1}e_n\,\bigr\rangle_{\BFZ\dbl z-\zeta\dbr}\;.
\]
Its stabilizer
\[
\Stab\Fq_0\es=\es\Bigl\{\left(\begin{array}{c|c}\GL_{n-1} & (z-\zeta) \\ \hline \raisebox{0em}[1.1em]{\LARGE $\ast$} & \ast \end{array}\right)\Bigr\}\es\subset\es\GL_n\bigl(\BFZ\dbl z-\zeta\dbr\bigr)\,.
\]
is a maximal parahoric subgroup with quotient $\CH=\BP^{n-1}_\BFZ$. The universal pair $(\CP,\CQ)$ of $\sigma$-modules over $\CO_{\CH_s^\an}\ancon$ satisfies $\CQ\supset\CP=\CF_{-1,n}$ and $\deg\CQ=0$ (see Example~\ref{Ex1.2.8} for the definition of $\CF_{-1,n}$). Hence at every point $x\in\CH_s^\an$, $\CQ_x$ is isoclinic of slope zero by Theorem~\ref{Thm1a} and Proposition~\ref{Prop0.7}. So $\CH^a_{\ul w,b}=\CH^{wa}_{\ul w,b}=\CH^\an_s=\BP^{n-1}_{\BF_{q^n}\dpl\zeta\dpr}$. The automorphism group $J_b\bigl(\BF_q\dpl z\dpr\bigr)$ of the $z$-isocrystal equals $\Delta\!\mal$ where $\Delta$ is the central skew field of invariant $\frac{1}{n}$ over $\BF_q\dpl z\dpr$. The space $\CH_{\ul w,b}^{wa}$ with its $J_b\bigl(\BF_q\dpl z\dpr\bigr)$-action is the period space studied by Hopkins and Gross \cite{HG1,HG2}. Note that they also studied the analogous situation in mixed characteristic and interpreted $\BP^{n-1}$ as a period space of filtered isocrystals; see also \cite[5.50]{RZ}. Our example provides the hitherto missing analogous interpretation of $\BP^{n-1}$ in equal characteristic.
\end{example}

\begin{example}\label{Ex8.2}
Let $n=2$ and fix a positive integer $d$. Let $D=\BF_q\dpl z\dpr^{\oplus2}$ and let 
\[
b\es=\es\left(\begin{array}{cl} 0 & z^{-d} \\ 1 & 0 \end{array}\right) \es\in\es\GL_2\bigl(\BF_q\dpl z\dpr\bigr)\,.
\]
Then $b$ satisfies a decency equation with the integer $s=2$. Fix the Hodge-Pink weights $\ul w=(d\geq 0)$. The lattice
\[
\Fq_0\es=\es\Bigl\langle(z-\zeta)^d\binom{1}{0}\,,\,\binom{0}{1}\Bigr\rangle_{\BFZ\dbl z-\zeta\dbr}\;.
\]
has stabilizer
\[
\Stab\Fq_0\es=\es\Bigl\{\left(\begin{array}{cc}\ast &  (z-\zeta)^d \\ \ast & \ast\end{array}\right)\Bigr\}\es\subset\es\GL_2\bigl(\BFZ\dbl z-\zeta\dbr\bigr)\,.
\]
The stabilizer corresponds to the Borel subgroup of lower triangular matrices inside the group $\wt G:=Restr_{\bigl(\BFZ\dbl z-\zeta\dbr/(z-\zeta)^d\bigr)/\BFZ} \GL_2$ and so $\CH=Restr \BP^1_\BFZ$ is the full $(d-1)$-jet bundle over the projective line. It has dimension $d$. 

\medskip

Let us first assume that $d$ is odd. Then $\CH^{wa}=\CH_{\ul w,b}^{wa}=\CH^\an_s=\CH^\an\otimes\BF_{q^2}\dpl\zeta\dpr$ since there are no $F_D$-invariant subspaces of $D$. The automorphism group $J_b\bigl(\BF_q\dpl z\dpr\bigr)$ of $(D,F_D)$ equals $\Delta\!\mal$ where $\Delta$ is the central skew field of invariant $-\frac{d}{2}$ over $\BF_q\dpl z\dpr$.
The universal pair of $\sigma$-modules over $\CO_{\CH^\an_s}\ancon$ satisfies $\CQ\subset\CP=\CF_{d,2}$. Let $x\in\CH^\an_s$ be a point and let $\BC=\ol{\kappa(x)}$ be the completion of the algebraic closure of $\kappa(x)$. The point $x$ is admissible if and only if $\CQ_x$ is isoclinic of slope zero. By Theorem~\ref{Thm1a} and Proposition~\ref{Prop0.7} this is equivalent to
\[
\CQ(-1)^F(\BC)\es=\es\{\,f\in\CP(-1)^F(\BC): f\in\Fq_x\otimes\BC\dbl z-\zeta\dbr\,\}\es=\es0\,.
\]
Now we distinguish the cases
\begin{description}
\item[$d=1$]: \es We have $\CP(-1)^F(\BC)\cong(\CF_{-1,2})^F(\BC)=0$ by reasons of degree, hence $\CH_{\ul w,b}^a=\CH_{\ul w,b}^{wa}=\CH^\an_s$. This situation is dual to the period space of Hopkins-Gross from the previous example.
\item[$d=3$]: \es We have $\DS\CP(-1)^F(\BC)\es\cong\es(\CF_{1,2})^F(\BC)\es=$
\[
\es z^\BZ\cdot\bigl\{\,f=\sum_{\nu\in\BZ}z^{-\nu}\binom{u^{q^{2\nu}}}{u^{q^{2\nu+1}}}:\es u\in\BC\,,\,|\zeta^{q^2}|\leq|u|\leq|\zeta|\,\bigr\}\,,
\]
and the point $x$ is not admissible if and only if $f\in\Fq_x$ for such an $f$.
Thus inside the $3$-dimensional space $\CH^{wa}$ the set $\CH^{wa}\setminus\CH^a$ is parametrized by the $\BC$-valued points of the ``curve'' (one dimensional compact Berkovich space) $\CM\bigl(\BFZ\langle\frac{u}{\zeta},\frac{\zeta^{q^2}}{u}\rangle\bigr)$. However, the map from this ``curve'' to $\CH^{wa}$ is not analytic but only continuous. By Theorem~\ref{Thm3.5} it does not hit any classical rigid analytic point of $\CH^{wa}$ (these are the ones whose residue field is a finite extension of the base field $\BF_{q^2}\dpl\zeta\dpr$). It does not even hit a point whose residue field satisfies condition (\ref{CondDoubleStar}) on page~\pageref{CondDoubleStar}. At first glance it might seem surprising that there are points in $\CH_{\ul w,b}^{wa}$ whose residue field, which is topologically finitely generated, does not satisfy (\ref{CondDoubleStar}). See however Example~\ref{ExAnalyticPoints}\ref{ExAnalyticPointsD} which shows that this occurs naturally.
\item[$d=5$]: \es We have $\DS\CP(-1)^F(\BC)\es\cong\es(\CF_{3,2})^F(\BC)\es=$
\begin{eqnarray*}
\es z^\BZ\cdot\bigl\{\,f=\sum_{\nu\in\BZ}z^{-3\nu}\binom{u_0^{q^{2\nu}}\,\:+\,\:zu_1^{q^{2\nu}}\,\:+\,\:z^2u_2^{q^{2\nu}}}{u_0^{q^{2\nu+1}}+zu_1^{q^{2\nu+1}}+z^2u_2^{q^{2\nu+1}}}:\es u_0,u_1,u_2\in\BC,\\
|u_0|,|u_1|,|u_2|\leq|\zeta|\quad\text{and}\quad|\zeta^{q^2}|\leq|u_i|\text{ for some }i\,\bigr\}\,.
\end{eqnarray*}
\end{description}
To summarize, for odd $d\geq3$ the points of $\CH^{wa}\setminus\CH^a$ inside the $d$-dimensional space $\CH^{wa}$ are parametrized by a $(d-2)$-dimensional compact space.

\medskip

Now let us turn to the case where $d=2c$ is even. Fix an element $\lambda\in\BF_{q^2}\setminus\BF_q$ and let 
\[
g\es=\es\left(\begin{array}{ll} z^{-c} & \lambda z^{-c} \\ 1 & \lambda^q \end{array}\right)\qquad\text{and}\qquad b'\es=\es g^{-1}bg^\sigma\es=\es\left(\begin{array}{ll} z^{-c} & 0 \\ 0 & z^{-c} \end{array}\right)\,.
\]
Then $b'$ satisfies a decency equation with the integer $s=1$. So we replace $b$ by $b'$. The automorphism group $J_b\bigl(\BF_q\dpl z\dpr\bigr)$ of $(D,F_D)$ equals $\GL_2\bigl(\BF_q\dpl z\dpr\bigr)$. It acts on $\CH_s$ through $\PGL_2\bigl(\BF_q\dpl z\dpr\bigr)$. 
The universal pair of $\sigma$-modules over $\CO_{\CH^\an_s}\ancon$ satisfies $\CQ\subset\CP=\CO(c)^{\oplus2}$ and $D=\sigma^\ast D=\CP(-c)^F(\CH_s^\an)$. 
Let $x\in\CH_s^\an$ be a point and set $\BC:=\ol{\kappa(x)}$. By the criteria formulated in (\ref{EqCondOnF}) and (\ref{EqCondNotWA}) in the proofs of Theorems~\ref{Thm2.11} and \ref{Thm2.11'} the point $x$ is

\medskip

not weakly admissible if and only if there exists an 
\[
h\in\bigl(\CP(-c)\otimes\BC\ancon\bigr)^F(\BC)\,,\,h\neq0\quad\text{with}\quad(z-\zeta)^{c-1}h\in\Fq_x\otimes\BC\dbl z-\zeta\dbr\,,
\]

not admissible if and only if there exists an
\[
f\in\bigl(\CP(-1)\otimes\BC\ancon\bigr)^F(\BC)\,,\,f\neq0\quad\text{with}\quad f\in\Fq_x\otimes\BC\dbl z-\zeta\dbr\,.
\]
Again we distinguish the cases
\begin{description}
\item[$d=2$]: Then $c=1$ and $\CH^a=\CH^{wa}$. The points which are not weakly admissible are precisely the points
\begin{eqnarray*}
& \Fq\es=\es\left(\begin{array}{c}g(\zeta)+g'(\zeta) (z-\zeta)\\[1mm] 1\end{array}\right)
\cdot\kappa(x)\dbl z-\zeta\dbr\;+\;(z-\zeta)^2\Fp  \qquad\text{for}\quad g\in\BF_q\dpl z\dpr \\[2mm]
& \text{and}\qquad\Fq\es=\es\left(\begin{array}{c}1 \\ 0\end{array}\right)
\cdot\kappa(x)\dbl z-\zeta\dbr\;+\;(z-\zeta)^2\Fp\,.
\end{eqnarray*}
Consider the projection of $\CH_s^\an$, which is the tangent bundle over $\BP^1_\BFZ$, to $\BP^1_\BFZ$. In each fiber over a point in $\BP^1_\BFZ\bigl(\BFZ\bigr)$ there is exactly one point which is not weakly admissible.
\item[$d=4$]: Then $\CP(-1)^F(\BC)=\CO(1)^F(\BC)^{\oplus 2}=\bigl(\BF_q\dpl z\dpr\cdot\{\,f_\alpha:\alpha\in\BC\,,\,|\zeta^q|<|\alpha|\leq|\zeta|\,\}\bigr)^{\oplus2}$ by Proposition~\ref{Prop0.5}. Let $f=\binom{g\,f_\alpha}{\es f_\beta}$ with $\alpha,\beta\in\BC\,,\,|\zeta^q|<|\alpha|,|\beta|\leq|\zeta|$ and $g\in\BF_q\dpl z\dpr$. If $\alpha=\beta$ then $f=f_\alpha\cdot h$ for $h=\binom{g}{1}\in\CP(-2)^F(\BC)=D$. So the points
\begin{eqnarray*}
& \Fq\es=\es\left(\begin{array}{c}g(z)+a_3 (z-\zeta)^3 \\[1mm] 1\end{array}\right)
\cdot\kappa(x)\dbl z-\zeta\dbr\;+\;(z-\zeta)^4\Fp \quad\text{for}\quad g\in\BF_q\dpl z\dpr\,,\,a_3\in\kappa(x)\\[2mm]
& \text{and}\quad \Fq\es=\es\left(\begin{array}{c}1\\ a_3 (z-\zeta)^3\end{array}\right)
\cdot\kappa(x)\dbl z-\zeta\dbr\;+\;(z-\zeta)^4\Fp\quad\text{for}\quad a_3\in\kappa(x)
\end{eqnarray*}
are not weakly admissible and the points
\[
\Fq\es=\es\binom{g\,f_\alpha}{\es f_\beta}\cdot\kappa(x)\dbl z-\zeta\dbr\;+\;(z-\zeta)^4\Fp
\]
for $\alpha,\beta\in\BC\,,\,|\zeta^q|<|\alpha|,|\beta|\leq|\zeta|\,,\,\alpha\neq\beta\,,\,g\in\BF_q\dpl z\dpr\mal$
are weakly admissible but not admissible.
\end{description}
Along these lines on can as well determine the sets $\CH_b^a$ and $\CH_b^{wa}$ for $d>4$.
\end{example}

\begin{example}\label{Ex8.3}
We let again $n=2$ and $D=\BF_q\dpl z\dpr^{\oplus2}$. Let $d$ be a positive integer and set
\[
b\es=\es\left(\begin{array}{lc} z^{-d} & 0 \\ 0 & 1 \end{array}\right) \es\in\es\GL_2\bigl(\BF_q\dpl z\dpr\bigr)\,.
\]
It is decent with the integer $s=1$. Fix the Hodge-Pink weights $\ul w=(d>0)$. As in Example~\ref{Ex8.2}, $\CH$ is the full $(d-1)$-jet bundle over $\BP^1_\BFZ$ with dimension $d$. 
The automorphism group $J_b\bigl(\BF_q\dpl z\dpr\bigr)$ of the $z$-isocrystal equals the maximal split torus of diagonal matrices in $\GL_2\bigl(\BF_q\dpl z\dpr\bigr)$. It acts on $\CH_s^\an$ through $\PGL_2\bigl(\BF_q\dpl z\dpr\bigr)$.
The universal pair of $\sigma$-modules over $\CO_{\CH^\an_s}\ancon$ satisfies $\CQ\subset\CP=\CO(d)\oplus\CO(0)$. Here $(D,F_D)$ is not isoclinic. It has precisely two $F_D$-invariant subspaces namely $D'=\BF_q\dpl z\dpr\binom{1}{0}$ and $D'=\BF_q\dpl z\dpr\binom{0}{1}$. On the second subspace $F_{D'}=\sigma$, so there $\pair(D',F_{D'},\Fq_{D'})=(\CP',\CQ')$ always satisfies $\CP'=\CO(0)$ and $\deg\CQ'\leq0$. This is in accordance with weak admissibility.

For the first subspace we obtain $\CP'=\CO(d)$ and $\deg\CQ'>0$ if and only if $(z-\zeta)^{d-1}\binom{1}{0}\;\in\;\Fq$. Therefore the points which are not weakly admissible are
\[
\Fq\es=\es\left(\begin{array}{c}1\\[1mm] a_1 (z-\zeta)+\ldots+a_{d-1} (z-\zeta)^{d-1}\end{array}\right)
\cdot\kappa(x)\dbl z-\zeta\dbr\;+\;(z-\zeta)^d\Fp\quad\text{for}\es a_i\in\kappa(x)\,.
\]
We may therefore identify $\CH_b^{wa}$ with affine $d$-space $\BA_\BFZ^d$ over $\BFZ$ by mapping the $B$-valued point 
$(a_0,\ldots,a_{d-1})\in\BA_\BFZ^d(B)$ with values in an affinoid $\BFZ$-algebra $B$ to 
\[
\Fq\es=\es\binom{a_0+\ldots+a_{d-1}(z-\zeta)^{d-1}}{1}\cdot B\dbl z-\zeta\dbr\;+\;(z-\zeta)^d\Fp\es\in\es\CH_b^{wa}(B)\,.
\]
Under this identification $J_b\bigl(\BF_q\dpl z\dpr\bigr)$ acts on $\BA^d_\BFZ$ as follows. The element $\left(\begin{array}{cc}g&0\\0&1\end{array}\right)\in J_b\bigl(\BF_q\dpl z\dpr\bigr)$ maps the point $(a_0,\ldots,a_{d-1})\in\BA^d_\BFZ$ to
\[
\Bigl(g(\zeta)a_0\es,\es g(\zeta)a_1+g'(\zeta)a_0\es,\es\ldots\es,\es g(\zeta)a_{d-1}+g'(\zeta)a_{d-2}+\ldots+\frac{g^{(d-1)}}{(d-1)!}(\zeta)\,a_0\Bigr)\,.
\]
We use the criterion (\ref{EqCondOnF}) in the proof of Theorem~\ref{Thm2.11} to describe the admissible locus. Let $\BC$ be a complete algebraically closed extension of $L$ and note that there is an isomorphism
\[
\CP(-1)^F(\BC)\es=\es\bigl(\CO(d-1)\oplus\CO(-1)\bigr)^F(\BC)\es\lbij\es\CO(d-1)^F(\BC)\,,\quad f\cdot{\TS\binom{1}{0}}\leftmapsto f\,.
\]
In particular, if $x\in\CH_s^\an$ is not admissible then $\binom{(z-\zeta)^{d-1}}{0}\in\Fq_x$ since $f$ can have a zero at $z=\zeta$ at most of order $d-1$. So $x$ is not even weakly admissible. Therefore $\CH_b^a=\CH_b^{wa}=\BA_\BFZ^d$. In the next section we will construct interesting \'etale coverings of $\CH_b^a$.
\end{example}

%%%%%%%%%%%%%%%%%%%%%%%%%%%%%%%%%%%%%%%%%%%%%%%%%%%%%%%%%%%%%%%%%%%%%%
%
%    The Conjecture of Rapoport and Zink
%
%%%%%%%%%%%%%%%%%%%%%%%%%%%%%%%%%%%%%%%%%%%%%%%%%%%%%%%%%%%%%%%%%%%%%%

\subsection{The Conjecture of Rapoport and Zink} \label{SectConjRZ}

In this section we want to formulate and prove the analog of the conjecture of Rapoport and Zink \cite[p.\ 29]{RZ}, \cite[p.\ 429]{Rapoport94} mentioned in the introduction. It states the existence of a universal local system of $\BF_q\dpl z\dpr$-vector spaces on $\CH_{\ul w,b}^a$. For the notion of local system see Definition~\ref{DefLocalSystem}.
We can summarize the results from Section~\ref{SectAdmLocusOpen} as follows.

\begin{theorem} \label{Thm2.15}
The sheaf $\CQ^F$  of $F$-invariants induces a canonical local system $\CV_\CQ$ of $\BF_q\dpl z\dpr$-vector spaces on $(\CH^a_b)^\rig$ with $\CQ^F=\invlim\CV_\CQ$.
\end{theorem}

\begin{proof}
Consider the covering of $\CH^a_b$ by the $\CM(B_i)$, the finite \'etale $B_i$-algebras $C_i$, and the $\sigma$-modules $N_i$ over $C_i\langle\frac{z}{\zeta}\rangle$ from Corollary~\ref{Cor2.12}. On $U_i:=\CM(C_i)$ consider the sheaf $N_i^F$ of $F$-invariants. By Proposition~\ref{Prop13}, $T_zN_i:=\Bigl(\bigl(N_i\otimes C_i\dbl z\dbr/(z^m)\bigr)^F\Bigr)_{m\in\BN_0}$ is a local system of $\BF_q\dbl z\dbr$-lattices.
From the fact that the $N_i\otimes C_i\con$ descend to $\CM(B_i)$ and glue on all of $\CH_b^a$ we obtain isomorphisms over $U_i\times_{\CH^a_b}U_j$
\[
\phi_{ij}:\es T_zN_i\otimes\BF_q\dpl z\dpr\es\isoto\es T_zN_j\otimes\BF_q\dpl z\dpr
\]
which satisfy the cocycle condition. Then
\[
\CV_\CQ\es:=\es\bigl(\{U_i\to\CH_b^a\}\,,\,T_zN_i\,,\,\phi_{ij}\bigr)
\]
is a local system of $\BF_q\dpl z\dpr$-vector spaces on $(\CH_b^a)^\rig$ with
\[
\invlim\CV_\CQ|_{U_i}\es=\es(\invlim T_zN_i)\otimes_{\BF_q\dbl z\dbr}\BF_q\dpl z\dpr\es=\es\CQ^F|_{U_i}\,.
\]
The local system $\CV_\CQ$ is uniquely determined since it actually only depends on the descended $\sigma$-module $\bigl(N_i\otimes C_i\con\bigr)_i$ over $\CH_b^a$ which is canonical by Corollary~\ref{Cor2.12}.
\end{proof}

From the theorem we derive the analog of the conjecture of Rapoport and Zink.
Namely, recall that the pair $(\CP,\CQ)$ and the sheaf $\CQ^F$ were constructed in Section~\ref{SectPeriodSpaces} out of the faithful representation $V=\BF_q\dpl z\dpr^{\oplus n}$ in $\ul{\rm Rep}_{\BF_q\dpl z\dpr}\GL_n$. Moreover the constructions of the pair $(\CP,\CQ)$ of $\sigma$-modules, the sheaf $\CQ^F$, and the local system $\CV_\CQ$ are given by tensor functors in Theorems~\ref{Thm1.7} and \ref{Thm2.15} and Proposition~\ref{PropTateModuleExact}.
Since any $\BF_q\dpl z\dpr$-rational representation of $\GL_n$ occurs as a subquotient of $\bigoplus_{i=1}^rV^{\otimes n_i}\otimes(V\dual)^{\otimes m_i}$ for suitable $r,n_i$ and $m_i$ this yields a tensor functor
\[
I:\es\ul{\rm Rep}_{\BF_q\dpl z\dpr}\GL_n\es\longto\es\PLoc_{\CH^a_b}\,.
\]
For every geometric point $\bar x$ of $(\CH^a_b)^\rig$ with underlying analytic point $x$ of $\CH^a_b$ we obtain a fiber functor $I_x$ which associates with a representation of $\GL_n$ the fiber in $x$ of the corresponding local system. Explicitly this means the following. 
With a representation $\rho:\GL_n\to\Aut(W)$ we associate the Hodge-Pink structure 
\[
H_{b,\gamma_x}(W)\es=\es\Bigl(W\otimes\BF_q^{\,\alg}\dpl z\dpr,\,\rho(b)\cdot\sigma,\,\rho(\gamma_x)\Fp_W\Bigr)
\]
where $\Fp_W:=W\otimes_{\BF_q\dpl z\dpr}\kappa(x)\dbl z-\zeta\dbr$ and $\gamma_x$ is any element of $\GL_n\bigl(\kappa(x)\dpl z-\zeta\dpr\bigr)$ with $\gamma_x\cdot\Fp_V=\Fq_x$.
We further consider the associated pair $(\CP_W,\CQ_{W,x})=\pair\bigl(H_{b,\gamma_x}(W)\bigr)$ of $\sigma$-modules over $\kappa(x)\ancon$. Then $\CQ_{W,x}$ is isoclinic of slope zero by Definition~\ref{DefAdmPair} and Theorem~\ref{Thm2.6}. By Corollary~\ref{Cor6.12b} there exists a $\sigma$-module $N_W$ over $\kappa(x)\langle\frac{z}{\zeta}\rangle$ with $N_W\otimes\kappa(x)\ancon\cong\CQ_{W,x}$. Its Tate module $T_zN_W\bigl(\kappa(x)^\sep\bigr)$ is a free $\BF_q\dbl z\dbr$-module by Proposition~\ref{Prop13}. Now the fiber functor $I_x$ associates with $W$ the $\BF_q\dpl z\dpr$-vector space $V_zN_W:=T_zN_W\bigl(\kappa(x)^\sep\bigr)\otimes_{\BF_q\dbl z\dbr}\BF_q\dpl z\dpr$. Note that $V_zN_W$ does not depend on the particular choice of the lattice $N_W$ by Corollary~\ref{Cor6.12b}.
Since the analogs in mixed characteristic of the pairs $(\CP,\CQ)$ are the $(\phi,\Gamma)$-modules over the Robba ring (also called $p$-adic differential equations) this gives the precise analog of the conjecture of Rapoport and Zink in view of \cite[Theorem C]{Berger1}.

\begin{theorem}[Analog of the Rapoport-Zink Conjecture] \mbox{ }\\
There exists a faithful $\BF_q\dpl z\dpr$-linear exact tensor functor
\[
I:\es\ul{\rm Rep}_{\BF_q\dpl z\dpr}\GL_n\es\longto\es\PLoc_{\CH^a_b}
\]
such that for every point $x\in\CH^a_b$ the induced fiber functor $I_x$ associates with the representation $W$ of $\GL_n$ the $\BF_q\dpl z\dpr$-vector space $V_zN_W$.
\qed
\end{theorem}

\medskip

Next we want to show that not only every point of $\CH^a_b$ is admissible but that the whole universal family $(D,F_D,\Fq)$ over $\CH^a_b$ arises from a rigidified local shtuka with good reduction. For an admissible formal scheme $X$ over $\Spf \BF_{q^s}\dbl\zeta\dbr$ we denote by $X^\rig$ the associated rigid analytic space over $\BF_{q^s}\dpl\zeta\dpr$ and by $X_0$ its special fiber $X\mod\zeta$. Recall the notion of \'etale covering space of $X^\rig$ from Definition~\ref{DefEtaleCoveringSpace}.

\begin{theorem}\label{Thm2.16}
There exists a quasi-paracompact admissible formal scheme $X$ over $\Spf\BF_{q^s}\dbl\zeta\dbr$, such that its associated rigid analytic space $X^\rig$ is an \'etale covering space $\pi:X^\rig\to(\CH^a_b)^\rig$, and there exists a bounded rigidified local shtuka $(M,F_M,\delta_M)$ over $X$ with constant $z$-isocrystal
\[
M\otimes_{\CO_X\dbl z\dbr}\CO_{X_0}\dbl z\dbr[z^{-1}]\es\cong\es(D,F_D)\otimes_{\BF_{q^s}\dpl z\dpr}\CO_{X_0}\dbl z\dbr[z^{-1}]
\]
over $X_0$ and $\pi^\ast\Fq\;=\;\sigma^\ast\delta_M\circ F_M^{-1}\bigl(M\otimes_{\CO_X\dbl z\dbr}\CO_{X^\rig}\dbl z-\zeta\dbr\bigr)$ inside
\[
\TS\sigma^\ast\delta_M\bigl(\sigma^\ast M\otimes_{\CO_X\dbl z\dbr}\CO_{X^\rig}\dbl z-\zeta\dbr[\frac{1}{z-\zeta}]\bigr)\es=\es \sigma^\ast D\otimes_{\BF_{q^s}\dpl z\dpr}\CO_{X^\rig}\dbl z-\zeta\dbr[\frac{1}{z-\zeta}]\es=\es\pi^\ast\Fp[\frac{1}{z-\zeta}]\,.
\]
Moreover, one may choose $X$ such that $X^\rig$ is the space of $\BF_q\dbl z\dbr$-lattices inside the local system $\CV_\CQ$ from Theorem~\ref{Thm2.15}.
\end{theorem}

\begin{proof}
By Theorem~\ref{ThmParacompact}, $(\CH^a_b)^\rig$ is connected%
% by \cite[Proposition 3.3.4]{Berkovich1}
.
We choose a geometric base point $\bar x\in(\CH^a_b)^\rig$ and consider the continuous representation $\rho:\pi_1^\et\bigl((\CH^a_b)^\rig,\bar x\bigr)\to\GL_n\bigl(\BF_q\dpl z\dpr\bigr)$ induced by Proposition~\ref{Prop2.13} from the local system $\CV_\CQ$ on $(\CH^a_b)^\rig$. It gives rise to a continuous representation of $\pi_1^\et\bigl((\CH^a_b)^\rig,\bar x\bigr)$ on the discrete set $\GL_n\bigl(\BF_q\dpl z\dpr\bigr)/\GL_n\bigl(\BF_q\dbl z\dbr\bigr)$ of lattices inside $\CV_\CQ$. By Theorem~\ref{ThmFundGpPrinciple} the latter corresponds to (that is, it is trivialized by) an \'etale covering space $\pi:X_L\to(\CH^a_b)^\rig$ of $(\CH^a_b)^\rig$.
Hence there is a sheaf $N$ of $\sigma$-modules over $\CO_{X_L}\langle\frac{z}{\zeta}\rangle$ contained in $\pi^\ast\CQ$ on which $F_\CQ$ is an isomorphism, such that $\pi^\ast\CQ=N\otimes_{\CO_{X_L}\langle\frac{z}{\zeta}\rangle}\CO_{X_L}\ancon$. Now the assertion follows from Proposition~\ref{Prop2.6b}.
\end{proof}

\begin{remark}
\begin{enumerate}
\item 
 This theorem is optimal in the following sense. The Tate module 
\[
T_zM\es=\es\Bigl(\bigl(M\otimes\CO_{X^\rig}\dbl z\dbr/(z^m)\bigr)^F\Bigr)_{m\in\BN_0}
\]
of $M$ is a local system of $\BF_q\dbl z\dbr$-lattices on $X^\rig$ that satisfies $T_zM\otimes\BF_q\dpl z\dpr\;\cong\;\pi^\ast\CV_\CQ$. Since in general the local system of $\BF_q\dpl z\dpr$-vector spaces $\CV_\CQ$ contains $\BF_q\dbl z\dbr$-lattices only locally we can expect the existence of $M$ only on an \'etale covering space of $(\CH^a_b)^\rig$ and not on $(\CH^a_b)^\rig$ itself.
\item 
The space $X^\rig$ is the rigid analytic fiber of the \emph{Rapoport-Zink space} of the local isoshtuka $(D,F_D)$ constructed in \cite[Theorem 5.3]{HV1}. The above theorem provides a $\zeta$-adic formal model of the rigid analytic fiber $X^\rig$ of this Rapoport-Zink space.
\end{enumerate}
\end{remark}

\vspace{2cm}\pagebreak

%%%%%%%%%%%%%%%%%%%%%%%%%%%%%%%%%%%%%%%%%%%%%%%%%%%%%%%%%%%%%%%%%%%%%%
%
%    Appendix
%
%%%%%%%%%%%%%%%%%%%%%%%%%%%%%%%%%%%%%%%%%%%%%%%%%%%%%%%%%%%%%%%%%%%%%%

\begin{appendix}

\setcounter{tocdepth}{0}

\section{Background from Rigid Analytic and Berkovich Geometry} 
\setcounter{equation}{0}

\begin{flushright}
\begin{tabular}{p{0.9\textwidth}r}
A.1\quad Rigid Analytic Spaces and Formal Schemes \dotfill & \pageref{AppRigFormal}\\
A.2\quad Berkovich Spaces \dotfill & \pageref{AppBerkovich}\\
A.3\quad {\'E}tale Sheaves on Rigid Analytic Spaces \dotfill & \pageref{AppEtaleSheaves}\\
A.4\quad The \'Etale Fundamental Group \dotfill & \pageref{AppFundamentalGroups}\\
A.5\quad An Approximation Lemma \dotfill & \pageref{AppApproxLemma}\\
A.6\quad An Extension Lemma \dotfill & \pageref{AppExtensionLemma}
\end{tabular}\\
\end{flushright}

\bigskip

In this whole appendix let again $R\supset\BF_q\dbl\zeta\dbr$ be a rank-$1$ valuation ring which is complete and separated with respect to the $\zeta$-adic topology. We do not assume that $R$ is noetherian. Let $L$ be the fraction field and $\Fm_R$ be the maximal ideal of $R$.

%%%%%%%%%%%%%%%%%%%%%%%%%%%%%%%%%%%%%%%%%%%%%%%%%%%%%%%%%%%%%%%%%%%%%%
%
%    Rigid Analytic Spaces and Formal Schemes
%
%%%%%%%%%%%%%%%%%%%%%%%%%%%%%%%%%%%%%%%%%%%%%%%%%%%%%%%%%%%%%%%%%%%%%%

\subsection*{A.1\quad Rigid Analytic Spaces and Formal Schemes} \label{AppRigFormal}
\addtocounter{subsection}{1}
\setcounter{theorem}{0}

We do not intend to give an introduction to rigid analytic geometry here. We only want to explain the terminology used in this article and to direct the reader to further literature.

The $L$-algebra
\[
L\langle\ul y\rangle\es:=\es L\langle
y_1,\ldots,y_n\rangle\es:=\es\bigl\{\,f=\sum_{\ul i\in\BN_0^{\,n}}a_{\ul
  i}\,y_1^{i_1}\cdots y_n^{i_n}:\es a_i\in L\,,\,|a_{\ul i}|\to0\es\text{for}\es i_1+\ldots+i_n\to\infty\,\bigr\}
\]
 is called the \emph{Tate algebra} in $n$ variables over $L$. It is noetherian and carries the \emph{Gau{\ss} norm} $|f|:=\sup\{|a_{\ul i}|:\ul i\in\BN_0^{\,n}\}$ with respect to which it is complete. The Gau{\ss} norm is an \emph{$L$-algebra norm}, that is, it extends the absolute value on $L$. An \emph{affinoid $L$-algebra} is an $L$-algebra $B$ which can be described as a quotient of a Tate algebra for some $n$. For every presentation $B=L\langle\ul y\rangle/\Fa$ the residue norm of the Gau{\ss} norm is a complete $L$-algebra norm  on $B$. All these norms induce the same topology on $B$. More generally any $L$-algebra norm on $B$ defining this topology is called an \emph{$L$-Banach norm} on $B$. 

Now let $B$ be an affinoid $L$-algebra. For any maximal ideal $\Fm\subset B$ the residue field $B/\Fm$ is a finite extension of $L$ by \cite[Corollary 6.1.2/3]{BGR}. In rigid analytic geometry one equips the set $\Spm B$ of all maximal ideals of $B$ with a Grothendieck topology and a structure sheaf. One calls $\Spm B$ an \emph{affinoid rigid analytic space}. Every $L$-algebra homomorphism $B\to C$ induces a map $\Spm C\to\Spm B$ and these are precisely the \emph{morphisms} between affinoid rigid analytic spaces.
In the Grothendieck topology on $\Spm B$ any finite covering of $\Spm B$ by affinoid subdomains is admissible, where an \emph{affinoid subdomain} of $\Spm B$ is an affinoid space $\Spm B'$ together with an $L$-algebra homomorphism $B\hookrightarrow B'$ which identifies the set $\Spm B'$ with a subset $U$ of $\Spm B$ and which is universal for morphisms $\Spm A\to\Spm B$ of affinoid rigid analytic spaces whose image lies in $U$.
Now a \emph{rigid analytic space} over $L$ is a set $X$ carrying a Grothendieck topology and a structure sheaf, such that $X$ possesses an admissible covering by affinoid rigid analytic spaces. Rigid analytic spaces were invented by Tate~\cite{Tate}. Their theory is carefully presented in \cite{BGR}, see also \cite{Bosch}.

\medskip

Rigid analytic spaces can be viewed as generic fibers of certain formal schemes over $\Spf R$. See \cite[I$_{\rm new}$ \S10]{EGA} for the notion of formal scheme. We let 
\[
R\langle\ul y\rangle\es:=\es R\langle y_1,\ldots,y_n\rangle\es:=\es\bigl\{\,\sum_{\ul i\in\BN_0^{\,n}}a_{\ul i}\,\ul y^{\ul i} \in L\langle\ul y\rangle:\es a_{\ul i}\in R \es\text{for all}\es \ul i\in\BN_0^{\,n}\,\bigr\}\,.
\]
An \emph{admissible formal $R$-algebra} is a flat $R$-algebra which is a quotient of some $R\langle\ul y\rangle$. An \emph{admissible formal scheme} over $\Spf R$ is a formal scheme over $\Spf R$ which locally is $R$-isomorphic to the formal spectrum of an admissible formal $R$-algebra. For any admissible formal $R$-algebra $B\open$ the $L$-algebra $B:=B\open\otimes_R L$ is affinoid. The assignment $\Spf B\open\mapsto\Spm B$ can be globalized to a functor
\[
\rig:\es\{\text{admissible formal schemes over $\Spf R$}\}\es\longto\es\{\text{rigid analytic spaces over $L$}\}\,.
\]
Under this functor all admissible formal blowing-ups are mapped to isomorphisms, where an \emph{admissible formal blowing-up} of an admissible formal scheme $X$ over $\Spf R$ is by definition the formal completion of a blowing-up in a coherent sheaf of ideals $\CJ\subset\CO_X$ such that locally on $X$, $\CJ$ contains some power of $\zeta$. After imposing some finiteness condition and localizing the category of admissible formal schemes over $\Spf R$ by the system of admissible formal blowing-ups the functor $\rig$ becomes an equivalence of categories (see Theorem~\ref{ThmFormalRigBerkovich} for the precise statement). These ideas were introduced by Raynaud~\cite{Raynaud}. See \cite{FRG1}, \cite{FRG2}, or \cite{Bosch} for a thorough introduction.

%%%%%%%%%%%%%%%%%%%%%%%%%%%%%%%%%%%%%%%%%%%%%%%%%%%%%%%%%%%%%%%%%%%%%%
%
%    Appendix Berkovich Spaces
%
%%%%%%%%%%%%%%%%%%%%%%%%%%%%%%%%%%%%%%%%%%%%%%%%%%%%%%%%%%%%%%%%%%%%%%

\subsection*{A.2\quad Berkovich Spaces} \label{AppBerkovich}
\addtocounter{subsection}{1}
\setcounter{theorem}{0}

There is yet another variant of rigid analytic geometry which remedies the problem that rigid analytic spaces are not topological spaces in the classical sense. The theory was developed by Berkovich~\cite{Berkovich1,Berkovich2}. Let $B$ be an affinoid $L$-algebra with $L$-Banach norm $|\,.\,|$. Berkovich calls these algebras \emph{strictly $L$-affinoid}.

\begin{definition} \label{DefAnalyticPoint} 
An \emph{analytic point\/} $x$ of $B$ is a semi-norm $|\,.\,|_x:B \to \BR_{\geq 0}$ which satisfies: 
\begin{enumerate} 
\item 
$|f+g|_x \leq \max\{\,|f|_x,|g|_x\,\}$ \quad for all $f,g \in B$, 
\item  
$|fg|_x = |f|_x \,|g|_x$ \quad for all $f,g \in B$, 
\item 
$|\lambda|_x = |\lambda|$ \quad for all $\lambda \in L$, 
\item 
$|\,.\,|_x:B \to \BR_{\geq 0}$ is continuous with respect to the norm $|\,.\,|$ on $B$. 
\end{enumerate} 
The set of all analytic points of $B$ is denoted~$\CM(B)$. 
On $\CM(B)$ one considers the coarsest topology such that for every $f\in B$ the map $\CM(B) \to \BR_{\geq0}$ given by $x \mapsto |f|_x$ is continuous. In this topology, $\CM(B)$ is a compact Hausdorff space; see~\cite[Theorem 1.2.1]{Berkovich1}. It is further equipped with a structure sheaf. Such a space is called a \emph{strictly $L$-affinoid space}.

Every morphism $\phi:B\to C$ of affinoid $L$-algebras is automatically continuous and hence induces a continuous morphism $\CM(\phi):\CM(C) \to \CM(B)$ by mapping the semi-norm $C \to \BR_{\geq 0}$ to the composition $B\to C\to \BR_{\geq 0}$. By definition the $\CM(\phi)$ are the \emph{morphisms} in the category of strictly $L$-affinoid spaces. In particular, for an affinoid subdomain $\Spm B'\subset \Spm B$ this morphism identifies $\CM(B')$ with a closed subset of $\CM(B)$. 
 \end{definition} 
 
For every analytic point $x\in\CM(B)$ we let $\ker|\,.\,|_x\,:=\,\{\,b\in B: |b|_x=0\,\}$. It is a prime ideal in $B$. We define the \emph{(complete) residue field of $x$} as the completion with respect to $|\,.\,|_x$ of the fraction field of $B/\ker|\,.\,|_x$. It will be denoted $\kappa(x)$. There is a natural continuous homomorphism $B \to \kappa(x)$ of $L$-algebras. Conversely let $K$ be a \emph{complete extension} of $L$, by which we mean a field extension of $L$ equipped with an absolute value $|\,.\,|:K\to\BR_{\geq0}$ which restricts on $L$ to the norm of $L$ such that $K$ is complete with respect to $|\,.\,|$. Any continuous $L$-algebra homomorphism $B\to K$ defines on $B$ a semi-norm which is an analytic point.

Although $\kappa(x)$ is topologically finitely generated it may be quite large as one can see from example \ref{ExAnalyticPointsD} below, which was communicated to us by V.\ Berkovich.

\begin{example} \label{ExAnalyticPoints}
We want to give some examples for analytic points in $\CM(B)$. 
\begin{enumerate}
\item \label{ExAnalyticPointsA}
Every classical point $x$ of the rigid analytic space $\Spm B$ defines an analytic point. Namely, $x$ corresponds to a maximal ideal $\Fm_x\subset B$ and induces the semi-norm  
\[ 
B \longrightarrow B/\Fm_x \xrightarrow{\enspace |\,.\,|_x\;} \BR_{\geq 0} 
\] 
where $|\,.\,|_x$ is the unique absolute value on the finite extension $B/\Fm_x$ of $L$ extending the absolute value of $L$; see \cite[Corollary 6.1.2/3]{BGR}. 
Because of its finiteness over $L$ the residue field $B/\Fm_x$ is already complete with respect to $|\,.\,|_x$ and we have $\kappa(x)=B/\Fm_x$.
\item \label{ExAnalyticPointsB}
Now let $B=L\langle y_1,\ldots,y_n\rangle$ be the Tate algebra in $n$ variables. By \cite[Proposition 5.1.2/1]{BGR} the Gau{\ss} norm $\DS\Bigl|\sum_{\ul i\in\BN_0^{\,n}}a_{\ul i}\,\ul y^{\ul i}\Bigr|\,:=\,\sup\{\,|a_{\ul i}|:\;\ul i\in\BN_0^{\,n}\,\}$ is an absolute value on $B$ and defines an analytic point in $\CM(B)$.
\item \label{ExAnalyticPointsC}
Again for $B=L\langle y_1,\ldots,y_n\rangle$ we let $b=(b_1,\ldots,b_n)\in (L^\alg)^n$ with $|b_i|\le1$. We expand each $f\in B$ around $b$
\[
f\es=\es\sum_{\ul i\in\BN_0^{\,n}} c_{\ul i}\, (y_1-b_1)^{i_1}\cdots(y_n-b_n)^{i_n}\qquad\text{with}\quad c_{\ul i}\in L(b_1,\ldots,b_n)\,.
\]
For every element $u=(u_1,\ldots,u_n)$ of the $n$-cube $[0,1]^n$ we obtain an analytic point $P=P(b,u)$ in $\CM(B)$ by setting
\[
|f|_{\SSC P}\es:=\es \sup\{\,|c_{\ul i}|\,u_1^{i_1}\cdots u_n^{i_n}:\es\ul i\in\BN_0^{\,n}\,\}\,.
\]
If $b$ is fixed it is easy to see that this defines a continuous and injective map $[0,1]^n\to\CM(B),\,u\mapsto P(b,u)$ sending $(0,\ldots,0)$ to the point corresponding to $b$ as considered in (a), and $(1,\ldots,1)$ to the point corresponding to the Gau{\ss} norm on $B$ as in (b).
\item \label{ExAnalyticPointsD}
Let $B=L\langle y\rangle$. We will exhibit an analytic point $x\in\CM(B)$ with $\kappa(x)=\ol L$ the completion of an algebraic closure of $L$.

Let $L_n\subset L^\sep$ for $n\in\BN_0$ be a sequence of finite extensions of $L$ with 
\[
L\es =\es L_0\es \subsetneq\es \ldots\es \subsetneq \es L_n\es \subsetneq \es L_{n+1}\es \subsetneq \es \ldots\es \subsetneq \es \bigcup_{n\in\BN_0} L_n\es = \es L^\sep\,.
\]
It is not hard to inductively choose a sequence $a_{n+1}\in L_{n+1}$ for $n\in\BN$ with $L_{n+1}=L_{n}(a_{n+1})$ and
\[
|a_{n+1}-a_{n}|\es <{\TS\frac{1}{n}}\qquad\text{and}\qquad|a_{n+1}-a_{n}|\es<\es|\gamma(a_m)-a_m|
\]
for all $\gamma\in\Gal(L^\sep/L_{m-1})$ with $\gamma(a_m)\ne a_m$, and for all $m=1,\ldots, n$.
Then the limit $a:=\lim_{n\to\infty}a_n$ exists in $\ol L$ and satisfies $|a-a_m|\es<\es|\gamma(a_m)-a_m|$ for all $\gamma\in\Gal(L^\sep/L_{m-1})$ with $\gamma(a_m)\ne a_m$ and for all $m$. Let $K:=\wh{L(a)}$ be the closure of $L(a)$ inside $\ol L$. 

We claim that $K=\ol L$.

\noindent
\emph{Proof.} Consider the action of $G_L:=\Gal(L^\sep/L)$ on $\ol L$ induced by continuity from the fact that $L^\sep$ is dense in $\ol L$; see \cite[3.4.1/6]{BGR}. We have $H:=\{\gamma\in G_L:\gamma|_K=\id_K\}=\Gal(K^\sep/K)$ for the separable closure $K^\sep$ of $K$ inside $\ol L$ since each element of $\Gal(K^\sep/K)$ acts continuously with respect to $|\,.\,|$. We set $M:=(L^\sep)^H\subset K$. Then $L\subset M=L^\sep\cap K$ and one easily sees that $M$ is henselian with respect to $|\,.\,|$ because $K$ is complete. By the Ax-Sen-Tate Theorem~\cite[p.\ 417]{Ax} we obtain $K\subset\ol L^H=\wh{M^\perf}$, the closure of the perfect hull of $M$. 

We show by induction that $L_n\subset M$ for all $n\in\BN$. So assume $L_{n-1}\subset M$. For each $n$ there is an element $b_n\in M^\perf$ with $|a-b_n|\le|a_{n+1}-a_n|$. In particular $|a_n-b_n|<|\gamma(a_n)-a_n|$ for all $\gamma\in\Gal(L^\sep/L_{n-1})$ with $\gamma(a_n)\ne a_n$. By Krasner's Lemma \cite[3.4.2/2]{BGR} this implies $a_n\in L_{n-1}(b_n)\subset M^\perf$. Since $a_n$ is separable over $L$ we find $a_n\in M$ and $L_n=L_{n-1}(a_n)\subset M$ as desired. We conclude that $L^\sep\subset M\subset K$ and $\ol L=\wh{L^\sep}=K$. \qed

Now if $B=L\langle y\rangle$ we obtain an analytic point $x\in\CM(B)$ with $\kappa(x)=\ol L$ by mapping $y$ to the element $a\in\ol L$ constructed above. 
\end{enumerate}
\end{example}

In \cite[\S1.2]{Berkovich2} Berkovich defines the category of \emph{strictly $L$-analytic spaces}. These spaces are topological spaces which admit an atlas with strictly $L$-affinoid charts. Berkovich calls such a space \emph{good} if  every point of it has a strictly $L$-affinoid neighborhood. The good strictly $L$-analytic spaces are the spaces studied in \cite{Berkovich1}, see \cite[\S1.5]{Berkovich2}.

\begin{definition}\label{DefBerkovichSpaces}
In this article we want to refer to good strictly $L$-analytic spaces as \emph{Berkovich spaces}. Also we here use the term \emph{affinoid covering} of a Berkovich space $X$ for a covering $\{U_i\}_i$ of $X$ by strictly $L$-affinoid subspaces $U_i\subset X$ such that the open interiors of the $U_i$ in $X$ still cover $X$.
\end{definition}

\begin{example}
Since our definition of \emph{affinoid covering} is non-standard we give an example. Let $X=\CM\bigl(L\langle y\rangle\bigr)$ be the closed unit disc over $L$. Then $X=\CM\bigl(L\langle y,\frac{\zeta}{y}\rangle\bigr)\cup \CM\bigl(L\langle\frac{y}{\zeta}\rangle\bigr)$, but this \emph{is not} an affinoid covering in our sense. Namely $\CM\bigl(L\langle y,\frac{\zeta}{y}\rangle\bigr)\open=\{|\zeta|<|y|\le1\}$ and $\CM\bigl(L\langle\frac{y}{\zeta}\rangle\bigr)\open=\{|y|<|\zeta|\}$ and these two open sets do not cover $X$. However, $X=\CM\bigl(L\langle y,\frac{\zeta^2}{y}\rangle\bigr)\cup \CM\bigl(L\langle\frac{y}{\zeta}\rangle\bigr)$ \emph{is} an affinoid covering in our sense, since $\CM\bigl(L\langle y,\frac{\zeta^2}{y}\rangle\bigr)\open=\{|\zeta|^2<|y|\le1\}$ and $\{|y|<|\zeta|\}$ cover $X$. 
\end{example}

The spaces $\CM(B)$ for affinoid $L$-algebras $B$ are examples for Berkovich spaces. 
Other examples for Berkovich spaces arise from schemes $Y$ which are locally of finite type over $L$, see \cite[\S3.4]{Berkovich1}. Namely if $Y=\BA_L^n$ is affine $n$-space over $L$ the associated Berkovich space $Y^\an=(\BA_L^n)^\an$ consists of all semi-norms $|\,.\,|_x$ on the polynomial ring $L[y_1,\ldots,y_n]$ as in Definition~\ref{DefAnalyticPoint} (a) - (c). The topology on $(\BA_L^n)^\an$ is the coarsest topology, such that for all polynomials $f\in L[y_1,\ldots,y_n]$ the map $(\BA_L^n)^\an\to\BR_{\ge0}, x\mapsto |f|_x$ is continuous. The space $(\BA_L^n)^\an$ is the union of the increasing sequence of compact polydiscs $\CM\bigl(L\langle \frac{y_1}{\zeta^{-m}},\ldots,\frac{y_n}{\zeta^{-m}}\rangle\bigr)$ of radii $(|\zeta^{-m}|,\ldots,|\zeta^{-m}|)$ for $m\in\BN_0$. If $Y\subset\BA_L^n$ is a closed (algebraic) subscheme of affine $n$-space with ideal sheaf $\CJ$, the ideal sheaf $\CJ\CO_{(\BA_L^n)^\an}$ defines a closed strictly $L$-analytic subspace $Y^\an$ of $(\BA_L^n)^\an$. Finally, if $Y$ is arbitrary and $\{Y_i\}_i$ is a covering of $Y$ by affine open subschemes then one can glue the associated strictly $L$-affinoid spaces $Y_i^\an$ to the Berkovich space $Y^\an$. Moreover $Y^\an$ is Hausdorff if and only if the scheme $Y$ is separated, see \cite[Theorems 3.4.1 and 3.4.8]{Berkovich1}.

\medskip

The relation between Berkovich spaces, rigid analytic spaces, and formal schemes is as follows.
With every strictly $L$-analytic space $X$ which is Hausdorff one can associate a quasi-separated rigid analytic space
\[
X^\rig\es:=\es\{\,x\in X:\es\kappa(x) \text{ is a finite extension of }L\,\},
\]
see \cite[\S1.6]{Berkovich2}. Recall that a rigid analytic space is called \emph{quasi-separated} if the intersection of any two affinoid subdomains is a finite union of affinoid subdomains. To describe the subcategories on which the functor $X\mapsto X^\rig$ is an equivalence we need the following terminology. A topological Hausdorff space is called \emph{paracompact} if every open covering $\{U_i\}_i$ has a locally finite refinement $\{V_j\}_j$, where \emph{locally finite} means that every point has a neighborhood which meets only finitely many of the $V_j$. On the other hand an admissible covering of a rigid analytic space is said to be \emph{of finite type} if every member of the covering meets only finitely many of the other members. A rigid analytic space over $L$ is called \emph{quasi-paracompact} if it possesses an admissible affinoid covering of finite type. Similarly we define the notions \emph{of finite type} and \emph{quasi-paracompact} also for (an open covering of) an admissible formal $R$-scheme.

\begin{theorem}\label{ThmFormalRigBerkovich}
The following three categories are equivalent:
\begin{enumerate}
\item 
the category of paracompact strictly $L$-analytic spaces,
\item 
the category of quasi-separated quasi-paracompact rigid analytic spaces over $L$, and
\item 
the category of quasi-paracompact admissible formal $R$-schemes, localized by admissible formal blowing-ups.
\end{enumerate}
\end{theorem}

\begin{proof}
It is shown in \cite[Theorem 1.6.1]{Berkovich2} that $X\mapsto X^\rig$ is an equivalence between (a) and (b). The equivalence of (c) with (b) is due to Raynaud, see \cite[Theorem 2.8/3]{Bosch} for a proof.
\end{proof}

Regarding paracompactness the following lemma is useful.

\begin{lemma} \label{LemmaParacompact}
Let $X$ be a Berkovich space over $\BFZ$. Assume that $X$ admits a countable affinoid covering (in the sense of Definition~\ref{DefBerkovichSpaces}). Then $X$ possesses a countable fundamental system of neighborhoods consisting of affinoid Berkovich subspaces. In particular if $X$ is Hausdorff every open subset of $X$ is a paracompact Berkovich space.
\end{lemma}

\begin{proof}
% BECAUSE $x\in X$ implies $\exists U\subset X$ with $x\in U\open$ the open interior of $U$ in $X$. If $x\in V\subset X$ open, then $x\in V\cap U\open\subset U$ open. Therefore $\exists W\subset V\cap U\open$, $W$ affinoid, $x\in W\open{}^U$ the open interior of $W$ in $U$. We have $W\open{}^U=\wt W\cap U$ for $\wt W\subset X$ open. Hence $x\in \wt W\cap V\cap U\open\subset W\open{}^U\subset W$ and $\wt W\cap V\cap U\open\subset X$ open.
It suffices to show that every member $U=\CM(B)$ of the countable affinoid covering of $X$ has a countable fundamental system of affinoid neighborhoods. Since $\BF_q(\zeta)$ is countable and dense in $\BFZ$, the ring $B$ possesses a countable dense subset $\{f_i\}_{i\in\BN_0}$. For rational numbers $r'\geq r$ we let $V(f_i;r,r')$ (respectively $V(f_i;r,\infty)$) be the strictly $\BFZ$-affinoid subset of $U$ on which $|\zeta|^{r'}\leq |f_i|\leq|\zeta|^r$ (respectively $|f_i|\leq|\zeta|^r$). Note that if $r=\frac{c}{d}$ and $r'=\frac{c'}{d}$ for integers $c,c',d$ then $V(f_i;r,r')\;=\;\CM\bigl(B\langle Y,Y'\rangle/(f_i^d-\zeta^c Y,f_i^d Y'-\zeta^{c'})\bigr)$. We claim that the countable collection
\[
\bigcap_{j=1}^\ell V(f_{i_j};r_j,r'_j) \qquad\text{for}\quad \ell\in\BN_0\,,\,i_j\in\BN_0\,,\,r_j\in\BQ\,,\,r'_j\in\BQ\cup\{\infty\}\,,\,r_j\leq r'_j
\]
is a fundamental system of neighborhoods for $U$. Indeed, let $V\subset U$ be any open subset and let $x_0\in V$. By definition of the topology on $U$ there exist $g_1,\ldots,g_\ell\in B$ and open intervals $I_1,\ldots,I_\ell\subset\BR_{\geq0}$ such that
\[
x_0\es\in\es\{\,x\in U:\es|g_j|_x\in I_j \text{ for all }j=1,\ldots,\ell\,\}\es\subset\es V\,.
\]
For all $j$ we can find numbers $r_j,r'_j\in\BQ$ with $r'_j\geq r_j$ and
\[
|g_j|_{x_0}\es\in\es\bigl(\,|\zeta|^{r'_j},|\zeta|^{r_j}\bigr)\es\subset\es\bigl[\,|\zeta|^{r'_j},|\zeta|^{r_j}\bigr]\es\subset\es I_j
\]
if $|g_j|_{x_0}>0$, respectively numbers $r_j\in\BQ$ and $r'_j=\infty$ with
\[
|g_j|_{x_0}\es\in\es\bigl[\,0\,,|\zeta|^{r_j}\bigr)\es\subset\es\bigl[\,0\,,|\zeta|^{r_j}\bigr]\es\subset\es I_j
\]
if $|g_j|_{x_0}=0$. We may further approximate $g_j$ by some $f_{i_j}$ such that $|f_{i_j}-g_j| <|\zeta|^{r'_j}$ if $r'_j\neq\infty$, respectively $|f_{i_j}-g_j| <|\zeta|^{r_j}$ if $r'_j=\infty$. Then 
\[
\bigcap_{j=1}^\ell V(f_{i_j},r_j,r'_j)
\]
is an affinoid neighborhood of $x_0$ contained in $V$ as desired.

If $X$ is Hausdorff it remains to show that any open subset $V$ of $X$ is paracompact. For every point $x\in V$ we choose an affinoid neighborhood $V_x\subset V$ from the countable fundamental system of neighborhoods for $X$ whose existence we just have established. Then $V$ is the union of the countably many compact sets $V_x$ and hence paracompact by \cite[Theorem I.9.5]{TG}.
\end{proof}

%%%%%%%%%%%%%%%%%%%%%%%%%%%%%%%%%%%%%%%%%%%%%%%%%%%%%%%%%%%%%%%%%%%%%%
%
%    Etale Sheaves on Berkovich Spaces
%
%%%%%%%%%%%%%%%%%%%%%%%%%%%%%%%%%%%%%%%%%%%%%%%%%%%%%%%%%%%%%%%%%%%%%%

\subsection*{A.3\quad {\'E}tale Sheaves on Rigid Analytic Spaces}\label{AppEtaleSheaves} 
\addtocounter{subsection}{1}
\setcounter{theorem}{0}

A morphism $f:Y\to X$ of rigid analytic spaces over $L$ is called \emph{{\'e}tale} if for every (classical) point $y\in Y$ the induced homomorphism of local rings $\CO_{X,f(y)} \to \CO_{Y,y}$ is flat and unramified. See~\cite[{\S} 3]{JP} for a thorough discussion of this notion. 
 
Let $X$ be a rigid analytic space over $L$. We recall the definition of the {\'e}tale site of $X$ from de Jong, van der Put  
\cite[{\S} 3.2]{JP}. The underlying category of the site $X_\et$ is the category of all {\'e}tale morphisms $f: Y\to X$ of rigid analytic spaces over $L$. A morphism from $f$ to $f'$ is a morphism $g:Y \to Y'$ such that $f'\circ g=f$. The morphism $g$ is automatically {\'e}tale. 
 
\begin{definition} \label{DefEtCov} 
A family of {\'e}tale morphisms $\{\,g_i:Z_i \to Y\,\}_{i\in I}$ is a \emph{covering for the {\'e}tale topology\/} if it has the following property: 
 
\medskip 
 
\noindent 
\hspace{0.05\textwidth} 
\parbox{0.9\textwidth}{ 
For every (some) choice of admissible affinoid covering $Z_i = \bigcup_j Z_{i,j}$ one has $Y = \bigcup_{i,j} g_i(Z_{i,j})$, and this is an admissible covering in the Grothendieck topology of $Y$.}  
\end{definition}

\medskip 
 
Clearly any admissible covering of $Y$ is a covering for the {\'e}tale topology. 
 
The property in Definition~\ref{DefEtCov} is local on $Y$ in the following sense: if $Y=\bigcup Y_l$ is an admissible affinoid covering, then $\{\,g_i:Z_i \to Y\,\}$ is a covering for the {\'e}tale topology if and only if for all $l$ the same is true for the covering $\{\,g_i:g_i^{-1}(Y_l) \to Y_l\,\}$. This implies that if $\{\,Z_i \to Y\,\}$ and $\{\,W_{i,j} \to Z_i\,\}$ for all $i$ are coverings for the {\'e}tale topology, then $\{\,W_{i,j} \to Y\,\}$ is a covering for the {\'e}tale topology.  
 
\smallskip 
 
The category $X_\et$ equipped with the family of  coverings for the {\'e}tale topology is thus a site, called the \emph{{\'e}tale site of $X$}. The sheaves on this site are called \emph{{\'e}tale sheaves on $X$}.

\begin{example}\label{ExsOfEtSh} 
The following are examples of {\'e}tale sheaves of abelian groups on $X$ which we need in this article ($f:Y\to X$ will denote a general object of $X_\et$): 
\begin{enumerate} 
\item 
The structure sheaf $\BG_a$ defined by $Y \mapsto \Gamma(Y,\CO_Y)$. 
\item 
For any quasi-coherent sheaf $\CF$ of $\CO_X$-modules on $X$ we define the {\'e}tale sheaf $W(\CF)$ on $X_\et$ by $Y\mapsto \Gamma(Y,f^\ast\CF)$, where $f^\ast$ denotes the pullback of quasi-coherent modules. In particular, one has $W(\CO_X)\cong\BG_a$. Any {\'e}tale sheaf $W(\CF)$ is a sheaf of $\BG_a$-modules. 
\item 
For any group or ring $A$ the constant {\'e}tale sheaf $\ul{A}_X$ is defined by $Y \mapsto \prod_{\pi_0(Y)}A$, where $\pi_0(Y)$ is the set of connected components of $Y$. (The restriction maps are the obvious ones.) If $X$ is clear from the context we also write $\ul{A}$ instead of $\ul{A}_X$. 
\item 
If $A$ is a ring, an \emph{{\'e}tale sheaf of $\ul{A}_X$-modules\/} is an {\'e}tale sheaf on $X$ which is a sheaf of modules over the sheaf of rings $\ul{A}_X$.  
\end{enumerate}
\end{example}

As usual \'etale sheaves possess stalks at geometric points. Since in this article we are only interested in quasi-separated quasi-paracompact rigid analytic spaces $X$ over $L$ whose associated strictly $L$-analytic space $X^\an$ is good, that is a Berkovich space in the sense of Definition~\ref{DefBerkovichSpaces}, we recall the definition only in this case. A \emph{geometric point} $\bar x$ of $X$ is a morphism $\bar x:\Spm K\to X\wh\otimes_L K$ for an algebraically closed complete extension $K$ of $L$. Here $X\wh\otimes_L K$ denotes the rigid analytic space over $K$ obtained from $X$ by base change; see \cite[\S9.3.6]{BGR}. The geometric point $\bar x$ can be viewed as a morphism $B\to K$ for a suitable affinoid subdomain $\Spm B\subset X$. Then the absolute value on $K$ defines an analytic point $x\in\CM(B)$ which we call the \emph{underlying analytic point} of $\bar x$. We may even choose $B$ such that $\CM(B)$ is an affinoid neighborhood of $x$ in $X^\an$. If $\CF$ is a sheaf on $X_\et$ and $\bar x$ is a geometric point of $X$ the \emph{stalk} $\CF_{\bar x}$ of $\CF$ at $\bar x$ is the inductive limit
\[
\CF_{\bar x}\es:=\es\dirlim[U] \CF(U)
\]
over all \emph{\'etale neighborhoods} $U$ of $\bar x$, that is pairs $(f,\bar y)$ consisting of \'etale morphisms $f:U\to X$ and commutative diagrams
\[
\xymatrix {
\Spm K \ar[r]^{\bar y} \ar[dr]_{\bar x} & U\wh\otimes_L K \ar[d]^{f} \\
& \es X\wh\otimes_L K\es.
}
\]

\medskip

Corresponding to the definitions of this section one can also define the \'etale site of a Berkovich space; see \cite[\S4.1]{Berkovich2}. For a paracompact Berkovich space $X^\an$ with associated rigid analytic space $X$ one obtains a morphism of sites $X\to X^\an$. The sheaves on $X_\et$ we meet in this article in fact come from sheaves on the \'etale site $X_\et^\an$. However we have preferred to formulate our results with rigid analytic spaces since we believe that they are more widely known than Berkovich spaces. Anyway Theorem~\ref{ThmFormalRigBerkovich} allows us to use both rigid analytic and Berkovich spaces whichever is more convenient.

%%%%%%%%%%%%%%%%%%%%%%%%%%%%%%%%%%%%%%%%%%%%%%%%%%%%%%%%%%%%% 
%
%     The \'Etale Fundamental Group
%
%%%%%%%%%%%%%%%%%%%%%%%%%%%%%%%%%%%%%%%%%%%%%%%%%%%%%%%%%%%%% 
 
\subsection*{A.4\quad The \'Etale Fundamental Group} \label{AppFundamentalGroups} 
\addtocounter{subsection}{1}
\setcounter{theorem}{0}

In this section we want to recall de Jong's \cite{dJ} definition of the \emph{{\'e}tale fundamental group\/} of a rigid analytic space $X$. Here again there is the corresponding definition for Berkovich spaces. We restrict ourselves in this section to the case where $X$ is quasi-separated and quasi-paracompact and its associated strictly $L$-analytic space $X^\an$ is good. We call the points of $X^\an$ \emph{analytic points} of $X$. Since $X^\an$ is assumed to be good every analytic point $x$ of $X$ has a fundamental system of neighborhoods consisting of affinoid subdomains $\Spm B\subset X$ such that $\CM(B)$ is a neighborhood of $x$ in $X^\an$.
 
\begin{definition} \label{DefEtaleCoveringSpace} 
Let $\pi:Y\to X$ be a morphism of rigid analytic spaces over $L$. 
$Y$ is called an \emph{{\'e}tale covering space of $X$\/} if $Y$ is quasi-paracompact and if for every analytic point $x$ of $X$ there exists an affinoid subspace $\Spm B\subset X$ such that $\CM(B)$ is a neighborhood of $x$ in $X^\an$ and such that $Y\times_X\Spm B$ is a disjoint union of affinoids mapping finite \'etale to $\Spm B$.
The category of {\'e}tale covering spaces of $X$ is denoted $\ul{\Cov}_X^\et$. 
\end{definition}

All {\'e}tale covering spaces are coverings for the {\'e}tale topology in the site $X_\et$; see Definition~\ref{DefEtCov}. The reader should not confuse the concepts of {\'e}tale covering spaces and coverings for the {\'e}tale topology. 
Note that the strictly $L$-analytic space associated with an \'etale covering space of $X$ is again a Berkovich space.

\medskip 
 
In order to define the \'etale fundamental group let $\bar x:\Spm\to X\otimes_LK$ be a geometric point of $X$. As in \cite{SGA1} consider the fiber functor at $\bar x$ 
\[ 
F_{\bar x}^\et: \es\ul{\Cov}_X^\et \es\longto \es\ul{\rm Sets}\es,\qquad F_{\bar x}^\et\bigl(\pi:Y\to X\bigr) \es:=\es\{\,\bar{y}:\Spm K\to Y\otimes_LK \es\text{with} \es \pi\circ\bar{y} = \bar{x}\,\} 
\] 
to the category of sets.
 
\begin{definition} \label{DefFundamentalGroup} 
We define the \emph{\'etale fundamental group} as the automorphism group of $F_{\bar x}^\et$
\[ 
\pi_1^\et(X,\bar{x}) \es := \es \Aut(F_{\bar x}^\et)\,. 
\] 
It is topologized as follows. The fundamental open neighborhoods of 1 are the 
stabilizers $\Stab_{X',\bar{x}'}$ in $\pi_1^\et(X,\bar{x})$ of arbitrary geometric points $\bar{x}'$ above $\bar{x}$ in arbitrary \'etale covering spaces $X'\in\ul{\Cov}_X^\et$.  
\end{definition}

The \'etale fundamental group classifies the \'etale covering spaces in the following sense. Let $\pi_1^\et(X,\bar{x})$-\ul{Sets} be the category of discrete sets endowed with a continuous left action of $\pi_1^\et(X,\bar{x})$. 
If $Y\in \ul{\Cov}_X^\et$ then $F_{\bar x}^\et(Y)$ naturally is an object of $\pi_1^\et(X,\bar{x})$-\ul{Sets}. The following theorem is due to de Jong \cite[Theorem 2.10 and \S5]{dJ}

\begin{theorem} \label{ThmFundGpPrinciple}
The fiber functor  
\[ 
F_{\bar x}^\et: \es\ul{\Cov}_X^\et \es\longto\es \pi_1^\et(X,\bar{x})\text{\rm -\ul{Sets}} 
\] 
is fully faithful, and extends to an equivalence of categories 
\[ 
\{\text{disjoint unions of objects of }\ul{\Cov}_X^\et\} \es\longto\es \pi_1^\et(X,\bar{x})\text{\rm -\ul{Sets}} 
\] 
Connected coverings correspond to $\pi_1^\et(X,\bar{x})$-orbits. 
\qed
\end{theorem} 
 
\medskip

We finally want to define the notion of local systems of $\BF_q\dpl z\dpr$-vector spaces.

\begin{definition}\label{DefLocalSystem}
We define a \emph{local system of $\BF_q\dbl z\dbr$-lattices} on $X$ as a projective system $\CF=(\CF_n,i_n)$ of sheaves $\CF_n$ of $\BF_q\dbl z\dbr/(z^n)$-modules on $X_\et$ such that $\CF_n$ is a locally constant free $\ul{\BF_q\dbl z \dbr/(z^n)}$-module of finite rank 
and $i_n$ induces an isomorphism of sheaves of $\ul{\BF_q\dbl z \dbr/(z^{n-1})}$-modules
\[
i_n\otimes\id:\es\CF_n\otimes_{\ul{\BF_q\dbl z \dbr/(z^n)}}\;\,\ul{\BF_q\dbl z \dbr/(z^{n-1})}\es\isoto\es\CF_{n-1}\,.
\]
(Of course locally constant means locally for the \'etale topology.) The category $\BLoc_X$ of local systems of $\BF_q\dbl z\dbr$-lattices with the obvious morphisms is an additive $\BF_q\dbl z\dbr$-linear rigid tensor category.
If $\bar x$ is a geometric point of $X$ we define the \emph{stalk} $\CF_{\bar x}$ of $\CF$ at $\bar x$ as 
\[
\CF_{\bar x}\es:=\es\invlim(\CF_{n,\bar x},\,i_n)\,.
\]
It is a finite free $\BF_q\dbl z\dbr$-module. Starting from $\BF_q\dbl z\dbr$-lattices one defines local systems of $\BF_q\dpl z\dpr$-vector spaces as in \cite[\S5]{dJ}.
In concrete terms this means that a \emph{local system of $\BF_q\dpl z\dpr$-vector spaces} on $X$ is given by the following data
\[
\CV\es=\es\bigl(\{U_i\to X\},\,\CF_i,\,\phi_{ij}\bigr)
\]
where
\begin{itemize}
\item
$\{U_i\to X\}$ is a covering for the \'etale topology on $X$,
\item
$\CF_i$ is a local system of $\BF_q\dbl z\dbr$-lattices over $U_i$ for each $i$,
\item
for each pair $i,j$, $\phi_{ij}$ is an invertible section over $U_i\times_XU_j$ of the sheaf 
\[
\Hom_{\BLoc_X}\bigl(\CF_i|_{U_i\times_X U_j}\,,\,\CF_j|_{U_i\times_X U_j}\bigr)\otimes_{\BF_q\dbl z\dbr}\BF_q\dpl z\dpr\,.
\]
\end{itemize}
These data are subject to the cocycle condition $pr_{jk}^\ast(\phi_{jk})\circ pr_{ij}^\ast(\phi_{ij})\;=\;pr_{ik}^\ast(\phi_{ik})$ on the triple product $U_i\times_X U_j\times_X U_k$. A refinement of the covering gives by definition an isomorphic object. Therefore morphisms $\CV\to\CV'$ need only be defined for systems given over the same covering $\{U_i\to X\}$. In this case (after possibly refining the covering $\{U_i\to X\}$ further) such a morphism is defined by a collection of sections $\phi_i$ of the sheaf $\Hom_{\BLoc_X}(\CF_i,\CF'_i)\otimes_{\BF_q\dbl z\dbr}\BF_q\dpl z\dpr$ over $U_i$ satisfying $\phi'_{ij}\circ pr_i^\ast(\phi_i)=pr_j^\ast(\phi_j)\circ\phi_{ij}$ over $U_i\times_X U_j$.

If $\bar x$ is a geometric point of $X$ we define the \emph{stalk} $\CV_{\bar x}$ of $\CV$ at $\bar x$ as follows. Choose a lift $\bar y$ of $\bar x$ in some $U_i$ and put 
\[
\CV_{\bar x}\es:=\es\CF_{i,\bar y}\otimes_{\BF_q\dbl z\dbr}\BF_q\dpl z\dpr\,.
\]
One easily verifies that $\CV_{\bar x}$ is a well defined finite dimensional $\BF_q\dpl z\dpr$-vector space.

The local systems of $\BF_q\dpl z\dpr$-vector spaces form a category $\PLoc_X$. It is an abelian $\BF_q\dpl z\dpr$-linear rigid tensor category. The theory of these local systems parallels the theory of local systems of $\BQ_p$-vector spaces developed in \cite{dJ}.
In particular there is the following description where for a topological group $G$ we denote by $\ul{\rm Rep}_{\BF_q\dpl z\dpr}G$ the category of continuous representations in finite dimensional $\BF_q\dpl z\dpr$-vector spaces.
\end{definition}

\begin{proposition} \label{Prop2.13}
(Compare \cite[Corollary 5.2]{dJ}.) For any geometric point $\bar x$ of $X$ there is a natural $\BF_q\dpl z\dpr$-linear tensor functor
\[
w_{\bar x}:\es\PLoc_X\es\to\es \ul{\rm Rep}_{\BF_q\dpl z\dpr}\bigl(\pi_1^\et(X,\bar x)\bigr)
\]
which assigns to a local system $\CV$ the $\pi_1^\et(X,\bar x)$-representation $\CV_{\bar x}$. It is an equivalence if $X$ is connected.
\qed
\end{proposition}

Correspondingly one can define the category $\PLoc_{X^\an}$ of local systems of $\BF_q\dpl z\dpr$-vector spaces on the Berkovich space $X^\an$ associated with $X$ as in \cite[\S4]{dJ}.  

\begin{proposition}
(Compare \cite[Proposition 5.1]{dJ}.) There is a natural tensor equivalence of categories $\PLoc_{X^\an}\to\PLoc_X$.
\qed
\end{proposition}

\begin{proposition} \label{Prop2.14}
(Compare \cite[Corollary 4.4]{dJ}.) A local system of $\BF_q\dpl z\dpr$-vector spaces $\CV$ on $X$ can always be given as $\CV\;=\;\bigl(\{U_i\to X\},\CF_i,\phi_{ij}\bigr)$ where the $U_i$ are affinoid subdomains of $X$ and the associated Berkovich spaces $U_i^\an$ form an affinoid covering of $X^\an$ in the sense of Definition~\ref{DefBerkovichSpaces}.
\qed
\end{proposition}

%%%%%%%%%%%%%%%%%%%%%%%%%%%%%%%%%%%%%%%%%%%%%%%%%%%%%%%%%%%%%%%%%%%%%%
%
%    An Approximation Lemma
%
%%%%%%%%%%%%%%%%%%%%%%%%%%%%%%%%%%%%%%%%%%%%%%%%%%%%%%%%%%%%%%%%%%%%%%

\subsection*{A.5\quad An Approximation Lemma} \label{AppApproxLemma}
\addtocounter{subsection}{1}
\setcounter{theorem}{0}

The following lemma generalizes the fact that $L^\sep$ lies dense in $L^\alg$.

\begin{lemma} \label{Lemma6.7'}
Let $B$ be an affinoid $L$-algebra and let $x\in\CM(B)$. Then to any given $b\in\wh{\kappa(x)^\alg}$ and $\epsilon\in\BR_{>0}$ there exists a connected affinoid neighborhood $\CM(B')$ of $x$, a finite \'etale
$B'$-algebra $C$ such that $\CM(C)$ contains a unique point $y$ above $x$, and an element $c\in C$ with $|b-c|_y<\epsilon$.
\end{lemma}

\begin{proof}
By \cite[3.4.1/6]{BGR}, $\kappa(x)^\sep$ is dense in $\wh{\kappa(x)^\alg}$. So we find an element $b'\in\kappa(x)^\sep$ with $|b-b'|_x<\epsilon$.
Consider the minimal polynomial $f=X^n+a'_{n-1}X^{n-1}+\ldots+a'_0$ of $b'$ over $\kappa(x)$. Since $\kappa(x)=\Quot\bigl(B/\ker|\,.\,|_x\bigr)\kompl$ we may approximate the $a'_i$ by 
elements $\bar a_i\in\Quot\bigl(B/\ker|\,.\,|_x\bigr)$ and replace them by the $\bar a_i$, neither changing the separability of $\kappa(x)(b')$ over $\kappa(x)$ nor the inequality $|b-b'|_x<\epsilon$.
We also may assume that $f$ remains irreducible over $\kappa(x)$ since otherwise we replace it by an appropriate irreducible factor and repeat the approximation.
Consider the smallest common denominator $\bar d$ in $B/\ker|\,.\,|_x$ of the $\bar a_i$ and lift it to an element $d\in B$. The Zariski closed subset on which $d$ vanishes does not contain $x$. 
Hence there is an affinoid neighborhood $\CM(B')$ of $x$ on which $d$ is invertible. We obtain $\bar a_i\in B'/\ker|\,.\,|_x$. Choose lifts $a_i\in B'$ of $\bar a_i$ and let $C$ be the finite $B'$-algebra
$B'[X]/(X^n+a_{n-1}X^{n-1}+\ldots+a_0)$. The Zariski closed subset of $\CM(B')$ above which $C$ is not \'etale does not contain $x$. So we may shrink $\CM(B')$ to a connected affinoid neighborhood of $x$ above which $C$ is finite \'etale.
By construction there is exactly one point in $\CM(C)$ which maps to $x$. Putting $c=X$ the lemma follows.
\end{proof}

For further reference let us record a consequence of \cite[Lemma 3.1.6]{JP}.

\begin{lemma}\label{Lemma6.7''}
Let $B$ be an affinoid $L$-algebra and let $x\in\CM(B)$. Consider a finite $B$-algebra $C$ such that there is precisely one point $y\in\CM(C)$ in the fiber over $x$. Let $\CM(C')$ be an affinoid neighborhood of $y$ in $\CM(C)$. Then there exists a connected affinoid neighborhood $\CM(B')$ of $x$ in $\CM(B)$ such that $\CM(B'\otimes_BC)$ is an affinoid neighborhood of $y$ contained in $\CM(C')$.
\qed
\end{lemma}

%%%%%%%%%%%%%%%%%%%%%%%%%%%%%%%%%%%%%%%%%%%%%%%%%%%%%%%%%%%%%%%%%%%%%%
%
%    Extension Lemma
%
%%%%%%%%%%%%%%%%%%%%%%%%%%%%%%%%%%%%%%%%%%%%%%%%%%%%%%%%%%%%%%%%%%%%%%

\subsection*{A.6\quad An Extension Lemma}\label{AppExtensionLemma}
\addtocounter{subsection}{1}
\setcounter{theorem}{0}

We introduce the following notation in the formal setting. Let $B\open$ be an admissible formal $R$-algebra (see \refAppRigFormal) and let $B=B\open\otimes_R L$ be the associated affinoid $L$-algebra. We equip $B$ with an $L$-Banach norm $|\,.\,|$ such that $B\open=\{b\in B:|b|\leq1\}$. We denote by $\ol B$ the $R/\Fm_R$-algebra $B\open/\{b\in B\open:|b|<1\}$. We define the $B\open$-algebras
\begin{tabbing}
$B\con[n]$\es\= \kill
$B\open\langle z\rangle$\>$=\es \DS\bigl\{\;\sum_{i=0}^\infty b_iz^i:\es b_i\in B\open,\;|b_i|\to0 \;(i\to\infty)\;\bigr\}$ \quad and\\[2mm]
$B\open\langle z,z^{-1}\rangle$\> $\DS =\es\bigl\{\,\sum_{i=-\infty}^\infty b_iz^i:\es b_i\in B\open,\;|b_i|\to0 \;(i\to\pm\infty)\;\bigr\}$\,.
\end{tabbing}
We denote by an over-bar the reduction modulo $\Fm_R$, for example $\ol{B\open\langle z\rangle}\;=\;\ol B[z]$.

Note that we do \emph{not assume that $B\open$ is noetherian}. For the notion of coherent sheaf in the non-noetherian case see \cite[0 \S5.3]{EGA}.

\begin{lemma}\label{Lemma2.17b}
Let $\CM$ be a  coherent sheaf on $\Spf B\open\langle z\rangle$. Then $\CM$ is locally free if and only if $\CM$ is flat over $\Spf B\open$ and for every maximal ideal $\Fm$ of $B\open$ the $(B\open/\Fm)[z]$-module $\CM\otimes_{B\open\langle z\rangle}(B\open\!/\Fm)[z]$ is free.
\end{lemma}

\begin{proof}
Since one direction is obvious we will now assume that $\CM$ is flat over $\Spf B\open$ and that $\CM\otimes_{B\open\langle z\rangle}(B\open\!/\Fm)[z]$ is locally free for each maximal ideal $\Fm$ of $B\open$. We must show that for any maximal ideal $\Fn$ of $B\open\langle z\rangle$ the $B\open\langle z\rangle_\Fn$-module $\CM_\Fn$ is free. Set $\Fm:=\Fn\cap B\open$. Since $\Fn$ is maximal we have $\Fm_R\subset\Fn$. Consider the ideal $\ol\Fn=\Fn/\Fm_R B\open\langle z\rangle\subset \ol B[z]$. Since $\ol B$ is of finite type over the field $R/\Fm_R$, it is a Jacobson ring. Hence $\ol\Fm:=\ol\Fn\cap \ol B=\Fm/\Fm_R B\open$ is a maximal ideal of $\ol B$ by \cite[Theorem 4.19]{Eisenbud} and so $\Fm$ is a maximal ideal of $B\open$.

By assumption the $B\open\langle z\rangle_\Fn/\Fm B\open\langle z\rangle_\Fn$-module $\CM_\Fn/\Fm\CM_\Fn$ is finite free. Lifting a basis to $\CM_\Fn$ we obtain an exact sequence of $B\open\langle z\rangle_\Fn$-modules
\[
\TS 0\to \CK\to B\open\langle z\rangle_\Fn^{\oplus\ell}\to\CM_\Fn\to 0\,,
\]
which is exact on the right by Nakayama. The kernel $\CK$ is coherent by \cite[0, Corollaire 5.3.4]{EGA}, because $B\open\langle z\rangle_\Fn$ is a coherent ring by \cite[Proposition 1.3]{FRG1}. Since $\CM_\Fn$ is flat over $B\open_\Fm$ the sequence remains exact after tensoring with $B\open_\Fm/\Fm B\open_\Fm$. In particular $\CK\otimes_{B\open_\Fm}B\open/\Fm=(0)$ and we find $\CK=\Fm \CK$, whence $\CK=0$ again by Nakayama. This shows that $\CM_\Fn$ is free as desired.
\end{proof}

The argument used to prove the following extension lemma is a refinement of Langton's method \cite[Proposition 6]{Langton}.

\begin{lemma}\label{Lemma2.17}
Let $X$ be a quasi-compact admissible formal scheme over $\Spf R$ with associated rigid analytic space $X_L$ over $L$. Let $\CN$ be a finitely generated free sheaf on $X\times_{\Spf R}\Spf R\langle z,z^{-1}\rangle$, let $\CM_L$ be a locally free sheaf on $X_L\times_L\Spm L\langle z\rangle$ of finite rank, and let
\[
\TS\phi:\es\CN\otimes_{\CO_X\langle z,z^{-1}\rangle}\CO_{X_L}\langle z,z^{-1}\rangle\es \isoto\es\CM_L\otimes_{\CO_{X_L}\langle z\rangle}\CO_{X_L}\langle z,z^{-1}\rangle
\]
be an isomorphism. Then there exists an admissible blowing-up $\pi:X'\to X$ such that $\CM_L$ and $\CN$ both come from a locally free sheaf $\CM'$ on $X'\times_{\Spf R}\Spf R\langle z\rangle$ of finite rank with the following property: For any morphism $\beta:\Spf C\open\to X'$ of admissible formal $\Spf R$-schemes we have
\[
\TS\Gamma\bigl(\Spf C\open\langle z\rangle\,,\,\beta^\ast\CM'\bigr)\es=\es\Gamma\bigl(\Spm (C\open\otimes_R L)\langle z\rangle\,,\,\beta^\ast\CM_L\bigr)\;\cap\;\Gamma\bigl(\Spf C\open\langle z,z^{-1}\rangle\,,\,\beta^\ast\pi^\ast\CN\bigr)
\]
inside $\Gamma\bigl(\Spm (C\open\otimes_R L)\langle z,z^{-1}\rangle\,,\,\beta^\ast\CM_L\bigr)$. This property uniquely determines $\CM'$.
\end{lemma}

\begin{proof}
In view of the uniqueness assertion the problem is local on $X$. We use Lemma~\ref{Lemma6.5} to produce an affinoid covering of $X_L$ over which $\CM_L$ is free. By \cite[Theorem5.5]{FRG2} this covering is induced by an open affine covering of $X$ after replacing $X$ by an admissible blowing-up of $X$. Thus we may assume that $X=\Spf B\open$ is affine and $\CM_L$ is associated with a free $B\langle z\rangle$-module $M_L$. We identify $\CN$ with its image under $\phi$. Let $\ul m=(m_1,\ldots,m_\ell)$ be a basis of $M_L$. Let $N$ be the $B\open\langle z,z^{-1}\rangle$-module associated with $\CN$. After multiplication with a power of $\zeta$ we may assume that the $m_i$ belong to $N$. We will show that the intersection $M=M_L\cap N$ inside $M_L\otimes B\langle z,z^{-1}\rangle$ is a coherent $B\open\langle z\rangle$-module with $M\otimes_{B\open\langle z\rangle}B\langle z\rangle=M_L$ and $M\otimes_{B\open\langle z\rangle}B\open\langle z,z^{-1}\rangle=N$.

Let $\ul n=(n_1,\ldots,n_\ell)$ be a basis of the $B\open\langle z,z^{-1}\rangle$-module $N$. The isomorphism between $M_L$ and $N$ yields a matrix $U\in\GL_\ell\bigl(B\langle z,z^{-1}\rangle\bigr)$ with $\ul n=\ul m\,U$, that is $n_j=\sum_i U_{ij}m_i$. We can write $U=z^{-a}U'+U''$ for $U''\in M_\ell\bigl(\zeta B\open\langle z,z^{-1}\rangle\bigr)$, $a\in\BZ$, and $U'\in M_\ell\bigl(B[z]\bigr)$. Since $m_i\in N$ we have
\[
U^{-1}\es\in\es M_\ell\bigl(B\open\langle z,z^{-1}\rangle\bigr)\qquad\text{and}\qquad U^{-1}U'\es=\es z^a(\Id_\ell-U^{-1}U'')\es\in\es\GL_\ell\bigl(B\open\langle z,z^{-1}\rangle\bigr)\,. 
\]
By \cite[Lemma 9.7.1/1]{BGR}, $\det U'\in B\langle z,z^{-1}\rangle\mal\;=\;B\mal\cdot B\open\langle z,z^{-1}\rangle\mal$. So we may write $\det U'=b\cdot f$ with $b\in B\mal$ and $f\in B\open\langle z,z^{-1}\rangle\mal$. We define the basis $\ul n'=(n'_1,\ldots,n'_\ell):=\ul n\,U^{-1}U'f^{-1}$ of $N$. Then $\ul n'=\ul m\,U'f^{-1}$. Considering the adjoint matrix $(U')^\ad$ of $U'$ one sees that
\[
\TS (U')^{-1}f\es=\es b^{-1}(U')^\ad\es=\es (U^{-1}U')^{-1}U^{-1}f\es\in\es M_\ell\bigl(B[z]\cap B\open\langle z,z^{-1}\rangle\bigr)\es=\es M_\ell\bigl(B\open[z]\bigr)\,.
\]
In particular $m_i\in\bigoplus_{j=1}^\ell B\open[z]\cdot n'_j$. Moreover, there exists an integer $e$ such that $U'\in M_\ell\bigl(\zeta^e B\open[z]\bigr)$.

Now let $m\in M=M_L\cap N$. Then $m=\ul m\cdot x=\ul n'\cdot y$ for vectors $x\in B\langle z\rangle^\ell$ and $y\in B\open\langle z,z^{-1}\rangle^\ell$. From the equations $x=U'f^{-1}\,y$ and $y=(U')^{-1}f\,x$ we derive $x\in\bigl(\zeta^e B\open\langle z\rangle\bigr)^\ell$ and $y\in B\open\langle z\rangle^\ell$. Thus $M$ equals the intersection
\begin{equation}\label{EqMIsIntersection}
M\es=\es\bigoplus_{i=1}^\ell B\open\langle{\TS z}\rangle \zeta^e\,m_i\es\cap\es \bigoplus_{j=1}^\ell B\open\langle{\TS z}\rangle\,n'_j
\end{equation}
inside $\bigoplus_{j=1}^\ell B\open\langle{\TS z}\rangle \zeta^e\,n'_j$. By \cite[Proposition 1.3]{FRG1}, $B\open\langle z\rangle$ is a coherent ring. Hence $M$ is a coherent $B\open\langle z\rangle$-module by \cite[0, Corollaire 5.3.6]{EGA}. Since $\ul m$ and $\ul n'\,f=\ul m\,U'$ belong to $M$ and form bases of $M_L$ and $N$ respectively, we find $M\otimes_{B\open\langle z\rangle}B\langle z\rangle=M_L$ and $M\otimes_{B\open\langle z\rangle}B\open\langle z,z^{-1}\rangle=N$. Let $\CM$ be the coherent sheaf on $\Spf B\open\langle z\rangle$ associated with $M$. By \cite[Theorem 4.1]{FRG2} there exists an admissible blowing-up $\pi:X'\to \Spf B\open$ such that $\CM':=\pi^\ast\CM/(\zeta\text{-torsion})$ is a coherent sheaf on $X'\times_{\Spf R}\Spf R\langle z\rangle$ which is flat over $X'$. To show that $\CM'$ is locally free we need the following lemma.

\begin{lemma}\label{Lemma2.17c}
Let $L'$ be a finite extension of $L$ with valuation ring $R'$. Let $C\open$ be an admissible formal $R'$-algebra and let $C:=C\open\otimes_{R'}L'$ be the associated affinoid $L'$-algebra. Let $\beta:\Spf C\open\to X'$ be a morphism of formal $\Spf R$-schemes and let $\beta_L:\Spm C\to X'_L$ be the associated morphism of rigid analytic spaces. Then
\[
\TS\Gamma\bigl(\Spf C\open\langle z\rangle\,,\,\beta^\ast\CM'\bigr)\es=\es M_L\otimes_{B\langle z\rangle}C\langle z\rangle\es\cap\es N\otimes_{B\open\langle z,z^{-1}\rangle}C\open\langle z,z^{-1}\rangle\,.
\]
\end{lemma}

\begin{proof}[Proof of Lemma~\ref{Lemma2.17c}]
We translate the equation (\ref{EqMIsIntersection}) describing $M$ into the following exact sequences of coherent $B\open\langle z\rangle$-modules
\[
\xymatrix  @C-0pc @R=0pc @M+0.3pc
{
0\ar[r] & \DS\bigoplus_{i=1}^\ell B\open\langle{\TS z}\rangle \zeta^e\,m_i \ar[r] & \DS\bigoplus_{j=1}^\ell B\open\langle{\TS z}\rangle \zeta^e\, n'_j \ar[r] & W_1 \ar[r] & 0\,, \\
0\ar[r] & \DS\bigoplus_{j=1}^\ell B\open\langle{\TS z}\rangle n'_j \ar[r] & \DS\bigoplus_{j=1}^\ell B\open\langle{\TS z}\rangle \zeta^e\, n'_j \ar[r] & W_2 \ar[r] & 0\,, \\
0\ar[r] & M \ar[r] & \DS\bigoplus_{j=1}^\ell B\open\langle{\TS z}\rangle \zeta^e\, n'_j \ar[r] & W_1\oplus W_2 \ar[r] & 0\,.
}
\]
We tensor with $C\open\langle z\rangle$ over $B\open\langle z\rangle$ and obtain the following commutative diagram with exact rows
\begin{equation}\label{EqW1Tensored}
\xymatrix @C-0.5pc @R-0.5pc 
{
0\ar@{-->}[r] & *!U(0.0) 
\objectbox{ \es\DS\bigoplus_{i=1}^\ell C\open\langle{\TS z}\rangle \zeta^e\,m_i \es}
\ar[r]\ar@{^{ (}->}[d] & 
*!U(0.05) 
\objectbox{ \DS\es\bigoplus_{j=1}^\ell C\open\langle{\TS z}\rangle \zeta^e\, n'_j\es} \ar[r]\ar@{^{ (}->}[d] & 
*!U(0.15) 
\objectbox{ \es W_1\otimes_{B\open\langle{\TS z}\rangle}C\open\langle{\TS z}\rangle \es} \ar[r] & 0 \\
0\ar[r] & *!U(0.0) 
\objectbox{ \es\DS\bigoplus_{i=1}^\ell C\langle{\TS z,z^{-1}}\rangle \zeta^e\,m_i\es} \ar[r]^\sim & 
*!U(0.05) 
\objectbox{ \DS\es\bigoplus_{j=1}^\ell C\langle{\TS z,z^{-1}}\rangle \zeta^e\, n'_j\es} \ar[r] & 0 \,.
}
\end{equation}
Here the injectivity of the vertical arrows shows that the upper row is also exact on the left. Similarly we find that the sequence
\begin{equation}\label{EqW2Tensored}
\xymatrix @M+0.3pc
{
0\ar[r] & \DS\bigoplus_{j=1}^\ell C\open\langle{\TS z}\rangle n'_j 
\ar[r] & \DS\bigoplus_{j=1}^\ell C\open\langle{\TS z}\rangle \zeta^e\, n'_j \ar[r] & *!U(0.05) \objectbox{\es W_2\otimes_{B\open\langle{\TS z}\rangle}C\open\langle{\TS z}\rangle\es} \ar[r] & 0 
}
\end{equation}
is exact. Consider the commutative diagram with exact first row
\begin{equation}\label{EqMTensored}
\xymatrix @C-1pc @R-0.5pc
{
0\ar[r] & V\ar[r] & **{!U(0.15) +<0.0pc,1pc> } \objectbox{\;M\otimes_{B\open\langle{\TS z}\rangle}C\open\langle{\TS z}\rangle\;} \ar[r]\ar[d]^f & 
*!U(0.05) \objectbox{\;\DS\bigoplus_{j=1}^\ell C\open\langle{\TS z}\rangle \zeta^e\, n'_j\;} \ar[r]\ar@{^{ (}->}[d] & 
*!U(0.25) \objectbox{\;(W_1\oplus W_2)\otimes_{B\open\langle{\TS z}\rangle}C\open\langle{\TS z}\rangle\;} \ar[r] & 0 \\
& & %*!U(0.4) 
**{+<0.0pc,1pc>} \objectbox{M\otimes_{B\open\langle{\TS z}\rangle}C\langle{\TS z}\rangle}
 \ar[r]\ar@{^{ (}->}[d]^g & \DS\bigoplus_{j=1}^\ell C\langle{\TS z}\rangle \zeta^e\, n'_j \ar@{^{ (}->}[d] \\
& & %*!U(0.4)
**{!R(0.2) +<0.2pc,1pc>} \objectbox{ M\otimes_{B\open\langle{\TS z}\rangle}C\langle{\TS z,z^{-1}}\rangle}
 \ar@{=}[r] & *!L(0.2) \objectbox{\;\DS\bigoplus_{j=1}^\ell C\langle{\TS z,z^{-1}}\rangle \zeta^e\, n'_j}
}
\end{equation}
The equality in the last row follows from the fact that $M\otimes_{B\open\langle z\rangle}B\langle z,z^{-1}\rangle\;=\;{\bigoplus_{j=1}^\ell}B\langle z,z^{-1}\rangle \zeta^e n'_j$. The vertical map $g$ is injective because $M\otimes_{B\open\langle{\TS z}\rangle}C\langle{\TS z}\rangle\;=\;M_L\otimes_{B\langle{\TS z}\rangle}C\langle{\TS z}\rangle$ is flat over $C\langle{\TS z}\rangle$. Let ${}_\zeta(\pi^\ast\CM)$ be the $\zeta$-torsion of $\pi^\ast\CM$. The sequence
\[
\xymatrix @M+0.3pc
{
0\ar[r] & \beta^\ast {}_\zeta(\pi^\ast\CM)\ar[r] & \beta^\ast\pi^\ast\CM \ar[r] & \beta^\ast\CM'\ar[r] & 0
}
\]
is exact on the left since $\CM'$ is flat over $X'$. We obtain the following commutative diagram with exact rows
\[
\xymatrix @C-1pc %@R-1pc %@M+0.3pc
{
0\ar[r] &\,\Gamma\bigl(\Spf C\open\langle z\rangle\,,\,\beta^\ast {}_\zeta(\pi^\ast\CM)\bigr)\, \ar[r] & 
**{!U(0.1) +<0pc,0.5pc>} \objectbox{\;M\otimes_{B\open\langle{\TS z}\rangle}C\open\langle{\TS z}\rangle\;} \ar[r]\ar[d]^f & **{+<0pc,0.5pc>} \objectbox{\;\Gamma\bigl(\Spf C\open\langle z\rangle\,,\,\beta^\ast\CM'\bigr)\;} \ar[r]\ar@{^{ (}->}[d] & 0 \\
& 0\ar[r] & **{!U(0.15) +<0pc,0.5pc>} \objectbox{\;M\otimes_{B\open\langle{\TS z}\rangle}C\langle{\TS z}\rangle\;} \ar[r] & **{!U(0.1) +<0pc,0.5pc>} \objectbox{\;\Gamma\bigl(\Spf C\open\langle z\rangle\,,\,\beta^\ast\CM'\bigr) \otimes_{C\open}C\;} \ar[r] & 0\,.
}
\]
The second row is exact since $\beta^\ast{}_\zeta(\pi^\ast\CM)\otimes_{C\open}C=0$. The vertical map on the right is injective since $\beta^\ast\CM'$ is flat over $C\open$. We conclude that 
\[
\TS V\es=\es \ker f\es=\es\Gamma\bigl(\Spf C\open\langle z\rangle\,,\,\beta^\ast {}_\zeta(\pi^\ast\CM)\bigr)
\]
and
\[
\TS\Gamma\bigl(\Spf C\open\langle z\rangle\,,\,\beta^\ast\CM'\bigr)\es=\es M\otimes_{B\open\langle{\TS z}\rangle}C\open\langle{\TS z}\rangle/V\es=\es M_L\otimes_{B\langle z\rangle}C\langle z\rangle\es\cap\es N\otimes_{B\open\langle z,z^{-1}\rangle}C\open\langle z,z^{-1}\rangle
\]
as claimed.
\end{proof}

\noindent
\emph{Proof of Lemma~\ref{Lemma2.17} continued}. To finish the proof we will show that $\CM'$ is locally free on $X'\times_{\Spf R}\Spf R\langle z\rangle$ using Lemma~\ref{Lemma2.17b}. Let $x\in X'$ be a closed point with residue field $\kappa(x)$. We have to show that $\CM'\otimes_{\CO_{X'}\langle z\rangle}\kappa(x)[z]$ is a free $\kappa(x)[z]$-module. By \cite[3.5 and 3.2]{FRG1} there exists a rank-$1$ valuation ring $R'$ extending $R$ whose fraction field $L'$ is a finite extension of $L$, and a morphism of formal $R$-schemes $\beta:\Spf R'\to X'$ which makes the residue field $k'$ of $R'$ into a finite extension of the residue field $\kappa(x)$ of $x$. By faithfully flat descent \cite[IV$_2$, Proposition 2.5.2]{EGA} it suffices to show that $\beta^\ast\CM'\otimes_{R'} k'$ is 
a free module over $R'\langle z\rangle\otimes_{R'}k'=k'[z]$. We put $M'_{R'}:=\Gamma\bigl(\Spf R'\langle z\rangle\,,\,\beta^\ast\CM'\bigr)$ and $N_{R'}:=\Gamma\bigl(\Spf R'\langle z,z^{-1}\rangle\,,\,\beta^\ast\pi^\ast\CN\bigr)$ and consider the restriction map
\[
M'_{R'}\es\longto\es\Gamma\bigl(\Spf R'\langle z,z^{-1}\rangle\,,\,\beta^\ast\CM'\bigr)\es=\es N_{R'}\,.
\]
We claim that the induced morphism $M'_{R'}/\Fm_{R'}M'_{R'}\to N_{R'}/\Fm_{R'}N_{R'}$ is injective. Indeed $M'_{R'}\cap\Fm_{R'}N_{R'}\;=\;\Fm_{R'}M'_{R'}$ since by Lemma~\ref{Lemma2.17c}, $M'_{R'}\;=\;N_{R'}\cap M_L\otimes_{B\langle z\rangle}L'\langle z\rangle$ and $M_L\otimes_{B\langle z\rangle}L'\langle z\rangle\;=\;\Fm_{R'}\cdot M_L\otimes_{B\langle z\rangle}L'\langle z\rangle$. Hence $M'_{R'}\otimes_{R'}k'\subset N_{R'}\otimes_{R'}k'$ and since $N'_{R'}\otimes_{R'}k'$ is a torsion free $k'[z]$-module, the same is true for $M'_{R'}\otimes_{R'}k'$. We conclude that $\beta^\ast\CM'\otimes_{R'}k'$ is a free $k'[z]$-module as desired.
\end{proof}

We apply Lemma~\ref{Lemma2.17} in Section~\ref{SectCriteriaWA} to extend $\sigma$-modules over rigid analytic spaces to formal models. For this purpose let us introduce the following $B\open$-algebra
\begin{tabbing}
$B\con[n]$\es\= \kill
$B\open\dbl z,z^{-1}\rangle$ \>$\DS =\es\bigl\{\;\sum_{i=-\infty}^\infty b_i z^i:\es b_i\in B\open,\; |b_i|\to0\es(i\to -\infty)\;\bigr\}$\,.
\end{tabbing}
Note that $B\open\dbl z\dbr\otimes_{B\open}B$ is strictly contained in the $B$-algebra $B\{z\}$ defined in Section~\ref{SectNotation} since for any element $\sum b_i z^i$ of $B\open\dbl z\dbr\otimes_{B\open}B$ the sequence $|b_i|$ is bounded. In contrast this is not true in $B\{z\}$. With this notation we may deduce from the proof of Lemma~\ref{Lemma2.17} also the following 

\begin{lemma}\label{Lemma2.18}
Keep the situation of Lemma~\ref{Lemma2.17}. If in addition we are given morphisms
\[
\begin{array}{lcl}
F_{\CM_L}: \quad\sigma^\ast\CM_L\otimes_{B\langle z\rangle}B\{z\} &\longto & \CM_L\otimes_{B\open\langle z\rangle}B\{z\} \qquad\text{and}\\[2mm]
F_\CN:\es \sigma^\ast\CN\otimes_{B\open\langle z,z^{-1}\rangle}B\open\dbl z,z^{-1}\rangle & \longto & \CN\otimes_{B\open\langle z,z^{-1}\rangle}B\open\dbl z,z^{-1}\rangle
\end{array}
\]
with $\phi\circ F_\CN=F_{\CM_L}\circ\sigma^\ast\phi$ then there exists a uniquely determined morphism
\[
F_{\CM'}:\sigma^\ast\CM'\otimes_{\CO_{X'}\langle z\rangle}\CO_{X'}\dbl z\dbr \es\longto\es \CM'\otimes_{\CO_{X'}\langle z\rangle}\CO_{X'}\dbl z\dbr
\]
of sheaves of $\CO_{X'}\dbl z\dbr$-modules on $X'$ extending $F_{\CM_L}$ and $F_\CN$.
\end{lemma}

We start the proof with the following 

\begin{lemma}\label{Lemma2.18c}
In the situation of Lemma~\ref{Lemma2.17c}, $\Gamma\bigl(\Spf C\open\langle z\rangle\,,\,\beta^\ast\CM'\bigr)\otimes_{C\open\langle z\rangle}C\open\dbl z\dbr$ equals the intersection
\[
\TS M_L\otimes_{B\langle z\rangle}\bigl(C\open\dbl z\dbr\otimes_{C\open}C\bigr)\es\cap\es N\otimes_{B\open\langle z,z^{-1}\rangle}C\open\dbl z,z^{-1}\rangle
\]
inside $M_L\otimes_{B\langle z\rangle}\bigl(C\open\dbl z,z^{-1}\rangle\otimes_{C\open}C\bigr)$.
\end{lemma}

\begin{proof}
This follows literally as in the proof of Lemma~\ref{Lemma2.17c}. One only replaces $C\open\langle z\rangle$ by $C\open\dbl z\dbr$, $C\langle z\rangle$ by $C\open\dbl z\dbr\otimes_{C\open}C$, and $C\langle z,z^{-1}\rangle$ by $C\open\dbl z,z^{-1}\rangle\otimes_{C\open}C$.
\end{proof}

\begin{proof}[Proof of Lemma~\ref{Lemma2.18}]
We continue with the notation introduced in the proof of Lemma~\ref{Lemma2.17}. In particular $M_L=\bigoplus_i B\langle z\rangle m_i$ and $N=\bigoplus_j B\open\langle z,z^{-1}\rangle n'_j$. The bases $\ul m=(m_1,\ldots,m_\ell)$ and $\ul n'=(n'_1,\ldots,n'_\ell)$ are related by $\ul n'=\ul m\;U'f^{-1}$ and $\ul m=\ul n'\,(U')^{-1}f$ where $f\in B\open\langle z,z^{-1}\rangle\mal$, $U'\in M_\ell\bigl(\zeta^eB\open[z]\bigr)$, and $(U')^{-1}f\in M_\ell\bigl(B\open[z]\bigr)$.

Let $F_{\CM_L}(\sigma^\ast\ul m)=\ul m\,\Phi$ and $F_\CN(\sigma^\ast\ul n')=\ul n'\,\Psi$ with matrices 
$\Phi\in M_\ell\bigl(B\{z\}\bigr)$ and $\Psi\in M_\ell\bigl(B\open\dbl z,z^{-1}\rangle\bigr)$. From the compatibility relation $\phi\circ F_\CN=\CF_{\CM_L}\circ\sigma^\ast\phi$ we deduce
\[
\Phi\;=\;U'f^{-1}\cdot \Psi\,\bigl((U')^{-1}f\bigr)^\sigma\es\in\es M_\ell\bigl(B\{z\}\cap\zeta^e B\open\dbl z,z^{-1}\rangle\bigr)\;\subset\;M_\ell\bigl(B\open\dbl z\dbr\otimes_{B\open}B\bigr)\,.
\]
So if $\beta:\Spf C\open\hookrightarrow X'$ is an affine open subset and $m\in\Gamma\bigl(\Spf C\open\langle z\rangle\,,\,\sigma^\ast\CM'\bigr)\otimes_{C\open\langle z\rangle}C\open\dbl z\dbr$ then $F_{\CM_L}(m)\;\in\;M_L\otimes_{B\langle z\rangle}\bigl(C\open\dbl z\dbr\otimes_{C\open}C\bigr)$ and $F_\CN(m)\;\in\;N\otimes_{B\open\langle z,z^{-1}\rangle}C\open\dbl z,z^{-1}\rangle$. Using Lemma~\ref{Lemma2.18c} we obtain the desired morphism 
\[
F_\CM:\es\Gamma\bigl(\Spf C\open\langle z\rangle\,,\,\sigma^\ast\CM'\bigr)\otimes_{C\open\langle z\rangle}C\open\dbl z\dbr\es\longto\es\Gamma\bigl(\Spf C\open\langle z\rangle\,,\,\CM'\bigr)\otimes_{C\open\langle z\rangle}C\open\dbl z\dbr\,.
\]
\end{proof}

\end{appendix}

%%%%%%%%%%%%%%%%%%%%%%%%%%%%%%%%%%%%%%%%%%%%%%%%%%%%%%%%%%%%%%%%%%%%%%
%
%    Bibliography
%
%%%%%%%%%%%%%%%%%%%%%%%%%%%%%%%%%%%%%%%%%%%%%%%%%%%%%%%%%%%%%%%%%%%%%%

\vfill

\noindent
Urs Hartl\\
University of Muenster\\
Institute of Mathematics\\
Einsteinstr.~62\\
D -- 48149 Muenster\\
Germany\\[1mm]
http:/\!/www.math.uni-muenster.de/u/urs.hartl/

\end{document}